\documentclass[thmsa,onecolumn,10pt,11pt]{article}
\usepackage{amssymb}
\usepackage{pstricks, pst-plot, pst-node, pst-tree}
\usepackage{amsfonts}
\usepackage{graphicx}
\usepackage{amsmath}

\newtheorem{theorem}{Theorem}[section]

\newtheorem{corollary}{Corollary}[section]

\newtheorem{definition}{Definition}

\newtheorem{lemma}{Lemma}[section]

\newtheorem{proposition}{Proposition}[section]

\newenvironment{proof}[1][Proof]{\textbf{#1.} }{\ \rule{0.5em}{0.5em}}
\begin{document}

\title{\textbf{Moments, cumulants and diagram formulae for non-linear
functionals of random measures }}
\author{Giovanni PECCATI\thanks{Equipe Modal'X, Universit\'{e}
Paris Ouest Nanterre La D\'{e}fense, 200 Avenue de la R\'{e}publique, 92000 Nanterre and LSTA, Universit\'{e} Paris VI, France. E-mail: \texttt{giovanni.peccati@gmail.com}} $\,$ and Murad S.\ TAQQU\thanks{%
Boston University, Departement of Mathematics, 111 Cummington Road, Boston
(MA), USA. E-mail: \texttt{murad@math.bu.edu}}}

\maketitle

\begin{abstract}
This survey provides a unified discussion of multiple integrals, moments,
cumulants and diagram formulae associated with functionals of completely
random measures. Our approach is combinatorial, as it is based on the
algebraic formalism of partition lattices and M\"{o}bius functions. Gaussian
and Poisson measures are treated in great detail. We also present several
combinatorial interpretations of some recent CLTs involving sequences of
random variables belonging to a fixed Wiener chaos.


\end{abstract}

\tableofcontents

\newpage

\section{Introduction}

\subsection{Overview}

The aim of this survey is to provide a unified treatment of moments and
cumulants associated with non-linear functionals of completely random
measures, such as Gaussian, Poisson or Gamma measures. We will focus on
multiple stochastic integrals, and we shall mainly adopt a combinatorial
point of view. In particular, our main inspiration is a truly remarkable
paper by Rota and Wallstrom \cite{RoWa}, building (among many others) on
earlier works by It\^{o} \cite{Ito}, Meyer \cite{Mey78, Mey92} and,
most importantly, Engel \cite{Engel} (see also Bitcheler \cite{Bitcheler},
Kussmaul \cite{Kussmaul}, Linde \cite{WLindeBook}, Masani \cite{Masani},
Neveu \cite{Neveu-1968} and Tsilevich and Vershik \cite{TV} for related
works). In particular, in \cite{RoWa} the authors point out a crucial
connection between the machinery of multiple stochastic integration and the
structure of the lattice of partitions of a finite set, with specific
emphasis on the role played by the associated M\"{o}bius function (see e.g.
\cite{Aigner}, as well as Section \ref{S : Lattice} below). As we will see
later on, the connection between multiple stochastic integration and
partitions is given by the natural isomorphism between the partitions of the
set $\left\{ 1,...,n\right\} $ and the \textsl{diagonal sets }associated
with the Cartesian product of $n$ measurable spaces (a diagonal set is just
a subset of the Cartesian product consisting of points that have two or more
coordinates equal).

The best description of the approach to stochastic integration followed in
the present survey is still given by the following sentences, taken from
\cite{RoWa}:

\begin{quotation}
The basic difficulty of stochastic integration is the following. We are
given a measure $\varphi $ on a set $S$, and we wish to extend such a
measure to the product set $S^{n}$. There is a well-known and established
way of carrying out such an extension, namely, taking the product measure.
While the product measure is adequate in most instances dealing with a
scalar valued measure, it turns out to be woefully inadequate when the
measure is vector-valued, or, in the case dealt with presently,
random-valued. The product measure of a nonatomic scalar measure will vanish
on sets supported by lower-dimensional linear subspaces of $S^{n}$. This is
not the case, however, for random measures. The problem therefore arises of
modifying the definition of product measure of a random measure in such a
way that the resulting measure will vanish on lower-dimensional subsets of $%
S^{n}$, or diagonal sets, as we call them.
\end{quotation}

As pointed out in \cite{RoWa}, as well as in Section \ref{S : MCA} below,
the combinatorics of partition lattices provide the correct framework in
order to define a satisfactory stochastic product measure.

As apparent from the title, in the subsequent sections a prominent role will
be played by moments and cumulants. In particular, the principal aims of our
work are the following:

\begin{description}
\item[--] \textbf{Put diagram formulae in a proper algebraic setting}\textsl{%
.} Diagram formulae are mnemonic devices, allowing to compute moments and
cumulants associated with one or more random variables. These tools have
been developed and applied in a variety of frameworks: see e.g. \cite{Shir, Sur}
for diagram formulae associated with general random variables;
see \cite{BrMa, ChaSlud, GiSu, Marinucci} for
non-linear functionals of Gaussian fields; see \cite{Surg1984} for
non-linear functionals of Poisson measures. They can be quite useful in the
obtention of Central Limit Theorem (CLTs) by means of the so-called \textsl{%
method of moments and cumulants }(see e.g. \cite{Major}).
Inspired by the works by McCullagh \cite{Mcc}, Rota and Shen \cite{Ro SHen}
and Speed \cite{Speed}, we shall show that all diagram formulae quoted above
can be put in a unified framework, based on the use of partitions of finite
sets. Although somewhat implicit in the previously quoted references, this
clear algebraic interpretation of diagrams is new. In particular, in Section %
\ref{S : DG} we will show that all diagrams encountered in the probabilistic
literature (such as Gaussian, non-flat and connected diagrams) admit a neat
translation in the combinatorial language of partition lattices.

\item[--] \textbf{Illustrate the Engel-Rota-Wallstrom theory}\textsl{.} We
shall show that the theory developed in \cite{Engel} and \cite{RoWa} allows
to recover several crucial results of stochastic analysis, such as
multiplication formulae for multiple Gaussian and Poisson integrals see
\cite{Kab, Nualart, Surg1984}. This extends the content of
\cite{RoWa}, which basically dealt with product measures. See also \cite%
{FarreJoUtz} for other results in this direction.

\item[--] \textbf{Enlight the combinatorial implications of new CLTs}\textsl{%
.} In a recent series of papers (see \cite{MaPeAb, NouNu, NouPeccPTRF, NouPecexact, NouPecReveillac, NuOrtiz, PNu05, Pecp, PecSoleTaqUtz, PeTaq2bleP, PeTaqMwi, PTu04}), a new set of tools has been developed, allowing
to deduce simple CLTs involving random variables having the form of multiple
stochastic integrals. All these results can be seen as simplifications of
the method of moments and cumulants. In Section \ref{S : SImpliCLT}, we will
illustrate these results from a combinatorial standpoint, by providing some
neat interpretations in terms of diagrams and graphs. In particular, we will
prove that in these limit theorems a fundamental role is played by the so-called
\textsl{circular diagrams}, that is, connected Gaussian diagrams whose edges
only connect subsequent rows.
\end{description}

We will develop the necessary combinatorial tools related to partitions,
diagram and graphs from first principles in Section \ref{S : Lattice} and
Section \ref{S : DG}. Section \ref{S : cum} provides a self-contained
treament of moments and cumulants from a combinatorial point of view.
Stochastic integration is introduced in Section \ref{S : MCA}. Section \ref%
{S : MF} and Section \ref{S : DF} deal, respectively, with product formulae
and diagram formulae. In Section \ref{S : IsonormalGP} one can find an
introduction to the concept of isonormal Gaussian process. Finally, Section %
\ref{S : SImpliCLT} deals with CLTs on Wiener chaos.

\subsection{Some related topics}

In this survey, we choose to follow a very precise path, namely starting with the basic properties of partition lattices and diagrams,
and develop from there as many as possible of the formulae associated with products, moments and cumulants in the theory of stochastic integration with respect to completely random measures. In order to keep the length of the present work within bounds, several crucial
topics are not included (or are just mentioned) in the discussion to follow. One remarkable
omission is of course a complete discussion of the connections between multiple stochastic integrals and orthogonal polynomials. This topic is partially treated in Section \ref{S : IsonormalGP} below, in the particular case of Gaussian processes. For recent references on more general stochastic processes (such as L\'{e}vy processes), see e.g. the monograph by Schoutens \cite{SchoutBook} and the two papers by Sol\'{e} and Utzet \cite{SolUtzAOP, SolUtzBer}. Other related (and missing) topics are detailed in the next list, whose entries are followed by a
brief discussion.

\begin{description}
\item[--] \textsl{Wick products}. Wick products are intimately related to
chaotic expansions. A complete treatment of this topic can be found e.g. in
Janson's book \cite{Janson}.

\item[--] \textsl{Malliavin calculus. }See the two monographs by Nualart
\cite{Nualart2, Nualart} for Malliavin calculus in a Gaussian
setting. A good introduction to Malliavin calculus for Poisson measures is contained in
the classic papers by Nualart and Vives \cite{NV}, Privault \cite{P1} and
Privault and Wu \cite{PrivaultWu}. A fundamental connection between Malliavin
operators and limit theorems has been first pointed out in \cite{NuOrtizSPA}%
. See \cite{NouNu, NouPec2007, NouPeccPTRF, PecSoleTaqUtz} for further developments.

\item[--] \textsl{Hu-Meyer formulae.} Hu-Meyer formulae connect
Stratonovich multiple integrals and multiple Wiener-It\^{o} integrals. See
\cite{Nualart} for a standard discussion of this topic in a Gaussian
setting. Hu-Meyer formulae for general L\'{e}vy processes can be naturally
obtained by means of the theory described in the present survey: see the
excellent paper by Farria \textit{et al. }\cite{FarreJoUtz} for a complete
treatment of this point.\textit{\ }

\item[--] \textsl{Stein's method}. Stein's method for normal and non-normal
approximation can be a very powerful tool, in order to obtain central and
non-central limit theorems for non-linear functionals of random fields. See
\cite{STeinbook} for a classic reference on the subject. See \cite%
{NouPeccPTRF, NouPecexact, NouPecReveillac} for several
limit theorems involving functionals of Gaussian fields, obtained by means
of Stein's method and Malliavin calculus. See \cite{PecSoleTaqUtz} for an
application of Stein's method to functionals of Poisson measures.

\item[--] \textsl{Free probability. }The properties of the lattice of
(non-crossing) partitions and the corresponding M\"{o}bius function are
crucial in free probability. See the monograph by Nica and Speicher \cite%
{Nica Speicher} for a valuable introduction to the combinatorial aspects of
free probability. See Anshelevich \cite{Ans1, Ans2} for some
instances of a \textquotedblleft free\textquotedblright\ theory of multiple
stochastic integration.
\end{description}

\section{The lattice of partitions of a finite set \label{S : Lattice}}

\setcounter{equation}{0}In this section we recall some combinatorial results
concerning the lattice of partitions of a finite set. These objects play an
important role in the obtention of the \textsl{diagram formulae} presented
in Section \ref{S : MCA}. The reader is referred to Stanley \cite[Chapter 3]%
{Stanley} and Aigner \cite{Aigner} for a detailed presentation of (finite)
partially ordered sets and M\"{o}bius inversion formulae.

\subsection{Partitions of a positive integer}

Given an integer $n\geq 1$, we define the set $\Lambda \left( n\right) $ of
\textsl{partitions}\textit{\ }of $n$ as the collection of all vectors of the
type $\lambda =\left( \lambda _{1},...,\lambda _{k}\right) $ ($k\geq 1$),
where:
\begin{equation}
\begin{array}{l}
\text{(i) \ }\lambda _{j}\text{\ is an integer for every \ }j=1,...,k, \\
\text{(ii) \ }\lambda _{1}\geq \lambda _{2}\geq \cdot \cdot \cdot \geq
\lambda _{k}\geq 1, \\
\text{(iii) }\lambda _{1}+\cdot \cdot \cdot +\lambda _{k}=n\text{.}%
\end{array}
\label{PartEnt}
\end{equation}%
We call $k$ the \textit{length }of $\lambda $. It is sometimes convenient to
write a partition $\lambda =\left( \lambda _{1},...,\lambda _{k}\right) \in
\Lambda \left( n\right) $ in the form $\lambda =\left(
1^{r_{1}}2^{r_{2}}\cdot \cdot \cdot n^{r_{n}}\right) $. This representation
(which encodes all information about $\lambda $) simply indicates that, for
every $i=1,...,n$, the vector $\lambda $ contains exactly $r_{i}$ ($\geq 0$)
components equal to $i$. Clearly, if $\lambda =\left( \lambda
_{1},...,\lambda _{k}\right) =\left( 1^{r_{1}}2^{r_{2}}\cdot \cdot \cdot
n^{r_{n}}\right) \in \Lambda \left( n\right) $, then
\begin{equation}\label{e : sumlabel}
 1r_{1}+\cdot \cdot
\cdot +nr_{n}=n
\end{equation}
and $r_{1}+\cdot \cdot \cdot +r_{n}=k$. We will sometimes
use the (more conventional) notation
\begin{equation*}
\lambda \vdash n\text{ \ \ instead of \ }\lambda \in \Lambda \left( n\right)
.
\end{equation*}

\bigskip

\textbf{Examples. }(i) If $n=5$, one can e.g. have $5=4+1$ or $5=1+1+1+1+1$.
In the first case the length is $k=2$, with $\lambda _{1}=4$ and $\lambda
_{2}=1$, and the partition is $\lambda =\left(
1^{1}2^{0}3^{0}4^{1}5^{0}\right) $. In the second case, the length is $k=5$
with $\lambda _{1}=...=\lambda _{5}=1$, and the partition is $\lambda
=\left( 1^{5}2^{0}3^{0}4^{0}5^{0}\right) .$

(ii) One can go easily from one representation to the other. Thus $\lambda
=\left( 1^{2}2^{3}3^{0}4^{2}\right) $ corresponds to%
\begin{equation*}
n=\left( 1\times 2\right) +\left( 2\times 3\right) +\left( 3\times 0\right)
+\left( 4\times 2\right) =16,
\end{equation*}%
that is, to the decomposition $16=4+4+2+2+2+1+1$, and thus to
\begin{equation*}
\lambda =\left( \lambda _{1},\lambda _{2},\lambda _{3},\lambda _{4},\lambda
_{5},\lambda _{6},\lambda _{7}\right) =\left( 4,4,2,2,2,1,1\right) \text{.}
\end{equation*}

\subsection{Partitions of a set}

Let $b$ denote a finite nonempty set and let
\begin{equation*}
\mathcal{P}\left( b\right) \text{ be the
set of \textsl{partitions }of } b.
\end{equation*}
By definition, an element $\pi $ of $%
\mathcal{P}\left( b\right) $ is a collection of nonempty and disjoint
subsets of $b$ (called \textsl{blocks}), such that their union equals $b$.
The symbol $\left\vert \pi \right\vert $ indicates the number of blocks (or
the \textsl{size}) of the partition $\pi $.

\bigskip

\textbf{Remark on notation. }For each pair $i,j\in b$ and for each $\pi \in
\mathcal{P}\left( b\right) $, we write $i\sim _{\pi }j$ whenever $i$ and $j$
belong to the same block of $\pi $.

\bigskip

We now define a partial ordering on $\mathcal{P}\left( b\right) $. For every
$\sigma ,\pi \in \mathcal{P}\left( b\right) $, we write $\sigma \leq \pi $
if, and only if, each block of $\sigma $ is contained in a block of $\pi $.
Borrowing from the terminology used in topology one also says that $\pi $ is
\textsl{coarser }than $\sigma $. It is clear that $\leq $ is a \textsl{%
partial ordering relation}, that is, $\leq $ is a \textsl{binary}, \textsl{%
reflexive}, \textsl{antisymmetric} and \textsl{transitive} relation on $%
\mathcal{P}\left( b\right) $ (see e.g. \cite[p. 98]{Stanley}). Also, $\leq $
induces on $\mathcal{P}\left( b\right) $ a \textsl{lattice }structure.
Recall that a lattice is a partially ordered set such that each pair of
elements has a least upper bound and a greatest lower bound (see \cite[p. 102%
]{Stanley}). In particular, the partition $\sigma \wedge \pi $, \textsl{meet}
of $\sigma ,\pi \in \mathcal{P}\left( b\right) $, is the partition of $b$
such that each block of $\sigma \wedge \pi $ is a nonempty intersection
between one block of $\sigma $ and one block of $\pi $. On the other hand,
the partition $\sigma \vee \pi $, \textsl{join} of $\sigma ,\pi \in \mathcal{%
P}\left( b\right) $, is the element of $\mathcal{P}\left( b\right) $ whose
blocks are constructed by taking the non-disjoint unions of the blocks of $%
\sigma $ and $\pi $, that is, by taking the union of those blocks that have
at least one element in common.

\bigskip

\textbf{Remarks. }(a) Whenever $\pi _{1}\leq \pi _{2}$, one has $\left\vert
\pi _{1}\right\vert \geq \left\vert \pi _{2}\right\vert $. In particular, $%
\left\vert \sigma \wedge \pi \right\vert \geq \left\vert \sigma \vee \pi
\right\vert $.

(b) The partition $\sigma \wedge \pi $ is the greatest lower bound associated
with the pair $\left( \sigma ,\pi \right) $. As such, $\sigma \wedge \pi $
is completely characterized by the property of being the unique element of $%
\mathcal{P}\left( b\right) $ such that: (i) $\sigma \wedge \pi \leq \sigma $%
, (ii)\ $\sigma \wedge \pi \leq \pi $, and (iii) $\rho \leq \sigma \wedge
\pi $ for every $\rho \in \mathcal{P}\left( b\right) $ such that $\rho \leq
\sigma ,\pi $.

(c) Analogously, the partition $\sigma \vee \pi $ is the least upper
bound associated with the pair $\left( \sigma ,\pi \right) $. It follows
that $\sigma \vee \pi $ is completely characterized by the property of being
the unique element of $\mathcal{P}\left( b\right) $ such that: (i) $\sigma
\leq \sigma \vee \pi $, (ii)\ $\pi \leq \sigma \vee \pi $, and (iii) $\sigma
\vee \pi \leq \rho $ for every $\rho \in \mathcal{P}\left( b\right) $ such
that $\sigma ,\pi \leq \rho $.

\bigskip

\textbf{Examples. }(i) Take $b=\left\{ 1,2,3,4,5\right\} $. If $\pi
=\left\{ \left\{ 1,2,3\right\} ,\left\{ 4,5\right\} \right\} $ and $\sigma
=\left\{ \left\{ 1,2\right\} ,\left\{ 3\right\} ,\left\{ 4,5\right\}
\right\} .$ Then, $\sigma \leq \pi $ (because each block of $\sigma $ is
contained in a block of $\pi $) and
\begin{equation*}
\sigma \wedge \pi =\sigma \text{ \ \ and \ \ }\sigma \vee \pi =\pi \text{.}
\end{equation*}%
A graphical representation of $\pi $, $\sigma $, $\sigma \wedge \pi $ and $%
\sigma \vee \pi $ is:%
\begin{eqnarray*}
\pi &=&\fbox{$1$ $2\ \ \ 3$}\fbox{$4$ $5$} \\
\sigma &=&\fbox{$1$ $2$ }\fbox{$3$}\fbox{$4$ $5$} \\
\sigma \wedge \pi &=&\fbox{$1$ $2$ }\fbox{$3$}\fbox{$4$ $5$} \\
\sigma \vee \pi &=&\fbox{$1$ $\ 2\ \ 3$}\fbox{$4$ $5$}
\end{eqnarray*}

(ii) If $\pi =\left\{ \left\{ 1,2,3\right\} ,\left\{ 4,5\right\} \right\} $
and $\sigma =\left\{ \left\{ 1,2\right\} ,\left\{ 3,4,5\right\} \right\} $,
then $\pi \ $and $\sigma $ are not ordered and
\begin{equation*}
\sigma \wedge \pi =\left\{ \left\{ 1,2\right\} ,\left\{ 3\right\} ,\left\{
4,5\right\} \right\} \text{ \ \ and \ \ }\sigma \vee \pi =\left\{ b\right\}
=\left\{ \left\{ 1,2,3,4,5\right\} \right\} \text{.}
\end{equation*}%
A graphical representation of $\pi $, $\sigma $, $\sigma \wedge \pi $ and $%
\sigma \vee \pi $ is:%
\begin{eqnarray*}
\pi &=&\fbox{$1$ $2\ \ \ 3$}\fbox{$4$ $5$} \\
\sigma &=&\fbox{$1$ $2$ }\fbox{$3$ $\ 4$ $5$} \\
\sigma \wedge \pi &=&\fbox{$1$ $2$ }\fbox{$3$}\fbox{$4$ $5$} \\
\sigma \vee \pi &=&\fbox{$1\,\,\,2\,\,\,3\,\,\,4\,\,\,5\,\!\,\,$}
\end{eqnarray*}

(iii) A convenient way to build $\sigma \vee \pi $ is to do it in successive
steps. Take the union of two blocks with a common element and let it be a
new block of $\pi $. See if it shares an element with another block of $%
\sigma $. If yes, repeat. For instance, suppose that $\pi
=\{\{1,2\},\{3\},\{4\}\}$ and $\sigma =\{\{1,3\},\{2,4\}\}$. Then, $\pi \ $%
and $\sigma $ are not ordered and
\begin{equation*}
\sigma \wedge \pi =\left\{ \left\{ 1\right\} ,\left\{ 2\right\} ,\left\{
3\right\} ,\left\{ 4\right\} \right\} \text{ \ \ and \ \ }\sigma \vee \pi
=\left\{ \left\{ 1,2,3,4\right\} \right\} \text{.}
\end{equation*}%
One now obtains $\sigma \vee \pi $ by noting that the element 2 is common to
$\left\{ 1,2\right\} \in \pi $ and $\left\{ 2,4\right\} \in \sigma $, and
the \textquotedblleft merged\textquotedblright\ block $\left\{ 1,2,4\right\}
$ shares the element 1 with the block $\left\{ 1,3\right\} \in \sigma $,
thus implying the conclusion. A graphical representation of $\pi $, $\sigma $%
, $\sigma \wedge \pi $ and $\sigma \vee \pi $ is:%
\begin{eqnarray*}
\pi &=&\fbox{$1$ $2$}\fbox{$3$}\fbox{ $\ 4$} \\
\sigma &=&\fbox{$1$ $3$}\fbox{$2$ \ \ \ $4$} \\
\sigma \wedge \pi &=&\fbox{$1$}\fbox{$2$}\fbox{$3$}\fbox{ $4$} \\
\sigma \vee \pi &=&\fbox{$1\,\,\,2\,\,\,3\,\ \ \,\,4$}
\end{eqnarray*}

\bigskip

\textbf{Remark on notation. }When displaying a partition $\pi $ of $\left\{
1,...,n\right\} $ ($n\geq 1$), the blocks $b_{1},...,b_{k}\in \pi $ \ will
always be listed in the following way: $b_{1}$ will always contain the
element $1$, and%
\begin{equation*}
\min \left\{ i:i\in b_{j}\right\} <\min \left\{ i:i\in b_{j+1}\right\} \text{%
, \ }j=1,...,k-1\text{.}
\end{equation*}%
Also, the elements within each block will be always listed in increasing
order. For instance, if $n=6$ and the partition $\pi $ involves the blocks $%
\left\{ 2\right\} ,\left\{ 4\right\} ,\left\{ 1,6\right\} $ and $\left\{
3,5\right\} $, we will write $\pi =\left\{ \left\{ 1,6\right\} ,\left\{
2\right\} ,\left\{ 3,5\right\} ,\left\{ 4\right\} \right\} $.\

\bigskip

The \textsl{maximal element }of\textit{\ }$\mathcal{P}\left( b\right) $ is
the trivial partition $\hat{1}=\left\{ b\right\} $. The \textsl{minimal
element}\textit{\ }of $\mathcal{P}\left( b\right) $ is the partition $\hat{0}
$, such that each block of $\hat{0}$ contains exactly one element of $b$.
Observe that $\left\vert \hat{1}\right\vert =1$ and $\left\vert \hat{0}%
\right\vert =\left\vert b\right\vert $, and also $\hat{0}\leq \hat{1}$. If $%
\sigma \leq \pi $, we write $\left[ \sigma ,\pi \right] $ to indicate the
\textsl{segment}\textit{\ }$\left\{ \rho \in \mathcal{P}\left( b\right)
:\sigma \leq \rho \leq \pi \right\} $, which is a subset of partitions of $b$%
. Plainly, $\mathcal{P}\left( b\right) =\left[ \hat{0},\hat{1}\right] $.

\subsection{Relations between partitions of a set and partitions of an
integer \label{SS : Class}}

We now focus on the notion of \textsl{class}, which associates with a
segment of partitions a partition of an integer. In particular, the \textsl{%
class}\textit{\ }of a segment $\left[ \sigma ,\pi \right] $ ($\sigma \leq
\pi $), denoted $\lambda \left( \sigma ,\pi \right) $, is defined as the
partition of the integer $\left\vert \sigma \right\vert $ given by
\begin{equation}
\lambda \left( \sigma ,\pi \right) =\left( 1^{r_{1}}2^{r_{2}}\cdot \cdot
\cdot \left\vert \sigma \right\vert ^{r_{\left\vert \sigma \right\vert
}}\right) \text{,}  \label{tool}
\end{equation}%
where $r_{i}$ indicates the number of blocks of $\pi $ that contain exactly $%
i$ blocks of $\sigma $. We stress that necessarily $\left\vert \sigma
\right\vert \geq \left\vert \pi \right\vert $, and also%
\begin{equation*}
\left\vert \sigma \right\vert =1r_{1}+2r_{2}+\cdot \cdot \cdot +\left\vert
\sigma \right\vert r_{\left\vert \sigma \right\vert }\text{ \ \ and \ \ }%
\left\vert \pi \right\vert =r_{1}+\cdot \cdot \cdot +r_{\left\vert \sigma
\right\vert }\text{.}
\end{equation*}
For instance, if $\pi =\left\{ \left\{ 1,2,3\right\} ,\left\{ 4,5\right\}
\right\} $ and $\sigma =\left\{ \left\{ 1,2\right\} ,\left\{ 3\right\}
,\left\{ 4,5\right\} \right\} $, then since $\left\{ 1,2\right\} $ and $%
\left\{ 3\right\} $ are contained in $\left\{ 1,2,3\right\} $ and $\left\{
4,5\right\} $ in $\left\{ 4,5\right\} $, we have $r_{1}=1,$ $r_{2}=1$, $%
r_{3}=0$, that is, $\lambda \left( \sigma ,\pi \right) =\left(
1^{1}2^{1}3^{0}\right) $, corresponding to the partition of the integer $%
3=2+1$. In view of (\ref{PartEnt}), one may suppress the terms $r_{i}=0$ in (%
\ref{tool}), and write for instance $\lambda \left( \sigma ,\pi \right)
=\left( 1^{1}2^{0}3^{2}\right) =\left( 1^{1}3^{2}\right) $ for the class of
the segment $\left[ \sigma ,\pi \right] $, associated with the two
partitions $\sigma =\left\{ \left\{ 1\right\} ,\left\{ 2\right\} ,\left\{
3\right\} ,\left\{ 4\right\} ,\left\{ 5\right\} ,\left\{ 6\right\} ,\left\{
7\right\} \right\} $ and $\pi =\left\{ \left\{ 1\right\} ,\left\{
2,3,4\right\} ,\left\{ 5,6,7\right\} \right\} $.

\bigskip

Now fix a set $b$ such that $\left\vert b\right\vert =n\geq 1$. Then, for a
fixed $\lambda =\left( 1^{r_{1}}2^{r_{2}}\cdot \cdot \cdot n^{r_{n}}\right)
\vdash n$, the number of partitions $\pi \in \mathcal{P}\left( b\right) $
such that $\lambda \left( \hat{0},\pi \right) =\lambda $ is given by
\begin{equation}
\Big[  \begin{array}{c}
                           n \\
                           \lambda
                         \end{array}\Big]
=\Big[  \begin{array}{c}
                           n \\
                           r_1,...,r_n
                         \end{array}\Big]=\frac{n!}{\left(
1!\right) ^{r_{1}}r_{1}!\left( 2!\right) ^{r_{2}}r_{2}!\cdot \cdot \cdot
\left( n!\right) ^{r_{n}}r_{n}!}  \label{Npart}
\end{equation}%
(see e.g. \cite{Stanley}). The requirement that $\lambda \left( \hat{0},\pi
\right) =\lambda =\left( 1^{r_{1}}2^{r_{2}}\cdot \cdot \cdot
n^{r_{n}}\right) $ simply means that, for each $i=1,...,n$, the partition $%
\pi $ must have exactly $r_{i}$ blocks containing $i$ elements of $b$. Recall that $r_1,...,r_n$
must satisfy (\ref{e : sumlabel}).

\medskip From now on, we let
\begin{equation}
\fbox{$\left[ n\right] =\left\{ 1,...,n\right\} ,$ \ \ $n\geq 1$.}
\label{e : bracket-n}
\end{equation}%
With this notation, the maximal and minimal element of the set $\mathcal{P}%
\left( \left[ n\right] \right) $ are given, respectively, by
\begin{equation}
\fbox{$\hat{1}=\left\{ \left[ n\right] \right\} =\left\{ \left\{
1,...,n\right\} \right\} $ \ \ and \ \ $\hat{0}=\left\{ \left\{ 1\right\}
,...,\left\{ n\right\} \right\} .$}  \label{e : 1 hat 0}
\end{equation}

\bigskip

\textbf{Examples. }(i) For any finite set $b$, one has always that%
\begin{equation*}
\lambda \left( \hat{0},\hat{1}\right) =\left( 1^{0}2^{0}\cdot \cdot \cdot
\left\vert b\right\vert ^{1}\right) ,
\end{equation*}%
because $\hat{1}$ has only one block, namely $b$, and that block contains $%
\left\vert b\right\vert $ blocks of $\hat{0}$.

(ii) Fix $k\geq 1$ and let $b$ be such that $\left\vert b\right\vert =n\geq
k+1$. Consider $\lambda =\left( 1^{r_{1}}2^{r_{2}}\cdot \cdot \cdot
n^{r_{n}}\right) \vdash n$ be such that $r_{k}=r_{n-k}=1$ and $r_{j}=0$ for
every $j\neq k,n-k$. For instance, if $n=5$ and $k=2$, then $\lambda =\left(
1^{0}2^{1}3^{1}4^{0}5^{0}\right) $. Then, each partition $\pi \in \mathcal{P}%
\left( b\right) $ such that $\lambda \left( \hat{0},\pi \right) =\lambda $
has only one block of $k$ elements and one block of $n-k$ elements. To
construct such a partition, it is sufficient to specify the block of $k$
elements. This implies that there exists a bijection between the set of
partitions $\pi \in \mathcal{P}\left( b\right) $ such that $\lambda \left(
\hat{0},\pi \right) =\lambda $ and the collection of the subsets of $b$
having exactly $k$ elements. In particular, (\ref{Npart}) gives
\begin{equation*}
\Big[  \begin{array}{c}
                         \small{n} \\
                           \small{\lambda}
                         \end{array}\Big]=\binom{n}{k}=n!/\left( k!\left( n-k\right) !\right).
\end{equation*}

(iii) Let $b=\left[ 7\right] =\left\{ 1,...,7\right\} $\textbf{\ }and $%
\lambda =\left( 1^{1}2^{3}3^{0}4^{0}5^{0}6^{0}7^{0}\right) $. Then, (\ref%
{Npart}) implies that there are exactly $\frac{7!}{3!\left( 2!\right) ^{3}}%
=105$ partitions $\pi \in \mathcal{P}\left( b\right) $, such that $\lambda
\left( \hat{0},\pi \right) =\lambda $. One of these partitions is $\left\{
\left\{ 1\right\} ,\left\{ 2,3\right\} ,\left\{ 4,5\right\} ,\left\{
6,7\right\} \right\} $. Another is $\left\{ \left\{ 1,7\right\} ,\left\{
2\right\} ,\left\{ 3,4\right\} ,\left\{ 5,6\right\} \right\} .$

(iv) Let $b=\left[ 5\right] =\left\{ 1,...,5\right\} $, $\sigma =\left\{
\left\{ 1\right\} ,\left\{ 2\right\} ,\left\{ 3\right\} ,\left\{ 4,5\right\}
\right\} $ and $\pi =\left\{
\left\{ 1,2,3\right\} ,\left\{ 4,5\right\} \right\} $. Then, $\sigma \leq \pi $ and the set of partitions defined by the
interval $\left[ \sigma ,\pi \right] $ is $\left\{ \sigma ,\pi ,\rho
_{1},\rho _{2},\rho _{3}\right\} $, where%
\begin{eqnarray*}
\rho _{1} &=&\left\{ \left\{ 1,2\right\} ,\left\{ 3\right\} ,\left\{
4,5\right\} \right\}  \\
\rho _{2} &=&\left\{ \left\{ 1,3\right\} ,\left\{ 2\right\} ,\left\{
4,5\right\} \right\}  \\
\rho _{3} &=&\left\{ \left\{ 1\right\} ,\left\{ 2,3\right\} ,\left\{
4,5\right\} \right\} .
\end{eqnarray*}%
The partitions $\rho _{1}$, $\rho _{2}$ and $\rho _{3}$ are not ordered
(i.e., for every $1\leq i\neq j\leq 3$, one cannot write $\rho _{i}\leq \rho
_{j}$), and are built by taking unions of blocks of $\sigma $ in such a way
that they are contained in blocks of $\pi $. Moreover, $\lambda \left(
\sigma ,\pi \right) =\left( 1^{1}2^{0}3^{1}4^{0}5^{0}\right) $, since there
is exactly one block of $\pi $ containing one block of $\sigma $, and one
block of $\pi $ containing three blocks of $\sigma $.

(v) This example is related to the techniques developed in Section 6.1. Fix $%
n\geq 2$, as well as a partition $\gamma =\left( \gamma _{1},...,\gamma
_{k}\right) \in \Lambda \left( n\right) $ such that $\gamma _{k}\geq 2$.
Recall that, by definition, one has that $\gamma _{1}\geq \gamma _{2}\geq
\cdot \cdot \cdot \geq \gamma _{k}$ and $\gamma _{1}+\cdot \cdot \cdot
+\gamma _{k}=n$. Now consider the segment $\left[ \hat{0},\pi \right] $,
where
\begin{eqnarray*}
\hat{0} &=&\left\{ \left\{ 1\right\} ,\left\{ 2\right\} ,...,\left\{
n\right\} \right\} \text{, \ and} \\
\pi &=&\left\{ \left\{ 1,...,\gamma _{1}\right\} ,\left\{ \gamma
_{1}+1,...,\gamma _{1}+\gamma _{2}\right\} ,...,\left\{ \gamma _{1}+\cdot
\cdot \cdot +\gamma _{k-1}+1,...,n\right\} \right\} .
\end{eqnarray*}%
Then, the $j$th block of $\pi $ contains exactly $\gamma _{j}$ blocks of $%
\hat{0}$, $j=1,...,k$, thus giving that the class $\lambda \left( \hat{0}%
,\pi \right) $ is such that $\lambda \left( \hat{0},\pi \right) =\left(
\gamma _{1}^{1}\gamma _{2}^{1}\cdot \cdot \cdot \gamma _{k}^{1}\right)
=\gamma $, after suppressing the indicators of the type $r^{0}$.

\subsection{M\"{o}bius functions and M\"{o}bius inversion formulae}

For $\sigma ,\pi \in \mathcal{P}\left( b\right) $, we denote by $\mu \left(
\sigma ,\pi \right) $ the \textsl{M\"{o}bius function}\textit{\ }associated
with the lattice $\mathcal{P}\left( b\right) $. It is defined as follows. If
$\sigma \nleq \pi $ (that is, if the relation $\sigma \leq \pi $ does not
hold), then $\mu \left( \sigma ,\pi \right) =0$. If $\sigma \leq \pi $, then
the quantity $\mu \left( \sigma ,\pi \right) $ depends only on the class $%
\lambda \left( \sigma ,\pi \right) $ of the segment $\left[ \sigma ,\pi %
\right] $, and is given by (see \cite{Aigner})%
\begin{eqnarray}
\mu \left( \sigma ,\pi \right) &=&\left( -1\right) ^{n-r}\left( 2!\right)
^{r_{3}}\left( 3!\right) ^{r_{4}}\cdot \cdot \cdot \left( \left( n-1\right)
!\right) ^{r_{n}}  \label{MobF} \\
&=&\left( -1\right) ^{n-r}\prod_{j=0}^{n-1}\left( j!\right) ^{r_{j+1}},
\label{MobF2}
\end{eqnarray}%
where $n=\left\vert \sigma \right\vert $, $r=\left\vert \pi \right\vert $,
and $\lambda \left( \sigma ,\pi \right) =\left( 1^{r_{1}}2^{r_{2}}\cdot
\cdot \cdot n^{r_{n}}\right) $ (that is, there are exactly $r_{i}$ blocks of
$\pi $ containing exactly $i$ blocks of $\sigma $). Since $0!=1!=1$,
expressions (\ref{MobF}) and (\ref{MobF2}) do not depend on the specific
values of $r_{1}$ (the number of blocks of $\pi $ containing exactly $1$
block of $\sigma $) and $r_{2}$ (the number of blocks of $\pi $ containing
exactly two blocks of $\sigma $).

\bigskip

\textbf{Examples. }(i) If $\left\vert b\right\vert =n\geq 1$ and $\sigma \in
\mathcal{P}\left( b\right) $ is such that $\left\vert \sigma \right\vert =k$ ( $%
\leq n$ ), then
\begin{equation}
\mu \left( \sigma ,\hat{1}\right) =\left( -1\right) ^{k-1}\left( k-1\right) !.
\label{mp1}
\end{equation}%
Indeed, in (\ref{MobF}) $r_{k}=1$, since $\hat{1}$ has a single block which
contains the $k$ blocks of $\sigma $. In particular, $\mu \left( \hat{0}%
,\left\{ b\right\} \right) $ $=\mu \left( \hat{0},\hat{1}\right) =\left(
-1\right) ^{n-1}$ $\left( n-1\right) !$.

(ii) For every $\pi \in \mathcal{P}\left( b\right) $, one has $\mu \left(
\pi ,\pi \right) =1$. Indeed, since (trivially) each element of $\pi $
contains exactly one element of $\pi $, one has $\lambda \left( \pi ,\pi
\right) =\left( 1^{\left\vert \pi \right\vert }2^{0}3^{0}\cdot \cdot \cdot
n^{0}\right) $.

\bigskip

The next result is crucial for the obtention of the combinatorial formulae
found in Section \ref{S : cum} and Section \ref{S : MCA} below. For every
pair of functions $G,F$, from $\mathcal{P}\left( b\right) $ into $\mathbb{C}$
and such that $\forall \sigma \in \mathcal{P}\left( b\right) $,%
\begin{equation}
G\left( \sigma \right) =\sum_{\hat{0}\leq \pi \leq \sigma }F\left( \pi
\right) \text{ \ \ (resp. \ }G\left( \sigma \right) =\sum_{\sigma \leq \pi
\leq \hat{1}}F\left( \pi \right) \text{)\ }  \label{MobPreInv}
\end{equation}%
one has the following \textsl{M\"{o}bius inversion formula}: $\forall \pi
\in \mathcal{P}\left( b\right) $,
\begin{equation}
F\left( \pi \right) =\sum_{\hat{0}\leq \sigma \leq \pi }\mu \left( \sigma
,\pi \right) G\left( \sigma \right) \text{ \ \ (resp. \ }F\left( \pi \right)
=\sum_{\pi \leq \sigma \leq \hat{1}}\mu \left( \pi ,\sigma \right) G\left(
\sigma \right) \text{),}  \label{MobInv}
\end{equation}%
where $\mu \left( \cdot ,\cdot \right) $ is the M\"{o}bius function given in
(\ref{MobF}). For a proof of (\ref{MobInv}), see e.g. \cite[Section 3.7]%
{Stanley} and \cite{Aigner}. To understand (\ref{MobInv}) as inversion formulae,
one can interpret the sum $\sum_{\hat{0}\leq \pi \leq \sigma }F\left( \pi
\right) $ as an integral of the type $\int_{\hat{0}}^{\sigma }F\left( \pi
\right) d\pi $ (and analogously for the other sums appearing in (\ref%
{MobPreInv}) and (\ref{MobInv})).

\bigskip

In general (see \cite[Section 3.7]{Stanley}), the M\"{o}bius function is
defined by recursion on any finite partially ordered set by the following
relations:%
\begin{equation}
\begin{array}{ll}
\mu \left( x,x\right) =1 & \forall x\in P, \\
\mu \left( x,y\right) =-\sum_{x\preceq z\prec y}\mu \left( x,z\right) \text{%
, } & \forall \text{ }x,y\in P:x\prec y, \\
\mu \left( x,y\right) =0 & \forall \text{ }x,y\in P:x\npreceq y,%
\end{array}
\label{mg}
\end{equation}%
where $P$ is a finite partially ordered set, with partial order $\preceq $,
and we write $x\prec y$ to indicate that $x\preceq y$ and $x\neq y$. For
instance, $P$ could be a set of subsets, with $\preceq $ equal to the
inclusion relation $\subseteq $. In our context $P=\mathcal{P}\left(
b\right) $, the set of partitions of $b$, and $\preceq $ is the partial
order $\leq $ considered above, so that (\ref{mg}) becomes
\begin{equation}
\begin{array}{ll}
\mu \left( \sigma,\sigma \right) =1 & \forall \sigma \in \mathcal{P}\left(
b\right) , \\
\mu \left( \sigma,\pi \right) =-\sum_{\sigma \leq \rho < \pi}\mu \left( \sigma,\rho \right) \text{%
, } & \forall \text{ }\sigma,\pi \in \mathcal{P}\left(
b\right) :\sigma < \pi, \\
\mu \left( \sigma, \pi \right) =0 & \forall \text{ }\sigma,\pi\in \mathcal{P}\left(
b\right) :\sigma \nleq \pi,%
\end{array}
\label{mgPART}
\end{equation}%
where we write $\sigma < \pi$ to indicate that $\sigma \leq \pi$ and $\sigma\neq \pi$ (and similarly for $\rho<\pi$). The recursion
formula (\ref{mg}) has the following consequence: for each $%
x\preceq y$,%
\begin{equation}
\sum_{x\preceq z\preceq y}\mu \left( z,y\right) =\sum_{x\preceq z\preceq
y}\mu \left( x,z\right) =\left\{
\begin{array}{ll}
0 & \text{if }x\neq y\text{,} \\
\mu \left( x,x\right) \text{(}=1\text{)} & \text{if }x=y\text{,}%
\end{array}%
\right.  \label{maliMob}
\end{equation}%
which will be used in the sequel. The second equality in (\ref{maliMob}) is
an immediate consequence of (\ref{mg}). To prove the first equality in (\ref%
{maliMob}), fix $x$ and write $G\left( z\right) =\mathbf{1}_{x\preceq z}$.
Since, trivially, $G\left( z\right) =\sum_{y\preceq z}\mathbf{1}_{y=x}$, one
can let $F\left( y\right) =\mathbf{1}_{y=x}$ in (\ref{MobPreInv}) and use
the inversion formula (\ref{MobInv}) to deduce that
\begin{equation*}
\mathbf{1}_{y=x}=\sum_{z\preceq y}\mu \left( z,y\right) G\left( z\right)
=\sum_{x\preceq z\preceq y}\mu \left( z,y\right) \text{,}
\end{equation*}%
which is equivalent to (\ref{maliMob}).

\bigskip

Now consider two finite partially ordered sets $P,Q$, whose order relations
are noted, respecively, $\preceq _{P}$ and $\preceq _{Q}$. The \textsl{%
lattice product }of $P$ and $Q$ is defined as the cartesian product $P\times
Q$, endowed with the following partial order relation: $\left( x,y\right) $ $%
\preceq _{P\times Q}$ $\left( x^{\prime },y^{\prime }\right) $ if, and only
if, $x\preceq _{P}x^{\prime }$ and $y\preceq _{Q}y^{\prime }$. Lattice
products of more than two partially ordered sets are defined analogously. We
say (see e.g. \cite[p. 98]{Stanley}) that $P$ and $Q$ are \textsl{isomorphic
}if there exists a bijection $\psi :P\rightarrow Q$ which is order-preserving
and such that the inverse of $\psi $ is also order-preserving; this
requirement on the bijection $\psi $ is equivalent to saying that, for every
$x,x^{\prime }\in P$,
\begin{equation}
x\preceq _{P}x^{\prime }\text{ if and only if }\psi \left( x\right)
\preceq _{Q}\psi \left( x^{\prime }\right) \text{.}  \label{iso}
\end{equation}%
Of course, two isomorphic partially ordered sets have the same cardinality.
The following result is quite useful for explicitly computing M\"{o}bius
functions. It states that the M\"{o}bius function is invariant under isomorphisms, and that the M\"{o}bius function of a lattice product
is the product of the associated M\"{o}bius functions. Point 1 is an immediate consequence of (\ref{mg}), for a proof of Point 2, see e.g. \cite[Section 3.8]{Stanley}.

\bigskip

\begin{proposition}
\label{P : latticeP}Let $P,Q$ be two partially ordered sets, and let $\mu
_{P}$ and $\mu _{Q}$ denote their M\"{o}bius functions. Then,

\begin{enumerate}
\item If $P$ and $Q$ are isomorphic, then $\mu _{P}\left( x,y\right) =\mu
_{Q}\left( \psi \left( x\right) ,\psi \left( y\right) \right) $ for every $%
x,y\in P$, where $\psi $ is the bijection appearing in (\ref{iso}).

\item The M\"{o}bius function associated with the partially ordered set $%
P\times Q$ is given by:%
\begin{equation*}
\mu _{P\times Q}\left[ \left( x,y\right) ,\left( x^{\prime },y^{\prime
}\right) \right] =\mu _{P}\left( x,x^{\prime }\right) \times \mu _{Q}\left(
y,y^{\prime }\right) \text{.}
\end{equation*}
\end{enumerate}
\end{proposition}

\bigskip

The next result is used in the proof of Theorem \ref{T : ProdRW}.

\bigskip

\begin{proposition}
\label{P : ISO}Let $b$ be a finite set, and let $\pi ,\sigma \in \mathcal{P}%
\left( b\right) $ be such that: (i) $\sigma \leq \pi $, and (ii) the segment
$\left[ \sigma ,\pi \right] $ has class $\left( \lambda _{1},...,\lambda
_{k}\right) \vdash |\sigma |$. Then, $\left[ \sigma ,\pi \right] $ is a
partially ordered set isomorphic to the lattice product of the $k$ sets $%
\mathcal{P}\left( \left[ \lambda _{i}\right] \right) $, $i=1,...,k$.
\end{proposition}

\begin{proof}
To prove the statement, we shall use the fact that each partition in $\left[
\sigma ,\pi \right] $ is obtained by taking unions of the blocks of $\sigma $
that are contained in the same block of $\pi $. Start by observing that $%
\left( \lambda _{1},...,\lambda _{k}\right) $ is the class of $\left[ \sigma
,\pi \right] $ if and only if for every $i=1,...,k$, there is a block $%
b_{i}\in \pi $ such that $b_{i}$ contains exactly $\lambda _{i}$ blocks of $%
\sigma $. In particular, $k=\left\vert \pi \right\vert $. We now construct a
bijection $\psi $, between $\left[ \sigma ,\pi \right] $ and the lattice
products of the $\mathcal{P}\left( \left[ \lambda _{i}\right] \right) $'s,
as follows.

\begin{description}
\item[i)] For $i=1,...,k$, write $b_{i,j}$, $j=1,...,\lambda _{i}$, to
indicate the blocks of $\sigma $ contained in $b_{i}$.

\item[ii)] For every partition $\rho \in \left[ \sigma ,\pi \right] $ and
every $i=1,...,k$, construct a partition $\zeta \left( i,\rho \right) $ of $%
\left[ \lambda _{i}\right] =\left\{ 1,...,\lambda _{i}\right\} $ by the
following rule: for every $j,l\in \left\{ 1,...,\lambda _{i}\right\} $, $%
j\sim _{\zeta \left( i,\rho \right) }l$ (that is, $j$ and $l$ belong to the
same block of $\zeta \left( i,\rho \right) $) if and only if the union $%
b_{i,j}\cup b_{i,l}$ is contained in a block of $\rho $.

\item[iii)] Define the application $\psi :\left[ \sigma ,\pi \right]
\rightarrow \mathcal{P}\left( \left[ \lambda _{1}\right] \right) \times
\cdot \cdot \cdot \times \mathcal{P}\left( \left[ \lambda _{k}\right]
\right) $ as%
\begin{equation}
\rho \mapsto \psi \left( \rho \right) :=\left( \zeta \left( 1,\rho \right)
,...,\zeta \left( k,\rho \right) \right) .  \label{zed}
\end{equation}
\end{description}

\noindent It is easily seen that the application $\psi $ in (\ref{zed}) is
indeed an order-preserving bijection, verifying (\ref{iso}) for $P=\left[
\sigma ,\pi \right] $ and $Q=\mathcal{P}\left( \left[ \lambda _{1}\right]
\right) \times \cdot \cdot \cdot \times \mathcal{P}\left( \left[ \lambda _{k}%
\right] \right) $.
\end{proof}

\section{Combinatorial expressions of cumulants and moments\label{S : cum}}

\setcounter{equation}{0}We recall here the definition of \textsl{cumulant},
and we present several of its properties. A detailed discussion of cumulants
is contained in the book by Shiryayev \cite{Shir}; see also the papers by
Rota and Shen \cite{Ro SHen}, Speed \cite{Speed} and Surgailis \cite{Sur}.
An analysis of cumulants involving different combinatorial structures can be
found in \cite[pp. 20-23]{Pit} and the references therein.

\subsection{Cumulants \label{SS : DefCum}}

For $n\geq 1$, we consider a vector of real-valued random variables $\mathbf{%
X}_{\left[ n\right] }=\left( X_{1},...,X_{n}\right) $ such that $\mathbb{E}%
\left\vert X_{j}\right\vert ^{n}<\infty $, $\forall j=1,...,n$. For every
subset $b=\left\{ j_{1},...,j_{k}\right\} \subseteq \left[ n\right] =\left\{
1,...,n\right\} $, we write%
\begin{equation}
\fbox{$\mathbf{X}_{b}=\left( X_{j_{1}},...,X_{j_{k}}\right) $}\text{ \ \ and
\ \ }\fbox{$\mathbf{X}^{b}=X_{j_{1}}\times \cdot \cdot \cdot \times
X_{j_{k}} $,}  \label{prod def}
\end{equation}%
where $\times $ denotes the usual product. For instance, $\forall m\leq n$,
\begin{equation*}
\mathbf{X}_{\left[ m\right] }=\left( X_{1},..,X_{m}\right) \text{ \ \ and \
\ }\mathbf{X}^{\left[ m\right] }=X_{1}\times \cdot \cdot \cdot \times X_{m}.
\end{equation*}%
For every $b=\left\{ j_{1},...,j_{k}\right\} \subseteq \left[ n\right] $ and
$\left( z_{1},...,z_{k}\right) \in \mathbb{R}^{k}$, we let $g_{\mathbf{X}%
_{b}}\left( z_{1},..,z_{k}\right) =\mathbb{E}\left[ \exp \left( i\sum_{\ell
=1}^{k}z_{\ell }X_{j_{\ell }}\right) \right] $. The \textsl{joint cumulant}%
\textit{\ }of the components of the vector $\mathbf{X}_{b}$ is defined as%
\begin{equation}
\chi \left( \mathbf{X}_{b}\right) =\left( -i\right) ^{k}\frac{\partial ^{k}}{%
\partial z_{1}\cdot \cdot \cdot \partial z_{k}}\log g_{\mathbf{X}_{b}}\left(
z_{1},..,z_{k}\right) \mid _{z_{1}=...=z_{k}=0}\text{,}  \label{cumDef}
\end{equation}%
thus%
\begin{equation*}
\chi \left( X_{1},...,X_{k}\right) =\left( -i\right) ^{k}\frac{\partial ^{k}%
}{\partial z_{1}\cdot \cdot \cdot \partial z_{k}}\log \mathbb{E}\left. \left[
\exp \left( i\sum_{\ell =1}^{k}z_{l}X_{l}\right) \right] \right\vert
_{z_{1}=...=z_{k}=0}.
\end{equation*}%
We recall the following facts.

\begin{description}
\item[(i)] The application $\mathbf{X}_{b}\mapsto \chi \left( \mathbf{X}%
_{b}\right) $ is \textsl{homogeneous}, that is, for every $\mathbf{h}=\left(
h_{1},...,h_{k}\right) \in \mathbb{R}^{k}$,
\begin{equation*}
\chi \left( h_{1}X_{j_{1}},...,h_{k}X_{j_{k}}\right) =(\Pi _{\ell
=1}^{k}h_{\ell })\times \chi \left( \mathbf{X}_{b}\right) \text{;}
\end{equation*}

\item[(ii)] The application $\mathbf{X}_{b}\mapsto \chi \left( \mathbf{X}%
_{b}\right) $ is invariant with respect to the permutations of $b$;

\item[(iii)] $\chi \left( \mathbf{X}_{b}\right) =0$, if the vector $\mathbf{X%
}_{b}$ has the form $\mathbf{X}_{b}=\mathbf{X}_{b^{\prime }}\cup \mathbf{X}%
_{b^{\prime \prime }}$, with $b^{\prime },b^{\prime \prime }\neq \varnothing
$, $b^{\prime }\cap b^{\prime \prime }=\varnothing $ and $\mathbf{X}%
_{b^{\prime }}\ $and $\mathbf{X}_{b^{\prime \prime }}$ independent;

\item[(iv)] if $\mathbf{Y}=\left\{ Y_{j}:j\in J\right\} $ is a Gaussian
family and if $\mathbf{X}_{\left[ n\right] }$ is a vector obtained by
juxtaposing $n\geq 3$ elements of $\mathbf{Y}$ (with possible repetitions),
then $\chi \left( \mathbf{X}_{\left[ n\right] }\right) =0$.
\end{description}

\bigskip

Properties (i) and (ii) follow immediately from (\ref{cumDef}). To see how
to deduce (iii) from (\ref{cumDef}), just observe that, if $\mathbf{X}_{b}$
has the structure described in (iii), then%
\begin{equation*}
\log g_{\mathbf{X}_{b}}\left( z_{1},..,z_{k}\right) =\log g_{\mathbf{X}%
_{b^{\prime }}}\left( z_{\ell }:j_{\ell }\in b^{\prime }\right) +\log g_{%
\mathbf{X}_{b^{\prime \prime }}}\left( z_{\ell }:j_{\ell }\in b^{\prime
\prime }\right)
\end{equation*}%
(by independence), so that
\begin{eqnarray*}
&&\frac{\partial ^{k}}{\partial z_{1}\cdot \cdot \cdot \partial z_{k}}\log
g_{\mathbf{X}_{b}}\left( z_{1},..,z_{k}\right) \\
&=&\frac{\partial ^{k}}{\partial z_{1}\cdot \cdot \cdot \partial z_{k}}\log
g_{\mathbf{X}_{b^{\prime }}}\left( z_{\ell }:j_{\ell }\in b^{\prime }\right)
+\frac{\partial ^{k}}{\partial z_{1}\cdot \cdot \cdot \partial z_{k}}\log g_{%
\mathbf{X}_{b^{\prime \prime }}}\left( z_{\ell }:j_{\ell }\in b^{\prime
\prime }\right) =0.
\end{eqnarray*}%
Finally, property (iv) is proved by using the fact that, if $\mathbf{X}_{%
\left[ n\right] }$ is obtained by juxtaposing $n\geq 3$ elements of a
Gaussian family (even with repetitions), then $\log g_{\mathbf{X}_{b}}\left(
z_{1},..,z_{k}\right) $ has necessarily the form $\sum_{l}a\left( l\right)
z_{l}+\sum_{i,j}b\left( i,j\right) z_{i}z_{j}$, where $a\left( k\right) $
and $b\left( i,j\right) $ are coefficients not depending on the $z_{l}$'s$.$
All the derivatives of order higher than 2 are then zero.

When $\left\vert b\right\vert =n$, one says that the cumulant $\chi \left(
\mathbf{X}_{b}\right) $, given by (\ref{cumDef}), \textsl{has order}\textit{%
\ }$n$. When $\mathbf{X}_{\left[ n\right] }=\left( X_{1},...,X_{n}\right) $
is such that $X_{j}=X$, $\forall j=1,...,n$, where $X$ is a random variable
in $L^{n}\left( \mathbb{P}\right) $, we write
\begin{equation}
\fbox{$\chi \left( \mathbf{X}_{\left[ n\right] }\right) =\chi _{n}\left(
X\right) $}  \label{e : cum-single}
\end{equation}%
and we say that $\chi _{n}\left( X\right) $ is the $n$th\textit{\ \textsl{%
cumulant} }(or the \textsl{cumulant of order}\textit{\ }$n$) of $X$. Of
course, in this case one has that
\begin{equation*}
\chi _{n}\left( X\right) =(-i)^n\left. \frac{\partial ^{n}}{\partial z^{n}}\log
g_{X}\left( z\right) \right\vert _{z=0}\text{,}
\end{equation*}%
where $g_{X}\left( z\right) =\mathbb{E}\left[ \exp \left( izX\right) \right]
$. Note that, if $X,Y\in L^{n}\left( \mathbb{P}\right) $ ($n\geq 1$) are
independent random variables, then (\ref{cumDef}) implies that
\begin{equation*}
\chi _{n}\left( X+Y\right) =\chi _{n}\left( X\right) +\chi _{n}\left(
Y\right) ,
\end{equation*}%
since $\chi _{n}\left( X+Y\right) $ involve the derivative of $\mathbb{E}%
\left[ \exp \left( i\left( X+Y\right) \Sigma _{j=1}^{n}z_{j}\right) \right] $
with respect to $z_{1},...,z_{n}$.

\subsection{Relations between moments and cumulants, and between cumulants}

We want to relate expectations of products of random variables, such as $%
\mathbb{E}\left[ X_{1}X_{2}X_{3}\right] $, to cumulants of vectors of random
variables, such as $\chi \left( X_{1},X_{2},X_{3}\right) $. Note the
disymmetry: moments involve products, while cumulants involve vectors. We
will have, for example, $\chi \left( X_{1},X_{2}\right) =\mathbb{E}\left[
X_{1}X_{2}\right] -\mathbb{E}\left[ X_{1}\right] \mathbb{E}\left[ X_{2}%
\right] $, and hence $\chi \left( X_{1},X_{2}\right) =\mathbf{Cov}\left(
X_{1},X_{2}\right) $, the covariance of the vector $\left(
X_{1},X_{2}\right) $. Conversely, we will have $\mathbb{E}\left[ X_{1}X_{2}%
\right] =\chi \left( X_{1}\right) \chi \left( X_{2}\right) +\chi \left(
X_{1},X_{2}\right) $. Thus, using the notation introduced above, we will
establish precise relations between objects of the type $\chi \left( \mathbf{%
X}_{b}\right) =\chi \left( X_{j}:j\in b\right) $ and $\mathbb{E}\left[
\mathbf{X}^{b}\right] =\mathbb{E}\left[ \Pi _{j\in b}X_{j}\right] $. We can
do this also for random variables that are products of other random
variables: for instance, to obtain $\chi \left( Y_{1}Y_{2},Y_{3}\right) $,
we apply the previous formula with $X_{1}=Y_{1}Y_{2}$ and $X_{2}=Y_{3}$, and
get $\chi \left( Y_{1}Y_{2},Y_{3}\right) =\mathbb{E}\left[ Y_{1}Y_{2}Y_{3}%
\right] $ $-$ $\mathbb{E}\left[ Y_{1}Y_{2}\right] \mathbb{E}\left[ Y_{3}%
\right] $. We shall also state a formula, due to Malyshev, which expresses $%
\chi \left( Y_{1}Y_{2},Y_{3}\right) $ in terms of other cumulants, namely in
this case
\begin{equation*}
\chi \left( Y_{1}Y_{2},Y_{3}\right) =\chi \left( Y_{1},Y_{3}\right) \chi
\left( Y_{2}\right) +\chi \left( Y_{1}\right) \chi \left( Y_{2},Y_{3}\right)
+\chi \left( Y_{1},Y_{2},Y_{3}\right) .
\end{equation*}%
The next result, first proved in \cite{LeoShy} (for Parts 1 and 2) and \cite%
{Maly} (for Part 3), contains three crucial relations, linking the cumulants
and the moments associated with a random vector $\mathbf{X}_{\left[ n\right]
}$. We use the properties of the lattices of partitions, as introduced in
the previous section.

\bigskip

\begin{proposition}[Leonov and Shiryayev \protect\cite{LeoShy} and Malyshev
\protect\cite{Maly}]
\label{P : LeoShy}For every $b\subseteq \left[ n\right] $,

\begin{enumerate}
\item
\begin{equation}
\mathbb{E}\left[ \mathbf{X}^{b}\right] =\sum_{\pi =\left\{
b_{1},...,b_{k}\right\} \in \mathcal{P}\left( b\right) }\chi \left( \mathbf{X%
}_{b_{1}}\right) \cdot \cdot \cdot \chi \left( \mathbf{X}_{b_{k}}\right) ;
\label{LS1}
\end{equation}

\item
\begin{equation}
\mathbb{\chi }\left( \mathbf{X}_{b}\right) =\sum_{\sigma =\left\{
a_{1},...,a_{r}\right\} \in \mathcal{P}\left( b\right) }\left( -1\right)
^{r-1}\left( r-1\right) !\mathbb{E}\left( \mathbf{X}^{a_{1}}\right) \cdot
\cdot \cdot \mathbb{E}\left( \mathbf{X}^{a_{r}}\right) ;  \label{LS2}
\end{equation}

\item $\forall \sigma =\left\{ b_{1},...,b_{k}\right\} \in \mathcal{P}\left(
b\right) $,
\begin{equation}
\chi \left( \mathbf{X}^{b_{1}},...,\mathbf{X}^{b_{k}}\right) =\sum
_{\substack{ \tau =\left\{ t_{1},...,t_{s}\right\} \in \mathcal{P}\left(
b\right) \text{ }  \\ \tau \vee \sigma =\hat{1}}}\chi \left( \mathbf{X}%
_{t_{1}}\right) \cdot \cdot \cdot \chi \left( \mathbf{X}_{t_{s}}\right)
\text{.}  \label{LS3}
\end{equation}
\end{enumerate}
\end{proposition}

\bigskip

\textbf{Remark. }To the best of our knowledge, our forthcoming proof of
equation (\ref{LS3}) (which is known as \textsl{Malyshev's formula}) is new.
As an illustration of (\ref{LS3}), consider the cumulant $\chi \left(
X_{1}X_{2},X_{3}\right) $, in which case one has $\sigma =\left\{
b_{1},b_{2}\right\} $, with $b_{1}=\left\{ 1,2\right\} $ and $b_{2}=\left\{
3\right\} $. There are three partitions $\tau \in \mathcal{P}\left( \left[ 3%
\right] \right) $ such that $\tau \vee \sigma =\hat{1}=\left\{ 1,2,3\right\}
$, namely $\tau _{1}=\hat{1}$, $\tau _{2}=\left\{ \left\{ 1,3\right\}
,\left\{ 2\right\} \right\} $ and $\tau _{3}=\left\{ \left\{ 1\right\}
,\left\{ 2,3\right\} \right\} $, from which it follows that $\chi \left(
X_{1}X_{2},X_{3}\right) =\chi \left( X_{1},X_{2},X_{3}\right) +\chi \left(
X_{1}X_{3}\right) \chi \left( X_{2}\right) +\chi \left( X_{1}\right) \chi
\left( X_{2}X_{3}\right) $.

\bigskip

\begin{proof}[Proof of Proposition \protect\ref{P : LeoShy}]
The proof of (\ref{LS1}) is obtained by differentiating the characteristic
function and its logarithm, and by identifying corresponding terms (see \cite%
{Shir}, \cite[Section 6]{Ro SHen} or \cite{Speed}). We now show how to
obtain (\ref{LS2}) and (\ref{LS3}) from (\ref{LS1}). Relation (\ref{LS1})
implies that, $\forall \sigma =\left\{ a_{1},...,a_{r}\right\} \in \mathcal{P%
}\left( b\right) $,
\begin{equation}
\prod_{j=1}^{r}\mathbb{E}\left[ \mathbf{X}^{a_{j}}\right] =\sum_{\substack{ %
\pi =\left\{ b_{1},...,b_{k}\right\} \leq \sigma  \\ \pi \in \mathcal{P}%
\left( b\right) }}\chi \left( \mathbf{X}_{b_{1}}\right) \cdot \cdot \cdot
\chi \left( \mathbf{X}_{b_{k}}\right) .  \label{F0}
\end{equation}%
We can therefore set $G\left( \sigma \right) =\prod_{j=1}^{r}\mathbb{E}\left[
\mathbf{X}^{a_{j}}\right] $ and $F\left( \pi \right) =\chi \left( \mathbf{X}%
_{b_{1}}\right) \cdot \cdot \cdot \chi \left( \mathbf{X}_{b_{k}}\right) $ in
(\ref{MobPreInv}) and (\ref{MobInv}), so as to deduce that, for every $\pi
=\left\{ b_{1},...,b_{k}\right\} \in \mathcal{P}\left( b\right) $,
\begin{equation}
\chi \left( \mathbf{X}_{b_{1}}\right) \cdot \cdot \cdot \chi \left( \mathbf{X%
}_{b_{k}}\right) =\sum_{\sigma =\left\{ a_{1},...,a_{r}\right\} \leq \pi
}\mu \left( \sigma ,\pi \right) \prod_{j=1}^{r}\mathbb{E}\left[ \mathbf{X}%
^{a_{j}}\right] .  \label{F}
\end{equation}%
Relation (\ref{LS2}) is therefore a particular case of (\ref{F}), obtained
by setting $\pi =\hat{1}$ and by using the equality $\mu \left( \sigma ,\hat{%
1}\right) =\left( -1\right) ^{\left\vert \sigma \right\vert -1}\left(
\left\vert \sigma \right\vert -1\right) !$, which is a consequence of (\ref%
{MobF}).

To deduce Malyshev's formula (\ref{LS3}) from (\ref{LS2}) and (\ref{F0}),
write $\mathbf{X}^{b_{j}}=Y_{j}$, $j=1,...,k$ (recall that the $\mathbf{X}%
^{b_{j}}$ are random variables defined in (\ref{prod def})), and apply (\ref%
{LS2}) to the vector $\mathbf{Y}=\left( Y_{1},...,Y_{k}\right) $ to obtain
that
\begin{eqnarray}
\chi \left( \mathbf{X}^{b_{1}},...,\mathbf{X}^{b_{k}}\right)  &=&\mathbb{%
\chi }\left( Y_{1},...,Y_{k}\right) =\chi \left( \mathbf{Y}\right)   \notag
\\
&=&\sum_{\beta =\left\{ p_{1},...,p_{r}\right\} \in \mathcal{P}\left( \left[
k\right] \right) }\left( -1\right) ^{r-1}\left( r-1\right) !\mathbb{E}\left(
\mathbf{Y}^{p_{1}}\right) \cdot \cdot \cdot \mathbb{E}\left( \mathbf{Y}%
^{p_{r}}\right) .  \label{e}
\end{eqnarray}%
Now write $\sigma =\left\{ b_{1},...,b_{k}\right\} $, and observe that $%
\sigma $ is a partition of the set $b$, while the partitions $\beta $ in (%
\ref{e}) are partitions of the first $k$ integers. Now fix $\beta \in
\left\{ p_{1},...,p_{r}\right\} \in \mathcal{P}\left( \left[ k\right]
\right) $. For $i=1,...r,$ take the union of the blocks $b_{j}\in \sigma $
having $j\in p_{i}$, and call this union $u_{i}$. One obtains therefore a
partition $\pi =\left\{ u_{1},...,u_{r}\right\} \in \mathcal{P}\left(
b\right) $ such that $\left\vert \pi \right\vert =\left\vert \beta
\right\vert =r$. Thanks to (\ref{MobF}) and (\ref{mp1}),%
\begin{equation}
\left( -1\right) ^{r-1}\left( r-1\right) !=\mu \left( \beta ,\hat{1}\right)
=\mu \left( \pi ,\hat{1}\right)   \label{q}
\end{equation}%
(note that the two M\"{o}bius functions appearing in (\ref{q}) are
associated with different lattices: indeed, $\mu \left( \beta ,\hat{1}%
\right) $ refers to $\mathcal{P}\left( \left[ k\right] \right) $, whereas $%
\mu \left( \pi ,\hat{1}\right) $ is associated with $\mathcal{P}\left(
b\right) $). With this notation, one has also that $\mathbb{E}\left( \mathbf{%
Y}^{p_{1}}\right) \cdot \cdot \cdot \mathbb{E}\left( \mathbf{Y}%
^{p_{r}}\right) =\mathbb{E}\left( \mathbf{X}^{u_{1}}\right) \cdot \cdot
\cdot \mathbb{E}\left( \mathbf{X}^{u_{r}}\right) $, so that, by (\ref{F0}),
\begin{equation}
\mathbb{E}\left( \mathbf{Y}^{p_{1}}\right) \cdot \cdot \cdot \mathbb{E}%
\left( \mathbf{Y}^{p_{r}}\right) =\mathbb{E}\left( \mathbf{X}^{u_{1}}\right)
\cdot \cdot \cdot \mathbb{E}\left( \mathbf{X}^{u_{r}}\right) =\sum
_{\substack{ \tau =\left\{ t_{1},...,t_{s}\right\} \leq \pi  \\ \tau \in
\mathcal{P}\left( b\right) }}\chi \left( \mathbf{X}_{t_{1}}\right) \cdot
\cdot \cdot \chi \left( \mathbf{X}_{t_{s}}\right) \text{.}  \label{ee}
\end{equation}%
By plugging (\ref{q}) and (\ref{ee}) into (\ref{e}) we obtain finally that%
\begin{eqnarray*}
\chi \left( \mathbf{X}^{b_{1}},...,\mathbf{X}^{b_{k}}\right)
&=&\sum_{\sigma \leq \pi \leq \hat{1}}\mu \left( \pi ,\hat{1}\right)
\sum_{\tau =\left\{ t_{1},...,t_{s}\right\} :\tau \leq \pi }\chi \left(
\mathbf{X}_{t_{1}}\right) \cdot \cdot \cdot \chi \left( \mathbf{X}%
_{t_{s}}\right)  \\
&=&\sum_{\tau \in \mathcal{P}\left( b\right) }\chi \left( \mathbf{X}%
_{t_{1}}\right) \cdot \cdot \cdot \chi \left( \mathbf{X}_{t_{s}}\right)
\sum_{\pi \in \left[ \tau \vee \sigma ,\hat{1}\right] }\mu \left( \pi ,\hat{1%
}\right) =\sum_{\tau :\tau \vee \sigma =\hat{1}}\chi \left( \mathbf{X}%
_{t_{1}}\right) \cdot \cdot \cdot \chi \left( \mathbf{X}_{t_{s}}\right)
\text{,}
\end{eqnarray*}%
where the last equality is a consequence of (\ref{maliMob}), since
\begin{equation}
\sum_{\pi \in \left[ \tau \vee \sigma ,\hat{1}\right] }\mu \left( \pi ,\hat{1%
}\right) =\left\{
\begin{array}{ll}
1 & \text{if }\tau \vee \sigma =\hat{1} \\
0 & \text{otherwise.}%
\end{array}%
\right.   \label{e : sum-mu}
\end{equation}
\end{proof}

\bigskip

For a single random variable $X$, one has (\ref{e : cum-single}): hence,
Proposition \ref{P : LeoShy} implies

\begin{corollary}
\label{C : SingleLeoShy} Let $X$ be a random variable such that $\mathbb{E}%
\left\vert X\right\vert ^{n}<\infty $. Then,
\begin{eqnarray}
\mathbb{E}\left[ X^{n}\right] &=&\sum_{\pi =\left\{ b_{1},...,b_{k}\right\}
\in \mathcal{P}\left( \left[ n\right] \right) }\chi _{\left\vert
b_{1}\right\vert }\left( X\right) \cdot \cdot \cdot \chi _{\left\vert
b_{k}\right\vert }\left( X\right)  \label{e : m-c2} \\
\mathbb{\chi }_{n}\left( X\right) &=&\sum_{\sigma =\left\{
a_{1},...,a_{r}\right\} \in \mathcal{P}\left( b\right) }\left( -1\right)
^{r-1}\left( r-1\right) !\mathbb{E}\left( X^{\left\vert a_{1}\right\vert
}\right) \cdot \cdot \cdot \mathbb{E}\left( X^{\left\vert a_{r}\right\vert
}\right)  \label{e : mc-2}
\end{eqnarray}
\end{corollary}

\bigskip

\textbf{Examples. }(i) Formula (\ref{LS2}), applied respectively to $%
b=\left\{ 1\right\} $ and to $b=\left\{ 1,2\right\} $, gives immediately the
well-known relations%
\begin{equation}
\chi \left( X\right) =\mathbb{E}\left( X\right) \text{ \ and \ }\chi \left(
X,Y\right) =\mathbb{E}\left( XY\right) -\mathbb{E}\left( X\right) \mathbb{E}%
\left( Y\right) =\mathbf{Cov}\left( X,Y\right) \text{.}  \label{easy}
\end{equation}

(ii) One has that
\begin{eqnarray*}
\chi \left( X_{1},X_{2},X_{3}\right)  &=&\mathbb{E}\left(
X_{1}X_{2}X_{3}\right) -\mathbb{E}\left( X_{1}X_{2}\right) \mathbb{E}\left(
X_{3}\right)  \\
&&-\mathbb{E}\left( X_{1}X_{3}\right) \mathbb{E}\left( X_{2}\right) -\mathbb{%
E}\left( X_{2}X_{3}\right) \mathbb{E}\left( X_{1}\right)  \\
&&+2\mathbb{E}\left( X_{1}\right) \mathbb{E}\left( X_{2}\right) \mathbb{E}%
\left( X_{3}\right) ,
\end{eqnarray*}%
so that, in particular,%
\begin{equation*}
\chi _{3}\left( X\right) =\mathbb{E}\left( X^{3}\right) -3\mathbb{E}\left(
X^{2}\right) \mathbb{E}\left( X\right) +2\mathbb{E}\left( X\right) ^{3}.
\end{equation*}

(iii) Let $\mathbf{G}_{\left[ n\right] }=\left( G_{1},...,G_{n}\right) $, $%
n\geq 3$, be a Gaussian vector such that $\mathbb{E}\left( G_{i}\right) =0$,
$i=1,...,n$. Then, for every $b\subseteq \left[ n\right] $ such that $%
\left\vert b\right\vert \geq 3$, we know from Section \ref{SS : DefCum} that%
\begin{equation*}
\chi \left( \mathbf{G}_{b}\right) =\chi \left( G_{i}:i\in b\right) =0.
\end{equation*}%
By applying this relation and formulae (\ref{LS1}) and (\ref{easy}) to $%
\mathbf{G}_{\left[ n\right] }$, one therefore obtains the well-known relation%
\begin{eqnarray}
&&\mathbb{E}\left[ G_{1}\times G_{2}\times \cdot \cdot \cdot \times G_{n}%
\right]  \label{Feynman} \\
&=&\left\{
\begin{array}{ll}
\sum_{\pi =\left\{ \left\{ i_{1},j_{1}\right\} ,...,\left\{
i_{k},j_{k}\right\} \right\} \in \mathcal{P}\left( \left[ n\right] \right) }%
\mathbb{E}\left( G_{i_{1}}G_{j_{1}}\right) \cdot \cdot \cdot \mathbb{E}%
\left( G_{i_{k}}G_{j_{k}}\right) , & n\text{ even} \\
0, & n\text{ odd.}%
\end{array}%
\right.  \notag
\end{eqnarray}%
Observe that, on the RHS of (\ref{Feynman}), the sum is taken over all
partitions $\pi $ such that each block of $\pi $ contains exactly two
elements.

\section{Diagrams and graphs\label{S : DG}}

In this section, we translate part of the notions presented in Section \ref%
{S : Lattice} into the language of diagrams and graphs, which are often used
in order to compute cumulants and moments of non-linear functionals of
random fields (see e.g. \cite{BrMa, ChaSlud, GiSu, Go
Taqqu, Marinucci, Surg1984, Sur}).

\subsection{Diagrams\label{SS : diagrams}}

Consider a finite set $b$. A \textsl{diagram }is a graphical representation
of a pair of partitions $\left( \pi ,\sigma \right) \subseteq \mathcal{P}%
\left( b\right) $, such that $\pi =\left\{ b_{1},...,b_{k}\right\} $ and $%
\sigma =\left\{ t_{1},...,t_{l}\right\} $. It is obtained as follows.

\begin{enumerate}
\item Order the elements of each block $b_{i}$, for $i=1,...,k$;

\item Associate with each block $b_{i}\in \pi $ a row of $\left\vert
b_{i}\right\vert $ \textsl{vertices} (represented as dots), in such a way
that the $j$th vertex of the $i$th row corresponds to the $j$th element of
the the block $b_{i}$;

\item For every $a=1,...,l$, draw a closed curve\ around the vertices
corresponding to the elements of the block $t_{a}\in \sigma $.
\end{enumerate}

We will denote by $\Gamma \left( \pi ,\sigma \right) $ the diagram of a pair
of partitions $\left( \pi ,\sigma \right) $.

\newpage

\textbf{Examples. }(i) If $b=\left[ 3\right] $ and $\pi =\sigma =\left\{
\left\{ 1,2\right\} ,\left\{ 3\right\} \right\} $, then $\Gamma \left( \pi
,\sigma \right) $ is represented in Fig. 1.

\begin{figure}[h]
\begin{center}
\psset{unit=0.7cm}
\begin{pspicture}(0,-1.5)(4.0,1.5)
\psframe[linewidth=0.02,dimen=outer](4.0,1.5)(0.0,-1.5)
\psdots[dotsize=0.15](0.6,0.88)
\psdots[dotsize=0.15](3.38,0.88)
\psdots[dotsize=0.15](0.6,-0.92)
\psellipse[linewidth=0.02,dimen=outer](1.99,0.89)(1.75,0.31)
\pscircle[linewidth=0.02,dimen=outer](0.61,-0.91){0.19}
\end{pspicture}
\caption{\sl A simple diagram}
\end{center}
\end{figure}
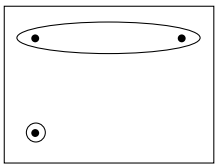

(ii) If $b=\left[ 8\right] $, and $\pi
=\left\{ \left\{ 1,2,3\right\} ,\left\{ 4,5\right\} ,\left\{ 6,7,8\right\}
\right\} $ and $\sigma =\left\{ \left\{ 1,4,6\right\} ,\left\{ 2,5\right\}
,\left\{ 3,7,8\right\} \right\} $, then $\Gamma \left( \pi ,\sigma \right) $
is represented in Fig. 2.
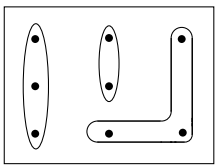
\begin{figure}[h]
\begin{center}
\psset{unit=0.7cm}
\begin{pspicture}(0,-1.5)(4.0,1.5)
\psframe[linewidth=0.02,dimen=outer](4.0,1.5)(0.0,-1.5)
\psdots[dotsize=0.15](0.6,0.88)
\psdots[dotsize=0.15](3.38,0.88)
\psdots[dotsize=0.15](0.6,-0.92)
\psdots[dotsize=0.15](0.6,-0.02)
\psdots[dotsize=0.15](2.0,0.9)
\psdots[dotsize=0.15](2.0,-0.92)
\psdots[dotsize=0.15](3.4,-0.9)
\psdots[dotsize=0.15](2.0,-0.02)
\psellipse[linewidth=0.02,dimen=outer](0.61,-0.02)(0.25,1.2)
\psellipse[linewidth=0.02,dimen=outer](2.0,0.41)(0.2,0.73)
\psline[linewidth=0.02cm](3.18,0.9)(3.18,-0.5)
\psline[linewidth=0.02cm](2.98,-0.68)(1.78,-0.68)
\psline[linewidth=0.02cm](3.58,-0.7)(3.58,0.9)
\psarc[linewidth=0.02](1.78,-0.88){0.2}{87.87891}{274.3987}
\psline[linewidth=0.02cm](1.78,-1.08)(3.0,-1.08)
\psarc[linewidth=0.02](3.0,-0.5){0.18}{261.8699}{6.340192}
\psline[linewidth=0.02cm](2.98,-1.08)(3.2,-1.08)
\psline[linewidth=0.02cm](3.58,-0.62)(3.58,-0.78)
\psline[linewidth=0.02cm](3.16,-1.08)(3.28,-1.08)
\psarc[linewidth=0.02](3.38,0.9){0.2}{0.0}{180.0}
\psline[linewidth=0.02cm](3.58,-0.72)(3.58,-0.94)
\psline[linewidth=0.02cm](3.24,-1.08)(3.44,-1.08)
\psarc[linewidth=0.02](3.44,-0.94){0.14}{270.0}{6.340192}
\end{pspicture}
\caption{\sl A diagram built from two three-block partitions}
\end{center}
\end{figure}

Hence, the rows in $\Gamma \left( \pi ,\sigma \right) $ indicate
the sets in $\pi $ and the curves indicate the sets in $\sigma $.

\bigskip

\textbf{Remarks. }(a) We use the terms \textquotedblleft
element\textquotedblright\ and \textquotedblleft vertex\textquotedblright\
interchangeably.

(b) Note that the diagram generated by the pair $\left( \pi ,\sigma \right) $
is different, in general, from the diagram generated by $\left( \sigma ,\pi
\right) $.

(c) Each diagram is a finite \textsl{hypergraph. }We recall that a finite
hypergraph is an object of the type $\left( V,E\right) $, where $V$ is a
finite set of vertices, and $E$ is a collection of (not necessarily
disjoint) nonempty subsets of $V$. The elements of $E$ are usually called
\textsl{edges}. In our setting, these are the blocks of $\sigma $.

(d) Note that, once a partition $\pi $ is specified, the diagram $\Gamma
\left( \pi ,\sigma \right) $\textsl{\ }encodes\textsl{\ }all the information
on $\sigma $.

\bigskip

Now fix a finite set $b$. In what follows, we will list and describe several
type of diagrams. They can be all characterized in terms of the lattice
structure of $\mathcal{P}\left( b\right) $, namely the partial ordering $%
\leq $ and the join and meet operations $\vee $ and $\wedge $, as described
in Section \ref{S : Lattice}. Recall that $\hat{1}=\left\{ b\right\} $, and $%
\hat{0}$ is the partition whose elements are the singletons of $b$.

\bigskip

\begin{description}
\item[Connected Diagrams.] The diagram $\Gamma \left( \pi ,\sigma \right) $
associated with two partitions $\left( \pi ,\sigma \right) $ is said to be
\textsl{connected} if $\pi \vee \sigma =\hat{1}$, that is, if the only
partition $\rho $ such that $\pi \leq \rho $ and $\sigma \leq \rho $ is the
maximal partition $\hat{1}$. The diagram appearing in Fig. 2 is connected,
whereas the one in Fig. 1 is not (indeed, in this case $\pi \vee \sigma =\pi
\vee \pi =\pi \neq \hat{1}$). Another example of a non-connected diagram
(see Fig. 3) is obtained by taking $b=\left[ 4\right] $, $\pi =\left\{
\left\{ 1,2\right\} ,\left\{ 3\right\} ,\left\{ 4\right\} \right\} $ and $%
\sigma =\left\{ \left\{ 1,2\right\} ,\left\{ 3,4\right\} \right\} $, so that
$\pi \leq \sigma $ (each block of $\pi $ is contained in a block of $\sigma $%
) and $\pi \vee \sigma =\sigma \neq \hat{1}$.

\begin{figure}[htbp]
\begin{center}
\psset{unit=0.7cm}
\begin{pspicture}(0,-1.5)(4.0,1.5)
\psframe[linewidth=0.02,dimen=outer](4.0,1.5)(0.0,-1.5)
\psdots[dotsize=0.15](0.6,0.88)
\psdots[dotsize=0.15](3.38,0.88)
\psdots[dotsize=0.15](0.6,-0.92)
\psellipse[linewidth=0.02,dimen=outer](1.99,0.89)(1.75,0.31)
\psdots[dotsize=0.15](0.6,-0.02)
\psellipse[linewidth=0.02,dimen=outer](0.61,-0.47)(0.27,0.75)
\end{pspicture}
\caption{\sl A non-connected diagram}
\end{center}
\end{figure}
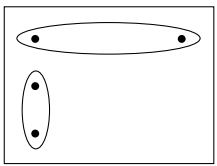

In other words, $\Gamma \left( \pi
,\sigma \right) $ \textsl{is connected if and only if the rows of the
diagram (the blocks of }$\pi $\textsl{) cannot be divided into two subsets,
each defining a separate diagram}. Fig. 4 shows that the diagram in Fig. 3
can be so divided.

\begin{figure}[htbp]
\begin{center}
\psset{unit=0.7cm}
\begin{pspicture}(0,-1.5)(4.0,1.5)
\psframe[linewidth=0.02,dimen=outer](4.0,1.5)(0.0,-1.5)
\psdots[dotsize=0.15](0.6,0.88)
\psdots[dotsize=0.15](3.38,0.88)
\psdots[dotsize=0.15](0.6,-0.92)
\psellipse[linewidth=0.02,dimen=outer](1.99,0.89)(1.75,0.31)
\psdots[dotsize=0.15](0.6,-0.02)
\psellipse[linewidth=0.02,dimen=outer](0.61,-0.47)(0.27,0.75)
\psline[linewidth=0.02cm](0.0,0.42)(3.98,0.42)
\end{pspicture}
\caption{\sl Dividing a non-connected diagram}
\end{center}
\end{figure}
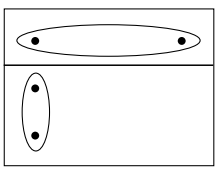

The diagram in Fig. 5, which has
the same partition $\pi $, but $\sigma =\left\{ \left\{ 1,3,4\right\}
,\left\{ 2\right\} \right\} $, is connected.

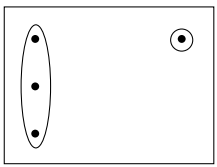
\begin{figure}[htbp]
\begin{center}
\psset{unit=0.7cm}
\begin{pspicture}(0,-1.5)(4.0,1.5)
\psframe[linewidth=0.02,dimen=outer](4.0,1.5)(0.0,-1.5)
\psdots[dotsize=0.15](0.6,0.88)
\psdots[dotsize=0.15](3.38,0.88)
\psdots[dotsize=0.15](0.6,-0.92)
\psdots[dotsize=0.15](0.6,-0.02)
\psellipse[linewidth=0.02,dimen=outer](0.61,-0.02)(0.29,1.18)
\pscircle[linewidth=0.02,dimen=outer](3.38,0.86){0.22}
\end{pspicture}
\caption{\sl A connected diagram}
\end{center}
\end{figure}

Note that we do not use\ the term
`connected' as one usually does in graph theory (indeed, the diagrams we
consider in this section are always \textsl{non-connected} \textsl{%
hypergraphs}, since their edges are disjoint by construction).

\item[Non-flat Diagrams.] The diagram $\Gamma \left( \pi ,\sigma \right) $
is \textsl{non-flat} if
\begin{equation*}
\pi \wedge \sigma =\hat{0},
\end{equation*}%
that is, if the only partition $\rho $ such that $\rho \leq \pi $ and $\rho
\leq \sigma $ is the minimal partition $\hat{0}$. It is easily seen that $%
\pi \wedge \sigma =\hat{0}$ if and only if for any two blocks $b_{j}\in
\pi $, $t_{a}\in \sigma $, the intersection $b_{j}\cap t_{a}$ either is
empty or contains exactly one element. Graphically, a non-flat graph is such
that the closed curves defining the blocks of $\sigma $ cannot join two
vertices in the same row (thus having a `flat' or `horizontal' portion). The
diagrams in Fig. 1-3 are all flat, whereas the diagram in Fig. 5 is
non-flat. Another non-flat diagram is given in Fig. 6, and is obtained by taking $b=\left[ 7%
\right] $, $\pi =\left\{ \left\{ 1,2,3\right\} ,\left\{ 4\right\} ,\left\{
5,6,7\right\} \right\} $ and $\sigma =\{\left\{ 1,4,5\right\} ,$ $\left\{
2,7\right\} ,$ $\left\{ 3,6\right\} \}$.

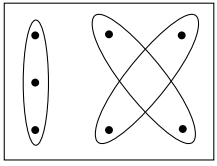
\begin{figure}[htbp]
\begin{center}
\psset{unit=0.7cm}
\begin{pspicture}(0,-1.86)(4.0,1.86)
\psframe[linewidth=0.02,dimen=outer](4.0,1.45)(0.0,-1.55)
\psdots[dotsize=0.15](0.6,0.83)
\psdots[dotsize=0.15](3.38,0.83)
\psdots[dotsize=0.15](0.6,-0.97)
\psdots[dotsize=0.15](0.6,-0.07)
\psdots[dotsize=0.15](2.0,0.85)
\psdots[dotsize=0.15](2.0,-0.97)
\psdots[dotsize=0.15](3.4,-0.95)
\psellipse[linewidth=0.02,dimen=outer](0.61,-0.07)(0.25,1.2)
\psbezier[linewidth=0.02](1.68,1.05)(1.68,0.25)(3.64,-1.85)(3.64,-1.05)(3.64,-0.25)(1.68,1.85)(1.68,1.05)
\psbezier[linewidth=0.02](3.7,1.05)(3.7,0.25)(1.74,-1.85)(1.74,-1.05)(1.74,-0.25)(3.7,1.85)(3.7,1.05)
\end{pspicture}
\caption{\sl A non-flat diagram}
\end{center}
\end{figure}

\item[Gaussian Diagrams.] We say that the diagram $\Gamma \left( \pi ,\sigma
\right) $ is \textsl{Gaussian}, whenever every block of $\sigma $ contains
exactly two elements. Plainly, Gaussian diagrams exists only if there is an
even number of vertices. When a diagram is Gaussian, one usually represents
the blocks of $\sigma $ not by closed curves, but by segments connecting two
vertices (which are viewed as the edges of the resulting graph). For
instance, a Gaussian (non-flat and connected) diagram is obtained in Fig. 7, where we have taken $%
b=\left[ 6\right] $, $\pi =\left\{ \left\{ 1,2,3\right\} ,\left\{ 4\right\}
,\left\{ 5,6\right\} \right\} $ and $\sigma =\left\{ \left\{ 1,4\right\}
,\left\{ 2,5\right\} ,\left\{ 3,6\right\} \right\} $.

\begin{figure}[htbp]
\begin{center}
\psset{unit=0.7cm}
\begin{pspicture}(0,-1.5)(4.0,1.5)
\psframe[linewidth=0.02,dimen=outer](4.0,1.5)(0.0,-1.5)
\psdots[dotsize=0.15](0.6,0.88)
\psdots[dotsize=0.15](3.38,0.88)
\psdots[dotsize=0.15](0.6,-0.92)
\psdots[dotsize=0.15](0.6,-0.02)
\psdots[dotsize=0.15](2.0,0.9)
\psdots[dotsize=0.15](2.0,-0.92)
\psline[linewidth=0.02cm](0.6,0.82)(0.6,0.02)
\psline[linewidth=0.02cm](0.6,-0.9)(2.0,0.9)
\psline[linewidth=0.02cm](1.98,-0.94)(3.38,0.86)
\end{pspicture}
\caption{\sl A Gaussian diagram}
\end{center}
\end{figure}
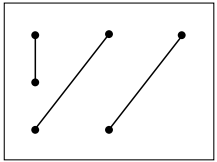

In the terminology of graph theory, a Gaussian diagram is a \textsl{non-connected
(non-directed) graph}. Since every vertex is connected with exactly another
vertex, one usually says that such a graph is a \textsl{perfect matching}.

\item[Circular (Gaussian) Diagrams.] Consider two partitions $\pi =\left\{
b_{1},...,b_{k}\right\} $ and $\sigma =\left\{ t_{1},...,t_{l}\right\} $
such that the blocks of $\sigma $ have size $\left\vert t_{a}\right\vert =2$
for every $a=1,...,l$. Then, the diagram $\Gamma \left( \pi ,\sigma \right) $
(which is Gaussian) is said to be \textsl{circular} if each row of $\Gamma
\left( \pi ,\sigma \right) $ is linked to both the previous and the next
row, with no other possible links except for the first and the last row,
which should also be linked together. This implies that the diagram is
connected. Formally, the diagram $\Gamma \left( \pi ,\sigma \right) $ is
circular (Gaussian) whenever the following properties hold (recall that $%
i\sim _{\sigma }j$ means that $i$ and $j$ belong to the same block of $%
\sigma $): (i) for every $p=2,...,k-1$ there exist $j_{1}\sim _{\sigma
}i_{1} $ and $j_{2}\sim _{\sigma }i_{2}$ such that $j_{1},j_{2}\in b_{p},$ $%
i_{1}\in b_{p-1}$ and $i_{2}\in b_{p+1}$, (ii) for every $p=2,...,k-1$ and
every $j\in b_{p}$, $j\sim _{\sigma }i$ implies that $i\in b_{p-1}\cup
b_{p+1}$, (iii) there exist $j_{1}\sim _{\sigma }i_{1}$ and $j_{2}\sim
_{\sigma }i_{2}$ such that $j_{1},j_{2}\in b_{k},$ $i_{1}\in b_{k-1}$ and $%
i_{2}\in b_{1}$, (iv) for every $j\in b_{k}$, $j\sim _{\sigma }i$ implies
that $i\in b_{1}\cup b_{k-1}$ (v) there exist $j_{1}\sim _{\sigma }i_{1}$
and $j_{2}\sim _{\sigma }i_{2}$ such that $j_{1},j_{2}\in b_{1},$ $i_{1}\in
b_{2}$ and $i_{2}\in b_{k}$,\ (vi) for every $j\in b_{1}$, $j\sim _{\sigma
}i $ implies that $i\in b_{k}\cup b_{2}$. For instance, a circular diagram
is obtained by taking $b=\left[ 10\right] $ and%
\begin{eqnarray*}
\pi &=&\left\{ \left\{ 1,2\right\} ,\left\{ 3,4\right\} ,\left\{ 5,6\right\}
\left\{ 7,8\right\} ,\left\{ 9,10\right\} \right\} \\
\sigma &=&\left\{ \left\{ 1,3\right\} ,\left\{ 2,9\right\} ,\left\{
4,6\right\} ,\left\{ 5,7\right\} ,\left\{ 8,10\right\} \right\} \text{,}
\end{eqnarray*}%
which implies that $\Gamma \left( \pi ,\sigma \right) $ is the diagram in Fig. 8.

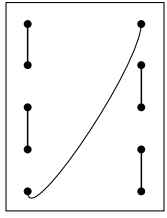
\begin{figure}[htbp]
\begin{center}
\psset{unit=0.7cm}
\begin{pspicture}(0,-2.21)(3.02,2.2)
\psframe[linewidth=0.02,dimen=outer](3.02,2.2)(0.0,-1.78)
\psdots[dotsize=0.15](0.42,1.78)
\psdots[dotsize=0.15](0.42,1.0)
\psdots[dotsize=0.15](0.42,0.2)
\psdots[dotsize=0.15](0.42,-0.6)
\psdots[dotsize=0.15](0.42,-1.4)
\psdots[dotsize=0.15](2.58,1.78)
\psdots[dotsize=0.15](2.58,1.0)
\psdots[dotsize=0.15](2.58,0.2)
\psdots[dotsize=0.15](2.58,-0.6)
\psdots[dotsize=0.15](2.58,-1.4)
\psline[linewidth=0.02cm](0.4,1.76)(0.38,1.76)
\psline[linewidth=0.02cm](0.42,1.76)(0.42,1.04)
\psline[linewidth=0.02cm](0.42,0.18)(0.42,-0.58)
\psline[linewidth=0.02cm](2.58,0.98)(2.58,0.22)
\psline[linewidth=0.02cm](2.58,-0.6)(2.58,-1.4)
\psbezier[linewidth=0.02](2.58,1.8)(2.58,1.0)(0.42,-2.2)(0.42,-1.4)
\end{pspicture}
\caption{\sl A circular diagram}
\end{center}
\end{figure}

Another example of a circular diagram is given in Fig. 9. It is obtained from $b=\left[
12\right] $ and%
\begin{eqnarray*}
\pi &=&\left\{ \left\{ 1,2,3\right\} ,\left\{ 4,5\right\} ,\left\{
6,7\right\} ,\left\{ 8,9\right\} ,\left\{ 10,11,12\right\} \right\} \\
\sigma &=&\left\{ \left\{ 1,4\right\} ,\left\{ 2,11\right\} ,\left\{
3,10\right\} ,\left\{ 5,7\right\} ,\left\{ 6,8\right\} ,\left\{ 9,12\right\}
\right\} \text{.}
\end{eqnarray*}

\begin{figure}[htbp]
\begin{center}
\psset{unit=0.7cm}
\begin{pspicture}(0,-2.21)(3.02,2.2)
\psframe[linewidth=0.02,dimen=outer](3.02,2.2)(0.0,-1.78)
\psdots[dotsize=0.15](0.42,1.78)
\psdots[dotsize=0.15](0.42,1.0)
\psdots[dotsize=0.15](0.42,0.2)
\psdots[dotsize=0.15](0.42,-0.6)
\psdots[dotsize=0.15](0.42,-1.4)
\psdots[dotsize=0.15](2.58,1.78)
\psdots[dotsize=0.15](2.58,1.0)
\psdots[dotsize=0.15](2.58,0.2)
\psdots[dotsize=0.15](2.58,-0.6)
\psdots[dotsize=0.15](2.58,-1.4)
\psline[linewidth=0.02cm](0.4,1.76)(0.4,1.76)
\psline[linewidth=0.02cm](0.42,1.76)(0.42,1.04)
\psline[linewidth=0.02cm](0.42,0.18)(0.42,-0.58)
\psline[linewidth=0.02cm](2.58,0.98)(2.58,0.22)
\psline[linewidth=0.02cm](2.58,-0.6)(2.58,-1.4)
\psdots[dotsize=0.15](1.5,1.8)
\psdots[dotsize=0.15](1.5,-1.4)
\psline[linewidth=0.02cm](1.48,1.82)(1.46,1.82)
\psline[linewidth=0.02cm](1.48,1.82)(1.48,1.8)
\psbezier[linewidth=0.02](2.58,1.82)(2.58,1.0)(0.42,-2.2)(0.42,-1.4)
\psline[linewidth=0.02cm](1.48,1.82)(1.48,1.82)
\psline[linewidth=0.02cm](1.48,1.8)(1.5,-1.36)
\end{pspicture}
\caption{\sl A circular diagram with rows of different size}
\end{center}
\end{figure}
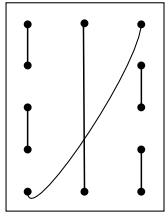

\end{description}

\bigskip

\textbf{Examples. }(i) Thanks to the previous discussion, one can
immediately reformulate Malyshev's formula (\ref{LS3}) as follows. For every
finite set $b$ and every $\sigma =\left\{ b_{1},...,b_{k}\right\} \in
\mathcal{P}\left( b\right) $,
\begin{equation}
\chi \left( \mathbf{X}^{b_{1}},...,\mathbf{X}^{b_{k}}\right) =\sum
_{\substack{ \tau =\left\{ t_{1},...,t_{s}\right\} \in \mathcal{P}\left(
b\right) \text{ }  \\ \Gamma \left( \sigma ,\tau \right) \text{ is connected}
}}\chi \left( \mathbf{X}_{t_{1}}\right) \cdot \cdot \cdot \chi \left(
\mathbf{X}_{t_{s}}\right) \text{.}  \label{Maly}
\end{equation}

(ii) Suppose that the random variables $X_{1},X_{2},X_{3}$ are such that $%
\mathbb{E}\left\vert X_{i}\right\vert ^{3}<\infty $, $i=1,2,3$. We have
already applied formula (\ref{Maly}) in order to compute the cumulant $\chi
\left( X_{1}X_{2},X_{3}\right) $. Here, we shall give a graphical
demonstration. Recall that, in this case, $b=\left[ 3\right] =\left\{
1,2,3\right\} $, and that the relevant partition is $\sigma =\left\{ \left\{
1,2\right\} ,\left\{ 3\right\} \right\} $. There are only three partitions $%
\tau _{1},\tau _{2},\tau _{3}\in \mathcal{P}\left( \left[ 3\right] \right) $
such that $\Gamma \left( \sigma ,\tau _{1}\right) ,$ $\Gamma \left( \sigma
,\tau _{2}\right) $ and $\Gamma \left( \sigma ,\tau _{3}\right) $ are
connected, namely $\tau _{1}=\hat{1}$, $\tau _{2}=\left\{ \left\{
1,3\right\} ,\left\{ 2\right\} \right\} $ and $\tau _{3}=\left\{ \left\{
1\right\} ,\left\{ 2,3\right\} \right\} $. The diagrams $\Gamma \left(
\sigma ,\tau _{1}\right) ,$ $\Gamma \left( \sigma ,\tau _{2}\right) $ and $%
\Gamma \left( \sigma ,\tau _{3}\right) $ are represented in Fig. 10. Relation (\ref{Maly}) thus implies that
\begin{eqnarray*}
\chi \left( X_{1}X_{2},X_{3}\right) &=&\chi \left( X_{1},X_{2},X_{3}\right)
+\chi \left( X_{1},X_{3}\right) \chi \left( X_{2}\right) +\chi \left(
X_{1}\right) \chi \left( X_{2},X_{3}\right) \\
&=&\chi \left( X_{1},X_{2},X_{3}\right) +\mathbf{Cov}\left(
X_{1},X_{3}\right) \mathbb{E}\left( X_{3}\right) +\mathbb{E}\left(
X_{1}\right) \mathbf{Cov}\left( X_{2},X_{3}\right) \text{,}
\end{eqnarray*}%
where we have used (\ref{easy}).

\begin{figure}[htbp]
\begin{center}
\psset{unit=0.7cm}
\begin{pspicture}(0,-2.0365882)(7.28,2.0365887)
\psframe[linewidth=0.02,dimen=outer](2.44,1.51)(0.0,-1.49)
\psframe[linewidth=0.02,dimen=outer](4.86,1.51)(2.42,-1.49)
\psframe[linewidth=0.02,dimen=outer](7.28,1.51)(4.84,-1.49)
\psdots[dotsize=0.15](0.42,1.11)
\psdots[dotsize=0.15](1.82,1.11)
\psdots[dotsize=0.15](0.42,-1.09)
\psdots[dotsize=0.15](2.82,1.11)
\psdots[dotsize=0.15](4.22,1.11)
\psdots[dotsize=0.15](5.42,1.11)
\psdots[dotsize=0.15](6.82,1.11)
\psdots[dotsize=0.15](2.82,-1.11)
\psdots[dotsize=0.15](5.42,-1.11)
\psellipse[linewidth=0.02,dimen=outer](2.83,-0.01)(0.25,1.32)
\pscircle[linewidth=0.02,dimen=outer](4.22,1.11){0.16}
\psdots[dotsize=0.15](5.42,1.11)
\pscircle[linewidth=0.02,dimen=outer](5.42,1.11){0.16}
\psbezier[linewidth=0.02](5.223007,-1.2309201)(5.306144,-2.0265887)(7.120129,0.43525168)(7.0369925,1.2309201)(6.953856,2.0265887)(5.1398706,-0.43525168)(5.223007,-1.2309201)
\psline[linewidth=0.02cm](0.2,1.21)(0.2,-1.09)
\psline[linewidth=0.02cm](0.6,-1.07)(0.6,0.73)
\psline[linewidth=0.02cm](0.74,0.91)(1.82,0.91)
\psline[linewidth=0.02cm](1.8,1.31)(0.34,1.31)
\psarc[linewidth=0.02](0.32,1.19){0.12}{81.869896}{186.3402}
\psline[linewidth=0.02cm](0.36,1.31)(0.32,1.31)
\psarc[linewidth=0.02](0.75,0.76){0.15}{85.601295}{186.3402}
\psline[linewidth=0.02cm](0.6,0.71)(0.6,0.75)
\psarc[linewidth=0.02](1.82,1.11){0.2}{263.65982}{96.340195}
\psarc[linewidth=0.02](0.4,-1.09){0.2}{180.0}{0.0}
\psline[linewidth=0.02cm](0.6,-1.01)(0.6,-1.09)
\end{pspicture}
\caption{\sl Computing cumulants by connected diagrams}
\end{center}
\end{figure}
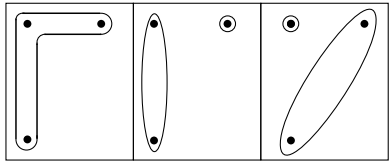

\subsection{Solving the equation $\protect\sigma \wedge \protect\pi =\hat{0}$
\label{SS : Solve}}

Let $\pi $ be a partition of $\left[ n\right] =\left\{ 1,...,n\right\} $.
One is often asked, as will be the case in Section \ref{SS : mult General},
to find all partitions $\sigma \in \mathcal{P}\left( \left[ n\right] \right)
$ such that
\begin{equation}
\sigma \wedge \pi =\hat{0},  \label{zero-chapeau}
\end{equation}%
where, as usual, $\hat{0}=\left\{ \left\{ 1\right\} ,...,\left\{ n\right\}
\right\} $, that is, $\hat{0}$ is the partition made up of singletons. The
use of diagrams provides an easy way to solve (\ref{zero-chapeau}), since (%
\ref{zero-chapeau}) holds if and only if the diagram $\Gamma \left( \pi
,\sigma \right) $ is non-flat. Hence, proceed as in Section \ref{SS :
diagrams}, by (\textbf{1}) ordering the blocks of $\pi $, (\textbf{2})
associating with each block of $\pi $ a row of the diagram, the number of
points in a row being equal to the number of elements in the block, and (%
\textbf{3}) drawing non-flat closed curves around the points of the diagram.

\smallskip

\textbf{Examples. }(i) Let $n=2$ and $\pi =\left\{ \left\{ 1\right\}
,\left\{ 2\right\} \right\} =\hat{0}.$ Then, $\sigma _{1}=\pi =\hat{0}$ and $%
\sigma _{2}=\hat{1}$ (as represented in Fig. 11) solve equation (\ref%
{zero-chapeau}). Note that $\mathcal{P}\left( \left[ 2\right] \right)
=\left\{ \sigma _{1},\sigma _{2}\right\} $.

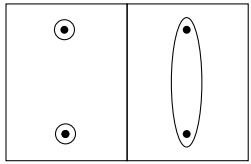
\begin{figure}[htbp]
\begin{center}
\psset{unit=0.7cm}
\begin{pspicture}(0,-1.5)(4.62,1.5)
\psframe[linewidth=0.02,dimen=outer](2.32,1.5)(0.0,-1.5)
\psframe[linewidth=0.02,dimen=outer](4.62,1.5)(2.3,-1.5)
\psdots[dotsize=0.15](1.12,1.0)
\psdots[dotsize=0.15](1.14,-0.98)
\psdots[dotsize=0.15](3.44,1.0)
\psdots[dotsize=0.15](3.44,-0.98)
\pscircle[linewidth=0.02,dimen=outer](1.12,1.0){0.2}
\psdots[dotsize=0.15](1.14,-0.98)
\pscircle[linewidth=0.02,dimen=outer](1.14,-0.98){0.2}
\psellipse[linewidth=0.02,dimen=outer](3.44,0.0)(0.3,1.24)
\end{pspicture}
\caption{\sl Solving $\sigma\wedge\pi = \hat{0}$ in the simplest case}
\end{center}
\end{figure}

(ii) Let $n=3$ and $\pi =\left\{ \left\{ 1,2\right\} ,\left\{ 3\right\}
\right\} $. Then, $\sigma _{1}=\hat{0}$, $\sigma _{2}=\left\{ \left\{
1,3\right\} ,\left\{ 2\right\} \right\} $ and $\sigma _{3}=\left\{ \left\{
1\right\} ,\left\{ 2,3\right\} \right\} $ (see Fig. 12) are the only
elements of $\mathcal{P}\left( \left[ 3\right] \right) $ solving (\ref{zero-chapeau}).

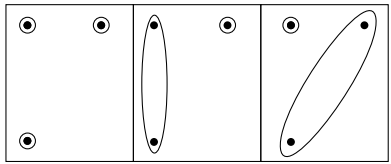
\begin{figure}[htbp]
\begin{center}
\psset{unit=0.7cm}
\begin{pspicture}(0,-2.0365882)(7.28,2.0365887)
\psframe[linewidth=0.02,dimen=outer](2.44,1.51)(0.0,-1.49)
\psframe[linewidth=0.02,dimen=outer](4.86,1.51)(2.42,-1.49)
\psframe[linewidth=0.02,dimen=outer](7.28,1.51)(4.84,-1.49)
\psdots[dotsize=0.15](0.42,1.11)
\psdots[dotsize=0.15](1.82,1.11)
\psdots[dotsize=0.15](0.42,-1.09)
\psdots[dotsize=0.15](2.82,1.11)
\psdots[dotsize=0.15](4.22,1.11)
\psdots[dotsize=0.15](5.42,1.11)
\psdots[dotsize=0.15](6.82,1.11)
\psdots[dotsize=0.15](2.82,-1.11)
\psdots[dotsize=0.15](5.42,-1.11)
\psellipse[linewidth=0.02,dimen=outer](2.83,-0.01)(0.25,1.32)
\pscircle[linewidth=0.02,dimen=outer](4.22,1.11){0.16}
\psdots[dotsize=0.15](5.42,1.11)
\pscircle[linewidth=0.02,dimen=outer](5.42,1.11){0.16}
\psbezier[linewidth=0.02](5.223007,-1.2309201)(5.306144,-2.0265887)(7.120129,0.43525168)(7.0369925,1.2309201)(6.953856,2.0265887)(5.1398706,-0.43525168)(5.223007,-1.2309201)
\psdots[dotsize=0.15](0.42,1.11)
\pscircle[linewidth=0.02,dimen=outer](0.42,1.11){0.16}
\psdots[dotsize=0.15](1.82,1.11)
\pscircle[linewidth=0.02,dimen=outer](1.82,1.11){0.16}
\psdots[dotsize=0.15](0.42,-1.09)
\pscircle[linewidth=0.02,dimen=outer](0.42,-1.09){0.16}
\end{pspicture}
\caption{\sl Solving $\sigma\wedge\pi = \hat{0}$ in a three-vertex diagram}
\end{center}
\end{figure}

(iii) Let $n=4$ and $\pi =\left\{ \left\{ 1,2\right\} ,\left\{ 3,4\right\}
\right\} $. Then, there are exactly seven $\sigma \in \mathcal{P}\left( %
\left[ 4\right] \right) $ solving (\ref{zero-chapeau}). They are all
represented in Fig. 13.

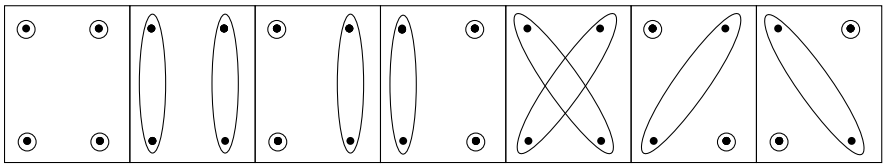
\begin{figure}[htbp]
\begin{center}
\psset{unit=0.7cm}
\begin{pspicture}(0,-1.9852024)(16.68,1.9852028)
\psframe[linewidth=0.02,dimen=outer](2.4,1.5)(0.0,-1.5)
\psdots[dotsize=0.15](0.42,1.06)
\psdots[dotsize=0.15](1.8,1.06)
\psdots[dotsize=0.15](0.44,-1.08)
\psdots[dotsize=0.15](1.82,-1.08)
\psframe[linewidth=0.02,dimen=outer](4.78,1.5)(2.38,-1.5)
\psdots[dotsize=0.15](2.8,1.06)
\psdots[dotsize=0.15](4.18,1.06)
\psdots[dotsize=0.15](2.82,-1.08)
\psdots[dotsize=0.15](4.2,-1.08)
\psframe[linewidth=0.02,dimen=outer](16.68,1.5)(14.28,-1.5)
\psdots[dotsize=0.15](14.7,1.06)
\psdots[dotsize=0.15](16.08,1.06)
\psdots[dotsize=0.15](14.72,-1.08)
\psdots[dotsize=0.15](16.1,-1.08)
\psframe[linewidth=0.02,dimen=outer](7.16,1.5)(4.76,-1.5)
\psdots[dotsize=0.15](5.18,1.06)
\psdots[dotsize=0.15](6.56,1.06)
\psdots[dotsize=0.15](5.2,-1.08)
\psdots[dotsize=0.15](6.58,-1.08)
\psframe[linewidth=0.02,dimen=outer](9.54,1.5)(7.14,-1.5)
\psdots[dotsize=0.15](7.56,1.06)
\psdots[dotsize=0.15](8.94,1.06)
\psdots[dotsize=0.15](7.58,-1.08)
\psdots[dotsize=0.15](8.96,-1.08)
\psframe[linewidth=0.02,dimen=outer](11.92,1.5)(9.52,-1.5)
\psdots[dotsize=0.15](9.94,1.06)
\psdots[dotsize=0.15](11.32,1.06)
\psdots[dotsize=0.15](9.96,-1.08)
\psdots[dotsize=0.15](11.34,-1.08)
\psframe[linewidth=0.02,dimen=outer](14.3,1.5)(11.9,-1.5)
\psdots[dotsize=0.15](12.32,1.06)
\psdots[dotsize=0.15](13.7,1.06)
\psdots[dotsize=0.15](12.34,-1.08)
\psdots[dotsize=0.15](13.72,-1.08)
\pscircle[linewidth=0.02,dimen=outer](0.42,1.04){0.18}
\psdots[dotsize=0.15](1.8,1.06)
\pscircle[linewidth=0.02,dimen=outer](1.8,1.04){0.18}
\psdots[dotsize=0.15](0.44,-1.08)
\psdots[dotsize=0.15](1.82,-1.08)
\pscircle[linewidth=0.02,dimen=outer](0.44,-1.1){0.18}
\psdots[dotsize=0.15](1.82,-1.08)
\pscircle[linewidth=0.02,dimen=outer](1.82,-1.1){0.18}
\psdots[dotsize=0.15](2.8,1.06)
\psdots[dotsize=0.15](2.8,1.06)
\psdots[dotsize=0.15](2.82,-1.08)
\psdots[dotsize=0.15](2.82,-1.08)
\psdots[dotsize=0.15](5.18,1.06)
\psdots[dotsize=0.15](5.2,-1.08)
\psdots[dotsize=0.15](5.18,1.06)
\psdots[dotsize=0.15](5.18,1.06)
\pscircle[linewidth=0.02,dimen=outer](5.18,1.04){0.18}
\psdots[dotsize=0.15](5.2,-1.08)
\psdots[dotsize=0.15](5.2,-1.08)
\pscircle[linewidth=0.02,dimen=outer](5.2,-1.1){0.18}
\psdots[dotsize=0.15](8.94,1.06)
\psdots[dotsize=0.15](8.96,-1.08)
\psdots[dotsize=0.15](8.94,1.06)
\psdots[dotsize=0.15](8.94,1.06)
\pscircle[linewidth=0.02,dimen=outer](8.94,1.04){0.18}
\psdots[dotsize=0.15](8.96,-1.08)
\psdots[dotsize=0.15](8.96,-1.08)
\pscircle[linewidth=0.02,dimen=outer](8.96,-1.1){0.18}
\psdots[dotsize=0.15](12.32,1.06)
\psdots[dotsize=0.15](12.32,1.06)
\psdots[dotsize=0.15](12.32,1.06)
\pscircle[linewidth=0.02,dimen=outer](12.32,1.04){0.18}
\psdots[dotsize=0.15](13.72,-1.08)
\psdots[dotsize=0.15](13.72,-1.08)
\psdots[dotsize=0.15](13.72,-1.08)
\pscircle[linewidth=0.02,dimen=outer](13.72,-1.1){0.18}
\psdots[dotsize=0.15](16.08,1.06)
\psdots[dotsize=0.15](16.08,1.06)
\psdots[dotsize=0.15](16.08,1.06)
\pscircle[linewidth=0.02,dimen=outer](16.08,1.04){0.18}
\psdots[dotsize=0.15](14.72,-1.08)
\psdots[dotsize=0.15](14.72,-1.08)
\psdots[dotsize=0.15](14.72,-1.08)
\pscircle[linewidth=0.02,dimen=outer](14.72,-1.1){0.18}
\psellipse[linewidth=0.02,dimen=outer](2.82,0.01)(0.26,1.33)
\psdots[dotsize=0.15](4.18,1.06)
\psdots[dotsize=0.15](4.18,1.06)
\psdots[dotsize=0.15](4.18,1.06)
\psellipse[linewidth=0.02,dimen=outer](4.2,0.01)(0.26,1.33)
\psdots[dotsize=0.15](6.56,1.06)
\psdots[dotsize=0.15](6.56,1.06)
\psdots[dotsize=0.15](6.56,1.06)
\psellipse[linewidth=0.02,dimen=outer](6.58,0.01)(0.26,1.33)
\psdots[dotsize=0.15](7.56,1.04)
\psdots[dotsize=0.15](7.56,1.04)
\psdots[dotsize=0.15](7.56,1.04)
\psellipse[linewidth=0.02,dimen=outer](7.58,-0.01)(0.26,1.33)
\psbezier[linewidth=0.02](9.688816,1.178577)(9.762214,0.38195118)(11.644583,-1.9552028)(11.571184,-1.158577)(11.497786,-0.36195117)(9.6154175,1.9752028)(9.688816,1.178577)
\psbezier[linewidth=0.02](14.448815,1.158577)(14.522214,0.36195117)(16.404583,-1.9752028)(16.331184,-1.178577)(16.257786,-0.38195118)(14.375418,1.9552028)(14.448815,1.158577)
\psbezier[linewidth=0.02](11.631185,1.178577)(11.557786,0.38195118)(9.675417,-1.9552028)(9.748816,-1.158577)(9.822214,-0.36195117)(11.704582,1.9752028)(11.631185,1.178577)
\psbezier[linewidth=0.02](13.991184,1.178577)(13.917786,0.38195118)(12.035418,-1.9552028)(12.108816,-1.158577)(12.182214,-0.36195117)(14.064583,1.9752028)(13.991184,1.178577)
\end{pspicture}
\caption{\sl The seven solutions of $\sigma\wedge\pi = \hat{0}$ in a four-vertex diagram}
\end{center}
\end{figure}

(iv) Let $n=4$ and $\pi =\left\{ \left\{ 1,2\right\} ,\left\{ 3\right\}
,\left\{ 4\right\} \right\} $. Then, there are ten $\sigma \in \mathcal{P}%
\left( \left[ 4\right] \right) $ that are solutions of (\ref{zero-chapeau}).
They are all represented in Fig. 14.

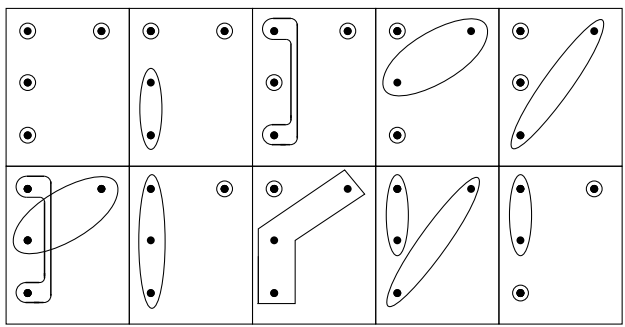
\begin{figure}[htbp]
%\begin{center}
\centering
\psset{unit=0.7cm}
\begin{pspicture}(0,-3.3513203)(11.72,3.3513203)
\psframe[linewidth=0.02,dimen=outer](2.36,2.95)(0.0,-0.07)
\psdots[dotsize=0.15](0.42,2.51)
\psdots[dotsize=0.15](1.82,2.51)
\psdots[dotsize=0.15](0.42,0.53)
\psframe[linewidth=0.02,dimen=outer](2.36,-0.05)(0.0,-3.07)
\psdots[dotsize=0.15](0.42,-0.49)
\psdots[dotsize=0.15](1.82,-0.49)
\psdots[dotsize=0.15](0.42,-2.47)
\psframe[linewidth=0.02,dimen=outer](11.72,2.95)(9.36,-0.07)
\psdots[dotsize=0.15](9.78,2.51)
\psdots[dotsize=0.15](11.18,2.51)
\psdots[dotsize=0.15](9.78,0.53)
\psframe[linewidth=0.02,dimen=outer](11.72,-0.05)(9.36,-3.07)
\psdots[dotsize=0.15](9.78,-0.49)
\psdots[dotsize=0.15](11.18,-0.49)
\psdots[dotsize=0.15](9.78,-2.47)
\psframe[linewidth=0.02,dimen=outer](4.7,2.95)(2.34,-0.07)
\psdots[dotsize=0.15](2.76,2.51)
\psdots[dotsize=0.15](4.16,2.51)
\psdots[dotsize=0.15](2.76,0.53)
\psframe[linewidth=0.02,dimen=outer](4.7,-0.05)(2.34,-3.07)
\psdots[dotsize=0.15](2.76,-0.49)
\psdots[dotsize=0.15](4.16,-0.49)
\psdots[dotsize=0.15](2.76,-2.47)
\psframe[linewidth=0.02,dimen=outer](7.04,2.95)(4.68,-0.07)
\psdots[dotsize=0.15](5.1,2.51)
\psdots[dotsize=0.15](6.5,2.51)
\psdots[dotsize=0.15](5.1,0.53)
\psframe[linewidth=0.02,dimen=outer](7.04,-0.05)(4.68,-3.07)
\psdots[dotsize=0.15](5.1,-0.49)
\psdots[dotsize=0.15](6.5,-0.49)
\psdots[dotsize=0.15](5.1,-2.47)
\psframe[linewidth=0.02,dimen=outer](9.38,2.95)(7.02,-0.07)
\psdots[dotsize=0.15](7.44,2.51)
\psdots[dotsize=0.15](8.84,2.51)
\psdots[dotsize=0.15](7.44,0.53)
\psframe[linewidth=0.02,dimen=outer](9.38,-0.05)(7.02,-3.07)
\psdots[dotsize=0.15](7.44,-0.49)
\psdots[dotsize=0.15](8.84,-0.49)
\psdots[dotsize=0.15](7.44,-2.47)
\psdots[dotsize=0.15](0.42,1.53)
\psdots[dotsize=0.15](2.76,2.51)
\psdots[dotsize=0.15](2.76,0.53)
\psdots[dotsize=0.15](2.76,1.53)
\psdots[dotsize=0.15](0.42,-0.49)
\psdots[dotsize=0.15](0.42,-2.47)
\psdots[dotsize=0.15](0.42,-1.47)
\psdots[dotsize=0.15](2.76,-0.49)
\psdots[dotsize=0.15](2.76,-2.47)
\psdots[dotsize=0.15](2.76,-1.47)
\psdots[dotsize=0.15](5.1,2.51)
\psdots[dotsize=0.15](5.1,0.53)
\psdots[dotsize=0.15](5.1,1.53)
\psdots[dotsize=0.15](5.1,-0.49)
\psdots[dotsize=0.15](5.1,-2.47)
\psdots[dotsize=0.15](5.1,-1.47)
\psdots[dotsize=0.15](7.44,2.51)
\psdots[dotsize=0.15](7.44,0.53)
\psdots[dotsize=0.15](7.44,1.53)
\psdots[dotsize=0.15](7.44,-0.49)
\psdots[dotsize=0.15](7.44,-2.47)
\psdots[dotsize=0.15](9.78,2.51)
\psdots[dotsize=0.15](9.78,0.53)
\psdots[dotsize=0.15](9.78,1.53)
\psdots[dotsize=0.15](9.78,-0.49)
\psdots[dotsize=0.15](9.78,-2.47)
\pscircle[linewidth=0.02,dimen=outer](0.42,2.51){0.16}
\psdots[dotsize=0.15](1.82,2.51)
\pscircle[linewidth=0.02,dimen=outer](1.82,2.51){0.16}
\psdots[dotsize=0.15](0.42,1.53)
\pscircle[linewidth=0.02,dimen=outer](0.42,1.53){0.16}
\psdots[dotsize=0.15](0.42,0.53)
\pscircle[linewidth=0.02,dimen=outer](0.42,0.53){0.16}
\psdots[dotsize=0.15](2.76,2.51)
\pscircle[linewidth=0.02,dimen=outer](2.76,2.51){0.16}
\psdots[dotsize=0.15](4.16,2.51)
\pscircle[linewidth=0.02,dimen=outer](4.16,2.51){0.16}
\psdots[dotsize=0.15](5.1,1.53)
\pscircle[linewidth=0.02,dimen=outer](5.1,1.53){0.16}
\psdots[dotsize=0.15](6.5,2.51)
\pscircle[linewidth=0.02,dimen=outer](6.5,2.51){0.16}
\psdots[dotsize=0.15](7.44,2.51)
\pscircle[linewidth=0.02,dimen=outer](7.44,2.51){0.16}
\psdots[dotsize=0.15](7.44,0.53)
\pscircle[linewidth=0.02,dimen=outer](7.44,0.53){0.16}
\psdots[dotsize=0.15](9.78,2.51)
\psdots[dotsize=0.15](9.78,1.53)
\pscircle[linewidth=0.02,dimen=outer](9.78,2.51){0.16}
\psdots[dotsize=0.15](9.78,1.53)
\pscircle[linewidth=0.02,dimen=outer](9.78,1.53){0.16}
\psdots[dotsize=0.15](4.16,-0.49)
\psdots[dotsize=0.15](4.16,-0.49)
\pscircle[linewidth=0.02,dimen=outer](4.16,-0.49){0.16}
\psdots[dotsize=0.15](5.1,-0.49)
\psdots[dotsize=0.15](5.1,-0.49)
\pscircle[linewidth=0.02,dimen=outer](5.1,-0.49){0.16}
\psdots[dotsize=0.15](11.18,-0.49)
\psdots[dotsize=0.15](11.18,-0.49)
\pscircle[linewidth=0.02,dimen=outer](11.18,-0.49){0.16}
\psdots[dotsize=0.15](9.78,-2.47)
\psdots[dotsize=0.15](9.78,-2.47)
\pscircle[linewidth=0.02,dimen=outer](9.78,-2.47){0.16}
\psellipse[linewidth=0.02,dimen=outer](2.76,1.03)(0.22,0.78)
\psdots[dotsize=0.15](7.44,-0.49)
\psellipse[linewidth=0.02,dimen=outer](7.44,-0.99)(0.22,0.78)
\psdots[dotsize=0.15](9.78,-0.49)
\psellipse[linewidth=0.02,dimen=outer](9.78,-0.99)(0.22,0.78)
\psdots[dotsize=0.15](7.44,-1.47)
\psdots[dotsize=0.15](9.78,-1.47)
\psellipse[linewidth=0.02,dimen=outer](2.78,-1.49)(0.26,1.28)
\psbezier[linewidth=0.02](9.147269,2.424366)(9.00792,1.6365958)(7.033382,0.8078638)(7.1727304,1.595634)(7.3120794,2.3834043)(9.286618,3.2121363)(9.147269,2.424366)
\psbezier[linewidth=0.02](2.1272693,-0.575634)(1.9879205,-1.3634042)(0.013381832,-2.1921363)(0.15273063,-1.404366)(0.29207942,-0.6165958)(2.2666183,0.21213622)(2.1272693,-0.575634)
\psbezier[linewidth=0.02](11.349141,2.5574944)(11.189088,1.7736684)(9.450807,-0.34132028)(9.610859,0.44250566)(9.770912,1.2263316)(11.509193,3.3413203)(11.349141,2.5574944)
\psbezier[linewidth=0.02](8.9891405,-0.44250566)(8.829088,-1.2263316)(7.0908065,-3.3413203)(7.2508593,-2.5574944)(7.410912,-1.7736684)(9.149194,0.34132028)(8.9891405,-0.44250566)
\psline[linewidth=0.02](4.8,-2.65)(4.8,-1.25)(6.44,-0.13)(6.82,-0.59)(5.5,-1.47)(5.5,-2.67)(4.78,-2.67)(4.78,-2.67)
\pscustom[linewidth=0.02]
{
\newpath
\moveto(4.8,-1.77)
\lineto(4.8,-1.82)
\curveto(4.8,-1.845)(4.8,-1.895)(4.8,-1.92)
\curveto(4.8,-1.945)(4.8,-2.0)(4.8,-2.03)
\curveto(4.8,-2.06)(4.8,-2.115)(4.8,-2.14)
\curveto(4.8,-2.165)(4.8,-2.22)(4.8,-2.25)
\curveto(4.8,-2.28)(4.8,-2.335)(4.8,-2.36)
\curveto(4.8,-2.385)(4.8,-2.435)(4.8,-2.46)
\curveto(4.8,-2.485)(4.8,-2.535)(4.8,-2.56)
\curveto(4.8,-2.585)(4.8,-2.625)(4.8,-2.64)
\curveto(4.8,-2.655)(4.8,-2.665)(4.8,-2.65)
}
\psline[linewidth=0.02](5.08,2.35)(5.34,2.35)(5.34,2.35)(5.34,2.35)
\psline[linewidth=0.02cm](5.06,0.73)(5.32,0.73)
\psline[linewidth=0.02cm](5.42,0.85)(5.42,2.25)
\psline[linewidth=0.02cm](5.06,0.35)(5.32,0.35)
\psline[linewidth=0.02cm](5.54,0.49)(5.54,2.65)
\psline[linewidth=0.02cm](5.08,2.75)(5.42,2.75)
\psarc[linewidth=0.02](5.08,2.55){0.2}{90.0}{270.0}
\psarc[linewidth=0.02](5.07,0.54){0.19}{90.0}{270.0}
\psarc[linewidth=0.02](5.41,2.62){0.13}{353.65982}{84.80557}
\psarc[linewidth=0.02](5.33,2.26){0.09}{353.65982}{84.80557}
\psline[linewidth=0.02cm](5.42,2.27)(5.42,2.25)
\psarc[linewidth=0.02](5.33,0.82){0.09}{258.69006}{9.462322}
\psline[linewidth=0.02cm](5.42,0.91)(5.42,0.81)
\psarc[linewidth=0.02](5.41,0.48){0.13}{258.69006}{0.0}
\psline[linewidth=0.02cm](5.22,0.35)(5.4,0.35)
\psline[linewidth=0.02cm](5.54,0.47)(5.54,0.51)
\psdots[dotsize=0.15](0.42,-0.49)
\psdots[dotsize=0.15](1.82,-0.49)
\psdots[dotsize=0.15](0.42,-2.47)
\psdots[dotsize=0.15](0.42,-0.49)
\psdots[dotsize=0.15](0.42,-2.47)
\psdots[dotsize=0.15](0.42,-1.47)
\psdots[dotsize=0.15](0.42,-1.47)
\psdots[dotsize=0.15](1.82,-0.49)
\psline[linewidth=0.02](0.4,-0.65)(0.66,-0.65)(0.66,-0.65)(0.66,-0.65)
\psline[linewidth=0.02cm](0.38,-2.27)(0.64,-2.27)
\psline[linewidth=0.02cm](0.74,-2.15)(0.74,-0.75)
\psline[linewidth=0.02cm](0.38,-2.65)(0.64,-2.65)
\psline[linewidth=0.02cm](0.86,-2.51)(0.86,-0.35)
\psline[linewidth=0.02cm](0.4,-0.25)(0.74,-0.25)
\psarc[linewidth=0.02](0.4,-0.45){0.2}{90.0}{270.0}
\psarc[linewidth=0.02](0.39,-2.46){0.19}{90.0}{270.0}
\psarc[linewidth=0.02](0.73,-0.38){0.13}{353.65982}{84.80557}
\psarc[linewidth=0.02](0.65,-0.74){0.09}{353.65982}{84.80557}
\psline[linewidth=0.02cm](0.74,-0.73)(0.74,-0.75)
\psarc[linewidth=0.02](0.65,-2.18){0.09}{258.69006}{9.462322}
\psline[linewidth=0.02cm](0.74,-2.09)(0.74,-2.19)
\psarc[linewidth=0.02](0.73,-2.52){0.13}{258.69006}{0.0}
\psline[linewidth=0.02cm](0.54,-2.65)(0.72,-2.65)
\psline[linewidth=0.02cm](0.86,-2.53)(0.86,-2.49)
\end{pspicture}
\caption{\sl The ten solutions of $\sigma\wedge\pi =\hat{0}$ in a three-row diagram}
%\end{center}
\end{figure}

\medskip

In what follows (see e.g. Theorem \ref{T : Diagrams} below), we will
sometimes be called to solve jointly the equations $\sigma \wedge \pi =\hat{0%
}$ and $\sigma \vee \pi =\hat{1}$, that is, given $\pi $, to find all
diagrams $\Gamma \left( \pi ,\sigma \right) $ that are non-flat ($\sigma
\wedge \pi =\hat{0}$) and connected ($\sigma \vee \pi =\hat{1}$). Having
found, as before, all those that are non-flat, one just has to choose among
them those that are connected, that is, the diagrams whose rows cannot be
divided into two subset, each defining a separate diagram. These are: the
second diagram in Fig. 11, the last two in Fig. 12, the last six in Fig. 13,
the sixth to ninth of Fig. 14. Again as an example, observe that the second
diagram in Fig. 14 is not connected: indeed, in this case, $\pi =\left\{
\left\{ 1,2\right\} ,\left\{ 3\right\} ,\left\{ 4\right\} \right\} $, $%
\sigma =\left\{ \left\{ 1\right\} ,\left\{ 2\right\} ,\left\{ 3,4\right\}
\right\} $, and $\pi \vee \sigma =\left\{ \left\{ 1,2\right\} ,\left\{
3,4\right\} \right\} \neq \hat{1}$.

\subsection{From Gaussian diagrams to multigraphs \label{ss : mgr}}

A \textsl{multigraph }is a graph in which (a) two vertices can be connected
by more than one edge, and (b) loops (that is, edges connecting one vertex
to itself) are allowed. Such objects are sometimes called \textquotedblleft
pseudographs\textquotedblright , but we will avoid this terminology. In what
follows, we show how a multigraph can be derived from a Gaussian diagram.
This representation of Gaussian diagrams can be used in the computation of
moments and cumulants (see \cite{Go Taqqu} or \cite{Marinucci}).

\bigskip

Fix a set $b$ and consider partitions $\pi ,\sigma \in \mathcal{P}\left(
b\right) $ such that $\Gamma \left( \pi ,\sigma \right) $ is Gaussian and $%
\pi =\left\{ b_{1},...,b_{k}\right\} $. Then, the multigraph $\hat{\Gamma}%
\left( \pi ,\sigma \right) $, with $k$ vertices and $\left\vert b\right\vert
/2$ edges, is obtained from $\Gamma \left( \pi ,\sigma \right) $ as follows.

\begin{enumerate}
\item Identify each row of $\Gamma \left( \pi ,\sigma \right) $ with a
vertex of $\hat{\Gamma}\left( \pi ,\sigma \right) $, in such a way that the $%
i$th row of $\Gamma \left( \pi ,\sigma \right) $ corresponds to the $i$th
vertex$\ v_{i}$ of $\hat{\Gamma}\left( \pi ,\sigma \right) $.

\item Draw an edge linking $v_{i}$ and $v_{j}$ for every pair $\left(
x,y\right) $ such that $x\in b_{i}$, $y\in b_{j}$ and $x\sim _{\sigma }y$.
\end{enumerate}

\bigskip

\textbf{Examples. }(i) The multigraph obtained from the diagram in Fig. 7 is given in Fig. 15.
\begin{figure}[htbp]
\begin{center}
\psset{unit=0.7cm}
\begin{pspicture}(0,-1.9976757)(4.01109,2.0023243)
\psframe[linewidth=0.02,dimen=outer](4.01109,2.0023243)(1.0110897,-1.9976757)
\psdots[dotsize=0.15](2.4510896,1.3823243)
\psdots[dotsize=0.15](2.4510896,0.02232436)
\psdots[dotsize=0.15](2.4310896,-1.3776757)
\psline[linewidth=0.02cm](2.4310896,1.3623244)(2.4110897,1.3623244)
\psline[linewidth=0.02cm](2.4510896,1.3823243)(2.4510896,0.02232436)
\rput{-66.40058}(1.1741974,1.7785621){\psarc[linewidth=0.02](1.9460475,-0.007889927){1.4779661}{355.56204}{135.0}}
\rput{-115.84482}(4.226715,2.5933304){\psarc[linewidth=0.02](2.9260476,-0.027889928){1.4779661}{225.0}{4.4379687}}
\end{pspicture}
\caption{\sl A multigraph with three vertices}
\end{center}
\end{figure}
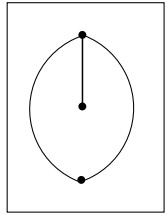

(ii) The multigraph associated with Fig. 8 is given in Fig. 16 (note that this graph has been
obtained from a circular diagram).

\begin{figure}[htbp]
\begin{center}
\psset{unit=0.7cm}
\begin{pspicture}(0,-2.0)(3.7414663,2.0)
\psframe[linewidth=0.02,dimen=outer](3.7414663,2.0)(0.7414663,-2.0)
\psdots[dotsize=0.15](2.1814663,1.28)
\psdots[dotsize=0.15](2.1814663,0.68)
\psdots[dotsize=0.15](2.1814663,0.1)
\psdots[dotsize=0.15](2.1814663,-0.5)
\psdots[dotsize=0.15](2.1814663,-1.1)
\psline[linewidth=0.02cm](2.1814663,1.26)(2.1814663,-1.08)
\rput{-68.76358}(0.9539096,1.6329645){\psarc[linewidth=0.02](1.6702104,0.11943224){1.2904245}{0.0}{138.17982}}
\end{pspicture}
\caption{\sl A multigraph built from a circular diagram}
\end{center}
\end{figure}
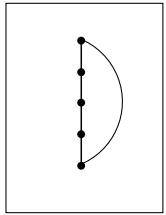

\bigskip

The following result is easily verified: it shows how the nature of a
Gaussian diagram can be deduced from its graph representation.

\bigskip

\begin{proposition}
\label{P : LOOP}Fix a finite set $b$, as well as a pair of partitions $%
\left( \pi ,\sigma \right) \subseteq \mathcal{P}\left( b\right) $ such that
the diagram $\Gamma \left( \pi ,\sigma \right) $ is Gaussian and $\left\vert
\pi \right\vert =k$. Then,

\begin{enumerate}
\item $\Gamma \left( \pi ,\sigma \right) $ is connected if and only if $%
\hat{\Gamma}\left( \pi ,\sigma \right) $ is a connected multigraph.

\item $\Gamma \left( \pi ,\sigma \right) $ is non-flat if and only if $%
\hat{\Gamma}\left( \pi ,\sigma \right) $ has no loops.

\item $\Gamma \left( \pi ,\sigma \right) $ is circular if and only if the
vertices $v_{1},...,v_{k}$ of $\hat{\Gamma}\left( \pi ,\sigma \right) $ are
such that: (i) there is an edge linking $v_{i}$ and $v_{i+1}$ for every $%
i=1,...,k-1$, and (ii) there is an edge linking $v_{k}$ and $v_{1}$.
\end{enumerate}
\end{proposition}

\newpage

As an illustration, in Fig. 17 we present the picture of a flat and non-connected
diagram (on the left), whose graph (on the right) is non-connected and
displays three loops.

\begin{figure}[htbp]
\begin{center}
\psset{unit=0.9cm}
\begin{pspicture}(0,-1.5)(4.78,1.5)
\psframe[linewidth=0.02,dimen=outer](2.4,1.5)(0.0,-1.5)
\psdots[dotsize=0.15](0.6,0.78)
\psdots[dotsize=0.15](1.2,0.78)
\psdots[dotsize=0.15](0.62,-0.02)
\psdots[dotsize=0.15](1.2,-0.02)
\psdots[dotsize=0.15](1.8,-0.02)
\psdots[dotsize=0.15](0.62,-0.8)
\psdots[dotsize=0.15](1.2,-0.8)
\psdots[dotsize=0.15](1.8,-0.8)
\psframe[linewidth=0.02,dimen=outer](4.78,1.5)(2.38,-1.5)
\psdots[dotsize=0.15](3.56,0.78)
\psdots[dotsize=0.15](3.58,-0.02)
\psdots[dotsize=0.15](3.58,-0.8)
\psline[linewidth=0.02cm](0.62,0.8)(1.18,0.8)
\psline[linewidth=0.02cm](0.62,0.0)(1.2,0.0)
\psline[linewidth=0.02cm](0.62,-0.78)(1.18,-0.78)
\psline[linewidth=0.02cm](1.78,0.02)(1.76,-0.02)
\psline[linewidth=0.02cm](1.8,0.0)(1.8,-0.76)
\psline[linewidth=0.02cm](3.56,-0.02)(3.54,0.0)
\psline[linewidth=0.02cm](3.58,0.0)(3.58,-0.76)
\pscircle[linewidth=0.02,dimen=outer](3.81,-0.79){0.23}
\pscircle[linewidth=0.02,dimen=outer](3.81,0.01){0.23}
\pscircle[linewidth=0.02,dimen=outer](3.79,0.79){0.23}
\end{pspicture}
\caption{\sl A non-connected flat diagram and its multigraph}
\end{center}
\end{figure}
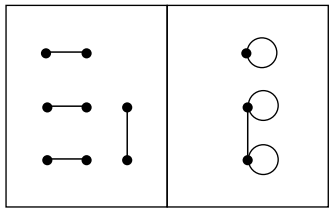
This situation corresponds to the case $b=\left[ 8\right] $,
\begin{eqnarray*}
\pi &=&\left\{ \left\{ 1,2\right\} ,\left\{ 3,4,5\right\} ,\left\{
6,7,8\right\} \right\} ,\text{ \ and} \\
\sigma &=&\left\{ \left\{ 1,2\right\} ,\left\{ 3,4\right\} ,\left\{
5,8\right\} ,\left\{ 6,7\right\} \right\} .
\end{eqnarray*}

\section{Completely random measures and Wiener-It\^{o} stochastic integrals
\label{S : MCA}}

\setcounter{equation}{0}We will now introduce the notion of a \textsl{%
completely random measure}\textit{\ }on a measurable space $\left( Z,%
\mathcal{Z}\right) $, as well as those of a \textsl{stochastic measure of
order}\textit{\ }$n\geq 2$, a \textsl{diagonal measure}\textit{\ }and a
\textsl{multiple (stochastic) Wiener-It\^{o} integral}. All these concepts
can be unified by means of the formalism introduced in Sections \ref{S :
Lattice}--\ref{S : DG}. We stress by now that the domain of the multiple
stochastic integrals defined in this section can be extended to more general
(and possibly random) classes of integrands. We refer the interested reader
to the paper by Kallenberg ans Szulga \cite{KaSz}, as well as to the monographs
by Kussmaul \cite{Kussmaul}, Kwapie\'{n} and Woyczy\'{n}ski \cite[Ch. 10]{KW}
and Linde \cite{WLindeBook}, for several results in this direction.

\subsection{Completely random measures\label{SS : CRM}}

Diagonals and subdiagonals play an important role in the context of multiple
integrals. The following definitions provide a convenient way to specify
them. In what follows, we will denote by $\left( Z,\mathcal{Z}\right) $ a
Polish space, where $\mathcal{Z}$ is the associated Borel $\sigma $-field.

\bigskip

\begin{definition}
\label{D : Diagonal sets}For every $n\geq 1$, we write $\left( Z^{n},%
\mathcal{Z}^{n}\right) =\left( Z^{\otimes n},\mathcal{Z}^{\otimes n}\right) $%
, with $Z^{1}=Z$. For every partition $\pi \in \mathcal{P}\left( \left[ n%
\right] \right) $ and every $B\in \mathcal{Z}^{n}$, we set%
\begin{equation}
\fbox{$Z_{\pi }^{n}\triangleq \left\{ \left( z_{1},...,z_{n}\right) \in
Z^{n}:z_{i}=z_{j}\text{ \ \ if and only if \ \ }i\sim _{\pi }j\right\} $}%
\text{ \ and\ \ }B_{\pi }\triangleq B\cap Z_{\pi }^{n}\text{.}
\label{base : Bipi}
\end{equation}
\end{definition}

\bigskip

Recall that $i\sim _{\pi }j$ means that the elements $i$ and $j$ belong to
the same block of the partition $\pi $. Relation (\ref{base : Bipi}) states
that the variables $z_{i}$ and $z_{j}$ should be equated if and only if $i$
and $j$ belong to the same block of $\pi $.

\bigskip

\textbf{Examples. }(i) Since $\hat{0}=\left\{ \left\{ 1\right\} ,...,\left\{
n\right\} \right\} $, no two elements can belong to the same block, and
therefore $B_{\hat{0}}$ coincides with the collection of all vectors $\left(
z_{1},...,z_{n}\right) \in B$ such that $z_{i}\neq z_{j}$, $\forall i\neq j$.

(ii) Since $\hat{1}=\left\{ \left\{ 1,...,n\right\} \right\} $, all elements
belong to the same block and therefore
\begin{equation*}
B_{\hat{1}}=\{\left(
z_{1},...,z_{n}\right) \in B: z_{1}=z_{2}=...=z_{n}\}.
\end{equation*}
A set such as $B_{\hat{1}}$ is said to be \textsl{purely diagonal}.

(iii) Suppose $n=3$ and $\pi =\left\{ \left\{ 1\right\} ,\left\{ 2,3\right\}
\right\} $. Then, $B_{\pi }=\{\left( z_{1},z_{2},z_{3}\right) \in
B:z_{2}=z_{3}$, $z_{1}\neq z_{2}\}$.

\bigskip

The following decomposition lemma (whose proof is immediate and left to the
reader) will be used a number of times.

\bigskip

\begin{lemma}
\label{L : decomposition}For every set $B\in \mathcal{Z}^{n}$,
\begin{equation*}
B=\cup _{\sigma \in \mathcal{P}\left( \left[ n\right] \right) }B_{\sigma
}=\cup _{\sigma \geq \hat{0}}B_{\sigma }\text{.}
\end{equation*}%
Moreover $B_{\pi }\cap B_{\sigma }=\varnothing $ if $\pi \neq \sigma $.
\end{lemma}

\bigskip

One has also that%
\begin{equation}
\left( A_{1}\times \cdot \cdot \cdot \times A_{n}\right) _{\hat{1}}=\underset%
{n\text{ times}}{(\underbrace{(\cap _{i=1}^{n}A_{i})\times \cdot \cdot \cdot
\times (\cap _{i=1}^{n}A_{i})})_{\hat{1}}}\text{;}  \label{Triv1}
\end{equation}%
indeed, since all coordinates are equal in the LHS of (\ref{Triv1}), their
common value must be contained in the intersection of the sets.

\bigskip

\textbf{Example. }As an illustration of (\ref{Triv1}), let $A_{1}=\left[ 0,1%
\right] $ and $A_{2}=\left[ 0,2\right] $ be intervals in $\mathbb{R}^{1}$,
and draw the rectangle $A_{1}\times A_{2}\in \mathbb{R}^{2}$. The set $%
\left( A_{1}\times A_{2}\right) _{\hat{1}}$ (that is, the subset of $%
A_{1}\times A_{2}$ composed of vectors whose coordinates are equal) is
therefore identical to the diagonal of the square $\left( A_{1}\cap
A_{2}\right) \times \left( A_{1}\cap A_{2}\right) =\left[ 0,1\right] ^{2}$.
The set $\left( A_{1}\times A_{2}\right) _{\hat{1}}$ can be visualized as
the thick diagonal segment in Fig. 18.

\begin{figure}[htbp]
\begin{center}
\psset{unit=0.8cm}
\begin{pspicture}(0,-1.97)(3.84,1.97)
\rput(1.004,-0.81){\psaxes[linewidth=0.02,ticksize=0.1358cm](0,0)(0,0)(2,2)}
\psline[linewidth=0.02cm](1.0,1.19)(2.0,1.19)
\psline[linewidth=0.02cm](2.0,1.19)(2.0,-0.69)
\psline[linewidth=0.02cm](1.0,0.19)(2.0,0.19)
\psframe[linewidth=0.02,dimen=outer](3.84,1.97)(0.0,-1.97)
\psline[linewidth=0.06cm](1.02,-0.79)(2.0,0.17)
\end{pspicture}
\caption{\sl \it A diagonal set}
\end{center}
\end{figure}
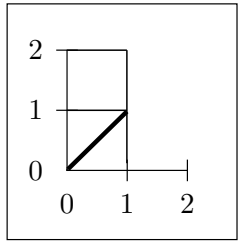

\bigskip

We shall now define a \textquotedblleft completely random
measure\textquotedblright\ $\varphi $, often called an \textquotedblleft
independently scattered random measure\textquotedblright . It has two
characteristics: it is a measure and it takes values in a space of random
variables. It will be denoted with its arguments as $\varphi \left( B,\omega
\right) $, where $B$ is a Borel set and $\omega $ is a point in the
underlying sample space $\Omega $. The \textquotedblleft
size\textquotedblright\ of $\varphi $ will be controlled by a non-random, $%
\sigma $-finite and non-atomic measure $\nu $, where $\nu \left( B\right) =%
\mathbb{E}\varphi \left( B\right) ^{2}$. The fact that $\nu $ is non-atomic
means that $\nu \left( \left\{ z\right\} \right) =0$ for every $z\in Z$. The
measure $\varphi $ will be used to define multiple integrals, where one
integrates either over a whole subset of $Z^{p}$, $p\geq 1$, or over a
subset \textquotedblleft without diagonals\textquotedblright . In the first
case we will need to suppose $\mathbb{E}\left\vert \varphi \left( B\right)
\right\vert ^{p}<\infty $; in the second case, one may suppose as little as $%
\mathbb{E}\varphi \left( B\right) ^{2}<\infty $. Since $p\geq 1$ will be
arbitrary, in order to deal with the first case we shall suppose that $%
\mathbb{E}\left\vert \varphi \left( B\right) \right\vert ^{p}<\infty $, $%
\forall p\geq 1$, that is, $\varphi \in \cap _{p\geq 1}L^{p}\left( \mathbb{P}%
\right) $. We now present the formal definition of $\varphi $.

\begin{definition}
\label{D : CRM}(1) Let $\nu $ be a positive, $\sigma $-finite and non-atomic
measure on $\left( Z,\mathcal{Z}\right) $, and let
\begin{equation}
\mathcal{Z}_{\nu }=\{B\in \mathcal{Z}:\nu \left( B\right) <\infty \}.
\label{ZetaNu}
\end{equation}%
A centered \textbf{completely random measure} (in $\cap _{p\geq
1}L^{p}\left( \mathbb{P}\right) $) on $\left( Z,\mathcal{Z}\right) $ with
\textbf{control measure} $\nu $, is a function $\varphi \left( \cdot ,\cdot
\right) $, from $\mathcal{Z}_{\nu }\times \Omega $ to $\mathbb{R}$, such that

\begin{description}
\item[(i)] For every fixed $B\in \mathcal{Z}_{\nu }$, the application $%
\omega \mapsto \varphi \left( B,\omega \right) $ is a random variable;

\item[(ii)] For every fixed $B\in \mathcal{Z}_{\nu }$, $\varphi \left(
B\right) \in \cap _{p\geq 1}L^{p}\left( \mathbb{P}\right) $;

\item[(iii)] For every fixed $B\in \mathcal{Z}_{\nu }$, $\mathbb{E}\left[
\varphi \left( B\right) \right] =0$;

\item[(iv)] $\varphi \left( \varnothing \right) =0$;

\item[(v)] For every collection of disjoint elements of $\mathcal{Z}_{\nu }$%
, $B_{1},...,B_{n}$, the variables $\varphi \left( B_{1}\right) ,$...$,$ $%
\varphi \left( B_{n}\right) $ are independent;

\item[(vi)] For every $B,C\in \mathcal{Z}_{\nu }$, $\mathbb{E}\left[ \varphi
\left( B\right) \varphi \left( C\right) \right] =\nu \left( B\cap C\right) .$
\end{description}

(2) When $\varphi \left( \cdot \right) $ verifies the properties (i) and
(iii)--(vi) above and $\varphi \left( B\right) \in L^{2}\left( \mathbb{P}%
\right) $, $\forall B\in \mathcal{Z}_{\nu }$ (so that it is not necessarily
true that $\varphi \left( B\right) \in L^{p}\left( \mathbb{P}\right) $, $%
p\geq 3$), we say that $\varphi $ is a \textbf{completely random measure in}
$L^{2}\left( \mathbb{P}\right) $.
\end{definition}

\bigskip

\textbf{Two crucial remarks. }(\textit{On additivity})\textbf{\ }(a) Let $%
B_{1},...,B_{n},...$ be a sequence of \textsl{disjoint} elements of $%
\mathcal{Z}_{\nu }$, and let $\varphi $ be a completely random measure on $%
\left( Z,\mathcal{Z}\right) $ with control $\nu $. Then, for every finite $%
N\geq 2$, one has that $\cup _{n=1}^{N}B_{n}\in \mathcal{Z}_{\nu }$, and, by
using Properties (iii), (v) and (vi) in Definition \ref{D : CRM}, one has
that
\begin{equation}
\mathbb{E}\left[ \left( \varphi \left( \cup _{n=1}^{N}B_{i}\right)
-\sum_{n=1}^{N}\varphi \left( B_{n}\right) \right) ^{2}\right] =\nu \left(
\cup _{n=1}^{N}B_{n}\right) -\sum_{n=1}^{N}\nu \left( B_{n}\right) =0\text{,}
\label{Crr1}
\end{equation}%
because $\nu $ is a measure, and therefore it is finitely additive. Relation
(\ref{Crr1}) implies in particular that%
\begin{equation}
\varphi \left( \cup _{n=1}^{N}B_{n}\right) =\sum_{n=1}^{N}\varphi \left(
B_{n}\right) \text{, \ \ a.s.-}\mathbb{P}\text{.}  \label{Crr2}
\end{equation}%
Now suppose that $\cup _{n=1}^{\infty }B_{n}\in \mathcal{Z}_{\nu }$. Then,
by (\ref{Crr2}) and again by virtue of Properties (iii), (v) and (vi) in
Definition \ref{D : CRM},
\begin{eqnarray*}
\mathbb{E}\left[ \left( \varphi \left( \cup _{n=1}^{\infty }B_{n}\right)
-\sum_{n=1}^{N}\varphi \left( B_{n}\right) \right) ^{2}\right]  &=&\mathbb{E}%
\left[ \left( \varphi \left( \cup _{n=1}^{\infty }B_{n}\right) -\varphi
\left( \cup _{n=1}^{N}B_{i}\right) \right) ^{2}\right]  \\
&=&\nu \left( \cup _{n=N+1}^{\infty }B_{n}\right) \underset{N\rightarrow
\infty }{\rightarrow }0\text{,}
\end{eqnarray*}%
because $\nu $ is $\sigma $-additive. This entails in turn that
\begin{equation}
\fbox{$\varphi \left( \cup _{n=1}^{\infty }B_{n}\right) =\sum_{n=1}^{\infty
}\varphi \left( B_{n}\right) $, \ \ a.s.-$\mathbb{P}$,}  \label{Crr3}
\end{equation}%
where the series on the RHS converges in $L^{2}\left( \mathbb{P}\right) $.
Relation (\ref{Crr3}) simply means that the application
\begin{equation*}
\mathcal{Z}_{\nu }\rightarrow L^{2}\left( \mathbb{P}\right) :B\mapsto
\varphi \left( B\right) \text{,}
\end{equation*}%
is $\sigma $-additive, and therefore that \textsl{every completely random
measure is a }$\sigma $-\textsl{additive measure with values in the Hilbert
space} $L^{2}\left( \mathbb{P}\right) $. See e.g. Engel \cite{Engel},
Kussmaul \cite{Kussmaul} or Linde \cite{WLindeBook} for further discussions
on vector-valued measures.

(b) In general, it is \textsl{not} true that, for a completely random
measure $\varphi $ and for a fixed $\omega \in \Omega $, the application%
\begin{equation*}
\mathcal{Z}_{\nu }\rightarrow \mathbb{R}:B\mapsto \varphi \left( B,\omega
\right) \text{ }
\end{equation*}%
is a $\sigma $-additive real-valued (signed) measure. The most remarkable
example of this phenomenon is given by Gaussian completely random measures.
See the discussion below for more details on this point.

\bigskip

\textbf{Remark on notation. }We consider the spaces $\left( Z,\mathcal{Z}%
\right) $ and $\left( Z^{n},\mathcal{Z}^{n}\right) =\left( Z^{\otimes n},%
\mathcal{Z}^{\otimes n}\right) $. Do not confuse the subset $Z_{\pi }^{n}$
in (\ref{base : Bipi}), where $\pi $ denotes a partition, with the $\sigma $%
-field $\mathcal{Z}_{\nu }^{n}$ in (\ref{ZetaNu}), where $\nu $ denotes a
control measure.

\bigskip

Now fix a completely random measure $\varphi $. For every $n\geq 2$ and
every rectangle $C=C_{1}\times \cdot \cdot \cdot \times C_{n}$, $C_{j}\in
\mathcal{Z}_{\nu }$, we define $\varphi ^{\left[ n\right] }\left( C\right)
\triangleq \varphi \left( C_{1}\right) \times \cdot \cdot \cdot \times
\varphi \left( C_{n}\right) $, so that the application $C\mapsto \varphi ^{%
\left[ n\right] }\left( C\right) $ defines a finitely additive application
on the ring of rectangular sets contained in $Z^{n}$, with values in the set
of $\sigma \left( \varphi \right) $-measurable random variables. In the next
definition we focus on those completely random measures such that the
application $\varphi ^{\left[ n\right] }$ admits a unique infinitely
additive (and square-integrable) extension on $\mathcal{Z}^{n}$. Here, the
infinite additivity is in the sense of the $L^{1}\left( \mathbb{P}\right) $
convergence. Note that we write $\varphi ^{\left[ n\right] }$ to emphasize
the dependence of $\varphi ^{\left[ n\right] }$ not only on $n$, but also on
the set $\left[ n\right] =\left\{ 1,...,n\right\} $, whose lattice of
partitions will be considered later.

\begin{definition}
\label{D : good}For $n\geq 2$, we write $\mathcal{Z}_{\nu }^{n}=\{C\in
\mathcal{Z}^{n}:$ $\nu ^{n}\left( C\right) <\infty \}$. A completely random
measure $\varphi $, verifying points (i)--(vi) of Definition \ref{D : CRM},
is said to be \textbf{good }if, for every fixed $n\geq 2$, there exists a
(unique) collection of random variables $\varphi ^{\left[ n\right] }=\left\{
\varphi ^{\left[ n\right] }\left( C\right) :C\in \mathcal{Z}^{n}\right\} $
such that

\begin{description}
\item[(i)] $\left\{ \varphi ^{\left[ n\right] }\left( C\right) :C\in
\mathcal{Z}_{\nu }^{n}\right\} \subseteq L^{2}\left( \mathbb{P}\right) $ ;

\item[(ii)] For every rectangle $C=C_{1}\times \cdot \cdot \cdot \times
C_{n} $, $C_{j}\in \mathcal{Z}_{\nu }$,%
\begin{equation}
\fbox{$\varphi ^{\left[ n\right] }\left( C\right) =\varphi \left(
C_{1}\right) \cdot \cdot \cdot \varphi \left( C_{n}\right) ;$}
\label{prodEng}
\end{equation}

\item[(iii)] $\varphi ^{\left[ n\right] }$ is a $\sigma $-additive random
measure in the following sense: if $C\in \mathcal{Z}_{\nu }^{n}$ is such
that $C=\cup _{j=1}^{\infty }C_{j}$, with the $\left\{ C_{j}\right\} $
disjoints, then
\begin{equation}
\varphi ^{\left[ n\right] }\left( C\right) =\sum_{j=1}^{\infty }\varphi ^{%
\left[ n\right] }\left( C_{j}\right) \text{, \ \ with convergence (at least)
in }L^{1}\left( \mathbb{P}\right) \text{.}  \label{sigmaEng}
\end{equation}
\end{description}
\end{definition}

\bigskip

Note that, in the case $n=1$, the $\sigma $-additivity of $\varphi $ follows
from point (vi) of Definition \ref{D : CRM} and the $\sigma $-additivity of $%
\nu $ (through $L^{2}$ convergence). In the case $n\geq 2$, the assumption
that the measure $\varphi $ is good\ implies $\sigma $-additivity in the
sense of (\ref{sigmaEng}).

\bigskip

\textbf{Remark. }The notion of a \textquotedblleft completely random
measure\textquotedblright\ can be traced back to Kingman's seminal paper
\cite{kingman67}. For further references on completely random measures, see
also the two surveys by Surgailis \cite{Sur} and \cite{SUr2} (note that, in
such references, completely random measures are called \textquotedblleft
independently scattered measures\textquotedblright ). The use of the term
\textquotedblleft good\textquotedblright , to indicate completely random
measures satisfying the requirements of Definition \ref{D : good}, is taken
from Rota and Wallstrom \cite{RoWa}. Existence of good measures is discussed
in Engel \cite{Engel} and Kwapie\'{n} and Woyczy\'{n}ski \cite[Ch. 10]{KW}.
For further generalizations of Engel's results the reader is referred e.g.
to \cite{KaSz} and \cite{RosWoy}.

\bigskip

\textbf{Examples.\ }\label{Ex. GaussPoiss}The following two examples of good
completely random measures will play a crucial role in the subsequent
sections.

(i) A \textsl{centered Gaussian random measure} with control $\nu $ is a
collection $G=\left\{ G\left( B\right) :B\in \mathcal{Z}_{\nu }\right\} $ of
jointly Gaussian random variables, centered and such that, for every $B,C\in
\mathcal{Z}_{\nu }$, $\mathbb{E}\left[ G\left( C\right) G\left( B\right) %
\right] $ $=$ $\nu \left( C\cap B\right) $. The family $G$ is clearly a
completely random measure. The fact that $G$ is also good is classic, and
can be seen as a special case of the main results in \cite{Engel}.

(ii) A \textsl{compensated Poisson measure with control} $\nu $ is a
completely random measure $\hat{N}=\{\hat{N}\left( B\right) :$ $B\in
\mathcal{Z}_{\nu }\}$, as in Definition \ref{D : CRM}, such that, $\forall
B\in \mathcal{Z}_{\nu }$, $\hat{N}\left( B\right) \overset{law}{=}N\left(
B\right) -\nu \left( B\right) $, where $N\left( B\right) $ is a Poisson
random variable with parameter $\nu \left( B\right) =\mathbb{E}N\left(
B\right) =\mathbb{E}N\left( B\right) ^{2}$. The fact that $\hat{N}$ is also
good derives once again from the main findings of \cite{Engel}. A more
direct proof of this last fact can be obtained by observing that, for almost
every $\omega $, $\hat{N}^{\left[ n\right] }\left( \cdot ,\omega \right) $
must necessarily coincide with the canonical product (signed) measure (on $%
\left( Z^{n},\mathcal{Z}^{n}\right) $) associated with the signed measure on
$\left( Z,\mathcal{Z}\right) $ given by $\hat{N}\left( \cdot ,\omega \right)
=N\left( \cdot ,\omega \right) -\nu \left( \cdot \right) $ (indeed, such a
canonical product measure satisfies necessarily (\ref{prodEng})). Note that
a direct proof of this type cannot be obtained in the Gaussian case. Indeed,
if $G$ is a Gaussian measure as in Point (i), one has that, for almost every
$\omega $, the mapping $B\mapsto G\left( B,\omega \right) $ \textsl{does not
define }a signed measure (see e.g. \cite[Ch. 1]{Janson})\textsl{.}

\subsection{Single integrals and infinite
divisibility\label{SS : infDIV}}

Let $\varphi $ be a completely random measure in the sense of Definition \ref%
{D : CRM}, with control measure $\nu $. Our aim in this section is twofolds:
(i) we shall define (single) Wiener-It\^{o} integrals with respect to $%
\varphi $, and (ii) we shall give a characterization of these integrals as
infinitely divisible random variables.

The fact that single Wiener-It\^{o} integrals are infinitely divisible
should not come as a surprise. Indeed, observe that, since $(Z,\mathcal{Z})$ is a Polish space and $\nu$ is non-atomic, the law of any random
variable of the type $\varphi \left( B\right) $, $B\in \mathcal{Z}_{\nu }$,
is infinitely divisible. Infinitely divisible laws are introduced in many textbooks, see e.g. Billingsley \cite{BillBook}. In particular, for every $B\in \mathcal{Z}_{\nu }$
there exists a unique pair $\left( c^{2}\left( B\right) ,\alpha _{B}\right) $
such that $c^{2}\left( B\right) \in \left[ 0,\infty \right) $ and $\alpha
_{B}$ is a measure on $\mathbb{R}$ satisfying
\begin{equation}
\alpha _{B}\left( \left\{ 0\right\} \right) =0\text{ \ \ and\ \ \ }\int_{%
\mathbb{R}}u^{2}\alpha _{B}\left( du\right) <\infty ,  \label{LK0}
\end{equation}%
and, for every $\lambda \in \mathbb{R}$,
\begin{equation}
\mathbb{E}\left[ \exp \left( i\lambda \varphi \left( B\right) \right) \right]
=\exp \left[ -\frac{c^{2}\left( B\right) \lambda ^{2}}{2}+\int_{\mathbb{R}%
}\left( \exp \left( i\lambda u\right) -1-i\lambda u\right) \alpha _{B}\left(
du\right) \right] .  \label{L-K}
\end{equation}%
The measure $\alpha_B$ is called a \textsl{L\'{e}vy measure}, and the components of the pair $(c^2(B), \alpha_B)$ are called the \textsl{L\'{e}vy-Khintchine exponent characteristics} associated with $\varphi \left( B\right) $. Also, the exponent
\begin{equation*}
-\frac{c^{2}\left( B\right) \lambda ^{2}}{2}+\int_{\mathbb{R}}\left( \exp
\left( i\lambda u\right) -1-i\lambda u\right) \alpha _{B}\left( du\right)
\end{equation*}%
is known as the \textsl{L\'{e}vy-Khintchine exponent}\textit{\ }associated
with $\varphi \left( B\right) $. Plainly, if $\varphi $ is Gaussian, then $%
\alpha _{B}=0$ for every $B\in \mathcal{Z}_{\nu }$ (the reader is referred
e.g. to \cite{Sato} for an exhaustive discussion of infinitely divisible
laws).

\bigskip

We now establish the existence of single Wiener-It\^{o} integrals with
respect to a completely random measure $\varphi $.

\begin{proposition}
\label{P : single MWI}Let $\varphi $ be a completely random measure in $%
L^{2}\left( \mathbb{P}\right) $, with $\sigma $-finite control measure $\nu $%
. Then, there exists a unique continuous linear operator $h\mapsto \varphi
\left( h\right) $, from $L^{2}\left( \nu \right) $ into $L^{2}\left( \mathbb{%
P}\right) $, such that
\begin{equation}
\varphi \left( h\right) =\sum_{j=1}^{m}c_{j}\varphi \left( B_{j}\right)
\label{mwww1}
\end{equation}%
for every elementary function of the type
\begin{equation}
h\left( z\right) =\sum_{j=1}^{m}c_{j}\mathbf{1}_{B_{j}}\left( z\right) ,
\label{mwww2}
\end{equation}%
where $c_{j}\in \mathbb{R}$ and the sets $B_{j}$ are in $\mathcal{Z}_{\nu }$
and disjoint.
\end{proposition}

\begin{proof}
In what follows, we call \textit{simple kernel} a kernel $h$ as in (\ref%
{mwww2}). For every simple kernel $h$, set $\varphi \left( h\right) $ to be
equal to (\ref{mwww1}). Then, by using Properties (iii), (v) and (vi) in
Definition \ref{D : CRM}, one has that, for every pair of simple kernels $%
h,h^{\prime }$,%
\begin{equation}\label{NickDrake}
\mathbb{E}\left[ \varphi \left( h\right) \varphi \left( h^{\prime }\right) %
\right] =\int_{Z}h\left( z\right) h^{\prime }\left( z\right) \nu \left(
dz\right) \text{.}
\end{equation}%
Since simple kernels are dense in $L^{2}\left( \nu \right) $, the proof is
completed by the following (standard) approximation argument. If $h\in L^2(\nu)$ and $\{h_n\}$
is a sequence of simple kernels converging to $h$,
then (\ref{NickDrake}) implies that $\{\varphi(h_n)\}$ is a Cauchy sequence in $L^2(\mathbb{P})$,
and one defines $\varphi(h)$ to be the $L^2(\mathbb{P})$ limit of $\varphi(h_n)$.
One easily verifies that the definition of $\varphi(h)$ does not depend on the chosen approximating sequence $\{h_n\}$. The
application $h\mapsto\varphi(h)$ is therefore well-defined, and
(by virtue of (\ref{NickDrake})) it is an isomorphism from $L^{2}\left( \nu \right) $ into $L^{2}\left( \mathbb{%
P}\right) $.
\end{proof}

\bigskip

The random variable $\varphi \left( h\right) $ is usually written as%
\begin{equation}
\int_{Z}h\left( z\right) \varphi \left( dz\right) \text{, \ \ }%
\int_{Z}hd\varphi \text{ \ \ or \ \ }I_{1}^{\varphi }\left( h\right) \text{,}
\label{WiIto1}
\end{equation}%
and it is called the \textsl{Wiener-It\^{o} stochastic integral} of $h$ with
respect to $\varphi $. By inspection of the previous proof, one sees that
Wiener-It\^{o} integrals verify the isometric relation%
 \begin{equation}
  \fbox{$\mathbb{E}\left[ \varphi \left( h\right) \varphi \left( g\right) \right]
=\int_{Z}h\left( z\right) g\left( z\right) \nu \left( dz\right) =\left(
g,h\right) _{L^{2}\left( \nu \right) }\text{, \ \ }\forall g,h\in
L^{2}\left( \nu \right) \text{.}$}   \label{ISO1}
\end{equation}
Observe also that \fbox{$\mathbb{E}\varphi \left( h\right) =0$.} If $B\in \mathcal{Z%
}_{\nu }$, we write interchangeably $\varphi \left( B\right) $ or $\varphi
\left( \mathbf{1}_{B}\right) $ (the two objects coincide, thanks to (\ref%
{mwww1})). For every $h\in L^{2}\left( \nu \right) $, the law of the random
variable $\varphi \left( h\right) $ is also infinitely divisible. The
following result provides a description of the L\'{e}vy-Khintchine exponent
of $\varphi \left( h\right) $. The proof is taken from \cite{PeTaqMultAOP}
and uses arguments and techniques developed in \cite{Raj Ros} (see also \cite%
[Section 5]{KW91}). Following the proof, we present an interpretation of the result.

\bigskip

\begin{proposition}
\label{P : LK}For every $B\in \mathcal{Z}_{\nu }$, let $\left( c^{2}\left(
B\right) ,\alpha _{B}\right) $ denote the pair such that $c^{2}\left(
B\right) \in \left[ 0,\infty \right) $, $\alpha _{B}$ verifies (\ref{LK0})
and%
\begin{equation}
\mathbb{E}\left[ \exp \left( i\lambda \varphi \left( B\right) \right) \right]
=\exp \left[ -\frac{c^{2}\left( B\right) \lambda ^{2}}{2}+\int_{\mathbb{R}%
}\left( \exp \left( i\lambda x\right) -1-i\lambda x\right) \alpha _{B}\left(
dx\right) \right] .  \label{LK4}
\end{equation}%
Then, the following holds

\begin{enumerate}
\item The application $B\mapsto c^{2}\left( B\right) $, from $\mathcal{Z}%
_{\nu }$ to $\left[ 0,\infty \right) $, extends to a unique $\sigma $%
-finite measure $c^{2}\left( dz\right) $ on $\left( Z,\mathcal{Z}\right) $,
such that $c^{2}\left( dz\right) \ll \nu \left( dz\right) .$

\item There exists a unique measure $\alpha $ on $\left( Z\times \mathbb{R},%
\mathcal{Z}\times \mathcal{B}\left( \mathbb{R}\right) \right) $ such that $%
\alpha \left( B\times C\right) =\alpha _{B}\left( C\right) $, for every $%
B\in \mathcal{Z}_{\nu }$ and $C\in \mathcal{B}\left( \mathbb{R}\right) $.

\item There exists a function $\rho _{\nu }:Z\times \mathcal{B}\left(
\mathbb{R}\right) \mapsto \left[ 0,\infty \right] $ such that (i) for every
$z\in Z$, $\rho _{\nu }\left( z,\cdot \right) $ is a L\'{e}vy measure%
\footnote{%
That is, $\rho _{\nu }\left( z,\left\{ 0\right\} \right) =0$ and $\int_{%
\mathbb{R}}\min \left( 1,x^{2}\right) \rho _{\nu }\left( z,dx\right)
<\infty $.} on $\left( \mathbb{R},\mathcal{B}\left( \mathbb{R}\right)
\right) $ satisfying $\int_{Z}x^{2}\rho _{\nu }\left( z,dx\right) <\infty $%
, (ii) for every $C\in \mathcal{B}\left( \mathbb{R}\right) $, $\rho _{\nu
}\left( \cdot ,C\right) $ is a Borel measurable function, (iii) for every
positive function $g\left( z,x\right) \in \mathcal{Z}\otimes \mathcal{B}%
\left( \mathbb{R}\right) $,
\begin{equation}
\int_{Z}\int_{\mathbb{R}}g\left( z,x\right) \rho _{\nu }\left( z,dx\right)
\nu \left( dz\right) =\int_{Z}\int_{\mathbb{R}}g\left( z,x\right) \alpha
\left( dz,dx\right) .  \label{LKdensity}
\end{equation}

\item For every $\left( \lambda ,z\right) \in \mathbb{R}\times Z$, define
\begin{equation}
K_{\nu }\left( \lambda ,z\right) =-\frac{\lambda ^{2}}{2}\sigma _{\nu
}^{2}\left( z\right) +\int_{\mathbb{R}}\left( e^{i\lambda x}-1-i\lambda
x\right) \rho _{\nu }\left( z,dx\right) \text{,}  \label{innerLK}
\end{equation}%
where $\sigma _{\nu }^{2}\left( z\right) =\frac{dc^{2}}{d\nu }\left(
z\right) $; then, for every $h\in L^{2}\left( \nu \right) $, $%
\int_{Z}\left\vert K_{\nu }\left( \lambda h\left( z\right) ,z\right)
\right\vert \nu \left( dz\right) <\infty $ and
\begin{eqnarray}
&&\mathbb{E}\left[ \exp \left( i\lambda \varphi \left( h\right) \right) %
\right]  \label{LK5*} \\
&=&\exp \left[ \int_{Z}K_{\nu }\left( \lambda h\left( z\right) ,z\right) \nu
\left( dz\right) \right]  \notag \\
&=&\exp \left[ -\frac{\lambda ^{2}}{2}\int_{Z}h^{2}\left( z\right) \sigma
_{\nu }^{2}\left( z\right) \nu \left( dz\right) +\int_{Z}\int_{\mathbb{R}%
}\left( e^{i\lambda h\left( z\right) x}-1-i\lambda h\left( z\right) x\right)
\rho _{\nu }\left( z,dx\right) \nu \left( dz\right) \right] .  \notag
\end{eqnarray}
\end{enumerate}
\end{proposition}

\begin{proof}
The proof follows from results contained in \cite[Section II]{Raj Ros}.
Point 1 is indeed a direct consequence of \cite[Proposition 2.1 (a)]{Raj Ros}%
. In particular, whenever $B\in \mathcal{Z}$ is such that $\nu \left(
B\right) =0$, then $\mathbb{E[}\varphi \left( B\right) ^{2}]=0$ (due to
Point (vi) of Definition \ref{D : CRM}) and therefore $c^{2}\left( B\right)
=0$, thus implying $c^{2}\ll \nu $. Point 2 follows from the first part of
the statement of \cite[Lemma 2.3]{Raj Ros}. To establish Point 3 define, as
in \cite[p. 456]{Raj Ros},%
\begin{equation*}
\gamma \left( B\right) =c^{2}\left( B\right) +\int_{\mathbb{R}}\min \left(
1,x^{2}\right) \alpha _{B}\left( dx\right) =c^{2}\left( B\right) +\int_{\mathbb{R}}\min \left(
1,x^{2}\right) \alpha\left(B, dx\right)  \text{,}
\end{equation*}%
whenever $B\in \mathcal{Z}_{\nu }$, and observe (see \cite[Definition 2.2]%
{Raj Ros}) that $\gamma \left( \cdot \right) $ can be canonically extended
to a $\sigma $-finite and positive measure on $\left( Z,\mathcal{Z}\right) $%
. Moreover, since $\nu \left( B\right) =0$ implies $\varphi \left( B\right)
=0$ \ a.s.-$\mathbb{P}$, the uniqueness of the L\'{e}vy-Khinchine
characteristics implies as before $\gamma \left( B\right) =0$, and therefore
$\gamma \left( dz\right) \ll \nu \left( dz\right) $. Observe also that, by
standard arguments, one can select a version of the density $\left( d\gamma
/d\nu \right) \left( z\right) $ such that $\left( d\gamma /d\nu \right)
\left( z\right) <\infty $ for every $z\in Z$. According to \cite[Lemma 2.3]%
{Raj Ros}, there exists a function $\rho :Z\times \mathcal{B}\left( \mathbb{R%
}\right) \mapsto \left[ 0,\infty \right] $, such that: (a) $\rho \left(
z,\cdot \right) $ is a L\'{e}vy measure on $\mathcal{B}\left( \mathbb{R}%
\right) $ for every $z\in Z$, (b) $\rho \left( \cdot ,C\right) $ is a Borel
measurable function for every $C\in \mathcal{B}\left( \mathbb{R}\right) $,
(c) for every positive function $g\left( z,x\right) \in \mathcal{Z}\otimes
\mathcal{B}\left( \mathbb{R}\right) $,
\begin{equation}
\int_{Z}\int_{\mathbb{R}}g\left( z,x\right) \rho \left( z,dx\right) \gamma
\left( dz\right) =\int_{Z}\int_{\mathbb{R}}g\left( z,x\right) \alpha \left(
dz,dx\right) .  \label{ffpp}
\end{equation}%
In particular, by using (\ref{ffpp}) in the case $g\left( z,x\right) =%
\mathbf{1}_{A}\left( z\right) x^{2}$ for $A\in \mathcal{Z}_{\mu }$,%
\begin{equation*}
\int_{A}\int_{\mathbb{R}}x^{2}\rho \left( z,dx\right) \gamma \left(
dz\right) =\int_{\mathbb{R}}x^{2}\alpha _{A}\left( dx\right) <\infty \text{,%
}
\end{equation*}%
since $\varphi \left( A\right) \in L^{2}\left( \mathbb{P}\right) $, and we
deduce that $\rho $ can be chosen in such a way that, for every $z\in Z$, $%
\int_{\mathbb{R}}x^{2}\rho \left( z,dx\right) <\infty $. Now define, for
every $z\in Z$ and $C\in \mathcal{B}\left( \mathbb{R}\right) $,
\begin{equation*}
\rho _{\nu }\left( z,C\right) =\frac{d\gamma }{d\nu }\left( z\right) \rho
\left( z,C\right) \text{,}
\end{equation*}%
and observe that, due to the previous discussion, the application $\rho
_{\nu }:Z\times \mathcal{B}\left( \mathbb{R}\right) \mapsto \left[ 0,\infty %
\right] $ trivially satisfies properties (i)-(iii) in the statement of Point
3, which is therefore proved. To prove Point 4, first define (as before) a
function $h\in L^{2}\left( \nu \right) $ to be \textit{simple }if $h\left(
z\right) =\sum_{i=1}^{n}a_{i}\mathbf{1}_{A_{i}}\left( z\right) $, where $%
a_{i}\in \mathbb{R}$, and $\left( A_{1},...,A_{n}\right) $ is a finite
collection of disjoint elements of $\mathcal{Z}_{\nu }$. Of course, the
class of simple functions (which is a linear space) is dense in $L^{2}\left(
\nu \right) $, and therefore for every $L^{2}\left( \nu \right) $ there
exists a sequence $h_{n}$, $n\geq 1$, of simple functions such that $%
\int_{Z}\left( h_{n}\left( z\right) -h\left( z\right) \right) ^{2}\nu \left(
dz\right) \rightarrow 0$. As a consequence, since $\nu $ is $\sigma $-finite
there exists a subsequence $n_{k}$ such that $h_{n_{k}}\left( z\right)
\rightarrow h\left( z\right) $ for $\nu $-a.e. $z\in Z$ (and therefore for $%
\gamma $-a.e. $z\in Z$) and moreover, for every $A\in \mathcal{Z}$, the
random sequence $\varphi \left( \mathbf{1}_{A}h_{n}\right) $ is a Cauchy
sequence in $L^{2}\left( \mathbb{P}\right) $, and hence it converges in
probability. In the terminology of \cite[p. 460]{Raj Ros}, this implies that
every $h\in L^{2}\left( \nu \right) $ is $\varphi $-integrable, and that,
for every $A\in \mathcal{Z}$, the random variable $\varphi \left( h\mathbf{1}%
_{A}\right) $, defined according to Proposition \ref{P : single MWI},
coincides with $\int_{A}h\left( z\right) \varphi \left( dz\right) $, i.e.
the integral of $h$ with respect to the restriction of $\varphi \left( \cdot
\right) $ to $A$, as defined in \cite[p. 460]{Raj Ros}. As a consequence, by
using a slight modification of \cite[Proposition 2.6]{Raj Ros}\footnote{%
The difference lies in the choice of the truncation.}, the function $K_{0}$
on $\mathbb{R}\times Z$ given by
\begin{equation*}
K_{0}\left( \lambda ,z\right) =-\frac{\lambda ^{2}}{2}\sigma _{0}^{2}\left(
z\right) +\int_{\mathbb{R}}\left( e^{i\lambda x}-1-i\lambda x\right) \rho
\left( z,dx\right) \text{,}
\end{equation*}%
where $\sigma _{0}^{2}\left( z\right) =\left( dc^{2}/d\gamma \right) \left(
z\right) $, is such that $\int_{Z}\left\vert K_{0}\left( \lambda h\left(
z\right) ,z\right) \right\vert \gamma \left( dz\right) <\infty $ for every $%
h\in L^{2}\left( \nu \right) $, and also
\begin{equation*}
\mathbb{E}\left[ \exp \left( i\lambda \varphi \left( h\right) \right) \right]
=\int_{Z}K_{0}\left( \lambda h\left( z\right) ,z\right) \gamma \left(
dz\right) .
\end{equation*}%
The fact that, by definition, $K_\nu$ in (\ref{innerLK}) verifies
\begin{equation*}
K_{\nu }\left( \lambda h\left( z\right) ,z\right) =K_{0}\left( \lambda
h\left( z\right) ,z\right) \frac{d\gamma }{d\nu }\left( z\right) \text{, \ \
}\forall z\in Z\text{, }\forall h\in L^{2}\left( \nu \right) \text{, }%
\forall \lambda \in \mathbb{R}\text{,}
\end{equation*}%
yields (\ref{LK5*}).
\end{proof}

\bigskip

{\bf Interpretation of Proposition \ref{P : LK}.} Let $B$ be a given set in $\mathcal{Z}_\nu$. The characteristic function of the random variable $\varphi(B)$ or $\varphi({\bf 1}_B)$ involves the L\'{e}vy characteristic  $(c^2(B),\alpha_B)$, where $c^2(B)$ is a non-negative constant and $\alpha_B(dx)$ is a L\'{e}vy measure on $\mathbb{R}$. We want now to view $B\in\mathcal{Z}_\nu$ as a ``variable'' and thus to extend $c^2(B)$ to a measure $c^2(dz)$ on $(Z,\mathcal{Z})$, and $\alpha_B(dx)$ to a measure $\alpha(dz,dx)$ on $\mathcal{Z}\otimes \mathcal{B}(\mathbb{R})$. Consider first $\alpha_B$. According to Proposition \ref{P : LK}, it is possible to extend it to a measure $\alpha(dz,dx)$ on $\mathcal{Z}\otimes \mathcal{B}(\mathbb{R})$, which can be expressed as
\begin{equation}\label{e : RR-ID}
\alpha(dz,dx) = \rho_\nu (z,dx)\nu(dz),
\end{equation}
where $\rho_\nu$ is a function on $Z\times \mathbb{R}$, with the property that $\rho_\nu(z,\cdot)$ is a L\'{e}vy measure for every $z\in Z$. In view of (\ref{e : RR-ID}), the measure $\alpha(dz,dx)$ is thus obtained as a ``mixture'' of the L\'{e}vy measures $\rho_\nu (z,\cdot)$ over the variable $z$, using the control measure $\nu$ as a mixing measure.
A similar approach is applied to the Gaussian part of the exponent in (\ref{LK4}), involving $c^2(B)$. The coefficient $c^2(B)$ can be extended to a measure $c^2(dz)$, and this measure can be moreover expressed as
\begin{equation}\label{e : RR-G}
c^2(dz) = \sigma^2_\nu(z) \nu(dz),
\end{equation}
where $\sigma^2_\nu$ is the density of $c^2$ with respect to $\nu$. This allows the to represent the characteristic function of the Wiener-It\^{o} integral $\varphi(h)$ as in (\ref{LK5*}). In that expression, the function $h(z)$ appears explicitly in the L\'{e}vy-Khinchine exponent as a factor to the argument $\lambda$ of the characteristic function.

\bigskip

\textbf{Examples. }(i) If $\varphi =G$ is a centered Gaussian measure with
control measure $\nu $, then $\alpha =0$ and $c^2 = \nu$ (therefore $\sigma^2_\nu =1$) and, for $h\in L^{2}\left( \nu
\right) $,%
\begin{equation*}
\mathbb{E}\left[ \exp \left( i\lambda G\left( h\right) \right) \right] =\exp %
\left[ -\frac{\lambda ^{2}}{2}\int_{Z}h^{2}\left( z\right) \nu \left(
dz\right) \right] .
\end{equation*}

(ii) If $\varphi =\hat{N}$ is a compensated Poisson measure with control
measure $\nu $, then $c^{2}\left( \cdot \right) =0$ and $\rho _{\nu }\left(
z,dx\right) =\delta _{1}\left( dx\right) $ for all $z\in Z$, where $\delta
_{1}$ is the Dirac mass at $x=1$. It follows that, for $h\in L^{2}\left( \nu
\right) $,%
\begin{equation*}
\mathbb{E}\left[ \exp \left( i\lambda \hat{N}\left( h\right) \right) \right]
=\int_{Z}\left( e^{i\lambda h\left( z\right) }-1-i\lambda h\left( z\right)
\right) \nu \left( dz\right) .
\end{equation*}

(iii) Let $\left( Z,\mathcal{Z}\right) $ be a measurable space, and let $%
\hat{N}$ be a centered Poisson random measure on $\mathbb{R\times }Z$
(endowed with the usual product $\sigma $-field) with $\sigma $-finite
control measure $\nu \left( du,dz\right) $. Define the measure $\mu $ on $%
\left( Z,\mathcal{Z}\right) $ by%
\begin{equation*}
\mu \left( B\right) =\int_{\mathbb{R}}\int_{Z}u^{2}\mathbf{1}_{B}\left(
z\right) \nu \left( du,dz\right) \text{.}
\end{equation*}%
Then, by setting $k_{B}\left( u,z\right) =u\mathbf{1}_{B}\left( z\right) $,
the mapping
\begin{equation}
B\mapsto \varphi \left( B\right) =\int_{\mathbb{R}}\int_{Z}k_{B}\left(
u,z\right) \hat{N}\left( du,dz\right) =\int_{\mathbb{R}}\int_{Z}u\mathbf{1}%
_{B}\left( z\right) \hat{N}\left( du,dz\right) ,  \label{cool}
\end{equation}%
where $B\in \mathcal{Z}_{\mu }=\left\{ B\in \mathcal{Z}:\mu \left( B\right)
<\infty \right\} $, is a completely random measure on $\left( Z,\mathcal{Z}%
\right) $, with control measure $\mu $. In particular, by setting $%
k_{B}\left( u,z\right) =u\mathbf{1}_{B}\left( z\right) $, one has that
\begin{eqnarray}
\mathbb{E}\left[ \exp \left( i\lambda \varphi \left( B\right) \right) \right]
&=&\mathbb{E}\left[ \exp \left( i\lambda \hat{N}\left( k_{B}\right) \right) %
\right]   \notag \\
&=&\exp \left[ \int_{\mathbb{R}}\int_{Z}\left( e^{i\lambda k_{B}\left(
u,z\right) }-1-i\lambda k_{B}\left( u,z\right) \right) \nu \left(
du,dz\right) \right]   \notag \\
&=&\exp \left[ \int_{\mathbb{R}}\int_{Z}\left( e^{i\lambda u\mathbf{1}%
_{B}\left( z\right) }-1-i\lambda u\mathbf{1}_{B}\left( z\right) \right) \nu
\left( du,dz\right) \right]   \notag \\
&=&\exp \left[ \int_{\mathbb{R}}\int_{Z}\left( e^{i\lambda u}-1-i\lambda
u\right) \mathbf{1}_{B}\left( z\right) \nu \left( du,dz\right) \right]
\notag \\
&=&\exp \left[ \int_{\mathbb{R}}\left( e^{i\lambda u}-1-i\lambda u\right)
\alpha _{B}\left( du\right) \right] \text{,}  \label{cool2}
\end{eqnarray}%
where $\alpha _{B}\left( du\right) =\int_{Z}\mathbf{1}_{B}\left( z\right)
\nu \left( du,dz\right) $ (compare with (\ref{L-K})).

(iv) Keep the framework of Point (iii). When the measure $\nu $ is a product measure of the type $\nu \left( du,dx\right) =\rho \left(
du\right) \beta \left( dx\right) $, where $\beta $ is $\sigma $-finite and $%
\rho \left( du\right) $ verifies $\rho \left( \left\{ 0\right\} \right) =0$
and $\int_{\mathbb{R}}u^{2}\rho \left( du\right) <\infty $ (and therefore $%
\alpha _{B}\left( du\right) =\beta \left( B\right) \rho \left( du\right) $),
one says that the completely random measure $\varphi $ in (\ref{cool}) is \textsl{homogeneous%
} (see e.g. \cite{PecPru}). In particular, for a homogeneous measure $%
\varphi $, relation (\ref{cool2}) gives%
\begin{equation}
\mathbb{E}\left[ \exp \left( i\lambda \varphi \left( B\right) \right) \right]
=\exp \left[ \beta \left( B\right) \int_{\mathbb{R}}\left( e^{i\lambda
u}-1-i\lambda u\right) \rho \left( du\right) \right] \text{.}  \label{c}
\end{equation}%

(v) From (\ref{c}) and the classic results on infinitely divisible random
variables (see e.g. \cite{Sato}), one deduces that a centered and
square-integrable random variable $Y$ is infinitely divisible if and only if
the following holds: there exists a homogeneous completely random measure $%
\varphi $ on some space $\left( Z,\mathcal{Z}\right) $, as well as an
independent centered standard Gaussian random variable $G$, such that
\begin{equation*}
Y\overset{law}{=}aG+\varphi \left( B\right) \text{, \ for some }a\geq 0\text{
and }B\in \mathcal{Z}\text{.}
\end{equation*}

(vi) Let the framework and notation of the previous Point (iii) prevail, and
assume moreover that: (1) $\left( Z,\mathcal{Z}\right) =\left( [0,\infty ),%
\mathcal{B}\left( [0,\infty )\right) \right) $, and (2) $\nu \left(
du,dx\right) =\rho \left( du\right) dx$, where $dx$ stands for the
restriction of the Lebesgue measure on $[0,\infty )$, and $\rho $ verifies $%
\rho \left( \left\{ 0\right\} \right) =0$ and $\int_{\mathbb{R}}u^{2}\rho
\left( du\right) <\infty $. Then, the process
\begin{equation}
t\mapsto \varphi \left( \left[ 0,t\right] \right) =\int_{\mathbb{R}%
}\int_{[0,t]}u\hat{N}\left( du,dz\right) \text{, \ \ }t\geq 0\text{,}
\label{zetat}
\end{equation}%
is a \textsl{centered and square-integrable L\'{e}vy process} (with no
Gaussian component) started from zero: in particular, the stochastic process
$t\mapsto \varphi \left( \left[ 0,t\right] \right) $ has independent and
stationary increments.

(vii) Conversely, every centered and square-integrable L\'{e}vy process $Z_{t}$
with no Gaussian component is such that $Z_{t}\overset{law}{=}\varphi
\left( \left[ 0,t\right] \right) $ (in the sense of stochastic processes)
for some $\varphi \left( \left[ 0,t\right] \right) $ defined as in (\ref%
{zetat}). To see this, just use the fact that, for every $t$,
\begin{equation*}
\mathbb{E}\left[ \exp \left( i\lambda Z_{t}\right) \right] =\exp \left[
t\int_{\mathbb{R}}\left( e^{i\lambda u}-1-i\lambda u\right) \rho \left(
du\right) \right] \text{,}
\end{equation*}%
where the L\'{e}vy measure verifies $\rho \left( \left\{ 0\right\} \right) =0
$ and $\int_{\mathbb{R}}u^{2}\rho \left( du\right) <\infty $, and observe that this last relation implies that $\varphi([0,t])$ and $Z_t$ have the same finite-dimensional distributions. This fact is
the starting point of the paper by Farras \textit{et al. }\cite{FarreJoUtz},
concerning Hu-Meyer formulae for L\'{e}vy processes.

\bigskip

\textbf{Remark. }Let $\left( Z,\mathcal{Z}\right) $ be a measurable space.
Point 4 in Proposition \ref{P : LK} implies that every centered completely
random measure $\varphi $ on $\left( Z,\mathcal{Z}\right) $ has the same law
as a random mapping of the type%
\begin{equation*}
B\mapsto G\left( B\right) +\int_{\mathbb{R}}\int_{Z}u\mathbf{1}_{B}\left(
z\right) \hat{N}\left( du,dz\right) \text{,}
\end{equation*}%
where $G$ and $\hat{N}$ are, respectively, a Gaussian measure on $Z$ and an
independent compensated Poisson measure on $\mathbb{R}\times Z$.

\bigskip

\textbf{Examples. }(i) (\textit{Gamma random measures}) Let $\left( Z,%
\mathcal{Z}\right) $ be a measurable space, and let $\hat{N}$ be a centered
Poisson random measure on $\mathbb{R\times }Z$ with $\sigma $-finite control
measure
\begin{equation*}
\nu \left( du,dz\right) =\frac{\exp \left( -u\right) }{u}\mathbf{1}%
_{u>0}du\beta \left( dz\right) \text{,}
\end{equation*}%
where $\beta \left( dz\right) $ is a $\sigma $-finite measure on $\left( Z,%
\mathcal{Z}\right) $. Now define the completely random measure $\varphi $
according to (\ref{cool}).\ By using (\ref{cool2}) and the fact that
\begin{equation*}
\alpha _{B}\left( du\right) =\beta \left( B\right) \frac{\exp \left(
-u\right) }{u}\mathbf{1}_{u>0}du,
\end{equation*}%
one infers that, for every $B\in \mathcal{Z}$ such that $\beta \left(
B\right) <\infty $ and every real $\lambda $,%
\begin{eqnarray*}
\mathbb{E}\left[ \exp \left( i\lambda \varphi \left( B\right) \right) \right]
&=&\exp \left[ \beta \left( B\right) \int_{0}^{\infty }\left( e^{i\lambda
u}-1-i\lambda u\right) \frac{\exp \left( -u\right) }{u}du\right]  \\
&=&\exp \left[ \beta \left( B\right) \int_{0}^{\infty }\left( e^{i\lambda
u}-1\right) \frac{\exp \left( -u\right) }{u}du\right] \exp \left( -i\lambda
\beta \left( B\right) \right)  \\
&=&\frac{1}{\left( 1-i\lambda \right) ^{\beta \left( B\right) }}\exp \left(
-i\lambda \beta \left( B\right) \right) \text{,}
\end{eqnarray*}%
thus yielding that $\varphi \left( B\right) $ is a centered Gamma random
variable, with unitary scale parameter and shape parameter $\beta \left(
B\right) $. The completely random measure $\varphi \left( B\right) $ has
control measure $\beta $, and it is called a (centered)\textsl{\ Gamma random
measure}. Note that $\varphi \left( B\right) +\beta \left( B\right) >0$,
a.s.-$\mathbb{P}$, whenever $0<\beta \left( B\right) <\infty $. See e.g.
\cite{Handa1, Handa2, JamesRYsurvey, TsVeYorJFA}, and
the references therein, for recent results on (non-centered) Gamma random
measures.

(ii) (\textit{Dirichlet processes}) Let the notation and assumptions of the
previous example prevail, and assume that $0<\beta \left( Z\right) <\infty $
(that is, $\beta $ is non-zero and finite). Then, $\varphi \left( Z\right)
+\beta \left( Z\right) >0$ and the mapping%
\begin{equation}
B\longmapsto \frac{\varphi \left( B\right) +\beta \left( B\right) }{\varphi
\left( Z\right) +\beta \left( Z\right) }  \label{dirichlet}
\end{equation}%
defines a random probability measure on $\left( Z,\mathcal{Z}\right) $,
known as \textsl{Dirichlet process with parameter }$\beta $. Since the
groundbreaking paper by Ferguson \cite{Ferg1973}, Dirichlet processes play a
fundamental role in Bayesian non-parametric statistics: see e.g.\ \cite%
{JamesLijPru, MSW} and the references therein. Note that (\ref%
{dirichlet}) \textsl{does not} define a completely random measure (the
independence over disjoint sets fails): however, as shown in \cite%
{PecBer2008dir}, one can develop a theory of (multiple) stochastic
integration with respect to general Dirichlet processes, by using some
appropriate approximations in terms of orthogonal $U$-statistics. See \cite%
{Pit} for a state of the art overview of Dirichlet processes in modern
probability.

\subsection{Multiple stochastic integrals of elementary functions\label{SS :
WISIelementary}}

We now fix a \underline{good} completely random measure $\varphi $, in the sense of
Definition \ref{D : good} of Section \ref{SS : CRM}, and consider what happens when $\varphi ^{\left[ n%
\right] }$ is applied not to $C\in \mathcal{Z}_{\nu }^{n}$ but to its
restriction $C_{\pi }$, where $\pi $ is a partition of $\left[ n\right]
=\{1,...,n\}$. The set $C_{\pi }$ is defined according to (\ref{base : Bipi}%
). We shall also apply $\varphi ^{\left[ n\right] }$ to the union $\cup
_{\sigma \geq \pi }C_{\sigma }$.\footnote{%
From here, and for the rest of the paper (for instance, in formula (\ref%
{messtoch2})), the expressions \textquotedblleft $\sigma \geq \pi $%
\textquotedblright\ and \textquotedblleft $\pi \leq \sigma $%
\textquotedblright\ are used interchangeably$.$} It will be convenient to
express the result in terms of $C$, and thus to view $\varphi ^{\left[ n%
\right] }\left( C_{\pi }\right) $, for example, not as the map $\varphi ^{%
\left[ n\right] }$ applied to $C_{\pi }$, but as a suitably restricted map
applied to $C$. This restricted map will be denoted St$_{\pi }^{\varphi ,%
\left[ n\right] }$, where \textquotedblleft St\textquotedblright\ stands for
\textquotedblleft Stochastic\textquotedblright . In this way, the
restriction is embodied in the map, that is, the measure, rather than in the
set.

Thus, fix a good completely random measure $\varphi $, as well as an integer
$n\geq 2$.

\bigskip

\begin{definition}
\label{Def : randomMEAS}For every $\pi \in \mathcal{P}\left( \left[ n\right]
\right) $, we define the two random measures:\footnote{%
Here, we use a slight variation of the notation introduced by Rota and
Wallstrom in \cite{RoWa}. In particular, Rota\ and Wallstrom write St$_{\pi
}^{\left[ n\right] }$ and $\varphi _{\pi }^{\left[ n\right] }$,
respectively, instead of St$_{\pi }^{\varphi ,\left[ n\right] }$ and St$%
_{\geq \pi }^{\varphi ,\left[ n\right] }$.}%
\begin{equation}
\fbox{${\rm St}${}$_{\pi }^{\varphi ,\left[ n\right] }\left( C\right)
\triangleq \varphi ^{\left[ n\right] }\left( C_{\pi }\right), $}\text{ \ \ }%
C\in \mathcal{Z}_{\nu }^{n}\text{,}  \label{mes stoch}
\end{equation}
and
\begin{equation}
\fbox{${\rm St}${}$_{\geq \pi }^{\varphi ,\left[ n\right] }\left(
C\right) \triangleq \varphi ^{\left[ n\right] }\left( \cup _{\sigma \geq \pi
}C_{\sigma }\right) =\sum_{\sigma \geq \pi }{\rm St}${}$_{\sigma
}^{\varphi ,\left[ n\right] }\left( C\right), $}\text{ \ \ }C\in \mathcal{Z}%
_{\nu }^{n}\text{,}  \label{messtoch2}
\end{equation}%
that are the restrictions of $\varphi ^{\left[ n\right] }$, respectively to
the sets $Z_{\pi }^{n}$ and $\cup _{\sigma \geq \pi }Z_{\sigma }^{n}$.
\end{definition}

\bigskip

In particular, one has the following relations:

\begin{itemize}
\item St$_{\geq \hat{0}}^{\varphi ,\left[ n\right] }=\varphi ^{\left[ n%
\right] }$, because the subscript \textquotedblleft\ $\geq \hat{0}$
\textquotedblright\ involves no restriction. Hence, St$_{\geq
\hat{0}}^{\varphi ,\left[ n\right] }$ charges the whole space, and therefore
coincides with $\varphi ^{\left[ n\right] }$ (see also Lemma \ref{L :
decomposition});

\item St$_{\hat{0}}^{\varphi ,\left[ n\right] }$ does not charge diagonals;

\item St$_{\hat{1}}^{\varphi ,\left[ n\right] }$ charges only the full
diagonal set $Z_{\hat{1}}^{n}$;

\item for every $\sigma \in \mathcal{P}\left( \left[ n\right] \right) $ and
every $C\in \mathcal{Z}_{\nu }^{n}$, St$_{\geq \sigma }^{\varphi ,\left[ n%
\right] }\left( C\right) =$ St$_{\geq \hat{0}}^{\varphi ,\left[ n\right]
}\left( C\cap Z_{\sigma }^{n}\right) $.
\end{itemize}

\bigskip

We also set
\begin{equation}
\text{St}_{\hat{1}}^{\varphi ,\left[ 1\right] }\left( C\right) =\text{St}_{%
\hat{0}}^{\varphi ,\left[ 1\right] }\left( C\right) =\varphi \left( C\right)
\text{, \ \ }C\in \mathcal{Z}_{\nu }\text{.}  \label{ovio}
\end{equation}%
Observe that (\ref{ovio}) is consistent with the trivial fact that the class
$\mathcal{P}\left( \left[ 1\right] \right) $ contains uniquely the trivial
partition $\left\{ \left\{ 1\right\} \right\} $, so that, in this case, $%
\hat{1}=\hat{0}=\left\{ \left\{ 1\right\} \right\} $.

\bigskip

We now define the class $\mathcal{E}\left( \nu ^{n}\right) $ of \textsl{%
elementary functions} on $Z^{n}$. This is the collection of all functions of
the type%
\begin{equation}
f\left( \mathbf{z}_{n}\right) =\sum_{j=1}^{m}k_{j}\mathbf{1}_{C_{j}}\left(
\mathbf{z}_{n}\right) ,  \label{PreElemenaires}
\end{equation}%
where $k_{j}\in \mathbb{R}$ and every $C_{j}\in \mathcal{Z}_{\nu }^{n}$ has
the form $C_{j}=C_{j}^{1}\times \cdot \cdot \cdot \times C_{j}^{n}$, $%
C_{j}^{\ell }\in \mathcal{Z}_{\nu }$ ($\ell =1,...,n$). For every $f\in
\mathcal{E}\left( \nu ^{n}\right) $ as above, we set
\begin{eqnarray}
\text{St}_{\pi }^{\varphi ,\left[ n\right] }\left( f\right) &=&\int_{Z^{n}}fd%
\text{St}_{\pi }^{\varphi ,\left[ n\right] }=\sum_{j=1}^{m}k_{j}\text{St}%
_{\pi }^{\varphi ,\left[ n\right] }\left( C_{j}\right)  \label{IntPi1} \\
\text{St}_{\geq \pi }^{\varphi ,\left[ n\right] }\left( f\right)
&=&\int_{Z^{n}}fd\text{St}_{\geq \pi }^{\varphi ,\left[ n\right]
}=\sum_{j=1}^{m}k_{j}\text{St}_{\geq \pi }^{\varphi ,\left[ n\right] }\left(
C_{j}\right) ,  \label{IntPi2}
\end{eqnarray}%
and we say that St$_{\pi }^{\varphi ,\left[ n\right] }\left( f\right) $
(resp. St$_{\geq \pi }^{\varphi ,\left[ n\right] }\left( f\right) $) is the
\textsl{stochastic integral}\textit{\ }of $f$ with respect to St$_{\pi
}^{\varphi ,\left[ n\right] }$ (resp. St$_{\geq \pi }^{\varphi ,\left[ n%
\right] }\left( f\right) $). For $C\in \mathcal{Z}_{\nu }^{n}$, we write
interchangeably St$_{\pi }^{\varphi ,\left[ n\right] }\left( C\right) $ and
St$_{\pi }^{\varphi ,\left[ n\right] }\left( \mathbf{1}_{C}\right) $ (resp.
St$_{\geq \pi }^{\varphi ,\left[ n\right] }\left( C\right) $ and St$_{\geq
\pi }^{\varphi ,\left[ n\right] }\left( \mathbf{1}_{C}\right) $). Note that (%
\ref{messtoch2}) yields the relation%
\begin{equation*}
\text{St}_{\geq \pi }^{\varphi ,\left[ n\right] }\left( f\right)
=\sum_{\sigma \geq \pi }\text{St}_{\sigma }^{\varphi ,\left[ n\right]
}\left( f\right) \text{.}
\end{equation*}%
We can therefore apply the M\"{o}bius formula (\ref{MobInv}) in order to
deduce the inverse relation%
\begin{equation}
\text{St}_{\pi }^{\varphi ,\left[ n\right] }\left( f\right) =\sum_{\sigma
\geq \pi }\mu \left( \pi ,\sigma \right) \text{St}_{\geq \sigma }^{\varphi ,%
\left[ n\right] }\left( f\right) ,  \label{RWMob}
\end{equation}%
(see also \cite[Proposition 1]{RoWa}).

\bigskip

\textbf{Remarks. }(i) The random variables St$_{\pi }^{\varphi ,\left[ n%
\right] }\left( f\right) $ and St$_{\geq \pi }^{\varphi ,\left[ n\right]
}\left( f\right) \ $are elements of $L^{2}\left( \mathbb{P}\right) $ for
every $f\in \mathcal{E}\left( \nu ^{n}\right) $. While here $f$ is an
elementary function, it is neither supposed that $f$ is symmetric nor that
it vanishes on the diagonals.

(ii) Because $f$ is elementary, the moments and cumulants of the integrals (%
\ref{IntPi1}) and (\ref{IntPi2}) are always defined. They will be computed
later via diagram formulae.

\subsection{Wiener-It\^{o} stochastic integrals \label{SS : WISI}}

We consider the extension of the integrals St$_{\pi }^{\varphi ,\left[ n%
\right] }\left( f\right) $ to non-elementary functions $f$ in the case $\pi =%
\hat{0}=\left\{ \left\{ 1\right\} ,...,\left\{ n\right\} \right\} .$ In view
of (\ref{base : Bipi}) and (\ref{mes stoch}), the random measure St$_{\hat{0}%
}^{\varphi ,\left[ n\right] }$ does not charge diagonals (see Definition \ref%
{D : Diagonal sets}, as well as the subsequent examples).

We start with a heuristic presentation. While relation (\ref{WiIto1})
involves a simple integral over $Z$, our goal here is to define integrals
over $Z^{n}$ with respect to St$_{\hat{0}}^{\varphi ,\left[ n\right] }$,
that is, multiple integrals of functions $f:Z^{n}\mapsto \mathbb{R}$, of the
form
\begin{equation*}
I_{n}^{\varphi }\left( f\right) =\int_{Z_{\hat{0}}^{n}}f\left(
z_{1},...,z_{n}\right) \varphi \left( dz_{1}\right) \cdot \cdot \cdot
\varphi \left( dz_{n}\right) \text{.}
\end{equation*}%
Since the integration is over $Z_{\hat{0}}^{n}$ we are excluding diagonals,
that is, we are asking that the support of the integrator is restricted to
the set of those $(z_{1},...,z_{n})$ such that $z_{i}\neq z_{j}$ for every $%
i\neq j$, $1\leq i,j\leq n$. To define the multiple integral, we approximate
the restriction of $f$ to $Z_{\hat{0}}$ by special elementary functions,
namely by finite linear combinations of indicator functions $\mathbf{1}%
_{C_{1}\times \cdot \cdot \cdot \times C_{n}}$, where the $C_{j}$'s are
disjoint sets in $\mathcal{Z}_{\nu }$. This will allow us to define the
extension by using isometry, that is, relations of the type%
\begin{eqnarray}
\mathbb{E}\left[ I_{n}^{\varphi }\left( f\right) ^{2}\right]
&=&n!\int_{Z^{n}}f\left( z_{1},...,z_{n}\right) ^{2}\nu \left( dx_{1}\right)
\cdot \cdot \cdot \nu \left( dx_{n}\right)  \notag \\
&=&n!\int_{Z_{\hat{0}}^{n}}f\left( z_{1},...,z_{n}\right) ^{2}\nu \left(
dx_{1}\right) \cdot \cdot \cdot \nu \left( dx_{n}\right) .  \label{vv}
\end{eqnarray}%
Note that the equality (\ref{vv}) is due to the fact that the control
measure $\nu $ is non-atomic, and therefore the associated product measure
never charges diagonals. It is enough, moreover, to suppose that $f$ is
symmetric, because if
\begin{equation}
\widetilde{f}\left( z_{1},...,z_{n}\right) =\frac{1}{n!}\sum_{w\in \mathfrak{%
S}_{n}}f\left( z_{w\left( 1\right) },...,z_{w\left( n\right) }\right)
\label{can sim}
\end{equation}%
is the canonical symmetrization of $f$ ($\mathfrak{S}_{n}$ is the group of
permutations of $\left[ n\right] $), then
\begin{equation}
\fbox{$I_{n}^{\varphi }\left( f\right) =I_{n}^{\varphi }(\widetilde{f})$.}  \label{qs}
\end{equation}%
This last equality is just a \textquotedblleft stochastic
equivalent\textquotedblright\ of the well known fact that integrals with
respect to deterministic symmetric measures are invariant with respect to
symmetrizations of the integrands. Indeed, an intuitive explanation of (\ref%
{qs}) can be obtained by writing
\begin{equation*}
I_{n}^{\varphi }\left( f\right) =\int_{Z^{n}}f\left[ \mathbf{1}_{Z_{\hat{0}%
}}d\varphi ^{\left[ n\right] }\right]
\end{equation*}%
and by observing that the set $Z_{\hat{0}}$ is symmetric\footnote{That is: $\left(
z_{1},...,z_{n}\right) \in Z_{\hat{0}}^{n}$ implies that $(z_{w\left(
1\right) },...,$ $z_{w\left( n\right) })$ $\in $ $Z_{\hat{0}}^{n}$ for every
$w\in \mathfrak{S}_{n}$}, so that $I_{n}^{\varphi }\left( f\right) $ appears
as an integral with respect to the symmetric stochastic measure $\mathbf{1}%
_{Z_{\hat{0}}}d\varphi ^{\left[ n\right] }$.

From now on, we will denote by $\mathcal{Z}_{s,\nu }^{n}=\mathcal{Z}_{s}^{n}$
(the dependence on $\nu $ is dropped, whenever there is no risk of
confusion) the \textsl{symmetric }$\sigma $\textsl{-field} generated by the
elements of $\mathcal{Z}_{\nu }^{n}$ of the type%
\begin{equation}
\widetilde{C}=\bigcup\limits_{w\in \mathfrak{S}_{n}}C_{w\left( 1\right)
}\times C_{w\left( 2\right) }\times \cdot \cdot \cdot \times C_{w\left(
n\right) }\text{,}  \label{Csym}
\end{equation}%
where $\left\{ C_{j}:j=1,...,n\right\} \subset \mathcal{Z}_{\nu }$ are
\textsl{pairwise} disjoint and $\mathfrak{S}_{n}$ is the group of the
permutations of $\left[ n\right] $.

\bigskip

\textbf{Remark. }One can easily show that $\mathcal{Z}_{s}^{n}$ is the $%
\sigma $-field generated by the symmetric functions on $Z^{n}$ that are
square-integrable with respect to $\nu ^{n}$, vanishing on every set $Z_{\pi
}^{n}$ such that $\pi \neq \hat{0}$, that is, on all diagonals of $Z^{n}$ of
the type $z_{i_{1}}=\cdot \cdot \cdot =z_{i_{j}}$, $1\leq i_{1}\leq \cdot
\cdot \cdot \leq i_{j}\leq n$.

\bigskip

By specializing (\ref{mes stoch})-(\ref{IntPi2}) to the case $\pi =\hat{0}$,
we obtain an intrinsic characterization of\textit{\ }\textsl{Wiener-It\^{o}
multiple stochastic integrals}, as well as of the concept of \textsl{%
stochastic measure of order}\textit{\ }$n\geq 2$. The key is the following
result, proved in \cite[p. 1268]{RoWa}.

\begin{proposition}
\label{P : ST0RW}Let $\varphi $ be a good completely random measure.

(A) For every $f\in \mathcal{E}\left( \nu ^{n}\right) $,
\begin{equation*}
{\rm St}\,\text{ \negthinspace \negthinspace \negthinspace }_{\hat{0}%
}^{\varphi ,\left[ n\right] }\left( f\right) ={\rm St}\,\text{
\negthinspace \negthinspace \negthinspace }_{\hat{0}}^{\varphi ,\left[ n%
\right] }\left( \widetilde{f}\right) ,
\end{equation*}%
where $\widetilde{f}$ is given in (\ref{can sim}). In other words, the
measure ${\rm St}\,_{\hat{0}}^{\varphi ,\left[ n\right] }$ is symmetric.

(B) The collection $\left\{ {\rm St}\,_{\hat{0}}^{\varphi ,\left[ n\right]
}\left( C\right) :C\in \mathcal{Z}_{\nu }^{n}\right\} $ is the unique
symmetric random measure on $\mathcal{Z}_{\nu }^{n}$ verifying the two
properties: (i) ${\rm St}\,_{\hat{0}}^{\varphi ,\left[ n\right] }\left(
C\right) =0$ for every $C\in \mathcal{Z}_{\nu }^{n}$ such that $C\subset
Z_{\pi }^{n}$ for some $\pi \neq \hat{0}$, and (ii)
\begin{equation}
{\rm St}\,{}_{\hat{0}}^{\varphi ,\left[ n\right] }\left( \widetilde{C}%
\right) ={\rm St}\,{}_{\hat{0}}^{\varphi ,\left[ n\right] }\left( \mathbf{1%
}_{\widetilde{C}}\right) =n!\varphi \left( C_{1}\right) \times \varphi
\left( C_{2}\right) \times \cdot \cdot \cdot \times \varphi \left(
C_{n}\right) \text{,}  \label{carST0}
\end{equation}%
for every set $\widetilde{C}$ as in (\ref{Csym}).
\end{proposition}

\bigskip

\textbf{Remark. }Note that ${\rm St}\,_{\hat{0}}^{\varphi ,\left[ n\right]
}$ is defined on the $\sigma $-field $\mathcal{Z}_{\nu }^{n}$, which also
contains non-symmetric sets. The measure ${\rm St}\,_{\hat{0}}^{\varphi ,%
\left[ n\right] }$ is \textquotedblleft symmetric\textquotedblright\ in the
sense that, for every set $C\in \mathcal{Z}_{\nu }^{n}$, the following
equality holds: ${\rm St}\,_{\hat{0}}^{\varphi ,\left[ n\right] }\left(
C\right) ={\rm St}\,_{\hat{0}}^{\varphi ,\left[ n\right] }\left(
C_{w}\right) $, a.s.-$\mathbb{P}$, where $w$ is a permutation of the set $%
\left[ n\right] $ and
\begin{equation*}
C_{w}=\left\{ \left( z_{1},...,z_{n}\right) \in Z^{n}:\left( z_{w\left(
1\right) },...,z_{w\left( n\right) }\right) \in C\right\} \text{.}
\end{equation*}

\bigskip

We denote by $L_{s}^{2}\left( \nu ^{n}\right) $ the Hilbert space of
symmetric and square integrable functions on $Z^{n}$ (with respect to $\nu
^{n}$). We also write $\mathcal{E}_{s,0}\left( \nu ^{n}\right) $ to indicate
the subset of $L_{s}^{2}\left( \nu ^{n}\right) $ composed of \textsl{%
elementary functions vanishing on diagonals}, that is, the functions of the
type $f=\sum_{j=1}^{m}k_{j}\mathbf{1}_{\widetilde{C}_{j}}$, where $k_{j}\in
\mathbb{R}$ and every $\widetilde{C}_{j}\subset Z_{\hat{0}}^{n}$ has the
form (\ref{Csym}). The index 0 in $\mathcal{E}_{s,0}\left( \nu ^{n}\right) $
refers to the fact that it is a set of functions which equals 0 on the
diagonals. Since $\nu $ is non-atomic, and $\nu ^{n}$ does not charge
diagonals, one easily deduces that $\mathcal{E}_{s,0}\left( \nu ^{n}\right) $
is dense in $L_{s}^{2}\left( \nu ^{n}\right) $. Moreover, the relation (\ref%
{carST0}) implies that, $\forall n,m\geq 2$,
\begin{equation}
\mathbb{E}\left[ {\rm St}\,{}_{\hat{0}}^{\varphi ,\left[ m\right] }\left(
f\right) {\rm St}\,{}_{\hat{0}}^{\varphi ,\left[ n\right] }\left( g\right) %
\right] =\delta _{n,m}\times n!\int_{Z^{n}}f\left( \mathbf{z}_{n}\right)
g\left( \mathbf{z}_{n}\right) \nu ^{n}\left( d\mathbf{z}_{n}\right) \text{,
\ \ }  \label{IsoMWI}
\end{equation}%
$\forall f\in \mathcal{E}_{s,0}\left( \nu ^{m}\right) $ and $\forall g\in
\mathcal{E}_{s,0}\left( \nu ^{n}\right) $, where $\delta _{n,m}=1$ if $n=m$,
and $=0$ otherwise. This immediately yields that, for every $n\geq 2$, the
linear operator $f\mapsto {\rm St}\,{}_{\hat{0}}^{\varphi ,\left[ n\right]
}\left( f\right) $, from $\mathcal{E}_{s,0}\left( \nu ^{n}\right) $ into $%
L^{2}\left( \mathbb{P}\right) $, can be uniquely extendend to a continuous
operator from $L_{s}^{2}\left( \nu ^{n}\right) $ into $L^{2}\left( \mathbb{P}%
\right) $. It is clear that these extensions also enjoy the orthogonality
and isometry properties given by (\ref{IsoMWI}).

\bigskip

\begin{definition}
\label{D : MWIIdefinition}For every $f\in L_{s}^{2}\left( \nu ^{n}\right) $,
the random variable ${\rm St}\,{}_{\hat{0}}^{\varphi ,\left[ n\right]
}\left( f\right) $ is the \textbf{multiple stochastic Wiener-It\^{o}
integral }(of order $n$) of $f$ with respect to $\varphi $. We also use the
classic notation%
\begin{equation}
\fbox{${\rm St}\,${}$_{\hat{0}}^{\varphi ,\left[ n\right] }\left( f\right)
=I_{n}^{\varphi }\left( f\right) ,$ \ \ $f\in L_{s}^{2}\left( \nu
^{n}\right) $.}  \label{MWIdef}
\end{equation}%
Note that
\begin{equation}
\mathbb{E}\left[ I_{m}^{\varphi }\left( f\right) I_{n}^{\varphi }\left(
g\right) \right] =\delta _{n,m}\times n!\int_{Z^{n}}f\left( \mathbf{z}%
_{n}\right) g\left( \mathbf{z}_{n}\right) \nu ^{n}\left( d\mathbf{z}%
_{n}\right) ,  \label{goo}
\end{equation}%
$\forall f\in L_{s}^{2}\left( \nu ^{m}\right) $ and $\forall g\in
L_{s}^{2}\left( \nu ^{n}\right) $. For $n\geq 2$, the random measure $%
\left\{ {\rm St}\,{}_{\hat{0}}^{\varphi ,\left[ n\right] }\left( C\right)
:C\in \mathcal{Z}_{\nu }^{n}\right\} $ is called the \textbf{stochastic
measure of order} $n$ associated with $\varphi $. When $f\in L^{2}\left( \nu
^{n}\right) $ (not necessarily symmetric), we set
\begin{equation}
I_{n}^{\varphi }\left( f\right) =I_{n}^{\varphi }(\widetilde{f})\text{,}
\label{Yo la}
\end{equation}%
where $\widetilde{f}$ is the symmetrization of $f$ given by (\ref{can sim}).
\end{definition}

\bigskip

We have supposed so far that $\varphi \left( C\right) \in L^{p}\left(
\mathbb{P}\right) $, $p\geq 3$, for every $C\in \mathcal{Z}_{\nu }$ (see
Definition \ref{D : CRM}). We shall now suppose that $\varphi \left(
C\right) \in L^{2}\left( \mathbb{P}\right) $, $C\in \mathcal{Z}_{\nu }$. In
this case, the notion of \textquotedblleft good measure\textquotedblright\
introduced in Definition \ref{D : good} and Proposition \ref{P : ST0RW} do
not apply since, in this case, $\varphi ^{\left[ n\right] }$ may not exist.
Indeed, consider (\ref{prodEng}) for example with $C_{1}=...=C_{n}\in
\mathcal{Z}_{\nu }$. Then, the quantity
\begin{equation*}
\mathbb{E}\left\vert \varphi ^{\left[ n\right] }\left( C\right) \right\vert
^{2}=\mathbb{E}\left\vert \varphi \left( C_{1}\right) \right\vert ^{2n}
\end{equation*}%
may be infinite (see also \cite{Engel}). It follows that, for $n\geq 2$, the
multiple Wiener-It\^{o} integral cannot be defined as a multiple integral
with respect to the restriction to $Z_{\hat{0}}^{n}$ of the
\textquotedblleft full stochastic measure\textquotedblright\ $\varphi ^{%
\left[ n\right] }$. Nontheless, one can always do as follows.

\bigskip

\begin{definition}
\label{DEF : IntMultB'}Let $\varphi $ be a completely random measure in $%
L^{2}\left( \mathbb{P}\right) $ (and not necessarily in $L^{p}\left( \mathbb{%
P}\right) $, $p\geq 3$), with non-atomic control measure $\nu $. For $n\geq
2 $, let%
\begin{equation}
I_{n}^{\varphi }\left( f\right) =n!\sum_{k=1}^{m}\gamma _{k}\times \left\{
\varphi \left( C_{1}^{\left( k\right) }\right) \varphi \left( C_{2}^{\left(
k\right) }\right) \cdot \cdot \cdot \varphi \left( C_{n}^{\left( k\right)
}\right) \right\} \text{,}  \label{IntMWI2}
\end{equation}%
for every simple function $f\in \sum_{k=1}^{p}\gamma _{k}\mathbf{1}_{%
\widetilde{C}^{\left( k\right) }}\in \mathcal{E}_{s,0}\left( \nu ^{n}\right)
$, where every $\widetilde{C}^{\left( k\right) }$ is as in (\ref{Csym}). It
is easily seen that the integrals $I_{n}^{\varphi }\left( f\right) $ defined
in (\ref{IntMWI2}) still verify the $L^{2}\left( \mathbb{P}\right) $
isometry property (\ref{goo}). Since the sets of the type $\widetilde{C}$
generate $\mathcal{Z}_{s}^{n}$, and $\nu $ is non-atomic, the operator $%
I_{n}^{\varphi }\left( \cdot \right) $ can be extended to a continuous
linear operator from $L_{s}^{2}\left( \nu ^{n}\right) $ into $L^{2}\left(
\mathbb{P}\right) $, such that (\ref{goo}) is verified. When $f\in
L^{2}\left( \nu ^{n}\right) $ (not necessarily symmetric), we set
\begin{equation*}
I_{n}^{\varphi }\left( f\right) =I_{n}^{\varphi }(\widetilde{f})\text{,}
\end{equation*}%
where $\widetilde{f}$ is given by (\ref{can sim}).
\end{definition}

\bigskip

\textbf{Remark. }Of course, if $\varphi \in \cap _{p\geq 1}L^{p}\left(
\mathbb{P}\right) $ (for example, when $\varphi $ is a Gaussian measure or a
compensated Poisson measure), then the definition of $I_{n}^{\varphi }$
obtained from\textit{\ }(\ref{IntMWI2}) coincides with the one given in (\ref%
{MWIdef}).

\subsection{Integral notation}

The following\ \textquotedblleft integral notation\textquotedblright\ is
somewhat cumbersome but quite suggestive: for every $n\geq 2$, every $\sigma
\in \mathcal{P}\left( \left[ n\right] \right) $ and every elementary
function $f\in \mathcal{E}\left( \nu ^{n}\right) $,%
\begin{equation*}
\fbox{St$_{\sigma }^{\varphi ,\left[ n\right] }\left( f\right)
=\int_{Z_{\sigma }^{n}}f\left( z_{1},...,z_{n}\right) \varphi \left(
dz_{1}\right) \cdot \cdot \cdot \varphi \left( dz_{n}\right) $,}
\end{equation*}%
and%
\begin{equation*}
\fbox{St$_{\geq \pi }^{\varphi ,\left[ n\right] }\left( f\right) =\int_{\cup
_{\sigma \geq \pi }Z_{\sigma }^{n}}$ $\ f\left( z_{1},...,z_{n}\right)
\varphi \left( dz_{1}\right) \cdot \cdot \cdot \varphi \left( dz_{n}\right)
. $}
\end{equation*}

For instance:

\begin{itemize}
\item if $n=2$, $f\left( z_{1},z_{2}\right) =f_{1}\left( z_{1}\right)
f_{2}\left( z_{2}\right) $ and $\sigma =\hat{0}=\left\{ \left\{ 1\right\}
,\left\{ 2\right\} \right\} $, then
\begin{equation*}
I_{2}\left( f\right) =\text{St}_{\sigma }^{\varphi ,\left[ 2\right] }\left(
f\right) =\int_{z_{1}\neq z_{2}}f_{1}\left( z_{1}\right) f_{2}\left(
z_{2}\right) \varphi \left( dz_{1}\right) \varphi \left( dz_{2}\right) \text{%
;}
\end{equation*}

\item if $n=2$, $f\left( z_{1},z_{2}\right) =f_{1}\left( z_{1}\right)
f_{2}\left( z_{2}\right) $ and $\sigma =\hat{1}=\left\{ 1,2\right\} $, then
\begin{equation*}
\text{St}_{\hat{1}}^{\varphi ,\left[ 2\right] }\left( f\right)
=\int_{z_{1}=z_{2}}f_{1}\left( z_{1}\right) f_{2}\left( z_{2}\right) \varphi
\left( dz_{1}\right) \varphi \left( dz_{2}\right) \text{;}
\end{equation*}

\item if $n=2$,
\begin{eqnarray*}
\text{St}_{\geq \hat{0}}^{\varphi ,\left[ 2\right] }\left( f\right) &=&\text{%
St}_{\left\{ \left\{ 1\right\} ,\left\{ 2\right\} \right\} }^{\varphi ,\left[
2\right] }\left( f\right) +\text{St}_{\left\{ \left\{ 1,2\right\} \right\}
}^{\varphi ,\left[ 2\right] }\left( f\right) \\
&=&\int_{z_{1}\neq z_{2}}f\left( z_{1},z_{2}\right) \varphi \left(
dz_{1}\right) \varphi \left( dz_{2}\right) +\int_{z_{1}=z_{2}}f\left(
z_{1},z_{2}\right) \varphi \left( dz_{1}\right) \varphi \left( dz_{2}\right)
\\
&=&\int_{Z^{2}}f\left( z_{1},z_{2}\right) \varphi \left( dz_{1}\right)
\varphi \left( dz_{2}\right) .
\end{eqnarray*}

\item if $n=3$, $f\left( z_{1},z_{2},z_{3}\right) =f_{1}\left(
z_{1},z_{2}\right) f_{2}\left( z_{3}\right) $ and $\sigma =\left\{ \left\{
1,2\right\} ,\left\{ 3\right\} \right\} $, then
\begin{equation*}
\text{St}_{\sigma }^{\varphi ,\left[ 3\right] }\left( f\right) =\int
_{\substack{ z_{3}\neq z_{1}  \\ z_{1}=z_{2}}}f_{1}\left( z_{1},z_{2}\right)
f_{2}\left( z_{3}\right) \varphi \left( dz_{1}\right) \varphi \left(
dz_{2}\right) \varphi \left( dz_{3}\right) ;
\end{equation*}

\item if $n=3$ and $f\left( z_{1},z_{2},z_{3}\right) =f_{1}\left(
z_{1},z_{2}\right) f_{2}\left( z_{3}\right) $ and $\sigma =\hat{1}=\left\{
\left\{ 1,2,3\right\} \right\} $, then
\begin{equation*}
\int_{z_{1}=z_{2}=z_{3}}f_{1}\left( z_{1},z_{2}\right) f_{2}\left(
z_{3}\right) \varphi \left( dz_{1}\right) \varphi \left( dz_{2}\right)
\varphi \left( dz_{3}\right) .
\end{equation*}
\end{itemize}

\subsection{Chaotic representation\label{SS : CH R}}

When $\varphi $ is a Gaussian measure or a compensated Poisson measure,
multiple stochastic Wiener-It\^{o} integrals play a crucial role, due to the
\textsl{chaotic representation property} enjoyed by $\varphi $. Indeed, when
$\varphi $ is Gaussian or compensated Poisson, one can show that every
functional $F\left( \varphi \right) \in L^{2}\left( \mathbb{P}\right) $ of $%
\varphi $, admits a unique \textsl{chaotic (Wiener-It\^{o}) decomposition}%
\begin{equation}
F\left( \varphi \right) =\mathbb{E}\left[ F\left( \varphi \right) \right]
+\sum_{n\geq 1}I_{n}^{\varphi }\left( f_{n}\right) \text{, \ \ \ }f_{n}\in
L_{s}^{2}\left( \nu ^{n}\right) \text{,}  \label{CHAOS!}
\end{equation}%
where the series converges in $L^{2}\left( \mathbb{P}\right) $ (see for
instance \cite{DMM5}, \cite{Janson} or \cite{Major}), and the kernels $%
\left\{ f_{n}\right\} $ are uniquely determined. Formula (\ref{CHAOS!})
implies that random variables of the type $I_{n}^{\varphi }\left(
f_{n}\right) $ are the basic \textquotedblleft building
blocks\textquotedblright\ of the space of square-integrable functionals of $%
\varphi $. In general, for a completely random measure $\varphi $, the
Hilbert space $C_{n}^{\varphi }=\left\{ I_{n}^{\varphi }\left( f\right)
:f\in L_{s}^{2}\left( \nu ^{n}\right) \right\} $, $n\geq 1$, is called the $%
n $th\textsl{\ Wiener chaos}\textit{\ }associated with $\varphi $. We set by
definition $C_{0}^{\varphi }=\mathbb{R}$ (that is, $C_{0}^{\varphi }$ is the
collection of all non-random constants) so that, in (\ref{CHAOS!}), $\mathbb{%
E}\left[ F\right] \in C_{0}^{\varphi }$. Observe that relation (\ref{CHAOS!}%
) can be reformulated in terms of Hilbert spaces as follows:%
\begin{equation*}
L^{2}\left( \mathbb{P},\sigma \left( \varphi \right) \right)
=\bigoplus\limits_{n=0}^{\infty }C_{n}^{\varphi }\text{,}
\end{equation*}%
where $\oplus $ indicates an orthogonal sum in the Hilbert space $%
L^{2}\left( \mathbb{P},\sigma \left( \varphi \right) \right) $.

\bigskip

\textbf{Remarks.\ \ }(i) If $\varphi =G$ is a Gaussian measure with
non-atomic control $\nu $, for every $p>2$ and every $n\geq 2$, there exists
a universal constant $c_{p,n}>0$, such that%
\begin{equation}
\mathbb{E}\left[ \left\vert I_{n}^{G}\left( f\right) \right\vert ^{p}\right]
^{1/p}\leq c_{n,p}\mathbb{E}\left[ I_{n}^{G}\left( f\right) ^{2}\right]
^{1/2}\text{,}  \label{GaussChaosCOntr}
\end{equation}%
$\forall $ $f\in L_{s}^{2}\left( \nu ^{n}\right) $ (see \cite[Ch. V]{Janson}%
). Moreover, on every finite sum of Wiener chaoses $\oplus
_{j=0}^{m}C_{j}^{G}$ and for every $p\geq 1$, the topology induced by $%
L^{p}\left( \mathbb{P}\right) $ convergence is equivalent to the $L^{0}$%
-topology induced by convergence in probability, that is, convergence in
probability is equivalent to convergence in $L^{p}$, for every $p\geq 1$
(see e.g. \cite{Sch}). We refer the reader to \cite{DMM5}, \cite{Janson} or
\cite[Ch. 1]{Nualart} for an exhaustive analysis of the properties of
multiple stochastic Wiener-It\^{o} integrals with respect to a Gaussian
measure $G$.

(ii) The \textquotedblleft chaotic representation
property\textquotedblright\ is enjoyed by other processes. One of the most
well-known examples is given by the class of \textit{normal martingales},
that is, real-valued martingales on $\mathbb{R}_{+}$ having a predictable
quadratic variation equal to $t$. See \cite{DMM5} and \cite{MPS} for a
complete discussion of this point.

\bigskip

When $\varphi $ is Gaussian or compensated Poisson, one can characterize the
measures St$_{\pi }^{\varphi ,\left[ n\right] }$, when $\pi \neq \hat{0}$,
that is, the effect of these measures on diagonals. The key fact is the
following elementary identity (corresponding to Proposition 2 in \cite{RoWa}%
).

\begin{proposition}
Let $\varphi $ be a good completely random measure. Then, for every $n\geq 2$%
, every $C_{1},...,C_{n}\subset Z_{\nu }$, and every partition $\pi \in
\mathcal{P}\left( \left[ n\right] \right) $,
\begin{eqnarray}
&&{\rm St}\,\text{\negthinspace }_{\geq \pi }^{\varphi ,\left[ n\right]
}\left( C_{1}\times \cdot \cdot \cdot \times C_{n}\right)  \notag \\
&=&\prod_{b=\left\{ i_{1},...,i_{\left\vert b\right\vert }\right\} \in \pi }%
\text{ }{\rm St}\,\text{\negthinspace }_{\hat{1}}^{\varphi ,\left[
\left\vert b\right\vert \right] }\left( C_{i_{1}}\times \cdot \cdot \cdot
\times C_{i_{\left\vert b\right\vert }}\right)  \label{bub} \\
&=&\prod_{b=\left\{ i_{1},...,i_{\left\vert b\right\vert }\right\} \in \pi }%
\text{ }{\rm St}\,\text{\negthinspace }_{\hat{1}}^{\varphi ,\left[
\left\vert b\right\vert \right] }(\underset{\left\vert b\right\vert \text{ \
times}}{\underbrace{(\cap _{k=1}^{\left\vert b\right\vert \text{%
\negthinspace }}C_{i_{k}})\times \cdot \cdot \cdot \times (\cap
_{k=1}^{\left\vert b\right\vert \text{\negthinspace }}C_{i_{k}})})}.
\label{bbbb}
\end{eqnarray}
\end{proposition}

\begin{proof}
Recall that ${\rm St}\,{}_{\geq \pi }^{\varphi ,\left[ n\right]
}=\sum_{\sigma \geq \pi }{\rm St}\,{}_{\sigma }^{\varphi ,\left[ n\right]
} $, by (\ref{messtoch2}). To prove the first equality, just observe that
both random measures on the LHS and the RHS of (\ref{bub}) are the
restriction of the product measure St$_{\geq \hat{0}}^{\varphi ,\left[ n%
\right] }$ to the union of the sets $Z_{\sigma }^{n}$ such that $\sigma \geq
\pi $. Equality (\ref{bbbb}) is an application of (\ref{Triv1}).
\end{proof}

\subsection{Computational rules}

We are now going to apply our setup to the Gaussian case ($\varphi =G$) and
to the Poisson case ($\varphi =\hat{N}=N-\nu $). We always suppose that the
control measure $\nu $ of either $G$ or $\hat{N}$ is non-atomic. Many of the
subsequent formulae can be understood intuitively, by applying the following
computational rules:

\begin{description}
\item[Gaussian case:]
\begin{equation}
G\left( dx\right) ^{2}=\nu \left( dx\right) \text{ and }G\left( dx\right)
^{n}=0\text{, \ for every }n\geq 3.  \label{e : Grule}
\end{equation}

\item[Poisson case:]
\begin{equation}
( \hat{N}\left( dx\right) ) ^{n}=\left( N\left( dx\right) \right)
^{n}=N\left( dx\right) \text{, \ for every }n\geq 2\text{.}
\label{e : Prule}
\end{equation}
\end{description}

\subsection{Multiple Gaussian stochastic integrals of elementary functions
\label{SS : MultipleWISIGauss}}

Suppose $\varphi $ is Gaussian. The next result (whose proof is sketched
below) can be deduced from \cite[Example G, p. 1272, and Proposition 2, 6
and 12]{RoWa}.

\begin{theorem}
\label{T : DiagGauss}Let $\varphi =G$ be a centered Gaussian completely
random measure with non-atomic control measure $\nu $. For every $n\geq 2$
and every $A\in \mathcal{Z}_{\nu }$%
\begin{equation}
{\rm St}\,\text{\negthinspace }_{\hat{1}}^{G,\left[ n\right] }\underset{n%
\text{ times}}{(\underbrace{A\times \cdot \cdot \cdot \times A}})\triangleq
\Delta _{n}^{G}\left( A\right) =\left\{
\begin{array}{lll}
0 &  & n\geq 3 \\
\nu \left( A\right) &  & n=2\text{,}%
\end{array}%
\right.  \label{GaussDia}
\end{equation}%
(the measure $\Delta _{n}^{G}\left( \cdot \right) $ is called the \textbf{%
diagonal measure}\textit{\ }of order $n$ associated with $G$). More
generally, for every $n\geq 2$, $\sigma \in \mathcal{P}\left( \left[ n\right]
\right) $ and $A_{1},...,A_{n}\in \mathcal{Z}_{\nu }$,
\begin{eqnarray}
&&{\rm St}\,\text{\negthinspace }_{\geq \sigma }^{G,\left[ n\right]
}\left( A_{1}\times \cdot \cdot \cdot \times A_{n}\right)  \notag \\
&=&\left\{
\begin{array}{lll}
0, &  & \text{if }\exists b\in \sigma :\left\vert b\right\vert \geq 3 \\
\prod_{b=\left\{ i,j\right\} \in \sigma }\nu \left( A_{i}\cap A_{j}\right)
\prod_{\ell =1}^{k}G\left( A_{j_{\ell }}\right) , &  & \text{otherwise,}%
\end{array}%
\right.  \label{mDG1}
\end{eqnarray}%
and
\begin{eqnarray}
&&{\rm St}\,\text{\negthinspace }_{\sigma }^{G,\left[ n\right] }\left(
A_{1}\times \cdot \cdot \cdot \times A_{n}\right)  \notag \\
&=&\left\{
\begin{array}{lll}
0, &  & \text{if }\exists b\in \sigma :\left\vert b\right\vert \geq 3 \\
\prod_{b=\left\{ i,j\right\} \in \sigma }\nu \left( A_{i}\cap A_{j}\right)
{\rm St}\,{}_{\hat{0}}^{G,\left[ k\right] }\left( A_{j_{1}}\times \cdot
\cdot \cdot \times A_{j_{k}}\right) , &  & \text{otherwise,}%
\end{array}%
\right.  \label{m-DG2}
\end{eqnarray}%
where $j_{1},...,j_{k}$ are the singletons contained in $\sigma \backslash
\left\{ b\in \sigma :\left\vert b\right\vert =2\right\} $.
\end{theorem}

\begin{proof} Relation (\ref{GaussDia}) is classic (for a proof, see e.g. \cite[%
Proposition 6]{RoWa}). Formula (\ref{mDG1}) is obtained by combining (\ref%
{GaussDia}) and (\ref{bub}). To prove (\ref{m-DG2}), suppose first that $%
\exists b\in \sigma $ such that $\left\vert b\right\vert \geq 3$. Then, by
using M\"{o}bius inversion (\ref{RWMob}),
\begin{equation*}
\text{St}_{\sigma }^{G,\left[ n\right] }\left( A_{1}\times \cdot \cdot \cdot
\times A_{n}\right) =\sum_{\sigma \leq \rho }\mu \left( \sigma ,\rho \right)
\text{St}_{\geq \rho }^{G,\left[ n\right] }\left( A_{1}\times \cdot \cdot
\cdot \times A_{n}\right) =0\text{,}
\end{equation*}%
where the last equality is due to (\ref{mDG1}) and to the fact that, if $%
\rho \geq \sigma $ and $\sigma $ contains a block with more than two
elements, then $\rho $ must also contain a block with more than two
elements. This proves the first line of (\ref{m-DG2}). Now suppose that all
the blocks of $\sigma $ have at most two elements, and observe that, by
Definition \ref{Def : randomMEAS} and (\ref{prodEng}),%
\begin{equation*}
\prod_{\ell =1}^{k}G\left( A_{j_{\ell }}\right) ={\rm St}\,\text{%
\negthinspace }_{\geq \hat{0}}^{G,\left[ k\right] } (A_{j_{1}}\times
\cdot \cdot \cdot \times A_{j_{k}}) .
\end{equation*}%
The proof is concluded by using the following relations:
\begin{eqnarray*}
 {\rm St}\,\text{\negthinspace }_{\sigma }^{G,\left[ n\right] }\left(
A_{1}\times \cdot \cdot \cdot \times A_{n}\right) &=& \sum_{\sigma \le \rho}
\mu(\sigma,\rho){\rm St}\,\text{\negthinspace }_{\ge \rho }^{G,\left[ n\right] }\left(
A_{1}\times \cdot \cdot \cdot \times A_{n}\right)\\
&=& \prod_{b=\left\{ i,j\right\} \in \sigma }\nu (A_i \cap A_j) \sum_{\sigma \le \rho} \mu(\sigma,\rho)\prod_{b=\left\{ r,l\right\} \in \rho \backslash \sigma }\nu (A_r \cap A_l)\times \\
&& \quad\quad \quad \times
{\rm St}\,\text{\negthinspace }_{\ge \hat{0}}^{G,\left[ m\right] } (A_{i_{1}}\times
\cdot \cdot \cdot \times A_{i_{m}}){\bf 1}_{\{\{i_1\},...,\{i_m\} \text{ are the singletons of } \rho\}},
\end{eqnarray*}
where we write $b = \{r,l\} \in \rho \backslash \sigma$ to indicate that the block $b$ is in $\rho$ and not in $\sigma$ (equivalently, $b$ is obtained by merging two singletons of $\sigma$). Indeed, by the previous discussion, one has that the partitions $\rho$ involved in the previous sums have uniquely blocks of size one or two, and moreover, by M\"{o}bius inversion,
\begin{eqnarray*}
&& \sum_{\sigma \le \rho} \mu(\sigma,\rho)\!\!\!\!
\prod_{b=\left\{ r,l\right\} \in \rho\backslash \sigma }\!\!\!\!\!\!\nu (A_r \cap A_l) \,\,\times
{\rm St}\,\text{\negthinspace }_{ \ge\hat{0}}^{G,\left[ m\right] } (A_{i_{1}}\times
\cdot \cdot \cdot \times A_{i_{m}}){\bf 1}_{\{\{i_1\},...,\{i_m\} \text{ are the singletons of } \rho\}}\\
&& =\sum_{\sigma^* \le \rho^*} \!\!\mu(\sigma^*,\rho^*)\!\!\!\!
\prod_{b=\left\{ r,l\right\} \in \rho^*  }\!\!\!\!\!\nu (A_r \cap A_l)\,\,\times
{\rm St}\,\text{\negthinspace }_{\ge\hat{0}}^{G,\left[ m\right] } (A_{i_{1}}\times
\cdot \cdot \cdot \times A_{i_{m}}){\bf 1}_{\{\{i_1\},...,\{i_m\} \text{ are the singletons of } \rho^*\}}\\
&& = \sum_{\hat{0} \le \rho^*} \mu(\hat{0},\rho^*)
{\rm St}\,\text{\negthinspace }_{\geq \rho^*}^{G,\left[ k\right] } (A_{j_{1}}\times
\cdot \cdot \cdot \times A_{j_{k}})\\
&& = {\rm St}\,\text{\negthinspace }_{\hat{0}}^{G,\left[ k\right] }(A_{j_{1}}\times
\cdot \cdot \cdot \times A_{j_{k}}),
\end{eqnarray*}
where $\sigma^*$ and $\rho^*$ indicate, respectively, the restriction of $\sigma$ and $\rho$ to $\{j_1,...,j_k\}$,
where $\{j_1\},...,\{j_k\}$ are the singletons of $\sigma$ (in particular, $\sigma^* = \hat{0}$).
Note that the fact that $\mu(\sigma,\rho) = \mu(\sigma^*, \rho^*) = \mu(\hat{0},\rho^*)$ is a consequence of (\ref{MobF})
and of the fact that $\rho$ has uniquely blocks of size one or two.\footnote{Thanks to F. Benaych-Georges for pointing out this argument.}
%The proof is concluded by using the elementary relations (recall that the
%sets $Z_{\sigma }^{n}$ and $Z_{\hat{0}}^{n}$ are defined according to (\ref%
%{base : Bipi})):
%\begin{eqnarray*}
%{\rm St}\,\text{\negthinspace }_{\sigma }^{G,\left[ n\right] }\left(
%A_{1}\times \cdot \cdot \cdot \times A_{n}\right) &=&{\rm St}\,\text{%
%\negthinspace }_{\geq \sigma }^{G,\left[ n\right] }\left( (A_{1}\times \cdot
%\cdot \cdot \times A_{n})\mathbf{1}_{Z_{\sigma }^{n}}\right) \\
%&=&\prod_{b=\left\{ i,j\right\} \in \sigma }\nu \left( A_{i}\cap
%A_{j}\right) {\rm St}\,\text{\negthinspace }_{\geq \hat{0}}^{G,\left[ k%
%\right] }\left( (A_{j_{1}}\times \cdot \cdot \cdot \times A_{j_{k}})\mathbf{1%
%}_{Z_{\hat{0}}^{k}}\right) \\
%&=&\prod_{b=\left\{ i,j\right\} \in \sigma }\nu \left( A_{i}\cap
%A_{j}\right) {\rm St}\,\text{\negthinspace }_{\hat{0}}^{G,\left[ k\right]
%}\left( A_{j_{1}}\times \cdot \cdot \cdot \times A_{j_{k}}\right) \text{.}
%\end{eqnarray*}%
%%Note that in the second line of the previous display we used relations (\ref{bub}) and (\ref{bbbb}).
\end{proof}

\bigskip

\textbf{Examples. }(i) One has ${\rm St}\,$\negthinspace $_{\geq \hat{0}%
}^{G,\left[ n\right] }\left( A_{1}\times \cdot \cdot \cdot \times
A_{n}\right) =G\left( A_{1}\right) \cdot \cdot \cdot G\left( A_{n}\right) $,
which follows from (\ref{mDG1}), since $\hat{0}=\left\{ \left\{ 1\right\}
,...,\left\{ n\right\} \right\} $, but also directly since the symbol
\textquotedblleft\ $\geq \hat{0}$ \textquotedblright\ entails no restriction
on the partition. In integral notation ($f$ is always supposed to be
elementary)%
\begin{equation*}
{\rm St}\,\!_{\geq \hat{0}}^{G,\left[ n\right] }\left( f\right)
=\int_{Z^{n}}f\left( z_{1},...,z_{n}\right) G\left( dz_{1}\right) \cdot
\cdot \cdot G\left( dz_{n}\right) .
\end{equation*}%
On the other hand, there is no way to \textquotedblleft
simplify\textquotedblright\ an object such as ${\rm St}\,\!_{\hat{0}}^{G,%
\left[ n\right] }\left( A_{1}\times \cdot \cdot \cdot \times A_{n}\right) $,
wich is expressed, in integral notation, as%
\begin{equation*}
{\rm St}\,\!_{\hat{0}}^{G,\left[ n\right] }\left( f\right)
=I_{n}^{G}\left( f\right) =\int_{z_{1}\neq \cdot \cdot \cdot \neq
z_{n}}f\left( z_{1},...,z_{n}\right) G\left( dz_{1}\right) \cdot \cdot \cdot
G\left( dz_{n}\right) .
\end{equation*}

(ii) For $n\geq 3$, one has
\begin{equation*}
{\rm St}\,\!_{\geq \hat{1}}^{G,\left[ n\right] }\left( A_{1}\times \cdot
\cdot \cdot \times A_{n}\right) ={\rm St}\,\!_{\hat{1}}^{G,\left[ n\right]
}\left( A_{1}\times \cdot \cdot \cdot \times A_{n}\right) =0\text{,}
\end{equation*}%
since the partition $\hat{1}$ contains a single block of size $\geq 2$. In
integral notation,
\begin{equation*}
{\rm St}\,\!_{\hat{1}}^{G,\left[ n\right] }\left( f\right)
=\int_{z_{1}=\cdot \cdot \cdot =z_{n}}f\left( z_{1},...,z_{n}\right) G\left(
dz_{1}\right) \cdot \cdot \cdot G\left( dz_{n}\right) =0.
\end{equation*}%
When $n=1$, however, one has that%
\begin{equation*}
{\rm St}\,\!_{\hat{1}}^{G,\left[ 1\right] }\left( f\right)
=\int_{Z}f\left( z\right) G\left( dz\right) \sim \mathcal{N}\left(
0,\int_{Z}f^{2}d\nu \right) .
\end{equation*}%
When $n=2$,
\begin{equation*}
{\rm St}\,\!_{\hat{1}}^{G,\left[ 2\right] }\left( f\right)
=\int_{Z}f\left( z,z\right) \nu \left( dz\right) \text{.}
\end{equation*}

(iii) Let $n=3$ and $\sigma =\left\{ \left\{ 1\right\} ,\left\{ 2,3\right\}
\right\} $, then ${\rm St}\,$\negthinspace $_{\hat{1}}^{G,\left[ 3\right]
}\left( A_{1}\times A_{2}\times A_{n}\right) =0$, and therefore
\begin{eqnarray*}
&&{\rm St}\,\text{\negthinspace }_{\geq \sigma }^{G,\left[ 3\right]
}\left( A_{1}\times A_{2}\times A_{3}\right) \\
&=&{\rm St}\,\text{\negthinspace }_{\sigma }^{G,\left[ 3\right] }\left(
A_{1}\times A_{2}\times A_{3}\right) +{\rm St}\,\text{\negthinspace }_{%
\hat{1}}^{G,\left[ n\right] }\left( A_{1}\times A_{2}\times A_{3}\right) \\
&=&{\rm St}\,\text{\negthinspace }_{\sigma }^{G,\left[ 3\right] }\left(
A_{1}\times A_{2}\times A_{3}\right) =G\left( A_{1}\right) \nu \left(
A_{2}\times A_{3}\right) \text{.}
\end{eqnarray*}%
In integral notation:
\begin{equation*}
{\rm St}\,\!_{\geq \sigma }^{G,\left[ 3\right] }\left( f\right) ={\rm
St}_{\sigma }^{G,\left[ 3\right] }\left( f\right) =\int_{Z}\int_{Z}f\left(
z,z,x\right) \nu \left( dz\right) G\left( dx\right) .
\end{equation*}

(iv) Let $n=6$, and $\sigma =\left\{ \left\{ 1,2\right\} ,\left\{ 3\right\}
,\left\{ 4\right\} ,\left\{ 5,6\right\} \right\} $. Then,
\begin{equation*}
{\rm St}\,\text{\negthinspace }_{\geq \sigma }^{G,\left[ 6\right] }\left(
A_{1}\times ...\times A_{6}\right) =\nu \left( A_{1}\cap A_{2}\right) \nu
\left( A_{5}\cap A_{6}\right) G\left( A_{3}\right) G\left( A_{4}\right)
\text{,}
\end{equation*}%
whereas
\begin{equation*}
{\rm St}\,\text{\negthinspace }_{\sigma }^{G,\left[ 6\right] }\left(
A_{1}\times ...\times A_{6}\right) =\nu \left( A_{1}\cap A_{2}\right) \nu
\left( A_{5}\cap A_{6}\right) {\rm St}\,\text{\negthinspace }_{\hat{0}}^{G,%
\left[ 2\right] }\left( A_{3}\times A_{4}\right) \text{.}
\end{equation*}%
These relations can be reformulated in integral notation as%
\begin{eqnarray*}
&&{\rm St}\,\text{\negthinspace }_{\geq \sigma }^{G,\left[ 6\right]
}\left( A_{1}\times ...\times A_{6}\right) \\
&=&\int_{Z}\int_{Z}\left\{ \int_{Z}\int_{Z}f\left( x,x,y,y,w,z\right) \nu
\left( dx\right) \nu \left( dy\right) \right\} G\left( dw\right) G\left(
dz\right) \\
&&{\rm St}\,\text{\negthinspace }_{\sigma }^{G,\left[ 6\right] }\left(
A_{1}\times ...\times A_{6}\right) \\
&=&\int_{w\neq z}\left\{ \int_{Z}\int_{Z}f\left( x,x,y,y,w,z\right) \nu
\left( dx\right) \nu \left( dy\right) \right\} G\left( dw\right) G\left(
dz\right) .
\end{eqnarray*}

(v) Let $n=6$, and $\sigma =\left\{ \left\{ 1,2\right\} ,\left\{ 3\right\}
,\left\{ 4\right\} ,\left\{ 5\right\} ,\left\{ 6\right\} \right\} $. Then,
\begin{equation*}
{\rm St}\,\!_{\geq \sigma }^{G,\left[ 6\right] }\left( A_{1}\times
...\times A_{6}\right) =\nu \left( A_{1}\cap A_{2}\right) G\left(
A_{3}\right) G\left( A_{4}\right) G\left( A_{5}\right) G(A_{6})
\end{equation*}%
and also%
\begin{equation*}
{\rm St}\,\text{\negthinspace }_{\sigma }^{G,\left[ 6\right] }\left(
A_{1}\times ...\times A_{6}\right) =\nu \left( A_{1}\cap A_{2}\right)
{\rm St}\,\text{\negthinspace }_{\hat{0}}^{G,\left[ 4\right] }\left(
A_{3}\times A_{4}\times A_{5}\times A_{6}\right) \text{.}
\end{equation*}

\subsection{Multiple stochastic Poisson integrals of elementary functions.}

A result analogous to Theorem \ref{T : DiagGauss} holds in the Poisson case.
To state this result in a proper way, we shall introduce the following
notation. Given $n\geq 2$ and $\sigma \in \mathcal{P}\left[ n\right] $, we
write%
\begin{equation*}
\mathbf{B}_{1}\left( \sigma \right) =\left\{ b\in \sigma :\left\vert
b\right\vert =1\right\} ,
\end{equation*}%
to denote the collection of singleton blocks of $\sigma $, and%
\begin{equation}
\mathbf{B}_{2}\left( \sigma \right) =\left\{ b=\left\{ i_{1},...,i_{\ell
}\right\} \in \sigma :\ell =\left\vert b\right\vert \geq 2\right\} \text{ }
\label{heavy not}
\end{equation}%
to denote the collection of the blocks of $\sigma $ containing \textsl{two
or more} elements. We shall also denote by $\mathbf{PB}_{2}\left( \sigma
\right) $ the set of all $2$-partitions of $\mathbf{B}_{2}\left( \sigma
\right) $, that is, $\mathbf{PB}_{2}\left( \sigma \right) $ is the
collection of all ordered pairs $\left( R_{1};R_{2}\right) $ of \textsl{%
non-empty} subsets of $\mathbf{B}_{2}\left( \sigma \right) $ such that $%
R_{1},R_{2}\subset \mathbf{B}_{2}\left( \sigma \right) $, $R_{1}\cap
R_{2}=\varnothing $, and $R_{1}\cup R_{2}=\mathbf{B}_{2}\left( \sigma
\right) $; whenever $\mathbf{B}_{2}\left( \sigma \right) =\varnothing $, one
sets $\mathbf{PB}_{2}\left( \sigma \right) =\varnothing $. We stress that $%
\mathbf{PB}_{2}\left( \sigma \right) $ is a partition of $\mathbf{B}%
_{2}\left( \sigma \right) $; the fact that $\mathbf{B}_{2}\left( \sigma
\right) $ is also a subset of the partition $\sigma $ should not create
confusion.

\bigskip

\textbf{Examples. }(i) Let $n=7$, and $\sigma =\left\{ \left\{ 1,2\right\}
,\left\{ 3,4\right\} ,\left\{ 5,6\right\} ,\left\{ 7\right\} \right\} $.
Then,
\begin{equation*}
\mathbf{B}_{2}\left( \sigma \right) =\left\{ \left\{ 1,2\right\} ,\left\{
3,4\right\} ,\left\{ 5,6\right\} \right\}
\end{equation*}%
and $\mathbf{PB}_{2}\left( \sigma \right) $ contains the six ordered pairs:%
\begin{eqnarray*}
&&(\left\{ \left\{ 1,2\right\} ,\left\{ 3,4\right\} \right\} ;\left\{
\left\{ 5,6\right\} \right\} )\text{ \ \ ; \ \ }(\left\{ \left\{ 5,6\right\}
\right\} ;\left\{ \left\{ 1,2\right\} ,\left\{ 3,4\right\} \right\} ) \\
&&(\left\{ \left\{ 1,2\right\} ,\left\{ 5,6\right\} \right\} ;\left\{
\left\{ 3,4\right\} \right\} )\text{ \ \ ; \ \ }(\left\{ \left\{ 3,4\right\}
\right\} ;\left\{ \left\{ 1,2\right\} ,\left\{ 5,6\right\} \right\} ) \\
&&(\left\{ \left\{ 3,4\right\} ,\left\{ 5,6\right\} \right\} ;\left\{
\left\{ 1,2\right\} \right\} )\text{ \ \ ; \ \ }(\left\{ \left\{ 1,2\right\}
\right\} ;\left\{ \left\{ 3,4\right\} ,\left\{ 5,6\right\} \right\} )\text{.}
\end{eqnarray*}%
For instance, the first ordered pair is made up of $R_{1}=\left\{ \left\{
1,2\right\} ,\left\{ 3,4\right\} \right\} $ and $R_{2}=\left\{ \left\{
5,6\right\} \right\} $, whose union is $\mathbf{B}_{2}\left( \sigma \right) $%
.

(ii) If $n=5$ and $\sigma =\left\{ \left\{ 1,2,3\right\} ,\left\{ 4\right\}
,\left\{ 5\right\} \right\} $, then $\mathbf{B}_{1}\left( \sigma \right)
=\left\{ \left\{ 4\right\} ,\left\{ 5\right\} \right\} $, $\mathbf{B}%
_{2}\left( \sigma \right) =\left\{ \left\{ 1,2,3\right\} \right\} $ and $%
\mathbf{PB}_{2}\left( \sigma \right) =\varnothing $.

(iii) If $n=7$ and $\sigma =\left\{ \left\{ 1,2,3\right\} ,\left\{
4,5\right\} ,\left\{ 6\right\} ,\left\{ 7\right\} \right\} $, then $\mathbf{B%
}_{1}\left( \sigma \right) =\left\{ \left\{ 6\right\} ,\left\{ 7\right\}
\right\} $ and $\mathbf{B}_{2}\left( \sigma \right) =\left\{ \left\{
1,2,3\right\} ,\left\{ 4,5\right\} \right\} $. Also, the set $\mathbf{PB}%
_{2}\left( \sigma \right) $ contains the two ordered pairs%
\begin{equation*}
\left( \left\{ 1,2,3\right\} ;\left\{ 4,5\right\} \right) \text{ \ \ and \ \
}\left( \left\{ 4,5\right\} ,\left\{ 1,2,3\right\} \right) \text{.}
\end{equation*}

\bigskip

We shall now suppose that $\varphi $ is a compensated Poisson measure.

\begin{theorem}
\label{T : PoissDiag}Let $\varphi =\hat{N}$ be a compensated Poisson measure
with non-atomic control measure $\nu $, and let $N\left( \cdot \right)
\triangleq \hat{N}\left( \cdot \right) +\nu \left( \cdot \right) $. For
every $n\geq 2$ and every $A\in \mathcal{Z}_{\nu }$,
\begin{equation}
{\rm St}\,\text{\negthinspace }_{\hat{1}}^{\hat{N},\left[ n\right] }%
\underset{n\text{ times}}{(\underbrace{A\times \cdot \cdot \cdot \times A}}%
)\triangleq \Delta _{n}^{\hat{N}}\left( A\right) =N\left( A\right)
\label{DiagPoiss}
\end{equation}%
($\Delta _{n}^{\hat{N}}\left( \cdot \right) $ is called the \textbf{diagonal
measure} of order $n$ associated with $\hat{N}$). Moreover, for every $%
A_{1},...,A_{n}\in \mathcal{Z}_{\nu }$,
\begin{equation}
{\rm St}\,\text{\negthinspace }_{\hat{1}}^{\hat{N},\left[ n\right]
}(A_{1}\times \cdot \cdot \cdot \times A_{n})=N\left( A_{1}\cap \cdot \cdot
\cdot \cap A_{n}\right) \text{.}  \label{VV}
\end{equation}%
More generally, for every $n\geq 2$, $\sigma \in \mathcal{P}\left( \left[ n%
\right] \right) $ and $A_{1},...,A_{n}\in \mathcal{Z}_{\nu }$,
\begin{eqnarray}
&&{\rm St}\,\text{\negthinspace }_{\geq \sigma }^{\hat{N},\left[ n\right]
}\left( A_{1}\times \cdot \cdot \cdot \times A_{n}\right)
\label{>sigmaPoiss} \\
&=&\prod_{b=\left\{ i_{1},...,i_{\ell}\right\} \in
{\bf B}_2(\sigma) }N\left( A_{i_{1}}\cap \cdot \cdot \cdot \cap A_{i_{\ell }}\right)
\prod_{a=1}^{k}\hat{N}\left( A_{j_{a}}\right) ,  \notag
\end{eqnarray}%
where $\left\{ \left\{ j_{1}\right\} ,...,\left\{ j_{k}\right\} \right\} =%
\mathbf{B}_{1}\left( \sigma \right) $, and also
\begin{eqnarray}
&&{\rm St}\,\!_{\sigma }^{\hat{N},\left[ n\right] }\left( A_{1}\times
\cdot \cdot \cdot \times A_{n}\right)  \label{=sigmaPoisson} \\
&=&\sum_{\left( R_{1};R_{2}\right) \in \mathbf{PB}_{2}\left( \sigma \right)
}\quad \prod_{b=\left\{ i_{1},...,i_{\ell }\right\} \in R_{1}}\nu \left(
A_{i_{1}}\cap \cdot \cdot \cdot \cap A_{i_{\ell }}\right) \times
\label{Milne} \\
&&\text{ \ }\times {\rm St}\,{}_{\hat{0}}^{\hat{N},\left[ \left\vert
R_{2}\right\vert +k\right] }\left( \underset{b=\left\{
e_{1},...,e_{m}\right\} \in R_{2}}{\mathrm{X}}(A_{e_{1}}\cap \cdot \cdot
\cdot \cap A_{e_{m}})\times A_{j_{1}}\times \cdot \cdot \cdot \times
A_{j_{k}}\right)  \label{AAAA} \\
&&\text{ \ \ \ }+{\rm St}\,{}_{\hat{0}}^{\hat{N},\left[ \left\vert \mathbf{%
B}_{2}\left( \sigma \right) \right\vert +k\right] }\left( \underset{%
b=\left\{ i_{1},...,i_{\ell }\right\} \in \mathbf{B}_{2}\left( \sigma
\right) }{\mathrm{X}}(A_{i_{1}}\cap \cdot \cdot \cdot \cap A_{i_{\ell
}})\times A_{j_{1}}\times \cdot \cdot \cdot \times A_{j_{k}}\right)
\label{AAA} \\
&&\text{ \ \ \ \ \ }+\prod_{b=\left\{ i_{1},...,i_{\ell }\right\} \in
\mathbf{B}_{2}\left( \sigma \right) }\nu \left( A_{i_{1}}\cap \cdot \cdot
\cdot \cap A_{i_{\ell }}\right) {\rm St}\,{}_{\hat{0}}^{\hat{N},\left[ k%
\right] }\left( A_{j_{1}}\times \cdot \cdot \cdot \times A_{j_{k}}\right) ,
\label{e : eg4}
\end{eqnarray}%
where $\left\{ \left\{ j_{1}\right\} ,...,\left\{ j_{k}\right\} \right\} =%
\mathbf{B}_{1}\left( \sigma \right) $ and where (by convention) $\Sigma
_{\varnothing }\equiv 0$, $\Pi _{\varnothing }\equiv 0$ and ${\rm St}\,{}_{%
\hat{0}}^{\hat{N},\left[ 0\right] }\equiv 1$. Also, $\left\vert
R_{2}\right\vert $ and $\left\vert \mathbf{B}_{2}\left( \sigma \right)
\right\vert $ stand, respectively, for the cardinality of $R_{2}$ and $%
\mathbf{B}_{2}\left( \sigma \right) $, and in formula (\ref{AAAA}) we used
the notation%
\begin{equation*}
\underset{b=\left\{ e_{1},...,e_{m}\right\} \in R_{2}}{\mathrm{X}}%
(A_{e_{1}}\cap \cdot \cdot \cdot \cap A_{e_{m}})=\underset{b\in R_{2}}{%
\mathrm{X}}(\cap _{e\in b}A_{e})=(\cap _{e\in b_{1}}A_{e})\times \cdot \cdot
\cdot \times (\cap _{e\in b_{\left\vert R_{2}\right\vert }}A_{e})\text{,}
\end{equation*}%
where $b_{1},....,b_{\left\vert R_{2}\right\vert }$ is some enumeration of $%
R_{2}$ (note that, due to the symmetry of ${\rm St}\,{}_{\hat{0}}^{\hat{N},%
\left[ \left\vert R_{2}\right\vert +k \right] }$, the choice of the
enumeration is immaterial). The summand appearing in formula (\ref{AAA}) is
defined via the same conventions.
\end{theorem}

\bigskip

\textbf{Remarks. }(a) When writing formula (\ref{AAA}), we implicitly use
the following convention: if $\mathbf{B}_{2}\left( \sigma \right)
=\varnothing $, then the symbol $\underset{b=\left\{ i_{1},...,i_{\ell
}\right\} \in \mathbf{B}_{2}\left( \sigma \right) }{\mathrm{X}}%
(A_{i_{1}}\cap \cdot \cdot \cdot \cap A_{i_{\ell }})$ is immaterial, and one
should read%
\begin{equation}
\underset{b=\left\{ i_{1},...,i_{\ell }\right\} \in \mathbf{B}_{2}\left(
\sigma \right) }{\mathrm{X}}(A_{i_{1}}\cap \cdot \cdot \cdot \cap A_{i_{\ell
}})\times A_{j_{1}}\times \cdot \cdot \cdot \times A_{j_{k}}=A_{j_{1}}\times
\cdot \cdot \cdot \times A_{j_{k}}=A_{1}\times \cdot \cdot \cdot \times A_{n}%
\text{,}  \label{bvb}
\end{equation}%
where the last equality follows from the fact that, in this case, $k=n$ and $%
\left\{ j_{1},...,j_{k}\right\} =\left\{ 1,...,n\right\} =\left[ n\right] $.
To see how the convention (\ref{bvb}) works, suppose that $\mathbf{B}%
_{2}\left( \sigma \right) =\varnothing $. Then, $\mathbf{PB}_{2}\left(
\sigma \right) =\varnothing $ and consequently, according to the conventions
stated in Theorem \ref{T : PoissDiag}, the lines (\ref{Milne})--(\ref{AAAA})
and (\ref{e : eg4}) are equal to zero (they correspond, respectively, to a
sum and a product over the empty set). The equality (\ref{=sigmaPoisson})
reads therefore%
\begin{equation}
{\rm St}\,\!_{\sigma }^{\hat{N},\left[ n\right] }\left( A_{1}\times \cdot
\cdot \cdot \times A_{n}\right) ={\rm St}\,{}_{\hat{0}}^{\hat{N},\left[ n%
\right] }\left( \underset{b=\left\{ i_{1},...,i_{\ell }\right\} \in \mathbf{B%
}_{2}\left( \sigma \right) }{\mathrm{X}}(A_{i_{1}}\cap \cdot \cdot \cdot
\cap A_{i_{\ell }})\times A_{j_{1}}\times \cdot \cdot \cdot \times
A_{j_{k}}\right) .  \label{vbv}
\end{equation}%
By using (\ref{bvb}) one obtains
\begin{eqnarray*}
&&{\rm St}\,{}_{\hat{0}}^{\hat{N},\left[ n\right] }\left( \underset{%
b=\left\{ i_{1},...,i_{\ell }\right\} \in \mathbf{B}_{2}\left( \sigma
\right) }{\mathrm{X}}(A_{i_{1}}\cap \cdot \cdot \cdot \cap A_{i_{\ell
}})\times A_{j_{1}}\times \cdot \cdot \cdot \times A_{j_{k}}\right) \\
&=&{\rm St}\,{}_{\hat{0}}^{\hat{N},\left[ n\right] }\left( A_{j_{1}}\times
\cdot \cdot \cdot \times A_{j_{k}}\right) ={\rm St}\,\!_{\sigma }^{\hat{N},%
\left[ n\right] }\left( A_{1}\times \cdot \cdot \cdot \times A_{n}\right)
\text{,}
\end{eqnarray*}%
entailing that, in this case, relation (\ref{=sigmaPoisson}) is equivalent
to the identity ${\rm St}\,{}_{\hat{0}}^{\hat{N},\left[ n\right] }=%
{\rm St}\,{}_{\hat{0}}^{\hat{N},\left[ n\right] }$.

(b) If $k=0$ (that is, if $\mathbf{B}_{1}\left( \sigma \right) $ equals the
empty set), then, according to the conventions stated in Theorem \ref{T :
PoissDiag}, one has ${\rm St}\,{}_{\hat{0}}^{\hat{N},\left[ k\right] }=%
{\rm St}\,{}_{\hat{0}}^{\hat{N},\left[ 0\right] }=1$. This yields that, in
this case, one should read line (\ref{e : eg4}) as follows:%
\begin{eqnarray*}
&&\prod_{b=\left\{ i_{1},...,i_{\ell }\right\} \ \in \mathbf{B}_{2}\left(
\sigma \right) }\nu \left( A_{i_{1}}\cap \cdot \cdot \cdot \cap A_{i_{\ell
}}\right) {\rm St}\,{}_{\hat{0}}^{\hat{N},\left[ k\right] }\left(
A_{j_{1}}\times \cdot \cdot \cdot \times A_{j_{k}}\right) \\
&=&\prod_{b=\left\{ i_{1},...,i_{\ell }\right\} \ \in \mathbf{B}_{2}\left(
\sigma \right) }\nu \left( A_{i_{1}}\cap \cdot \cdot \cdot \cap A_{i_{\ell
}}\right) \text{.}
\end{eqnarray*}

(c) If $k=0$, one should also read line (\ref{AAAA}) as%
\begin{eqnarray*}
&&{\rm St}{}_{\hat{0}}^{\hat{N},\left[ \left\vert R_{2}\right\vert +k%
\right] }\left( \underset{b=\left\{ e_{1},...,e_{m}\right\} \in R_{2}}{%
\mathrm{X}}(A_{e_{1}}\cap \cdot \cdot \cdot \cap A_{e_{m}})\times
A_{j_{1}}\times \cdot \cdot \cdot \times A_{j_{k}}\right) \\
&=&{\rm St}{}_{\hat{0}}^{\hat{N},\left[ \left\vert R_{2}\right\vert %
\right] }\left( \underset{b=\left\{ e_{1},...,e_{m}\right\} \in R_{2}}{%
\mathrm{X}}(A_{e_{1}}\cap \cdot \cdot \cdot \cap A_{e_{m}})\right) ,
\end{eqnarray*}%
and line (\ref{AAA}) as%
\begin{eqnarray*}
&&{\rm St}{}_{\hat{0}}^{\hat{N},\left[ \left\vert \mathbf{B}_{2}\left(
\sigma \right) \right\vert +k\right] }\left( \underset{b=\left\{
i_{1},...,i_{\ell }\right\} \in \mathbf{B}_{2}\left( \sigma \right) }{%
\mathrm{X}}(A_{i_{1}}\cap \cdot \cdot \cdot \cap A_{i_{\ell }})\times
A_{j_{1}}\times \cdot \cdot \cdot \times A_{j_{k}}\right) \\
&=&{\rm St}{}_{\hat{0}}^{\hat{N},\left[ \left\vert \mathbf{B}_{2}\left(
\sigma \right) \right\vert \right] }\left( \underset{b=\left\{
i_{1},...,i_{\ell }\right\} \in \mathbf{B}_{2}\left( \sigma \right) }{%
\mathrm{X}}(A_{i_{1}}\cap \cdot \cdot \cdot \cap A_{i_{\ell }})\right) .
\end{eqnarray*}

\bigskip

\begin{proof}[Proof of Theorem \protect\ref{T : PoissDiag}]
To see that (\ref{VV}) must necessarily hold, use the fact that $\nu $ is
non-atomic by assumption. Therefore,
\begin{eqnarray*}
&&{\rm St}^{\hat{N},[n]}_{\hat{1}}(A_1\times\cdot\cdot\cdot\times A_n) = {\rm St}^{\hat{N},[1]}_{\hat{0}}(A_1\cap\cdot\cdot\cdot\cap A_n) +\nu(A_1\cap\cdot\cdot\cdot\cap A_n) \\
&& = \hat{N}(A_1\cap\cdot\cdot\cdot\cap A_n) +\nu(A_1\cap\cdot\cdot\cdot\cap A_n)=N(A_1\cap\cdot\cdot\cdot\cap A_n).
\end{eqnarray*}
%the application
%\begin{equation*}
%(A_{1}\times \cdot \cdot \cdot \times A_{n})\mapsto N\left( A_{1}\cap \cdot
%\cdot \cdot \cap A_{n}\right)
%\end{equation*}%
%coincides a.s. with the restriction to the set $Z_{\hat{1}}^{n}$ (defined
%according to (\ref{base : Bipi})) of the signed measure defined by the
%relation%
%\begin{equation*}
%(A_{1}\times \cdot \cdot \cdot \times A_{n})\mapsto \hat{N}\left(
%A_{1}\right) \times \cdot \cdot \cdot \times \hat{N}\left( A_{n}\right)
%\text{. }
%\end{equation*}%
Observe also that (\ref{VV}) implies (\ref{DiagPoiss}). Equation (\ref{>sigmaPoiss}) is an
immediate consequence of (\ref{bub}), (\ref{DiagPoiss}) and (\ref{VV}). To prove (\ref{=sigmaPoisson}), use (\ref{>sigmaPoiss}) and the relation $N=%
\hat{N}+\nu $, to write
\begin{eqnarray*}
&&{\rm St}\,\text{\negthinspace }_{\geq \sigma }^{\hat{N},\left[ n\right]
}\left( A_{1}\times \cdot \cdot \cdot \times A_{n}\right) \\
&=&\prod_{\ell =2}^{n}\prod_{b=\left\{ i_{1},...,i_{\ell }\right\} \in
\sigma } \left[ \hat{N}\left( A_{i_{1}}\cap \cdot \cdot \cdot \cap
A_{i_{\ell }}\right) +\nu \left( A_{i_{1}}\cap \cdot \cdot \cdot \cap
A_{i_{\ell }}\right) ^{^{^{{}}}}\right] \prod_{a=1}^{k}\hat{N}\left(
A_{j_{a}}\right),
\end{eqnarray*}
and the last expression equals
\begin{eqnarray*}
&&\sum_{\left( R_{1};R_{2}\right) \in \mathbf{PB}_{2}\left( \sigma \right)
}\prod_{b=\left\{ i_{1},...,i_{\ell }\right\} \ \in R_{1}}\nu \left(
A_{i_{1}}\cap \cdot \cdot \cdot \cap A_{i_{\ell }}\right) \times \\
&&\times {\rm St}\,{}_{\geq \hat{0}}^{\hat{N},\left[ \left\vert
R_{2}\right\vert +k \right] }\left( \underset{b=\left\{ e_{1},...,e_{m}\right\}
\in R_{2}}{\mathrm{X}}(A_{e_{1}}\cap \cdot \cdot \cdot \cap A_{e_{m}})\times
A_{j_{1}}\times \cdot \cdot \cdot \times A_{j_{k}}\right) \\
&&+{\rm St}\,{}_{\geq \hat{0}}^{\hat{N},\left[\left\vert \mathbf{B}%
_{2}\left( \sigma \right) \right\vert+k \right] }\left( \underset{b=\left\{
i_{1},...,i_{\ell }\right\} \in \mathbf{B}_{2}\left( \sigma \right) }{%
\mathrm{X}}(A_{i_{1}}\cap \cdot \cdot \cdot \cap A_{i_{\ell }})\times
A_{j_{1}}\times \cdot \cdot \cdot \times A_{j_{k}}\right) \\
&&+\prod_{b=\left\{ i_{1},...,i_{\ell }\right\} \ \in \mathbf{B}_{2}\left(
\sigma \right) }\nu \left( A_{i_{1}}\cap \cdot \cdot \cdot \cap A_{i_{\ell
}}\right) {\rm St}\,{}_{\geq \hat{0}}^{\hat{N},\left[ k\right] }\left(
A_{j_{1}}\times \cdot \cdot \cdot \times A_{j_{k}}\right) \text{,}
\end{eqnarray*}%
since ${\rm St}\,{}_{\geq \hat{0}}^{\hat{N},\left[ n\right] }\left(
A_{1}\times \cdot \cdot \cdot \times A_{n}\right) =\hat{N}\left(
A_{1}\right) \cdot \cdot \cdot \hat{N}\left( A_{n}\right) $. The term before
last in the displayed equation corresponds to $R_{1}=\varnothing $, $R_{2}=%
\mathbf{B}_{2}\left( \sigma \right) $, and the last term to $R_{1}=\mathbf{B}%
_{2}\left( \sigma \right) $ and $R_{2}=\varnothing $. By definition, these
two cases are not involved in $\mathbf{PB}_{2}\left( \sigma \right) $. The
last displayed equation yields%
\begin{eqnarray*}
&&{\rm St}\,\text{\negthinspace }_{\sigma }^{\hat{N},\left[ n\right]
}\left( A_{1}\times \cdot \cdot \cdot \times A_{n}\right) ={\rm St}\,\text{%
\negthinspace }_{\geq \sigma }^{\hat{N},\left[ n\right] }\left( (A_{1}\times
\cdot \cdot \cdot \times A_{n})\mathbf{1}_{Z_{\sigma }^{n}}\right) \\
&=&\sum_{\left( R_{1};R_{2}\right) \in \mathbf{PB}_{2}\left( \sigma \right)
}\quad \prod_{b=\left\{ i_{1},...,i_{\ell }\right\} \ \in R_{1}}\nu \left(
A_{i_{1}}\cap \cdot \cdot \cdot \cap A_{i_{\ell }}\right) \times \\
&&\times {\rm St}\,{}_{\geq \hat{0}}^{\hat{N},\left[ \left\vert
R_{2}\right\vert+k \right] }\left( \left[ \underset{b=\left\{
e_{1},...,e_{m}\right\} \in R_{2}}{\mathrm{X}}(A_{e_{1}}\cap \cdot \cdot
\cdot \cap A_{e_{m}})\times A_{j_{1}}\times \cdot \cdot \cdot \times
A_{j_{k}}\right] \mathbf{1}_{Z_{\hat{0}}^{\left\vert R_{2}\right\vert
+k}}\right) \\
&&+{\rm St}\,{}_{\geq \hat{0}}^{\hat{N},\left[ \left\vert \mathbf{B}%
_{2}\left( \sigma \right) \right\vert +k \right] }\left( \left[ \underset{%
b=\left\{ i_{1},...,i_{\ell }\right\} \in \mathbf{B}_{2}\left( \sigma
\right) }{\mathrm{X}}(A_{i_{1}}\cap \cdot \cdot \cdot \cap A_{i_{\ell
}})\times A_{j_{1}}\times \cdot \cdot \cdot \times A_{j_{k}}\right] \mathbf{1%
}_{Z_{\hat{0}}^{\left\vert \mathbf{B}_{2}\left( \sigma \right) \right\vert +k
}}\right) \\
&&+\prod_{b=\left\{ i_{1},...,i_{\ell }\right\} \ \in \mathbf{B}_{2}\left(
\sigma \right) }\nu \left( A_{i_{1}}\cap \cdot \cdot \cdot \cap A_{i_{\ell
}}\right) {\rm St}\,{}_{\geq \hat{0}}^{\hat{N},\left[ k\right] }\left( %
\left[ A_{j_{1}}\times \cdot \cdot \cdot \times A_{j_{k}}\right] \mathbf{1}%
_{Z_{\hat{0}}^{k}}\right) \text{.}
\end{eqnarray*}%
Since, by definition,
\begin{eqnarray*}
&&{\rm St}\,{}_{\geq \hat{0}}^{\hat{N},\left[ \left\vert
R_{2}\right\vert +k \right] }\left( \left[ \underset{b=\left\{
e_{1},...,e_{m}\right\} \in R_{2}}{\mathrm{X}}(A_{e_{1}}\cap \cdot \cdot
\cdot \cap A_{e_{m}})\times A_{j_{1}}\times \cdot \cdot \cdot \times
A_{j_{k}}\right] \mathbf{1}_{Z_{\hat{0}}^{\left\vert R_{2}\right\vert +k
}}\right) \\
&&={\rm St}\,{}_{\hat{0}}^{\hat{N},\left[ \left\vert R_{2}\right\vert +k%
\right] }\left( \underset{b=\left\{ e_{1},...,e_{m}\right\} \in R_{2}}{%
\mathrm{X}}(A_{e_{1}}\cap \cdot \cdot \cdot \cap A_{e_{m}})\times
A_{j_{1}}\times \cdot \cdot \cdot \times A_{j_{k}}\right) \text{,}
\end{eqnarray*}
and
\begin{eqnarray*}
&&{\rm St}\,{}_{\geq \hat{0}}^{\hat{N},\left[ \left\vert \mathbf{B}%
_{2}\left( \sigma \right) \right\vert +k \right] }\left( \left[ \underset{%
b=\left\{ i_{1},...,i_{\ell }\right\} \in \mathbf{B}_{2}\left( \sigma
\right) }{\mathrm{X}}(A_{i_{1}}\cap \cdot \cdot \cdot \cap A_{i_{\ell
}})\times A_{j_{1}}\times \cdot \cdot \cdot \times A_{j_{k}}\right] \mathbf{1%
}_{Z_{\hat{0}}^{\left\vert \mathbf{B}_{2}\left( \sigma \right) \right\vert +k
}}\right) \\
&&={\rm St}\,{}_{\hat{0}}^{\hat{N},\left[ \left\vert \mathbf{B}%
_{2}\left( \sigma \right) \right\vert +k \right] }\left( \underset{b=\left\{
i_{1},...,i_{\ell }\right\} \in \mathbf{B}_{2}\left( \sigma \right) }{%
\mathrm{X}}(A_{i_{1}}\cap \cdot \cdot \cdot \cap A_{i_{\ell }})\times
A_{j_{1}}\times \cdot \cdot \cdot \times A_{j_{k}}\right),
\end{eqnarray*}
and
\begin{eqnarray*}
{\rm St}\,{}_{\geq \hat{0}}^{\hat{N},\left[ k\right] }\left( \left[
A_{j_{1}}\times \cdot \cdot \cdot \times A_{j_{k}}\right] \mathbf{1}_{Z_{%
\hat{0}}^{k}}\right) = {\rm St}\,{}_{\hat{0}}^{\hat{N},\left[ k\right] }\left( A_{j_{1}}\times
\cdot \cdot \cdot \times A_{j_{k}}\right) ,
\end{eqnarray*}
one obtains immediately the desired conclusion.
\end{proof}

\bigskip

\textbf{Remark on integral notation. }It is instructive to express the
results of Theorem \ref{T : PoissDiag} in integral notation. With $f\left(
z_{1},...,z_{n}\right) =g\left( z_{1}\right) \cdot \cdot \cdot g\left(
z_{n}\right) $, formula (\ref{>sigmaPoiss}) becomes
\begin{eqnarray*}
&&\int_{\cup _{\pi \geq \sigma }Z_{\pi }^{n}}g\left( z_{1}\right) \cdot
\cdot \cdot g\left( z_{n}\right) \hat{N}\left( dz_{1}\right) \cdot \cdot
\cdot \hat{N}\left( dz_{n}\right) \\
&=&\left( \prod_{b\in \sigma ,\left\vert b\right\vert \geq 2}\int_{Z}g\left(
z\right) ^{\left\vert b\right\vert }N\left( dz\right) \right) \times \left(
\int_{Z}g\left( z\right) \hat{N}\left( dz\right) \right) ^{k}\text{,}
\end{eqnarray*}%
where $k=\left\vert \mathbf{B}_{1}\left( \sigma \right) \right\vert .$ Again
with $f\left( z_{1},..,z_{n}\right) =g\left( z_{1}\right) \cdot \cdot \cdot
g\left( z_{n}\right) $, (\ref{=sigmaPoisson})--(\ref{e : eg4}) become
\begin{eqnarray*}
&&\int_{Z_{\sigma }^{n}}g\left( z_{1}\right) \cdot \cdot \cdot g\left(
z_{n}\right) \hat{N}\left( dz_{1}\right) \cdot \cdot \cdot \hat{N}\left(
dz_{n}\right) \\
&=&\sum_{\left( R_{1};R_{2}\right) \in \mathbf{PB}_{2}\left( \sigma \right)
}\prod_{b\in R_{1}}\int_{Z}g\left( z\right) ^{\left\vert b\right\vert }\nu
\left( dz\right) \times \mathbf{1}_{\left\{ R_{2}=\left\{
b_{1},...,b_{\left\vert R_{2}\right\vert }\right\} \right\} }\times \\
&&\text{ \ \ \ \ \ \ \ }\times \int_{z_{1}\neq \cdot \cdot \cdot \neq
z_{\left\vert R_{2}\right\vert +k}}g\left( z_{1}\right) ^{\left\vert
b_{1}\right\vert }\cdot \cdot \cdot g\left( z_{\left\vert R_{2}\right\vert
}\right) ^{\left\vert b_{\left\vert R_{2}\right\vert }\right\vert }g\left(
z_{\left\vert R_{2}\right\vert +1}\right) \cdot \cdot \cdot g\left(
z_{\left\vert R_{2}\right\vert +k}\right) \\
&&\text{ \ \ \ \ \ \ \ \ \ \ \ \ \ \ \ \ \ \ \ \ \ \ \ \ \ \ \ \ \ \ \ \ \ \
\ \ \ \ \ \ \ \ \ \ \ \ \ \ \ \ \ \ \ \ \ \ \ \ \ \ \ \ \ \ \ \ \ \ \ \ \ \
\ \ \ \ \ \ \ \ \ \ \ \ \ \ }\hat{N}\left( dz_{1}\right) \cdot \cdot \cdot
\hat{N}\left( dz_{\left\vert R_{2}\right\vert +k}\right) \\
&&+\mathbf{1}_{\left\{ \mathbf{B}_{2}\left( \sigma \right) =\left\{
b_{1},...,b_{\left\vert \mathbf{B}_{2}\left( \sigma \right) \right\vert
}\right\} \right\} }\times \\
&&\text{ }\times \int_{z_{1}\neq \cdot \cdot \cdot \neq z_{\left\vert
\mathbf{B}_{2}\left( \sigma \right) \right\vert +k}}g\left( z_{1}\right)
^{\left\vert b_{1}\right\vert }\cdot \cdot \cdot g\left( z_{\left\vert
\mathbf{B}_{2}\left( \sigma \right) \right\vert }\right) ^{\left\vert
b_{\left\vert \mathbf{B}_{2}\left( \sigma \right) \right\vert }\right\vert
}g\left( z_{\left\vert \mathbf{B}_{2}\left( \sigma \right) \right\vert
+1}\right) \cdot \cdot \cdot g\left( z_{\left\vert \mathbf{B}_{2}\left(
\sigma \right) \right\vert +k}\right) \\
&&\text{ \ \ \ \ \ \ \ \ \ \ \ \ \ \ \ \ \ \ \ \ \ \ \ \ \ \  \ \ \ \ \ \
\ \ \ \ \ \ \ \ \ \ \ \ \ \ \ \ \ \ \ \ \ \ \ \ \ \ \ \ \ \ \ \ \ \ \ \ \ \
\ \ \ \ \ \ \ \ \ \ \ \ \ \ \ }\hat{N}\left( dz_{1}\right) \cdot \cdot \cdot
\hat{N}\left( dz_{\left\vert R_{2}\right\vert +k}\right) \\
&&\text{ \ \ \ \ \ }+\prod_{b\in \mathbf{B}_{2}\left( \sigma \right)
}\int_{Z}g\left( z\right) ^{\left\vert b\right\vert }\nu \left( dz\right)
\times \int_{z_{1}\neq \cdot \cdot \cdot \neq z_{k}}g\left( z_{1}\right)
\cdot \cdot \cdot g\left( z_{k}\right) \hat{N}\left( dz_{1}\right) \cdot
\cdot \cdot \hat{N}\left( dz_{k}\right) ,
\end{eqnarray*}

where $k=\left\vert \mathbf{B}_{1}\left( \sigma \right) \right\vert $.

\bigskip

\textbf{Examples. }The examples below apply to a compensated Poisson measure
$\hat{N}$, and\ should be compared with those discussed after Theorem \ref{T
: DiagGauss}. We suppose throughout that $n\geq 2$.

(i) When $\sigma =\hat{0}=\left\{ \left\{ 1\right\} ,...,\left\{ n\right\}
\right\} $ one has, as in the Gaussian case,
\begin{equation*}
{\rm St}\,\text{\negthinspace }_{\geq \hat{0}}^{\hat{N},\left[ n\right]
}\left( A_{1}\times \cdot \cdot \cdot \times A_{n}\right) =\hat{N}\left(
A_{1}\right) \cdot \cdot \cdot \hat{N}\left( A_{n}\right)
\end{equation*}%
because the symbol \textquotedblleft\ $\geq \hat{0}$ \textquotedblright\
entails no restriction on the considered partitions. In integral notation,
this becomes
\begin{equation*}
{\rm St}\,\text{\negthinspace }_{\geq \hat{0}}^{\hat{N},\left[ n\right]
}\left( f\right) =\int_{Z^{n}}f\left( z_{1},...,z_{n}\right) \hat{N}\left(
dz_{1}\right) \cdot \cdot \cdot \hat{N}\left( dz_{n}\right) .
\end{equation*}%
The case of equality (\ref{=sigmaPoisson}) has already been discussed:
indeed, since $\sigma =\hat{0}$, and according to the conventions stated
therein, one has that $\mathbf{B}_{2}\left( \sigma \right) =\mathbf{PB}%
_{2}\left( \sigma \right) =\varnothing $ and therefore (\ref{=sigmaPoisson})
becomes an identity given by (\ref{e : eg4}), namely ${\rm St}\,$%
\negthinspace $_{\hat{0}}^{\hat{N},\left[ n\right] }\left( \cdot \right) =%
{\rm St}\,$\negthinspace $_{\hat{0}}^{\hat{N},\left[ n\right] }\left(
\cdot \right) $. Observe that, in integral notation, ${\rm St}\,$%
\negthinspace $_{\hat{0}}^{\hat{N},\left[ n\right] }\left( \cdot \right) $
is expressed as%
\begin{equation*}
{\rm St}\,\!_{\hat{0}}^{\hat{N},\left[ n\right] }\left( f\right)
=\int_{z_{1}\neq \cdot \cdot \cdot \neq z_{n}}f\left( z_{1},...,z_{n}\right)
\hat{N}\left( dz_{1}\right) \cdot \cdot \cdot \hat{N}\left( dz_{n}\right)
\text{.}
\end{equation*}

(ii) Suppose now $\sigma =\hat{1}=\left\{ \left\{ 1,...,n\right\} \right\} $%
. Then, (\ref{=sigmaPoisson}) reduces to (\ref{VV}). To see this, note that $%
k=0$ and that $\mathbf{B}_{2}\left( \hat{1}\right) $ contains only the block
$\left\{ 1,...,n\right\} $, so that $\mathbf{PB}_{2}\left( \sigma \right)
=\varnothing $. Hence the sum appearing in (\ref{Milne}) vanishes and one
has
\begin{eqnarray*}
&&{\rm St}\,\!_{\hat{1}}^{\hat{N},\left[ n\right] }\left( A_{1}\times
\cdot \cdot \cdot \times A_{n}\right) ={\rm St}\,{}_{\hat{0}}^{\hat{N},%
\left[ 1\right] }\left( A_{1}\cap \cdot \cdot \cdot \cap A_{n}\right) +\nu
\left( A_{1}\cap \cdot \cdot \cdot \cap A_{n}\right) \\
&=&\hat{N}\left( A_{1}\cap \cdot \cdot \cdot \cap A_{n}\right) +\nu \left(
A_{1}\cap \cdot \cdot \cdot \cap A_{n}\right) =N\left( A_{1}\cap \cdot \cdot
\cdot \cap A_{n}\right) .
\end{eqnarray*}%
In integral notation,
\begin{eqnarray*}
{\rm St}\,\!_{\hat{1}}^{\hat{N},\left[ n\right] }\left( f\right)
&=&\int_{z_{1}=\cdot \cdot \cdot =z_{n}}f\left( z_{1},...,z_{n}\right) \hat{N%
}\left( dz_{1}\right) \cdot \cdot \cdot \hat{N}\left( dz_{n}\right) \\
&=&\int_{Z}f\left( z,...,z\right) N\left( dz\right) .
\end{eqnarray*}%
This last relation makes sense heuristically, in view of the computational
rule
\begin{equation*}
\left( \hat{N}\left( dx\right) \right) ^{2}=\left( N\left( dx\right) \right)
^{2}-2N\left( dx\right) \nu \left( dx\right) +\left( \nu \left( dx\right)
\right) ^{2}=N\left( dx\right) ,
\end{equation*}%
since $\left( N\left( dx\right) \right) ^{2}=N\left( dx\right) $ and $\nu $
is non-atomic.

(iii) Let $n=3$ and $\sigma =\left\{ \left\{ 1\right\} ,\left\{ 2,3\right\}
\right\} $, so that $\mathbf{B}_{2}\left( \sigma \right) =$ $\left\{ \left\{
2,3\right\} \right\} $ and $\mathbf{PB}_{2}\left( \sigma \right)
=\varnothing $. According to (\ref{>sigmaPoiss}),%
\begin{equation}
{\rm St}\,\text{\negthinspace }_{\geq \sigma }^{\hat{N},\left[ 3\right]
}\left( A_{1}\times A_{2}\times A_{3}\right) =\hat{N}\left( A_{1}\right)
N\left( A_{2}\cap A_{3}\right) .  \label{e : N31}
\end{equation}%
On the other hand, (\ref{=sigmaPoisson}) yields
\begin{equation}
{\rm St}\,\!_{\sigma }^{\hat{N},\left[ 3\right] }\left( A_{1}\times
A_{2}\times A_{3}\right) ={\rm St}\,\!_{\hat{0}}^{\hat{N},\left[ 2\right]
}\left( A_{1}\times (A_{2}\cap A_{3})\right) +\hat{N}\left( A_{1}\right) \nu
\left( A_{2}\cap A_{3}\right) .  \label{e : N33}
\end{equation}%
In integral form, relation (\ref{e : N31}) becomes
\begin{equation*}
{\rm St}\,\text{\negthinspace }_{\geq \sigma }^{\hat{N},\left[ 3\right]
}\left( f\right) =\int_{Z}\int_{Z}f\left( z_{1},z_{2},z_{2}\right) \hat{N}%
\left( dz_{1}\right) N\left( dz_{2}\right) \text{,}
\end{equation*}%
and (\ref{e : N33}) becomes
\begin{eqnarray*}
{\rm St}\,\text{\negthinspace }_{\sigma }^{\hat{N},\left[ 3\right] }\left(
f\right) &=&\int_{z_{1}\neq z_{2},\text{ }z_{2}=z_{3}}f\left(
z_{1},z_{2},z_{3}\right) \hat{N}\left( dz_{1}\right) \hat{N}\left(
dz_{2}\right) \hat{N}\left( dz_{3}\right) \\
&=&\int_{z_{1}\neq z_{2}}f\left( z_{1},z_{2},z_{2}\right) \hat{N}\left(
dz_{1}\right) \hat{N}\left( dz_{2}\right) \\
&&+\int_{Z^{2}}f\left( z_{1},z_{2},z_{2}\right) \hat{N}\left( dz_{1}\right)
\nu \left( dz_{2}\right) .
\end{eqnarray*}%
This last relation makes sense heuristically, by noting that
\begin{eqnarray*}
\hat{N}\left( dz_{1}\right) \hat{N}\left( dz_{2}\right) \hat{N}\left(
dz_{2}\right) &=&\hat{N}\left( dz_{1}\right) N\left( dz_{2}\right) \\
&=&\hat{N}\left( dz_{1}\right) \hat{N}\left( dz_{2}\right) +\hat{N}\left(
dz_{1}\right) \nu \left( dz_{2}\right) .
\end{eqnarray*}%
We also stress that ${\rm St}\,$\negthinspace $_{\sigma }^{\hat{N},\left[ 3%
\right] }\left( f\right) $ can be also be expressed as
\begin{equation}
{\rm St}\,\text{\negthinspace }_{\sigma }^{\hat{N},\left[ 3\right] }\left(
f\right) =I_{2}^{\hat{N}}\left( g_{1}\right) +I_{1}^{\hat{N}}\left(
g_{2}\right) ,  \label{e : I21}
\end{equation}%
where $g_{1}\left( x,y\right) =f\left( x,y,y\right) $ and $g_{2}\left(
x\right) =\int f\left( x,y,y\right) \nu \left( dy\right) $. The form (\ref{e
: I21}) will be needed later. Since, by (\ref{VV}), ${\rm St}\,$%
\negthinspace $_{\hat{1}}^{\hat{N},\left[ 3\right] }\left( A_{1}\times
A_{2}\times A_{3}\right) =N\left( A_{1}\cap A_{2}\cap A_{3}\right) $ and
since, for our $\sigma $, one has ${\rm St}\,$\negthinspace $_{\geq \sigma
}^{\hat{N},\left[ 3\right] }={\rm St}\,$\negthinspace $_{\sigma }^{\hat{N},%
\left[ 3\right] }+{\rm St}\,$\negthinspace $_{\hat{1}}^{\hat{N},\left[ 3%
\right] }$, one also deduces the relation%
\begin{equation*}
\hat{N}\left( A_{1}\right) N\left( A_{2}\cap A_{3}\right) ={\rm St}\,\!_{%
\hat{0}}^{\hat{N},\left[ 2\right] }\left( A_{1}\times (A_{2}\cap
A_{3})\right) +N\left( A_{1}\cap A_{2}\cap A_{3}\right) +\hat{N}\left(
A_{1}\right) \nu \left( A_{2}\cap A_{3}\right) ,
\end{equation*}%
or, equivalently, since $\hat{N}=N+\nu $,%
\begin{equation}
\hat{N}\left( A_{1}\right) \hat{N}\left( A_{2}\cap A_{3}\right) ={\rm St}\,%
\!_{\hat{0}}^{\hat{N},\left[ 2\right] }\left( A_{1}\times (A_{2}\cap
A_{3})\right) +\nu \left( A_{1}\cap A_{2}\cap A_{3}\right) +\hat{N}\left(
A_{1}\cap A_{2}\cap A_{3}\right) \text{.}  \label{sq}
\end{equation}

We will see that (\ref{sq}) is consistent with the multiplication formulae
of next section.

(iv) Let $n=6$, and $\sigma =\left\{ \left\{ 1,2\right\} ,\left\{ 3\right\}
,\left\{ 4\right\} ,\left\{ 5,6\right\} \right\} $, so that%
\begin{equation*}
\mathbf{B}_{2}\left( \sigma \right) =\left\{ \left\{ 1,2\right\} ,\left\{
5,6\right\} \right\}
\end{equation*}%
and the class $\mathbf{PB}_{2}\left( \sigma \right) \ $contains the two pairs%
\begin{equation*}
\left( \left\{ 1,2\right\} ;\left\{ 5,6\right\} \right) \text{ \ \ and \ \ }%
\left( \left\{ 5,6\right\} ;\left\{ 1,2\right\} \right) \text{.}
\end{equation*}%
First, (\ref{>sigmaPoiss}) gives
\begin{equation*}
{\rm St}\,\text{\negthinspace }_{\geq \sigma }^{\hat{N},\left[ 6\right]
}\left( A_{1}\times ...\times A_{6}\right) =N\left( A_{1}\cap A_{2}\right)
N\left( A_{5}\cap A_{6}\right) \hat{N}\left( A_{3}\right) \hat{N}\left(
A_{4}\right) \text{.}
\end{equation*}%
Moreover, we deduce from (\ref{=sigmaPoisson}) that%
\begin{eqnarray*}
{\rm St}\,\text{\negthinspace }_{\sigma }^{\hat{N},\left[ 6\right] }\left(
A_{1}\times ...\times A_{6}\right) &=&\nu \left( A_{1}\cap A_{2}\right)
{\rm St}\,\text{\negthinspace }_{\hat{0}}^{\hat{N},\left[ 3\right] }\left(
\left( A_{5}\cap A_{6}\right) \times A_{3}\times A_{4}\right) \\
&&+\nu \left( A_{5}\cap A_{6}\right) {\rm St}\,\text{\negthinspace }_{\hat{%
0}}^{\hat{N},\left[ 3\right] }\left( \left( A_{1}\cap A_{2}\right) \times
A_{3}\times A_{4}\right) \\
&&+{\rm St}\,\text{\negthinspace }_{\hat{0}}^{\hat{N},\left[ 4\right]
}\left( \left( A_{1}\cap A_{2}\right) \times \left( A_{5}\cap A_{6}\right)
\times A_{3}\times A_{4}\right) \\
&&+\nu \left( A_{1}\cap A_{2}\right) \nu \left( A_{5}\cap A_{6}\right)
{\rm St}\,\text{\negthinspace }_{\hat{0}}^{\hat{N},\left[ 2\right] }\left(
A_{3}\times A_{4}\right) \text{.}
\end{eqnarray*}%
The last displayed equation becomes in integral form%
\begin{eqnarray*}
&&{\rm St}\,\text{\negthinspace }_{\sigma }^{\hat{N},\left[ 6\right]
}\left( f\right) \\
&=&\int_{\substack{ z_{1}=z_{2},\text{ }z_{5}=z_{6}\text{, }z_{3}\neq z_{4}
\\ z_{3}\neq z_{1}\text{, }z_{4}\neq z_{1}  \\ z_{3}\neq z_{5}\text{,\ }%
z_{4}\neq z_{5}}}f\left( z_{1},...,z_{6}\right) \hat{N}\left( dz_{1}\right)
\cdot \cdot \cdot \hat{N}\left( dz_{6}\right) \\
&=&\int_{w,x\neq y\neq z}f\left( w,w,x,y,z,z\right) \nu \left( dw\right)
\hat{N}\left( dx\right) \hat{N}\left( dy\right) \hat{N}\left( dz\right) \\
&&+\int_{w\neq x\neq y,z}f\left( w,w,x,y,z,z\right) \hat{N}\left( dw\right)
\hat{N}\left( dx\right) \hat{N}\left( dy\right) \nu \left( dz\right) \\
&&+\int_{w\neq x\neq y\neq z}f\left( w,w,x,y,z,z\right) \hat{N}\left(
dw\right) \hat{N}\left( dx\right) \hat{N}\left( dy\right) \hat{N}\left(
dz\right) \\
&&+\int_{w,x\neq y,z}f\left( w,w,x,y,z,z\right) \nu \left( dw\right) \hat{N}%
\left( dx\right) \hat{N}\left( dy\right) \nu\left( dz\right) .
\end{eqnarray*}%
Indeed, let us denote the RHS of the last expression as $(I)+(II)+(III)+(IV)$%
. For $\left( I\right) $ and $\left( II\right) $, we use (\ref{Milne})--(\ref%
{AAAA}) with $R_{1}=\left\{ \left\{ 1,2\right\} \right\} $ and $%
R_{2}=\left\{ \left\{ 5,6\right\} \right\} $ and $k=2$, which corresponds to
the number of singletons $\left\{ 3\right\} ,$ $\left\{ 4\right\} $. For $%
\left( III\right) $, we use (\ref{AAA}), with $k+\left\vert \mathbf{B}%
_{2}\left( \sigma \right) \right\vert =2+2=4$, and $\mathbf{B}_{2}\left(
\sigma \right) =\left\{ \left\{ 1,2\right\} ,\left\{ 5,6\right\} \right\} .$
For (\ref{e : eg4}) we use $k=2$ and our $\mathbf{B}_{2}\left( \sigma
\right) $.

(v) Let $n=6$, and $\sigma =\left\{ \left\{ 1,2\right\} ,\left\{ 3\right\}
,\left\{ 4\right\} ,\left\{ 5\right\} ,\left\{ 6\right\} \right\} $. Here, $%
\mathbf{B}_{2}\left( \sigma \right) =$ $\left\{ \left\{ 1,2\right\} \right\}
$ and the class $\mathbf{PB}_{2}\left( \sigma \right) $ is empty. Then,
\begin{equation*}
{\rm St}\,\text{\negthinspace }_{\geq \sigma }^{\hat{N},\left[ 6\right]
}\left( A_{1}\times ...\times A_{6}\right) =N \left( A_{1}\cap
A_{2}\right) \hat{N}\left( A_{3}\right) \hat{N}\left( A_{4}\right) \hat{N}%
\left( A_{5}\right) \hat{N}\left( A_{6}\right) \text{,}
\end{equation*}%
and also%
\begin{eqnarray*}
{\rm St}\,\text{\negthinspace }_{\sigma }^{\hat{N},\left[ 6\right] }\left(
A_{1}\times ...\times A_{6}\right) &=&\nu \left( A_{1}\cap A_{2}\right)
{\rm St}\,\text{\negthinspace }_{\hat{0}}^{\hat{N},\left[ 4\right] }\left(
A_{3}\times A_{4}\times A_{5}\times A_{6}\right) \\
&&+{\rm St}\,\text{\negthinspace }_{\hat{0}}^{\hat{N},\left[ 4\right]
}\left( \left( A_{1}\cap A_{2}\right) \times A_{3}\times A_{4}\times
A_{5}\times A_{6}\right) \text{.}
\end{eqnarray*}%
In integral form,%
\begin{eqnarray*}
&&{\rm St}\,\text{\negthinspace }_{\sigma }^{\hat{N},\left[ 6\right]
}\left( f\right) \\
&=&\int_{\substack{ z_{1}=z_{2},  \\ z_{1}\neq z_{j}\text{, }j=3,...,6  \\ %
z_{i}\neq z_{j}\text{,\ }3\leq i\neq j\leq 6}}f\left(
z_{1},z_{2},z_{3},z_{4},z_{5},z_{6}\right) \prod_{j=1}^{6}\hat{N}\left(
dz_{j}\right) \\
&=&\int_{\substack{ z_{1}\neq z_{j}\text{, }j=3,...,6  \\ z_{i}\neq z_{j}%
\text{,\ }3\leq i\neq j\leq 6}}f\left(
z_{1},z_{1},z_{3},z_{4},z_{5},z_{6}\right) \nu \left( dz_{1}\right)
\prod_{j=3}^{6}\hat{N}\left( dz_{j}\right) \\
&&+\int_{\substack{ z_{1}\neq z_{j}\text{, }j=3,...,6  \\ z_{i}\neq z_{j}%
\text{,\ }3\leq i\neq j\leq 6}}f\left(
z_{1},z_{1},z_{3},z_{4},z_{5},z_{6}\right) \hat{N}\left( dz_{1}\right)
\prod_{j=3}^{6}\hat{N}\left( dz_{j}\right) .
\end{eqnarray*}

\bigskip

\begin{corollary}
\label{C : meanN} Suppose that the Assumptions of Theorem \ref{T : PoissDiag}
hold. Fix $\sigma \in \mathcal{P}\left( \left[ n\right] \right) $ and assume
that $\left\vert b\right\vert \geq 2$, for every $b\in \sigma .$ Then,%
\begin{equation}
\mathbb{E}\left[ {\rm St}\,\text{\negthinspace }_{\sigma }^{\hat{N},\left[
n\right] }\left( A_{1}\times \cdot \cdot \cdot \times A_{n}\right) \right]
=\prod_{m=2}^{n}\prod_{b\in \left\{ j_{1},...,j_{m}\right\} \in \sigma }\nu
\left( A_{j_{1}}\cap \cdot \cdot \cdot \cap A_{j_{m}}\right) .
\label{e : meanN}
\end{equation}
\end{corollary}

\begin{proof}
Use (\ref{=sigmaPoisson})--(\ref{e : eg4}) and note that, by assumption, $%
k=0 $ (the partition $\sigma $ does not contain any singleton $\left\{
j\right\} $). It follows that the sum in (\ref{Milne}) vanishes, and one is
left with%
\begin{eqnarray*}
&&\mathbb{E}\left[ {\rm St}\,\text{\negthinspace }_{\sigma }^{\hat{N},%
\left[ n\right] }\left( A_{1}\times \cdot \cdot \cdot \times A_{n}\right) %
\right] \\
&=&\mathbb{E}\left[ {\rm St}\,{}_{\hat{0}}^{\hat{N},\left[ k+\left\vert
\mathbf{B}_{2}\left( \sigma \right) \right\vert \right] }\left( \underset{%
b=\left\{ i_{1},...,i_{\ell }\right\} \in \mathbf{B}_{2}\left( \sigma
\right) }{\mathrm{X}}(A_{i_{1}}\cap \cdot \cdot \cdot \cap A_{i_{\ell
}})\times A_{j_{1}}\times \cdot \cdot \cdot \times A_{j_{k}}\right) \right]
\\
&&+\prod_{b=\left\{ i_{1},...,i_{\ell }\right\} \in \mathbf{B}_{2}\left(
\sigma \right) }\nu \left( A_{i_{1}}\cap \cdot \cdot \cdot \cap A_{i_{\ell
}}\right) \\
&=&\prod_{b=\left\{ i_{1},...,i_{\ell }\right\} \in \mathbf{B}_{2}\left(
\sigma \right) }\nu \left( A_{i_{1}}\cap \cdot \cdot \cdot \cap A_{i_{\ell
}}\right) \text{,}
\end{eqnarray*}%
which is equal to the RHS of (\ref{e : meanN}).
\end{proof}

\section{Multiplication formulae \label{S : MF}}

\subsection{The general case \label{SS : mult General}}

\setcounter{equation}{0}The forthcoming Theorem \ref{T : ProdRW} applies to
every good completely random measure $\varphi $. It gives a universal
combinatorial rule, according to which every product of multiple stochastic
integrals can be represented as a sum over diagonal measures that are
indexed by non-flat diagrams (as defined in Section \ref{SS : diagrams}). We
will see that product formulae are crucial in order to deduce explicit
expressions for the cumulants and the moments of multiple integrals. As
discussed later in this section, Theorem \ref{T : ProdRW} contains (as
special cases) two celebrated \textsl{product formulae} for integrals with
respect to Gaussian and Poisson random measures. We provide two proofs of
Theorem \ref{T : ProdRW}: the first one is new and it is based on a decomposition of
partially diagonal sets; the second consists in a slight variation of the
combinatorial arguments displayed in the proofs of \cite[Th. 3 and Th. 4]%
{RoWa}, and is included for the sake of completeness. The theorem is
formulated for simple kernels to ensure that the integrals are always
defined, in particular the quantity ${\rm St}\,\!_{\sigma }^{\varphi ,%
\left[ n\right] }$, which appears on the RHS\ of (\ref{productformula}).

\begin{theorem}[Rota et Wallstrom]
\label{T : ProdRW}Let $\varphi $ be a good completely random measure with
non-atomic control $\nu $. For $n_{1},n_{2},...,n_{k}\geq 1$, we write $%
n=n_{1}+\cdot \cdot \cdot +n_{k}$, and we denote by $\pi ^{\ast }$ the
partition of $\left[ n\right] $ given by
\begin{equation}
\pi ^{\ast }=\left\{ \left\{ 1,...,n_{1}\right\} ,\left\{
n_{1}+1,...,n_{1}+n_{2}\right\} ,...,\left\{
n_{1}+...+n_{k-1}+1,...,n\right\} \right\} .  \label{pistar}
\end{equation}%
Then, if the kernels $f_{1},...,f_{k}$ are such that $f_{j}\in \mathcal{E}%
_{s,0}\left( \nu ^{n_{j}}\right) $ $(j=1,...,k)$, one has that%
\begin{equation}
\fbox{$\prod_{j=1}^{k}$I$_{n_{j}}^{\varphi }\left( f_{j}\right) $=$%
\prod_{j=1}^{k}{\rm St}\,$\negthinspace $_{\hat{0}}^{\varphi ,\left[ n_{j}%
\right] }\left( f_{j}\right) $=$\sum_{\sigma \in \mathcal{P}\left( \left[ n%
\right] \right) :\sigma \wedge \pi ^{\ast }=\hat{0}}{\rm St}\,$%
\negthinspace $_{\sigma }^{\varphi ,\left[ n\right] }\left( f_{1}\otimes
_{0}f_{2}\otimes _{0}\cdot \cdot \cdot \otimes _{0}f_{k}\right) $,}
\label{productformula}
\end{equation}%
where, by definition, the function in $n$ variables $f_{1}\otimes
_{0}f_{2}\otimes _{0}\cdot \cdot \cdot \otimes _{0}f_{k}\in \mathcal{E}%
\left( \nu ^{n}\right) $ is defined as%
\begin{equation}
f_{1}\otimes _{0}f_{2}\otimes _{0}\cdot \cdot \cdot \otimes _{0}f_{k}\left(
x_{1},x_{2},...,x_{n}\right) =\prod_{j=1}^{k}f_{j}\left(
x_{n_{1}+...+n_{j-1}+1},...,x_{n_{1}+...+n_{j}}\right) \text{, \ \ (}n_{0}=0%
\text{)}.  \label{0multContr}
\end{equation}
\end{theorem}

\begin{proof}
(\textit{First proof}) From the discussion of the previous section, one
deduces that
\begin{equation*}
\prod_{j=1}^{k}{\rm St}\,\!_{\hat{0}}^{\varphi ,\left[ n_{j}\right]
}\left( f_{j}\right) ={\rm St}\,\!_{\geq \hat{0}}^{\varphi ,\left[ n\right]
}\left[ \left( f_{1}\otimes _{0}\cdot \cdot \cdot \otimes _{0}f_{k}\right)
\mathbf{1}_{A^{\ast }}\right] \text{,}
\end{equation*}%
where
\begin{equation*}
A^{\ast }=\left\{ \left( z_{1},...,z_{n}\right) \in \mathcal{Z}%
^{n}:z_{i}\neq z_{j}\text{, }\forall i\neq j\text{ such that }i\sim _{\pi
^{\ast }}j\right\} \text{,}
\end{equation*}%
that is, $A^{\ast }$ is obtained by excluding all diagonals within each block of $\pi ^{\ast }$.
We shall prove that
\begin{equation}
A^{\ast }=\bigcup\limits_{\sigma \in \mathcal{P}\left( \left[ n\right]
\right) :\sigma \wedge \pi ^{\ast }=\hat{0}}Z_{\sigma }^{n}\text{.}
\label{uh}
\end{equation}%
Suppose first $\sigma $ is such that $\sigma \wedge \pi ^{\ast }=\hat{0}$,
that is, the meet of $\sigma $ and $\pi ^{\ast }$ is given by the
singletons. For every $\left( z_{1},...,z_{n}\right) \in Z_{\sigma }^{n}$
the following implication holds: if $i\neq j$ and $i\sim _{\pi ^{\ast }}j$,
then $i$ and $j$ are in two different blocks of $\sigma $, and therefore $%
z_{i}\neq z_{j}$. This implies that $Z_{\sigma }^{n}\subset A^{\ast }$. For
the converse, take $\left( z_{1},...,z_{n}\right) \in A^{\ast }$, and
construct a partition $\sigma \in \mathcal{P}\left( \left[ n\right] \right) $
by the following rule: $i\sim _{\sigma }j$ if and only if $z_{i}=z_{j}$.
For every pair $i\neq j$ such that $i\sim _{\pi ^{\ast }}j$, one has (by
definition of $A^{\ast }$) $z_{i}\neq z_{j}$, so that $\sigma \wedge \pi
^{\ast }=\hat{0}$, and hence (\ref{uh}). To conclude the proof of the
theorem, just use the additivity of ${\rm St}\,\!_{\geq \hat{0}}^{\varphi ,%
\left[ n\right] }$ to write
\begin{eqnarray*}
{\rm St}\,\!_{\geq \hat{0}}^{\varphi ,\left[ n\right] }\left[ \left(
f_{1}\otimes _{0}\cdot \cdot \cdot \otimes _{0}f_{k}\right) \mathbf{1}%
_{A^{\ast }}\right] &=&\sum_{\sigma \wedge \pi ^{\ast }=\hat{0}}{\rm St}\,%
\!_{\geq \hat{0}}^{\varphi ,\left[ n\right] }\left[ \left( f_{1}\otimes
_{0}\cdot \cdot \cdot \otimes _{0}f_{k}\right) \mathbf{1}_{Z_{\sigma }^{n}}%
\right] \\
&=&\sum_{\sigma \wedge \pi ^{\ast }=\hat{0}}{\rm St}\,\!_{\sigma
}^{\varphi ,\left[ n\right] }\left( f_{1}\otimes _{0}\cdot \cdot \cdot
\otimes _{0}f_{k}\right) ,
\end{eqnarray*}%
by using the relation ${\rm St}\,\!_{\geq \hat{0}}^{\varphi ,\left[ n%
\right] }\left[ \left( \cdot \right) \mathbf{1}_{Z_{\sigma }^{n}}\right] =%
{\rm St}\,\!_{ \sigma }^{\varphi ,\left[ n\right] }\left[ \cdot \right]
$.

(\textit{Second proof -- see \cite{RoWa}}) This proof uses Proposition \ref%
{P : latticeP} and Proposition \ref{P : ISO}. To simplify the discussion
(and without loss of generality) we can assume that $n_{1}\geq n_{2}\geq
\cdot \cdot \cdot \geq n_{k}$. For $j=1,...,k$ we have that%
\begin{equation*}
{\rm St}\,\!_{\hat{0}}^{\varphi ,\left[ n_{j}\right] }\left( f_{j}\right)
=\sum_{\sigma \in \mathcal{P}\left( \left[ n_{j}\right] \right) }\mu \left(
\hat{0},\sigma \right) {\rm St}\,\!_{\geq \sigma }^{\varphi ,\left[ n_{j}%
\right] }\left( f_{j}\right) \text{,}
\end{equation*}%
where we have used (\ref{RWMob}) with $\pi =\hat{0}$. From this relation one
obtains
\begin{eqnarray}
&&\prod_{j=1}^{k}{\rm St}\,\!_{\hat{0}}^{\varphi ,\left[ n_{j}\right]
}\left( f_{j}\right) =\sum_{\substack{\sigma _{1}\in \mathcal{P}\left( %
\left[ n_{1}\right] \right)}}\cdot\cdot\cdot \sum_{\substack{\sigma _{k}\in \mathcal{P}\left( %
\left[ n_{k}\right] \right) }} \prod_{j=1}^{k}\mu \left( \hat{0},\sigma
_{j}\right){\rm St}\,\!_{\geq \sigma _{j}}^{\varphi ,\left[
n_{j}\right] }\left( f_{j}\right)  \notag \\
&&=\sum_{\rho \in \mathcal{P}\left( \left[ n\right] \right) :\rho \leq \pi
^{\ast }}\mu \left( \hat{0},\rho \right) {\rm St}\,\!_{\geq \rho
}^{\varphi ,\left[ n\right] }\left( f_{1}\otimes _{0}\cdot \cdot \cdot
\otimes _{0}f_{k}\right) ,  \label{qq}
\end{eqnarray}%
where $\pi ^{\ast }$ is defined in (\ref{pistar}). To prove equality (\ref%
{qq}), recall the definition of \textquotedblleft class\textquotedblright\
in Section \ref{SS : Class}, as well as the last example in that section.
Observe that the segment $\left[ \hat{0},\pi ^{\ast }\right] $ has class $%
\left( n_{1},...,n_{k}\right) $, thus yielding (thanks to Proposition \ref{P
: ISO}) that $\left[ \hat{0},\pi ^{\ast }\right] $ is isomorphic to the
lattice product of the $\mathcal{P}\left( \left[ n_{j}\right] \right) $'s.
This implies that each vector
\begin{equation*}
\left( \sigma _{1},...,\sigma _{k}\right) \in \mathcal{P}\left( \left[ n_{1}%
\right] \right) \times \cdot \cdot \cdot \times \mathcal{P}\left( \left[
n_{k}\right] \right)
\end{equation*}%
has indeed the form
\begin{equation*}
\left( \sigma _{1},...,\sigma _{k}\right) =\psi ^{-1}\left( \rho \right)
\text{ }
\end{equation*}%
for a unique $\rho \in \left[ \hat{0},\pi ^{\ast }\right] $, where $\psi $
is a bijection defined as in (\ref{zed}). Now use, in order, Part 2 and Part
1 of Proposition \ref{P : latticeP} to deduce that
\begin{equation}
\prod_{j=1}^{k}\mu \left( \hat{0},\sigma _{j}\right) =\mu \left( \hat{0}%
,\left( \sigma _{1},...,\sigma _{k}\right) \right) =\mu \left( \hat{0},\psi
^{-1}\left( \rho \right) \right) =\mu \left( \hat{0},\rho \right) \text{.}
\label{*}
\end{equation}%
Observe that
\begin{equation*}
\hat{0}=\left\{ \left\{1\right\} ,...,\left\{ n_{j}\right\}
\right\} \text{ \ \ in \ \ }\mu \left( \hat{0},\sigma _{j}\right) \text{,}
\end{equation*}%
whereas
\begin{equation*}
\hat{0}=\left\{ \left\{ 1\right\} ,...,\left\{ n\right\} \right\} \text{ \ \
in \ \ }\mu ( \hat{0},\rho)
\text{.}
\end{equation*}%
Also, one has the relation
\begin{equation}
\prod_{j=1}^{k}{\rm St}\,\!_{\geq \sigma _{j}}^{\varphi ,\left[ n_{j}%
\right] }\left( f_{j}\right) ={\rm St}\,\!_{\geq \rho }^{\varphi ,\left[ n%
\right] }\left( f_{1}\otimes _{0}\cdot \cdot \cdot \otimes _{0}f_{k}\right)
\text{,}  \label{**}
\end{equation}%
(by the definition of ${\rm St}\,\!_{\geq \rho }^{\varphi ,\left[ n\right]
}$ as the measure charging all the diagonals contained in the diagonals
associated with the blocks of the $\sigma _{j}$ ($j=1,...,k$)). Then, (\ref%
{*}) and (\ref{**}) yield immediately (\ref{qq}).\ To conclude the proof,
write%
\begin{eqnarray*}
&&\sum_{\rho \in \mathcal{P}\left( \left[ n\right] \right) :\rho \leq \pi
^{\ast }}\mu \left( \hat{0},\rho \right) {\rm St}\,\!_{\geq \rho
}^{\varphi ,\left[ n\right] }\left( f_{1}\otimes _{0}\cdot \cdot \cdot
\otimes _{0}f_{k}\right) \\
&=&\sum_{\rho \in \mathcal{P}\left( \left[ n\right] \right) :\rho \leq \pi
^{\ast }}\mu \left( \hat{0},\rho \right) \sum_{\gamma \geq \rho }{\rm St}\,%
\!_{\gamma }^{\varphi ,\left[ n\right] }\left( f_{1}\otimes _{0}\cdot \cdot
\cdot \otimes _{0}f_{k}\right) \\
&=&\sum_{\gamma \in \mathcal{P}\left( \left[ n\right] \right) }{\rm St}\,%
\!_{\gamma }^{\varphi ,\left[ n\right] }\left( f_{1}\otimes _{0}\cdot \cdot
\cdot \otimes _{0}f_{k}\right) \sum_{\hat{0}\leq \rho \leq \pi ^{\ast
}\wedge \gamma }\mu \left( \hat{0},\rho \right) \text{.}
\end{eqnarray*}%
Since, by (\ref{maliMob}),
\begin{equation*}
\sum_{\hat{0}\leq \rho \leq \pi ^{\ast }\wedge \gamma }\mu \left( \hat{0}%
,\rho \right) =\left\{
\begin{array}{cc}
0 & \text{if \ }\pi ^{\ast }\wedge \gamma \neq \hat{0} \\
1 & \text{if \ }\pi ^{\ast }\wedge \gamma =\hat{0}.%
\end{array}%
\right. \text{,}
\end{equation*}%
relation (\ref{productformula}) is obtained.
\end{proof}

\bigskip

\textbf{Remark.} The RHS of (\ref{productformula}) can also be reformulated
in terms of diagrams and in terms of graphs, as follows:%
\begin{equation*}
\sum_{\sigma \in \mathcal{P}\left( \left[ n\right] \right) :\Gamma \left(
\pi ^{\ast },\sigma \right) \text{ is non-flat}}{\rm St}\,\!_{\sigma
}^{\varphi ,\left[ n\right] }\left( f_{1}\otimes _{0}f_{2}\otimes _{0}\cdot
\cdot \cdot \otimes _{0}f_{k}\right) \text{,}
\end{equation*}%
where $\Gamma \left( \pi ^{\ast },\sigma \right) $ is the diagram of $\left(
\pi ^{\ast },\sigma \right) $, as defined in Section \ref{SS : diagrams}, or,
whenever every $\Gamma \left( \pi ^{\ast },\sigma \right) $ involved in the previous sum is Gaussian,
\begin{equation*}
\sum_{\sigma \in \mathcal{P}\left( \left[ n\right] \right) :\hat{\Gamma}%
\left( \pi ^{\ast },\sigma \right) \text{ has no loops}}{\rm St}\,%
\!_{\sigma }^{\varphi ,\left[ n\right] }\left( f_{1}\otimes _{0}f_{2}\otimes
_{0}\cdot \cdot \cdot \otimes _{0}f_{k}\right) \text{.}
\end{equation*}%
where $\hat{\Gamma}\left( \pi ^{\ast },\sigma \right) $ is the graph of $%
\left( \pi ^{\ast },\sigma \right) $ defined in Section \ref{ss : mgr}. This
is because, thanks to Proposition \ref{P : LOOP}, the relation $\pi ^{\ast
}\wedge \sigma =\hat{0}$ indicates that $\Gamma \left( \pi ^{\ast },\sigma
\right) $ is non-flat or, equivalently in the case of Gaussian diagrams, that $\hat{\Gamma}\left( \pi ^{\ast
},\sigma \right) $ has no loops.

\bigskip

\textbf{Examples. }(i) Set $k=2$ and $n_{1}=n_{2}=1$ in Theorem \ref{T :
ProdRW}. Then, $n=2$, $\mathcal{P}\left( \left[ 2\right] \right) =\left\{
\hat{0},\hat{1}\right\} $ and $\pi ^{\ast }=\left\{ \left\{ 1\right\}
,\left\{ 2\right\} \right\} =\hat{0}$. Since $\hat{0}\wedge \hat{1}=\hat{0}$%
, (\ref{productformula}) gives immediately that, for every pair of
elementary functions $f_{1},f_{2}$,%
\begin{eqnarray}
I_{1}^{\varphi }\left( f_{1}\right) \times I_{1}^{\varphi }\left(
f_{2}\right) &=&{\rm St}\,\!_{\hat{0}}^{\varphi ,\left[ 2\right] }\left(
f_{1}\otimes _{0}f_{2}\right) +{\rm St}\,\!_{\hat{1}}^{\varphi ,\left[ 2%
\right] }\left( f_{1}\otimes _{0}f_{2}\right)  \notag \\
&=&I_{2}^{\varphi }\left( f_{1}\otimes _{0}f_{2}\right) +{\rm St}\,\!_{%
\hat{1}}^{\varphi ,\left[ 2\right] }\left( f_{1}\otimes _{0}f_{2}\right) ,
\label{simpleMult}
\end{eqnarray}%
Note that, if $\varphi =G$ is Gaussian, then relation (\ref{GaussDia})
yields that
\begin{equation*}
{\rm St}\,\!_{\hat{1}}^{G,\left[ 2\right] }\left( f_{1}\otimes
_{0}f_{2}\right) =\int_{Z}f_{1}\left( z\right) f_{2}\left( z\right) \nu
\left( dz\right) ,
\end{equation*}%
so that, in integral notation,
\begin{equation*}
I_{1}^{G}\left( f_{1}\right) \times I_{1}^{G}\left( f_{2}\right) =\int
\int_{z_{1}\neq z_{2}}f_{1}\left( z_{1}\right) f_{2}\left( z_{2}\right)
G\left( dz_{1}\right) G\left( dz_{2}\right) +\int_{Z}f_{1}\left( z\right)
f_{2}\left( z\right) \nu \left( dz\right) \text{.}
\end{equation*}%
On the other hand, if $\varphi $ is compensated Poisson, then ${\rm St}\,%
\!_{\hat{1}}^{\varphi ,\left[ 2\right] }\left( f_{1}\otimes _{0}f_{2}\right)
=\int_{Z}f_{1}\left( z\right) f_{2}\left( z\right) N\left( dz\right) $, so
that, by using the relation $N=\hat{N}+\nu $, (\ref{simpleMult}) reads%
\begin{eqnarray*}
I_{1}^{\hat{N}}\left( f_{1}\right) \times I_{1}^{\hat{N}}\left( f_{2}\right)
&=&I_{2}^{\hat{N}}\left( f_{1}\otimes _{0}f_{2}\right) +\int_{z}f_{1}\left(
z\right) f_{2}\left( z\right) \hat{N}\left( dz\right) +\int_{z}f_{1}\left(
z\right) f_{2}\left( z\right) \nu \left( dz\right) \\
&=&\int \int_{z_{1}\neq z_{2}}f_{1}\left( z_{1}\right) f_{2}\left(
z_{2}\right) \hat{N}\left( dz_{1}\right) \hat{N}\left( dz_{2}\right) \\
&&+\int_{z}f_{1}\left( z\right) f_{2}\left( z\right) \hat{N}\left( dz\right)
+\int_{z}f_{1}\left( z\right) f_{2}\left( z\right) \nu \left( dz\right) \\
&=&I_{2}^{\hat{N}}\left( f_{1}\otimes _{0}f_{2}\right) +I_{1}^{\hat{N}%
}\left( f_{1}f_{2}\right) +\int_{z}f_{1}\left( z\right) f_{2}\left( z\right)
\nu \left( dz\right)
\end{eqnarray*}

(ii) Consider the case $k=2$, $n_{1}=2$ and $n_{2}=1$. Then, $n=3$, and $\pi
^{\ast }=\left\{ \left\{ 1,2\right\} ,\left\{ 3\right\} \right\} $. There
are three elements $\sigma _{1},\sigma _{2}$, $\sigma _{3}\in \mathcal{P}%
\left( \left[ 3\right] \right) $ such that $\sigma _{i}\wedge \pi ^{\ast }=%
\hat{0}$, namely $\sigma _{1}=\hat{0}$, $\sigma _{2}=\left\{ \left\{
1,3\right\} ,\left\{ 2\right\} \right\} $ and $\sigma _{3}=\left\{ \left\{
1\right\} ,\left\{ 2,3\right\} \right\} .$ Then, (\ref{productformula})
gives that, for every pair $f_{1}\in \mathcal{E}_{s,0}\left( \nu ^{2}\right)
,f_{2}\in \mathcal{E}\left( \nu \right) $,%
\begin{eqnarray*}
I_{2}^{\varphi }\left( f_{1}\right) \times I_{1}^{\varphi }\left(
f_{2}\right) &=&{\rm St}\,\!_{\hat{0}}^{\varphi ,\left[ 3\right] }\left(
f_{1}\otimes _{0}f_{2}\right) +{\rm St}\,\!_{\sigma _{2}}^{\varphi ,\left[
3\right] }\left( f_{1}\otimes _{0}f_{2}\right) +{\rm St}\,\!_{\sigma
_{3}}^{\varphi ,\left[ 3\right] }\left( f_{1}\otimes _{0}f_{2}\right) . \\
&=&{\rm St}\,\!_{\hat{0}}^{\varphi ,\left[ 3\right] }\left( f_{1}\otimes
_{0}f_{2}\right) +2{\rm St}\,\!_{\sigma _{2}}^{\varphi ,\left[ 3\right]
}\left( f_{1}\otimes _{0}f_{2}\right) ,
\end{eqnarray*}%
where we have used the fact that, by the symmetry of $f_{1}$, ${\rm St}\,%
\!_{\sigma _{2}}^{\varphi ,\left[ 3\right] }\left( f_{1}\otimes
_{0}f_{2}\right) ={\rm St}\,\!_{\sigma _{3}}^{\varphi ,\left[ 3\right]
}\left( f_{1}\otimes _{0}f_{2}\right) $.

When $\varphi =G$ is a Gaussian measure, one can use (\ref{m-DG2}) applied
to $\sigma _{2}$ to deduce that
\begin{equation*}
{\rm St}\,\!_{\sigma _{2}}^{G,\left[ 3\right] }\left( f_{1}\otimes
_{0}f_{2}\right) =I_{1}^{G}\left[ \int_{Z}f_{1}\left( \cdot ,z\right)
f_{2}\left( z\right) \nu \left( dz\right) \right] \text{,}
\end{equation*}%
or, more informally,
\begin{equation*}
{\rm St}\,\!_{\sigma _{2}}^{G,\left[ 3\right] }\left( f_{1}\otimes
_{0}f_{2}\right) =\int_{Z}\int_{Z}f_{1}\left( z^{\prime },z\right)
f_{2}\left( z^{\prime }\right) \nu \left( dz\right) G\left( dz^{\prime
}\right) ,
\end{equation*}%
so that one gets%
\begin{eqnarray*}
I_{2}^{G}\left( f_{1}\right) \times I_{1}^{G}\left( f_{2}\right) &=&{\rm
St}\!_{\hat{0}}^{G,\left[ 3\right] }\left( f_{1}\otimes _{0}f_{2}\right)
+2I_{1}^{G}\left[ \int_{Z}f_{1}\left( \cdot ,z\right) f_{2}\left( z\right)
\nu \left( dz\right) \right] \\
&=&\int \int \int_{z_{1}\neq z_{2}\neq z_{3}}f_{1}\left( z_{1},z_{2}\right)
f_{2}\left( z_{3}\right) G\left( dz_{1}\right) G\left( dz_{2}\right) G\left(
dz_{3}\right) \\
&&+2\int_{Z}\int_{Z}f_{1}\left( z^{\prime },z\right) f_{2}\left( z\right)
\nu \left( dz\right) G\left( dz^{\prime }\right) .
\end{eqnarray*}

When $\varphi =\hat{N}$ is compensated Poisson, as shown in (\ref{e : I21}),
formula (\ref{=sigmaPoisson}), applied to $\sigma _{2}$, yields%
\begin{equation*}
{\rm St}\,\!_{\sigma _{2}}^{\hat{N},\left[ 3\right] }\left( f_{1}\otimes
_{0}f_{2}\right) =I_{1}^{\hat{N}}\left[ \int_{Z}f_{1}\left( \cdot ,z\right)
f_{2}\left( z\right) \nu \left( dz\right) \right] +I_{2}^{\hat{N}}\left[
f_{1}\otimes _{1}^{0}f_{2}\right] \text{,}
\end{equation*}%
where $f_{1}\otimes _{1}^{0}f_{2}\left( z^{\prime },z\right) =f_{1}\left(
z^{\prime },z\right) f_{2}\left( z\right) $.

(iii) Consider the case $k=3$, $n_{1}=n_{2}=n_{3}=1$. Then, $n=3$, and $\pi
^{\ast }=\left\{ \left\{ 1\right\} ,\left\{ 2\right\} ,\left\{ 3\right\}
\right\} =\hat{0}$. For every $\sigma \in \mathcal{P}\left( \left[ 3\right]
\right) $ one has that $\sigma \wedge \pi ^{\ast }=\hat{0}.$ Note also that $%
\mathcal{P}\left( \left[ 3\right] \right) =\left\{ \hat{0},\rho _{1},\rho
_{2},\rho _{3},\hat{1}\right\} $, where
\begin{equation*}
\rho _{1}=\left\{ \left\{ 1,2\right\} ,\left\{ 3\right\} \right\} \text{, \ }%
\rho _{2}=\left\{ \left\{ 1,3\right\} ,\left\{ 2\right\} \right\} \text{, \ }%
\rho _{3}=\left\{ \left\{ 1\right\} ,\left\{ 2,3\right\} \right\} \text{,}
\end{equation*}%
so that (\ref{productformula}) gives that, for every $f_{1},f_{2},f_{3}\in
\mathcal{E}\left( \nu \right) $,%
\begin{eqnarray*}
I_{1}^{\varphi }\left( f_{1}\right) I_{1}^{\varphi }\left( f_{2}\right)
I_{1}^{\varphi }\left( f_{3}\right) &=&{\rm St}\,\!_{\hat{0}}^{\varphi ,%
\left[ 3\right] }\left( f_{1}\otimes _{0}f_{2}\otimes _{0}f_{3}\right) +%
{\rm St}\,\!_{\rho _{1}}^{\varphi ,\left[ 3\right] }\left( f_{1}\otimes
_{0}f_{2}\otimes _{0}f_{3}\right) \\
&&+{\rm St}\,\!_{\rho _{2}}^{\varphi ,\left[ 3\right] }\left( f_{1}\otimes
_{0}f_{2}\otimes _{0}f_{3}\right) +{\rm St}\,\!_{\rho _{3}}^{\varphi ,%
\left[ 3\right] }\left( f_{1}\otimes _{0}f_{2}\otimes _{0}f_{3}\right) \\
&&+{\rm St}\,\!_{\hat{1}}^{\varphi ,\left[ 3\right] }\left( f_{1}\otimes
_{0}f_{2}\otimes _{0}f_{3}\right) .
\end{eqnarray*}%
In particular, by taking $f_{1}=f_{2}=f_{3}=f$ and by symmetry,
\begin{equation}
I_{1}^{\varphi }\left( f\right) ^{3}={\rm St}\,\!_{\hat{0}}^{\varphi ,%
\left[ 3\right] }\left( f\otimes _{0}f\otimes _{0}f\right) +{\rm St}\,\!_{%
\hat{1}}^{\varphi ,\left[ 3\right] }\left( f\otimes _{0}f\otimes
_{0}f\right) +3{\rm St}\,\!_{\rho _{1}}^{\varphi ,\left[ 3\right] }\left(
f\otimes _{0}f\otimes _{0}f\right) .  \label{formulaj}
\end{equation}
When $\varphi =G$ is Gaussian, then ${\rm St}\,\!_{\hat{1}}^{G,\left[ 3%
\right] }=0$ by (\ref{GaussDia}) and ${\rm St}\,\!_{\rho _{1}}^{\varphi ,%
\left[ 3\right] }\left( f\otimes _{0}f\otimes _{0}f\right) =\left\Vert
f\right\Vert ^{2}I_{1}^{G}\left( f\right) $, so that (\ref{formulaj}) becomes%
\begin{eqnarray*}
I_{1}^{G}\left( f\right) ^{3} &=&I_{3}^{G}\left( f\otimes _{0}f\otimes
_{0}f\right) +3\left\Vert f\right\Vert ^{2}I_{1}^{G}\left( f\right) \\
&=&\int \int \int_{z_{1}\neq z_{2}\neq z_{3}}f\left( z_{1}\right) f\left(
z_{2}\right) f\left( z_{3}\right) G\left( dz_{1}\right) G\left(
dz_{2}\right) G\left( dz_{3}\right) \\
&&+3\int_{Z}f\left( z\right) ^{2}\nu \left( dz\right) \times \int_{Z}f\left(
z\right) G\left( dz\right) .
\end{eqnarray*}
When $\varphi =\hat{N}$ is compensated Poisson, from (\ref{DiagPoiss}), then
\begin{equation*}
{\rm St}\,\!_{\hat{1}}^{\hat{N},\left[ 3\right] }\left( f\otimes
_{0}f\otimes _{0}f\right) =\int f\left( z\right) ^{3}N\left( dz\right)
=I_{1}^{\hat{N}}\left( f^{3}\right) +\int_{Z}f^{3}\left( z\right) \nu \left(
dz\right)
\end{equation*}%
by (\ref{e : I21}), and also
\begin{equation*}
{\rm St}\,\!_{\rho _{1}}^{\hat{N},\left[ 3\right] }\left( f\otimes
_{0}f\otimes _{0}f\right) =\left\Vert f\right\Vert ^{2}I_{1}^{\hat{N}}\left(
f\right) +I_{2}^{\hat{N}}\left( f^{2}\otimes _{0}f\right) \text{,}
\end{equation*}%
where $f^{2}\otimes _{0}f\left( z,z^{\prime }\right) =f^{2}\left( z\right)
f\left( z^{\prime }\right) $, so that (\ref{formulaj}) becomes%
\begin{eqnarray*}
I_{1}^{\hat{N}}\left( f\right) ^{3} &=&I_{3}^{\hat{N}}\left( f\otimes
_{0}f\otimes _{0}f\right) +I_{1}^{\hat{N}}\left( f^{3}\right)
+\int_{Z}f^{3}\left( z\right) \nu \left( dz\right) \\
&&+3\left\Vert f\right\Vert ^{2}I_{1}^{\hat{N}}\left( f\right) +I_{2}^{\hat{N%
}}\left( f^{2}\otimes _{0}f\right) \\
&=&\int \int \int_{z_{1}\neq z_{2}\neq z_{3}}f\left( z_{1}\right) f\left(
z_{2}\right) f\left( z_{3}\right) \hat{N}\left( dz_{1}\right) \hat{N}\left(
dz_{2}\right) \hat{N}\left( dz_{3}\right) \\
&&+\int_{Z}f\left( z\right) ^{3}\left( \hat{N}+\nu \right) \left( dz\right)
+3\int_{Z}f\left( z\right) ^{2}\nu \left( dz\right) \times \int_{Z}f\left(
z\right) \hat{N}\left( dz\right) \\
&&+\int \int_{z_{1}\neq z_{2}}f^{2}\left( z_{1}\right) f\left( z_{2}\right)
\hat{N}\left( dz_{1}\right) \hat{N}\left( dz_{2}\right) \text{. }
\end{eqnarray*}

\bigskip

General applications to the Gaussian and Poisson cases are discussed,
respectively, in Subsection \ref{SS ProdG} and Subsection \ref{SS ProdP}.

\subsection{Contractions}

As anticipated, the statement of Theorem \ref{T : ProdRW} contains two
well-known \textsl{multiplication formulae}, associated with the Gaussian
and Poisson cases. In order to state these results, we shall start with a
standard definition of the \textsl{contraction kernels }associated with two
symmetric functions $f$ and $g$\textsl{. }Roughly speaking, given $f\in
L_{s}^{2}\left( \nu ^{p}\right) $ and $g\in L_{s}^{2}\left( \nu ^{q}\right) $%
, the contraction of $f$ and $g$ on $Z^{p+q-r-l}$ ($r=0,...,q\wedge p$ and $%
l=1,...,r$), noted $f\star _{r}^{l}g$, is obtained by reducing the number of
variables in the tensor product $f\left( x_{1},...,x_{p}\right) g\left(
x_{p+1},...,x_{p+q}\right) $ as follows: $r$ variables are identified, and
of these, $l$ are integrated out with respect to $\nu $. The formal
definition of $f\star _{r}^{l}g$ is given below.

\bigskip

\begin{definition}
Let $\nu $ be a $\sigma $-finite measure on $\left( Z,\mathcal{Z}\right) $.
For every $q,p\geq 1$, $f\in L^{2}\left( \nu ^{p}\right) $, $g\in
L^{2}\left( \nu ^{q}\right) $ (not necessarily symmetric), $r=0,...,q\wedge
p $ and $l=1,...,r$, the \textsl{contraction}\textit{\ (of index }$\left(
r,l\right) $) \textit{of }$f$ and $g$ on $Z^{p+q-r-l}$, is the function $%
f\star _{r}^{l}g$ of $p+q-r-l$ variables defined as follows:%
\begin{eqnarray}
&&f\star _{r}^{l}g(\gamma _{1},\ldots ,\gamma _{r-l},t_{1},\ldots
,t_{p-r},s_{1},\ldots ,s_{q-r})  \notag \\
&=&\int_{Z^{l}}f(z_{1},\ldots ,z_{l},\gamma _{1},\ldots ,\gamma
_{r-l},t_{1},\ldots ,t_{p-r})\times  \label{preCTR} \\
&&\text{ \ \ \ \ \ \ \ \ \ \ \ \ }\times g(z_{1},\ldots ,z_{l},\gamma
_{1},\ldots ,\gamma _{r-l},s_{1},\ldots ,s_{q-r})\nu ^{l}\left(
dz_{1}...dz_{l}\right) \text{.} \notag
\end{eqnarray}
and, for $l=0$,
\begin{eqnarray}
&&f\star _{r}^{0}g(\gamma _{1},\ldots ,\gamma _{r},t_{1},\ldots
,t_{p-r},s_{1},\ldots ,s_{q-r})  \label{Ocont} \\
&=&f(\gamma _{1},\ldots ,\gamma _{r},t_{1},\ldots ,t_{p-r})g(\gamma
_{1},\ldots ,\gamma _{r},s_{1},\ldots ,s_{q-r}),  \notag
\end{eqnarray}%
so that
\begin{equation*}
f\star _{0}^{0}g(t_{1},\ldots ,t_{p},s_{1},\ldots ,s_{q})=f(t_{1},\ldots
,t_{p})g(s_{1},\ldots ,s_{q}).
\end{equation*}%
For instance, if $p=q=2$, one gets%
\begin{eqnarray}
f\star _{1}^{0}g\left( \gamma ,t,s\right) &=&f\left( \gamma ,t\right)
g\left( \gamma ,s\right) \text{, \ \ }f\star _{1}^{1}g\left( t,s\right)
=\int_{Z}f\left( z,t\right) g\left( z,s\right) \nu \left( dz\right)
\label{exCont1} \\
f\star _{2}^{1}g\left( \gamma \right) &=&\int_{Z}f\left( z,\gamma \right)
g\left( z,\gamma \right) \nu \left( dz\right) \text{, \ }  \label{exCont2} \\
\text{\ }f\star _{2}^{2}g &=&\int_{Z}\int_{Z}f\left( z_{1},z_{2}\right)
g\left( z_{1},z_{2}\right) \nu \left( dz_{1}\right) \nu \left( dz_{2}\right)
.  \label{ex}
\end{eqnarray}%
One also has
\begin{eqnarray}
&&f\star _{r}^{r}g\left( x_{1},...,x_{p+q-2r}\right)  \label{contr} \\
&=&\int_{Z^{r}}f\left( z_{1},...,z_{r},x_{1},...,x_{p-r}\right) g\left(
z_{1},...,z_{r},x_{p-r+1},...,x_{p+q-2r}\right) \nu \left( dz_{1}\right)
\cdot \cdot \cdot \nu \left( dz_{r}\right) \text{,}  \notag
\end{eqnarray}%
but, in analogy with (\ref{0multContr}), we set $\star _{r}^{r}=\otimes _{r}$%
, and consequently write%
\begin{equation}
f\star _{r}^{r}g\left( x_{1},...,x_{p+q-2r}\right) =f\otimes _{r}g\left(
x_{1},...,x_{p+q-2r}\right) \text{,}  \label{croceContr}
\end{equation}%
so that, in particular,
\begin{equation*}
f\star _{0}^{0}g=f\otimes _{0}g.
\end{equation*}
\end{definition}

\bigskip

The following elementary result is proved by using the Cauchy-Schwarz
inequality. It ensures that the contractions of the type (\ref{croceContr})
are still square-integrable kernels.

\begin{lemma}
Let $f\in L^{2}\left( \nu ^{p}\right) $ and $g\in L^{2}\left( \nu
^{q}\right) $. Then, for every $r=0,...,p\wedge q$, one has that $f\otimes
_{r}g\in L^{2}\left( \nu ^{p+q-2r}\right) $.
\end{lemma}

\begin{proof}
Just write
\begin{eqnarray*}
&&\int_{Z^{p+q-2r}}\left( f\otimes _{r}g\right) ^{2}d\nu ^{p+q-2r} \\
&=&\int_{Z^{p+q-2r}}\left( \int_{Z^{r}}f\left(
a_{1},...,a_{r},z_{1},...,z_{p-r}\right) \right. \\
&&\text{ \ \ }\left. g\left( a_{1},...,a_{r},z_{p-r+1},...,z_{p+q-r}\right)
^{^{^{{}}}}\nu ^{r}\left( da_{1},...da_{r}\right) \right) ^{2}\nu
^{p+q-2r}(dz_{1},...,dz_{p+q-r}) \\
&\leq &\left\Vert f\right\Vert _{L^{2}\left( \nu ^{p}\right) }^{2}\times
\left\Vert g\right\Vert _{L^{2}\left( \nu ^{q}\right) }^{2}.
\end{eqnarray*}
\end{proof}

\subsection{Symmetrization of contractions}

Suppose that $f\in L^{2}\left( \nu ^{p}\right) $ and $g\in L^{2}\left( \nu
^{q}\right) $, and let \textquotedblleft $\ \ \widetilde{}$ \
\textquotedblright\ denote symmetrization. Then $f=\widetilde{f}$ and $g=%
\widetilde{g}$. However, in general, the fact that $f$ and $g$ are symmetric
\textsl{does not }imply that the contraction $f\otimes _{r}g$ is symmetric.
For instance, if $p=q=1$,%
\begin{equation*}
\widetilde{f\otimes _{0}g}\left( s,t\right) =\frac{1}{2}\left[ f\left(
s\right) g\left( t\right) +g\left( s\right) f\left( t\right) \right] \text{;}
\end{equation*}%
if $p=q=2$%
\begin{equation*}
\widetilde{f\otimes _{1}g}\left( s,t\right) =\frac{1}{2}\int_{Z}\left[
f\left( x,s\right) g\left( x,t\right) +g\left( x,s\right) f\left( x,t\right) %
\right] \nu \left( dx\right) .
\end{equation*}%
In general, due to the symmetry of $f$ and $g$, for every $p,q\geq 1$ and
every $r=0,...,p\wedge q$ one has the relation
\begin{eqnarray*}
\widetilde{f\otimes _{r}g}\left( t_{1},...,t_{p+q-2r}\right) &=&\frac{1}{%
\binom{p+q-2r}{p-r}}\times \\
&&\times \sum_{1\leq i_{1}<\cdot \cdot \cdot <i_{p-r}\leq
p+q-2r}\int_{Z^{r}}f\left( \mathbf{t}_{\left( i_{1},...,i_{p-r}\right) },%
\mathbf{a}_{r}\right) g\left( \mathbf{t}_{\left( i_{1},...,i_{p-r}\right)
^{c}},\mathbf{a}_{r}\right) \nu ^{r}\left( d\mathbf{a}_{r}\right) ,
\end{eqnarray*}%
where we used the shorthand notation%
\begin{eqnarray*}
\mathbf{t}_{\left( i_{1},...,i_{p-r}\right) } &=&\left(
t_{i_{1}},...,t_{i_{p-r}}\right) \\
\mathbf{t}_{\left( i_{1},...,i_{p-r}\right) ^{c}} &=&\left(
t_{1},...,t_{p+q-2r}\right) \backslash \left(
t_{i_{1}},...,t_{i_{p-r}}\right) \\
\mathbf{a}_{r} &=&\left( a_{1},...,a_{r}\right) \\
\nu ^{r}\left( d\mathbf{a}_{r}\right) &=&\nu ^{r}\left(
da_{1},...,da_{r}\right) \text{.}
\end{eqnarray*}

Using the definition (\ref{preCTR}), one has also that $\widetilde{f\star
_{r}^{l}g}$ indicates the symmetrization of $f\star _{r}^{l}g$, where $l<r$.
For instance, if $p=3$, $q=2$, $r=2$ and $l=1$, one has that
\begin{equation*}
f\star _{r}^{l}g\left( s,t\right) =f\star _{2}^{1}g\left( s,t\right)
=\int_{Z}f\left( z,s,t\right) g\left( z,s\right) \nu \left( dz\right) \text{,%
}
\end{equation*}%
and consequently, since $f$ is symmetric,
\begin{equation*}
\widetilde{f\star _{r}^{l}g}\left( s,t\right) =\widetilde{f\star _{2}^{1}g}%
\left( s,t\right) =\frac{1}{2}\int_{Z}\left[ f\left( z,s,t\right) g\left(
z,s\right) +f\left( z,s,t\right) g\left( z,t\right) \right] \nu \left(
dz\right) \text{.}
\end{equation*}

\subsection{The product of two integrals in the Gaussian case \label{SS
ProdG}}

The main result of this section is the following general formula for
products of Gaussian multiple integrals.

\begin{proposition}
\label{P : 2bleprodGauss}Let $\varphi =G$ be a centered Gaussian measure
with $\sigma $-finite and non-atomic control measure $\nu $. Then, for every
$q,p\geq 1$, $f\in L_{s}^{2}\left( \nu ^{p}\right) $ and $g\in
L_{s}^{2}\left( \nu ^{q}\right) $,
\begin{equation}
I_{p}^{G}\left( f\right) I_{q}^{G}\left( g\right) =\sum_{r=0}^{p\wedge q}r!%
\dbinom{p}{r}\dbinom{q}{r}I_{p+q-2r}^{G}\left( \widetilde{f\otimes _{r}g}%
\right) \text{,}  \label{ProdGauss}
\end{equation}%
where the symbol ( $\widetilde{}$ ) indicates a symmetrization, the
contraction $f\otimes _{r}g$ is defined in (\ref{croceContr}), and for $%
p=q=r $, we write
\begin{eqnarray*}
&&I_{0}^{G}\left( \widetilde{f\otimes _{p}g}\right) \\
&=&f\otimes _{r}g=\int_{Z^{p}}f\left( z_{1},...,z_{p}\right) g\left(
z_{1},...,z_{p}\right) \nu \left( dz_{1}\right) \cdot \cdot \cdot \nu \left(
dz_{p}\right) \\
&=&(f,g)_{L^{2}\left( \nu ^{p}\right) }.
\end{eqnarray*}
\end{proposition}

\textbf{Remark. }Since, in general, one has that $I_{n}^{G}( \tilde{%
h}) =I_{n}^{G}\left( h\right) $ (see formula (\ref{Yo la})), one could
dispense with the symmetrization \ \textquotedblleft\ \ $\widetilde{}$\ \
\textquotedblright\ in formula (\ref{ProdGauss}).

\begin{proof}[Proof of Proposition \protect\ref{P : 2bleprodGauss}]
We start by assuming that $f\in \mathcal{E}_{s,0}\left( \nu ^{p}\right) $
and $g\in \mathcal{E}_{s,0}\left( \nu ^{q}\right) $, and we denote by $\pi
^{\ast }$ the partition of the set $\left[ p+q\right] =\left\{
1,...,p+q\right\} $ given by
\begin{equation*}
\pi ^{\ast }=\left\{ \left\{ 1,...,p\right\} ,\left\{ p+1,...,p+q\right\}
\right\} .
\end{equation*}%
According to formula (\ref{productformula})
\begin{equation*}
I_{p}^{G}\left( f\right) I_{q}^{G}\left( g\right) =\sum_{\sigma \in \mathcal{%
P}\left( \left[ p+q\right] \right) :\sigma \wedge \pi ^{\ast }=\hat{0}}%
{\rm St}\,\!_{\sigma }^{G,\left[ n\right] }\left( f\otimes _{0}g\right) .
\end{equation*}%
Every partition $\sigma \in \mathcal{P}\left( \left[ p+q\right] \right) $
such that $\sigma \wedge \pi ^{\ast }=\hat{0}$ is necessarily composed of $r$
($0\leq r\leq p\wedge q$) two-elements blocks of the type $\left\{
i,j\right\} $ where $i\in \left\{ 1,...,p\right\} $ and $j\in \left\{
p+1,...,p+q\right\} $, and $p+q-2r$ singletons. Moreover, for every fixed $%
r\in \left\{ 0,...,p\wedge q\right\} $, there are exactly $r!\dbinom{p}{r}%
\dbinom{q}{r}$ partitions of this type. To see this, observe that, to build
such a partition, one should first select one of the $\dbinom{p}{r}$ subsets
of size $r$ of $\left\{ 1,...,p\right\} $, say $A_{r}$, and a one of the $%
\dbinom{q}{r}$ subset of size $r$ of $\left\{ p+1,...,p+q\right\} $, say $%
B_{r}$, and then choose one of the $r!$ bijections between $A_{r}$ and $%
B_{r} $. When $r=0$, and therefore $\sigma =\hat{0}$, one obtains
immediately
\begin{equation*}
{\rm St}\,\!_{\sigma }^{G,\left[ p+q\right] }\left( f\otimes _{0}g\right) =%
{\rm St}\,\!_{\hat{0}}^{G,\left[ p+q\right] }\left( \widetilde{f\otimes
_{0}g}\right) =I_{p+q}^{G}\left( \widetilde{f\otimes _{0}g}\right) ,
\end{equation*}%
where the first equality is a consequence of the symmetry of ${\rm St}\,%
\!_{\hat{0}}^{G,\left[ n\right] }$ (see Proposition \ref{P : ST0RW}). On the
other hand, we claim that every partition $\sigma \in \mathcal{P}\left( %
\left[ p+q\right] \right) $ such that $\sigma \wedge \pi ^{\ast }=\hat{0}$
and $\sigma $ contains $r\geq 1$ two-elements blocks of the type $b=\{i,j\}$ (with $i\in\{1,...,p\}$ and $j\in\{p+1,...,p+q\}$) and $p+q-2r$
singletons, satisfies also
\begin{equation}
{\rm St}\,\!_{\sigma }^{G,\left[ p+q\right] }\left( f\otimes _{0}g\right) =%
{\rm St}\,\!_{\hat{0}}^{G,\left[ p+q-2r\right] }\left( f\otimes
_{r}g\right) =I_{p+q-2r}^{G}\left( \widetilde{f\otimes _{r}g}\right).
\label{string}
\end{equation}%
We give a proof of (\ref{string}). Consider first the (not necessarily
symmetric) functions
\begin{equation*}
f%
%TCIMACRO{\U{b0}}%
%BeginExpansion
{{}^\circ}%
%EndExpansion
=\mathbf{1}_{A_{1}\times \cdot \cdot \cdot \times A_{p}}\text{ \ and \ \ }g%
%TCIMACRO{\U{b0}}%
%BeginExpansion
{{}^\circ}%
%EndExpansion
=\mathbf{1}_{A_{p+1}\times \cdot \cdot \cdot \times A_{p+q}},
\end{equation*}%
where $A_{l}\in \mathcal{Z}_{\nu }$, $l=1,...,p+q$. Then, one may use (\ref%
{m-DG2}), in order to obtain
\begin{eqnarray*}
{\rm St}\,\!_{\sigma }^{G,\left[ p+q\right] }\left( f%
%TCIMACRO{\U{b0}}%
%BeginExpansion
{{}^\circ}%
%EndExpansion
\otimes _{0}g%
%TCIMACRO{\U{b0}}%
%BeginExpansion
{{}^\circ}%
%EndExpansion
\right) &=&\prod_{b=\left\{ i,j\right\} \in \sigma }\nu \left( A_{i}\cap
A_{j}\right) {\rm St}\,{}_{\hat{0}}^{G,\left[ p+q-2r\right] }\left(
A_{j_{1}}\times \cdot \cdot \cdot \times A_{j_{p+q-2r}}\right) \\
&=&{\rm St}\,{}_{\hat{0}}^{G,\left[ p+q-2r\right] }\left( \prod_{b=\left\{
i,j\right\} \in \sigma }\nu \left( A_{i}\cap A_{j}\right) \mathbf{1}%
_{A_{j_{1}}\times \cdot \cdot \cdot \times A_{j_{p+q-2r}}}\right) ,
\end{eqnarray*}%
where $\left\{ \left\{ j_{1}\right\} ,...,\left\{ j_{p+q-2r}\right\}
\right\} $ are the singletons of $\sigma $ (by the symmetry of ${\rm St}\,%
\!_{\hat{0}}^{G,\left[ p+q-2r\right] }$ we can always suppose, here and in
the following, that the singletons $\left\{ j_{1}\right\} ,...,\left\{
j_{p-r}\right\} $ are contained in $\left\{ 1,...,p\right\} $ and that the
singletons $\left\{ j_{p-r+1}\right\} ,...,\left\{ j_{p+q-2r}\right\} $ are
in $\left\{ p+1,...,p+q\right\} $). For every pair of permutations $w\in
\mathfrak{S}_{\left[ p\right] }$ and $w^{\prime }\in \mathfrak{S}%
_{[p+1,p+q]} $ (the group of permutations of the set $[p+1,p+q]=\left\{
p+1,...,p+q\right\} $), we define $f%
%TCIMACRO{\U{b0}}%
%BeginExpansion
{{}^\circ}%
%EndExpansion
^{,w}=$ $\mathbf{1}_{A_{1}^{w}\times \cdot \cdot \cdot \times A_{p}^{w}}$
and $g%
%TCIMACRO{\U{b0}}%
%BeginExpansion
{{}^\circ}%
%EndExpansion
^{,w^{\prime }}=\mathbf{1}_{A_{p+1}^{w^{\prime }}\times \cdot \cdot \cdot
\times A_{p+q}^{w^{\prime }}}$, where $A_{j}^{w}=A_{w\left( j\right) }$, $%
j=1,...,p$ (and analogously for $w^{\prime }$). In this way,
\begin{eqnarray}
&&{\rm St}\,\!_{\sigma }^{G,\left[ p+q\right] }\left( f%
%TCIMACRO{\U{b0}}%
%BeginExpansion
{{}^\circ}%
%EndExpansion
^{,w}\otimes _{0}g%
%TCIMACRO{\U{b0}}%
%BeginExpansion
{{}^\circ}%
%EndExpansion
^{,w^{\prime }}\right)  \label{OzerO} \\
&=&{\rm St}\,{}_{\hat{0}}^{G,\left[ p+q-2r\right] }\left( \prod_{b=\left\{
i,j\right\} \in \sigma }\nu \left( A_{i}^{w}\cap A_{j}^{w'}\right) \mathbf{1}%
_{A_{j_{1}}^{w}\times \cdot \cdot \cdot \times A_{j_{p-r}}^{w}\times
A_{j_{p-r+1}}^{w'}\times \cdot \cdot \cdot \times A_{j_{p+q-2r}}^{w^{\prime
}}}\right)  \notag \\
&=&{\rm St}\,{}_{\hat{0}}^{G,\left[ p+q-2r\right] }\left( \prod_{b=\left\{
i,j\right\} \in \sigma }\nu \left( A_{i}^{w}\cap A_{j}^{w'}\right) \widetilde{%
\mathbf{1}_{A_{j_{1}}^{w}\times \cdot \cdot \cdot \times A_{j_{p-r}}^{w}\times
A_{j_{p-r+1}}^{w'}\times \cdot \cdot \cdot \times A_{j_{p+q-2r}}^{w^{\prime }}}}%
\right) .  \notag
\end{eqnarray}%
Now write
\begin{equation*}
f=\sum_{w\in \mathfrak{G}_{\left[ p\right] }}f%
%TCIMACRO{\U{b0}}%
%BeginExpansion
{{}^\circ}%
%EndExpansion
^{,w}\text{ \ and \ }g=\sum_{w\in \mathfrak{G}_{\left[ p+1,p+q\right] }}g%
%TCIMACRO{\U{b0}}%
%BeginExpansion
{{}^\circ}%
%EndExpansion
^{,w},
\end{equation*}%
and observe that (by using (\ref{croceContr}))
\begin{equation*}
\sum_{w\in \mathfrak{G}_{\left[ p\right] }}\sum_{w'\in \mathfrak{G}_{\left[
p+1,p+q\right] }}\prod_{b=\left\{ i,j\right\} \in \sigma }\nu \left(
A_{i}^{w}\cap A_{j}^{w'}\right) \widetilde{\mathbf{1}_{A_{j_{1}}^{w}\times
\cdot \cdot \cdot \times A_{j_{p-r}}^{w}\times A_{j_{p-r+1}}^{w'}\times \cdot
\cdot \cdot \times A_{j_{p+q-2r}}^{w^{\prime }}}}=\widetilde{f\otimes _{r}g}.
\end{equation*}%
Since (\ref{OzerO}) gives
\begin{eqnarray*}
&&{\rm St}\,\!_{\sigma }^{G,\left[ p+q\right] }\left( f\otimes _{0}g\right)
\\
&=&{\rm St}\,{}_{\hat{0}}^{G,\left[ p+q-2r\right] }\left( \sum_{w\in
\mathfrak{S}_{\left[ p\right] }}\sum_{w'\in \mathfrak{S}_{\left[ p+1,p+q%
\right] }}\prod_{b=\left\{ i,j\right\} \in \sigma }\nu \left( A_{i}^{w}\cap
A_{j}^{w'}\right) \widetilde{\mathbf{1}_{A_{j_{1}}^{w}\times \cdot \cdot
\cdot \times A_{j_{p-r}}^{w}\times A_{j_{p-r+1}}^{w'}\times \cdot \cdot \cdot
\times A_{j_{p+q-2r}}^{w^{\prime }}}}\right) ,
\end{eqnarray*}%
we obtain (\ref{string}), so that, in particular, (\ref{ProdGauss}) is
proved for symmetric simple functions vanishing on diagonals. The general
result is obtained by using the fact that the linear spaces $\mathcal{E}%
_{s,0}\left( \nu ^{p}\right) $ and $\mathcal{E}_{s,0}\left( \nu ^{q}\right) $
are dense, respectively, in $L_{s}^{2}\left( \nu ^{p}\right) $ and $%
L_{s}^{2}\left( \nu ^{q}\right) $. Indeed, to conclude the proof it is
sufficient to observe that, if $\left\{ f_{k}\right\} \subset \mathcal{E}%
_{s,0}\left( \nu ^{p}\right) $ and $\left\{ g_{k}\right\} \subset \mathcal{E}%
_{s,0}\left( \nu ^{q}\right) $ are such that $f_{k}\rightarrow f$ in $%
L_{s}^{2}\left( \nu ^{p}\right) $ and $g_{k}\rightarrow g$ in $%
L_{s}^{2}\left( \nu ^{q}\right) $, then, for instance by Cauchy-Schwarz, $%
I_{p}^{G}\left( f_{k}\right) I_{q}^{G}\left( g_{k}\right) \rightarrow
I_{p}^{G}\left( f\right) I_{q}^{G}\left( g\right) $ in any norm $L^{s}\left(
\mathbb{P}\right) $, $s\geq 1$ (use e.g. (\ref{GaussChaosCOntr})), and also
\begin{equation*}
\widetilde{f_{k}\otimes _{r}g_{k}}\rightarrow \widetilde{f\otimes _{r}g}
\end{equation*}%
in $L_{s}^{2}\left( \nu ^{p+q-2r}\right) $, so that $I_{p+q-2r}^{G}\left(
\widetilde{f_{k}\otimes _{r}g_{k}}\right) \rightarrow I_{p+q-2r}^{G}\left(
\widetilde{f\otimes _{r}g}\right) $ in $L^{2}\left( \mathbb{P}\right) $.
\end{proof}

\bigskip

Other proofs of Proposition \ref{P : 2bleprodGauss} can be found e.g. in
\cite{Major}, \cite{DMM5} or \cite[Proposition 1.1.3]{Nualart}.

\bigskip

\textbf{Examples. }(i) When $p=q=1$, one obtains
\begin{equation*}
I_{1}^{G}\left( f\right) I_{1}^{G}\left( g\right) =I_{2}^{G}\left(
\widetilde{f\otimes _{0}g}\right) +I_{0}^{G}\left( \widetilde{f\otimes _{1}g}%
\right) =I_{2}^{G}\left( \widetilde{f\otimes _{0}g}\right) +\left\langle
f,g\right\rangle _{L^{2}\left( \nu \right) }\text{,}
\end{equation*}%
which is consistent with (\ref{simpleMult}).

(ii) When $p=q=2$, one obtains
\begin{equation*}
I_{2}^{G}\left( f\right) I_{2}^{G}\left( g\right) =I_{4}^{G}\left(
\widetilde{f\otimes _{0}g}\right) +4I_{2}^{G}\left( \widetilde{f\otimes _{1}g%
}\right) +\left\langle f,g\right\rangle _{L^{2}\left( \nu \right) }.
\end{equation*}

(iii) When $p=3$ and $q=2$, one obtains
\begin{equation*}
I_{3}^{G}\left( f\right) I_{2}^{G}\left( g\right) =I_{5}^{G}\left(
\widetilde{f\otimes _{0}g}\right) +6I_{3}^{G}\left( \widetilde{f\otimes _{1}g%
}\right) +6I_{1}^{G}\left( f\otimes _{2}g\right) ,
\end{equation*}%
where $f\otimes _{2}g\left( z\right) =\int_{Z^{2}}f\left( z,x,y\right)
g\left( x,y\right) \nu \left( dx\right) \nu \left( dy\right) $.

\subsection{The product of two integrals in the Poisson case \label{SS ProdP}%
}

We now focus on the product of two Poisson integrals.

\begin{proposition}
\label{P : 2bleProdPoiss}Let $\varphi =\hat{N}$ be a compensated Poisson
measure, with $\sigma $-finite and non-atomic control measure $\nu $. Then,
for every $q,p\geq 1$, $f\in \mathcal{E}_{s,0}\left( \nu ^{p}\right) $ and $%
g\in \mathcal{E}_{s,0}\left( \nu ^{q}\right) $,%
\begin{equation}
I_{p}^{\widehat{N}}(f)I_{q}^{\widehat{N}}(g)=\sum_{r=0}^{p\wedge q}r!\dbinom{%
p}{r}\dbinom{q}{r}\sum_{l=0}^{r}\dbinom{r}{l}I_{p+q-r-l}^{\widehat{N}}(%
\widetilde{f\star _{r}^{l}g}).  \label{PoissProdct}
\end{equation}

Formula (\ref{PoissProdct}) continues to hold for functions $f\in
L_{s}^{2}\left( \nu ^{p}\right) $ and $g\in L_{s}^{2}\left( \nu ^{q}\right) $
such that $f\star _{r}^{l}g\in L^{2}\left( \nu ^{q+p-r-l}\right) $, $\forall
r=0,...,p\wedge q$, $\forall l=0,...,r.$
\end{proposition}

\noindent \begin{proof}[Sketch of the proof] We shall only prove formula (\ref{PoissProdct}) in the simple case where $%
p=q=2$. The generalization to general indices $p,q\geq 1$ (left to the
reader) does not present any particular additional difficulty, except for
the need of a rather heavy notation. We shall therefore prove that
\begin{eqnarray}
I_{2}^{\widehat{N}}(f)I_{2}^{\widehat{N}}(g) &=&\sum_{r=0}^{2}r!\dbinom{2}{r}%
\dbinom{2}{r}\sum_{l=0}^{r}\dbinom{r}{l}I_{4-r-l}^{\widehat{N}}(\widetilde{%
f\star _{r}^{l}g})  \label{-a} \\
&=&I_{4}^{\widehat{N}}(\widetilde{f\star _{0}^{0}g})  \label{a} \\
&&+4\left[ I_{3}^{\widehat{N}}(\widetilde{f\star _{1}^{0}g})+I_{2}^{\widehat{%
N}}(\widetilde{f\star _{1}^{1}g})\right]  \label{aa} \\
&&+2\left[ I_{2}^{\widehat{N}}(\widetilde{f\star _{2}^{0}g})+2I_{1}^{%
\widehat{N}}(\widetilde{f\star _{2}^{1}g})+\left\langle f,g\right\rangle
_{L^{2}\left( \nu ^{2}\right) }\right] .  \label{aaa}
\end{eqnarray}%
Moreover, by linearity, we can also assume that
\begin{equation*}
f=\mathbf{1}_{A_{1}\times A_{2}}+\mathbf{1}_{A_{2}\times A_{1}}\text{ \ and
\ }g=\mathbf{1}_{B_{1}\times B_{2}}+\mathbf{1}_{B_{2}\times B_{1}}\text{,}
\end{equation*}%
where $A_{1}\cap A_{2}=B_{1}\cap B_{2}=\varnothing $. Denote by $\pi ^{\ast
} $ the partition of $\left[ 4\right] =\left\{ 1,...,4\right\} $ given by
\begin{equation*}
\pi ^{\ast }=\left\{ \left\{ 1,2\right\} ,\left\{ 3,4\right\} \right\} \text{%
,}
\end{equation*}%
and apply the general result (\ref{productformula}) to deduce that%
\begin{equation*}
I_{2}^{\widehat{N}}(f)I_{2}^{\widehat{N}}(g)=\sum_{\sigma \in \mathcal{P}%
\left( \left[ 4\right] \right) :\sigma \wedge \pi ^{\ast }=\hat{0}}{\rm%
St}\!_{\sigma }^{\hat{N},\left[ n\right] }\left( f\otimes _{0}g\right) .
\end{equation*}%
We shall prove that
\begin{equation*}
\sum_{\sigma \in \mathcal{P}\left( \left[ 4\right] \right) :\sigma \wedge
\pi ^{\ast }=\hat{0}}{\rm St}\,\!_{\sigma }^{\hat{N},\left[ n\right]
}\left( f\otimes _{0}g\right) =(\ref{a})+(\ref{aa})+(\ref{aaa}).
\end{equation*}%
To see this, observe that the class
\begin{equation*}
\left\{ \sigma \in \mathcal{P}\left( \left[ 4\right] \right) :\sigma \wedge
\pi ^{\ast }=\hat{0}\right\}
\end{equation*}%
contains exactly 7 elements, that is:

\begin{description}
\item[(I)] the trivial partition $\hat{0}$, containing only singletons;

\item[(II)] four partitions $\sigma _{1},...,\sigma _{4}$ containing one
block of two elements and two singletons, namely
\begin{eqnarray*}
\sigma _{1} &=&\left\{ \left\{ 1,3\right\} ,\left\{ 2\right\} ,\left\{
4\right\} \right\} \text{, \ \ }\sigma _{2}=\left\{ \left\{ 1,4\right\}
,\left\{ 2\right\} ,\left\{ 3\right\} \right\} \\
\sigma _{3} &=&\left\{ \left\{ 1\right\} ,\left\{ 2,3\right\} ,\left\{
4\right\} \right\} \text{ \ and\ \ }\sigma _{4}=\left\{ \left\{ 1\right\}
,\left\{ 2,4\right\} ,\left\{ 3\right\} \right\} \text{;}
\end{eqnarray*}

\item[(III)] two partitions $\sigma _{5},\sigma _{6}$ composed of two blocks
of two elements, namely%
\begin{equation*}
\sigma _{5}=\left\{ \left\{ 1,3\right\} ,\left\{ 2,4\right\} \right\} \text{
\ and \ }\sigma _{6}=\left\{ \left\{ 1,4\right\} ,\left\{ 2,3\right\}
\right\} .
\end{equation*}
\end{description}

By definition, one has that
\begin{equation*}
{\rm St}\,\!_{\hat{0}}^{\hat{N},\left[ n\right] }\left( f\otimes
_{0}g\right) =I_{4}^{\widehat{N}}(\widetilde{f\star _{0}^{0}g}),
\end{equation*}%
giving (\ref{a}). Now consider the partition $\sigma _{1}$, as defined in
Point {\bf (II)} above. By using the notation (\ref{heavy not}), one has that $%
\mathbf{B}_{2}\left( \sigma _{1}\right) =\left\{ \left\{ 1,3\right\}
\right\} $, $\left\vert \mathbf{B}_{2}\left( \sigma _{1}\right) \right\vert
=1$ and $\mathbf{PB}_{2}\left( \sigma _{1}\right) =\varnothing $. It follows
from formula (\ref{=sigmaPoisson}) that%
\begin{eqnarray*}
&&{\rm St}\,\!_{\sigma _{1}}^{\hat{N},\left[ 4\right] }\left( f\otimes
_{0}g\right) ={\rm St}\,\!_{\sigma _{1}}^{\hat{N},\left[ 4\right] }\left(
\left( \mathbf{1}_{A_{1}\times A_{2}}+\mathbf{1}_{A_{2}\times A_{1}}\right)
\otimes _{0}\left( \mathbf{1}_{B_{1}\times B_{2}}+\mathbf{1}_{B_{2}\times
B_{1}}\right) \right) \\
&=&{\rm St}\,{}_{\hat{0}}^{\hat{N},\left[ 3\right] }\left( \mathbf{1}%
_{\left( A_{1}\cap B_{1}\right) \times A_{2}\times B_{2}}+\mathbf{1}_{\left(
A_{1}\cap B_{2}\right) \times A_{2}\times B_{1}}+\mathbf{1}_{\left(
A_{2}\cap B_{1}\right) \times A_{1}\times B_{2}}+\mathbf{1}_{\left(
A_{2}\cap B_{2}\right) \times A_{1}\times B_{1}}\right) \\
&&+\nu \left( A_{1}\cap B_{1}\right) {\rm St}\,{}_{\hat{0}}^{\hat{N},\left[
2\right] }\left( \mathbf{1}_{A_{2}\times B_{2}}\right) +\nu \left( A_{1}\cap
B_{2}\right) {\rm St}\,{}_{\hat{0}}^{\hat{N},\left[ 2\right] }\left(
\mathbf{1}_{A_{2}\times B_{1}}\right) \\
&&+\nu \left( A_{2}\cap B_{1}\right) {\rm St}\,{}_{\hat{0}}^{\hat{N},\left[
2\right] }\left( \mathbf{1}_{A_{1}\times B_{2}}\right) +\nu \left( A_{2}\cap
B_{2}\right) {\rm St}\,{}_{\hat{0}}^{\hat{N},\left[ 2\right] }\left(
\mathbf{1}_{A_{1}\times B_{1}}\right) .
\end{eqnarray*}%
Observe that
\begin{eqnarray*}
&&{\rm St}\,{}_{\hat{0}}^{\hat{N},\left[ 3\right] }\left( \mathbf{1}%
_{\left( A_{1}\cap B_{1}\right) \times A_{2}\times B_{2}}+\mathbf{1}_{\left(
A_{1}\cap B_{2}\right) \times A_{2}\times B_{1}}+\mathbf{1}_{\left(
A_{2}\cap B_{1}\right) \times A_{1}\times B_{2}}+\mathbf{1}_{\left(
A_{2}\cap B_{2}\right) \times A_{1}\times B_{1}}\right) \\
&=&I_{3}^{\widehat{N}}( \widetilde{f\star _{1}^{0}g}) \text{,}
\end{eqnarray*}%
and moreover,%
\begin{eqnarray*}
&&\nu \left( A_{1}\cap B_{1}\right) {\rm St}\,{}_{\hat{0}}^{\hat{N},\left[
2\right] }\left( \mathbf{1}_{A_{2}\times B_{2}}\right) +\nu \left( A_{1}\cap
B_{2}\right) {\rm St}\,{}_{\hat{0}}^{\hat{N},\left[ 2\right] }\left(
\mathbf{1}_{A_{2}\times B_{1}}\right) + \\
&&\nu \left( A_{2}\cap B_{1}\right) {\rm St}\,{}_{\hat{0}}^{\hat{N},\left[
2\right] }\left( \mathbf{1}_{A_{1}\times B_{2}}\right) +\nu \left( A_{2}\cap
B_{2}\right) {\rm St}\,{}_{\hat{0}}^{\hat{N},\left[ 2\right] }\left(
\mathbf{1}_{A_{1}\times B_{1}}\right) \\
&=&I_{2}^{\widehat{N}}(f\star _{1}^{1}g)=I_{2}^{\widehat{N}}(
\widetilde{f\star _{1}^{1}g}) \text{.}
\end{eqnarray*}%
By repeating exactly same argument, one sees immediately that%
\begin{equation*}
{\rm St}\,\!_{\sigma _{1}}^{\hat{N},\left[ 4\right] }\left( f\otimes
_{0}g\right) ={\rm St}\,\!_{\sigma _{i}}^{\hat{N},\left[ 4\right] }\left(
f\otimes _{0}g\right) \text{,}
\end{equation*}%
for every for $i=2,3,4$ (the partitions $\sigma _{i}$ being defined as in
Point (II) above) so that the quantity%
\begin{equation*}
\sum_{i=1,...,4}{\rm St}\,\!_{\sigma _{i}}^{\hat{N},\left[ n\right]
}\left( f\otimes _{0}g\right)
\end{equation*}%
equals necessarily the expression appearing in (\ref{aa})$.$ Now we focus on
the partition $\sigma _{5}$ appearing in Point (III). Plainly (by using once
again the notation introduced in (\ref{heavy not})), $\mathbf{B}_{2}\left(
\sigma _{5}\right) =\left\{ \left\{ 1,3\right\} ,\left\{ 2,4\right\}
\right\} $, $\left\vert \mathbf{B}_{2}\left( \sigma _{5}\right) \right\vert
=2$, and the set $\mathbf{PB}_{2}\left( \sigma _{5}\right) $ contains two
elements, namely%
\begin{equation*}
\left( \left\{ \left\{ 1,3\right\} \right\} ;\left\{ \left\{ 2,4\right\}
\right\} \right) \text{ \ \ and \ \ }\left( \left\{ \left\{ 2,4\right\}
\right\} ;\left\{ \left\{ 1,3\right\} \right\} \right)
\end{equation*}%
(note that we write $\left\{ \left\{ 1,3\right\} \right\} $ (with two
accolades), since the elements of $\mathbf{PB}_{2}\left( \sigma _{5}\right) $
are pairs of collections of blocks of $\sigma _{2}$, so that $\left\{
\left\{ 1,3\right\} \right\} $ is indeed the singleton whose only element is
$\left\{ 1,3\right\} $). We can now apply formula (\ref{=sigmaPoisson}) to
deduce that%
\begin{eqnarray}
&&{\rm St}\,\!_{\sigma _{5}}^{\hat{N},\left[ 4\right] }\left( f\otimes
_{0}g\right) ={\rm St}\,\!_{\sigma _{5}}^{\hat{N},\left[ 4\right] }\left(
\left( \mathbf{1}_{A_{1}\times A_{2}}+\mathbf{1}_{A_{2}\times A_{1}}\right)
\otimes _{0}\left( \mathbf{1}_{B_{1}\times B_{2}}+\mathbf{1}_{B_{2}\times
B_{1}}\right) \right)  \label{z} \\
&=&2\left[ \nu \left( A_{1}\cap B_{1}\right) {\rm St}\,{}_{\hat{0}}^{\hat{N%
},\left[ 1\right] }\left( \mathbf{1}_{A_{2}\cap B_{2}}\right) +\nu \left(
A_{1}\cap B_{2}\right) {\rm St}\,{}_{\hat{0}}^{\hat{N},\left[ 1\right]
}\left( \mathbf{1}_{A_{2}\cap B_{1}}\right) \right.  \notag \\
&&\left. +\nu \left( A_{2}\cap B_{1}\right) {\rm St}\,{}_{\hat{0}}^{\hat{N}%
,\left[ 1\right] }\left( \mathbf{1}_{A_{1}\cap B_{2}}\right) +\nu \left(
A_{2}\cap B_{2}\right) {\rm St}\,{}_{\hat{0}}^{\hat{N},\left[ 1\right]
}\left( \mathbf{1}_{A_{1}\cap B_{1}}\right) \right]  \notag \\
&&+{\rm St}\,{}_{\hat{0}}^{\hat{N},\left[ 2\right] }\left( \mathbf{1}%
_{\left( A_{1}\cap B_{1}\right) \times \left( A_{2}\cap B_{2}\right) }+%
\mathbf{1}_{\left( A_{1}\cap B_{2}\right) \times \left( A_{2}\cap
B_{1}\right) }+\mathbf{1}_{\left( A_{2}\cap B_{1}\right) \times \left(
A_{1}\cap B_{2}\right) }+\mathbf{1}_{\left( A_{2}\cap B_{2}\right) \times
\left( A_{1}\cap B_{1}\right) }\right)  \notag \\
&&+\nu \left( A_{1}\cap B_{1}\right) \nu \left( A_{2}\cap B_{2}\right) +\nu
\left( A_{1}\cap B_{2}\right) \nu \left( A_{2}\cap B_{1}\right)  \notag \\
&&+\nu \left( A_{2}\cap B_{1}\right) \nu \left( A_{1}\cap B_{2}\right) +\nu
\left( A_{2}\cap B_{2}\right) \nu \left( A_{1}\cap B_{1}\right) \text{.}
\notag
\end{eqnarray}

One easily verifies that%
\begin{eqnarray}
&&2\left[ \nu \left( A_{1}\cap B_{1}\right) {\rm St}\,{}_{\hat{0}}^{\hat{N}%
,\left[ 1\right] }\left( \mathbf{1}_{A_{2}\cap B_{2}}\right) +\nu \left(
A_{1}\cap B_{2}\right) {\rm St}\,{}_{\hat{0}}^{\hat{N},\left[ 1\right]
}\left( \mathbf{1}_{A_{2}\cap B_{1}}\right) \right.  \notag \\
&&\left. +\nu \left( A_{2}\cap B_{1}\right) {\rm St}\,{}_{\hat{0}}^{\hat{N}%
,\left[ 1\right] }\left( \mathbf{1}_{A_{1}\cap B_{2}}\right) +\nu \left(
A_{2}\cap B_{2}\right) {\rm St}\,{}_{\hat{0}}^{\hat{N},\left[ 1\right]
}\left( \mathbf{1}_{A_{1}\cap B_{1}}\right) \right]  \notag \\
&=&2I_{1}^{\widehat{N}}(f\star _{2}^{1}g)=2I_{1}^{\widehat{N}}(\widetilde{%
f\star _{2}^{1}g}),  \label{zz}
\end{eqnarray}%
and moreover%
\begin{eqnarray}
&&\mathbf{1}_{\left( A_{1}\cap B_{1}\right) \times \left( A_{2}\cap
B_{2}\right) }+\mathbf{1}_{\left( A_{1}\cap B_{2}\right) \times \left(
A_{2}\cap B_{1}\right) }+\mathbf{1}_{\left( A_{2}\cap B_{1}\right) \times
\left( A_{1}\cap B_{2}\right) }+\mathbf{1}_{\left( A_{2}\cap B_{2}\right)
\times \left( A_{1}\cap B_{1}\right) }  \notag \\
&=&\widetilde{f\star _{2}^{0}g}  \label{zzz}
\end{eqnarray}%
and%
\begin{eqnarray}
\left\langle f,g\right\rangle _{L^{2}\left( \nu ^{2}\right) } &=&\nu \left(
A_{1}\cap B_{1}\right) \nu \left( A_{2}\cap B_{2}\right) +\nu \left(
A_{1}\cap B_{2}\right) \nu \left( A_{2}\cap B_{1}\right)  \label{zzzz} \\
&&+\nu \left( A_{2}\cap B_{1}\right) \nu \left( A_{1}\cap B_{2}\right) +\nu
\left( A_{2}\cap B_{2}\right) \nu \left( A_{1}\cap B_{1}\right) \text{.}
\notag
\end{eqnarray}%
Since, trivially, ${\rm St}\,\!_{\sigma _{5}}^{\hat{N},\left[ 4\right]
}\left( f\otimes _{0}g\right) ={\rm St}\,\!_{\sigma _{6}}^{\hat{N},\left[ 4%
\right] }\left( f\otimes _{0}g\right) $, we deduce immediately from (\ref{z}%
)--(\ref{zzzz}) that the sum ${\rm St}\,\!_{\sigma _{5}}^{\hat{N},\left[ 4%
\right] }\left( f\otimes _{0}g\right) +{\rm St}\,\!_{\sigma _{6}}^{\hat{N},%
\left[ 4\right] }\left( f\otimes _{0}g\right) $ equals the expression
appearing in (\ref{aaa}). This proves the first part of the Proposition. The
last assertion in the statement can be proved by a density argument, similar
to the one used in order to conclude the proof of Proposition \ref{P :
2bleprodGauss}.
\end{proof}

\bigskip

Other proofs of Proposition \ref{P : 2bleProdPoiss} can be found for
instance in \cite{Kab, Surg1984, CT}.

\bigskip

\textbf{Examples. }(i) When $p=q=1$, one obtains
\begin{equation*}
I_{1}^{\hat{N}}\left( f\right) I_{1}^{\hat{N}}\left( g\right) =I_{2}^{\hat{N}%
}\left( \widetilde{f\otimes _{0}g}\right) +I_{1}^{\hat{N}}\left( \widetilde{%
f\star _{1}^{0}g}\right) +\left\langle f,g\right\rangle _{L^{2}\left( \nu
\right) }\text{.}
\end{equation*}

(ii) When $p=2,$ and $q=1$, one has%
\begin{eqnarray*}
I_{2}^{\hat{N}}\left( f\right) I_{1}^{\hat{N}}\left( g\right) &=&I_{3}^{\hat{%
N}}\left( \widetilde{f\otimes _{0}g}\right) +2I_{2}^{\hat{N}}\left(
\widetilde{f\star _{1}^{0}g}\right) +2I_{1}^{\hat{N}}\left( \widetilde{%
f\star _{1}^{1}g}\right) \\
&=&\int \int \int_{z_{1}\neq z_{2}\neq z_{3}}f\left( z_{1},z_{2}\right)
g\left( z_{3}\right) \hat{N}\left( dz_{1}\right) \hat{N}\left( dz_{2}\right)
\hat{N}\left( dz_{3}\right) \\
&&+2\int \int_{z_{1}\neq z_{2}}f\left( z_{1},z_{2}\right) g\left(
z_{1}\right) \hat{N}\left( dz_{1}\right) \hat{N}\left( dz_{2}\right) \\
&&+2\int \left( \int f\left( z_{1},x\right) g\left( z_{1}\right) \nu \left(
dx\right) \right) \hat{N}\left( dz_{1}\right) .
\end{eqnarray*}

\section{Diagram formulae \label{S : DF}}

\setcounter{equation}{0}We now want a general formula for computing
cumulants and expectations of products of multiple integrals, that is,
formulae for objects of the type
\begin{equation*}
\mathbb{E}\left[ I_{n_{1}}^{\varphi }\left( f_{1}\right) \times \cdot \cdot
\cdot \times I_{n_{k}}^{\varphi }\left( f_{k}\right) \right] \text{ \ \ and
\ \ }\chi \left( I_{n_{1}}^{\varphi }\left( f_{1}\right)
,...,I_{n_{k}}^{\varphi }\left( f_{k}\right) \right) \text{.}
\end{equation*}

\subsection{Formulae for moments and cumulants}

As in the previous sections, we shall focus on completely random measures $%
\varphi $ that are also \textsl{good }(in the sense of Definition \ref{D :
good}), so that moments and cumulants are well-defined. As usual, we shall
assume that the control measure $\nu $ is non-atomic. This last assumption
is not enough, however, because while the measure $\nu \left( A\right) =%
\mathbb{E}\varphi \left( A\right) ^{2}$ may be non-atomic, for some $n\geq 2$
the mean measure (concentrated on the \textquotedblleft full
diagonal\textquotedblright )
\begin{equation*}
\left\langle \Delta _{n}^{\varphi }\left( A\right) \right\rangle \triangleq
\mathbb{E}\left[ {\rm St}\,\!_{\hat{1}}^{\varphi ,\left[ n\right] }\left(
A\right) \right] =\mathbb{E}\left[ \varphi ^{\otimes n}\left\{ \left(
z_{1},...,z_{n}\right) \in A:z_{1}=\cdot \cdot \cdot =z_{n}\right\} \right]
\end{equation*}%
may be atomic. We shall therefore assume that $\varphi $ is
\textquotedblleft multiplicative\textquotedblright , that is, that this
phenomenon does not take place for any $n\geq 2$.

Proceeding formally, let $\varphi $ be a good random measure on $Z$, and fix
$n\geq 2$. Recall that $\mathcal{Z}_{\nu }^{n}$ denotes the collection of
all sets $B$ in $\mathcal{Z}^{\otimes n}$ such that $\nu ^{\otimes n}\left(
B\right) =\nu ^{n}\left( B\right) <\infty $ (see (\ref{ZetaNu})). As
before, for every partition $\pi \in \mathcal{P}\left( \left[ n\right]
\right) $, the class $\mathcal{Z}_{\pi }^{n}$ is the collection of all $\pi $%
-diagonal elements of $\mathcal{Z}^{n}$ (see (\ref{base : Bipi})). Recall
also that ${\rm St}\,\!_{\pi }^{\varphi ,\left[ n\right] }$ is the
restriction of the measure $\varphi ^{\otimes n}=\varphi ^{\left[ n\right] }$
on $Z_{\pi }^{n}$ (see (\ref{mes stoch})). Now let
\begin{eqnarray}
\left\langle \text{St}_{\pi }^{\varphi ,\left[ n\right] }\right\rangle
\left( C\right) &=&\mathbb{E}\left[ \text{St}_{\pi }^{\varphi ,\left[ n%
\right] }\left( C\right) \right] \text{, \ \ }C\in \mathcal{Z}_{\nu }^{n}%
\text{,}  \label{expD1} \\
\Delta _{1}^{\varphi }\left( A\right) &=&\varphi \left( A\right) \text{,}
\label{expD1/2} \\
\Delta _{n}^{\varphi }\left( A\right) &=&\text{St}_{\hat{1}}^{\varphi ,\left[
n\right] }\underset{n\text{ times}}{(\underbrace{A\times \cdot \cdot \cdot
\times A})}\text{, \ \ }A\in \mathcal{Z}_{\nu }\text{,}  \label{expD2} \\
\left\langle \Delta _{n}^{\varphi }\right\rangle \left( A\right) &=&\mathbb{E%
}\left[ \Delta _{n}^{\varphi }\left( A\right) \right] \text{, \ \ }A\in
\mathcal{Z}_{\nu }\text{.}  \label{expD3}
\end{eqnarray}

Thus, $\Delta _{n}^{\varphi }\left( A\right) $ denotes the random measure
concentrated on the full diagonal $z_{1}=...=z_{n}$ of the $n$tuple product $%
A\times \cdot \cdot \cdot \times A$, and $\left\langle \cdot \right\rangle $
denotes expectation.

\bigskip

\begin{definition}
\label{D : Multipl}We say that the good completely random measure $\varphi $
is \textbf{multiplicative}\textit{\ }if the deterministic measure $A\mapsto
\left\langle \Delta _{n}^{\varphi }\right\rangle \left( A\right) $ is
non-atomic for every $n\geq 2$. We show in the examples below that a
Gaussian or compensated Poisson measure, with non-atomic control measure $%
\nu $, is always multiplicative.
\end{definition}

\bigskip

The term \textquotedblleft multiplicative\textquotedblright\ (which we take
from \cite{RoWa})\ originates from the fact that $\varphi $ is
multiplicative (in the sense of the previous definition) if and only if
for every partition $\pi $ the non-random measure $\left\langle \text{St}%
_{\pi }^{\varphi ,\left[ n\right] }\right\rangle \left( \cdot \right) $ can
be written as a product measure. In particular (see Proposition 8 in \cite%
{RoWa}), the completely random measure $\varphi $ is multiplicative if and only if for every $\pi \in \mathcal{P}\left( \left[ n\right] \right) $ and
every $A_{1},...,A_{n}\in \mathcal{Z}_{\nu }$,
\begin{equation}
\left\langle \text{St}_{\pi }^{\varphi ,\left[ n\right] }\right\rangle
\left( A_{1}\times \cdot \cdot \cdot \times A_{n}\right) =\prod_{b\in \pi
}\left\langle \text{St}_{\hat{1}}^{\varphi ,\left[ \left\vert b\right\vert %
\right] }\right\rangle \left( \underset{j\in b}{\mathrm{X}}A_{j}\right)
\text{,}  \label{MultExp}
\end{equation}%
where, for every $b=\left\{ j_{1},...,j_{k}\right\} \in \pi $, we used once
again the notation $\underset{j\in b}{\mathrm{X}}A_{j}\triangleq
A_{j_{1}}\times \cdot \cdot \cdot \times A_{j_{k}}$. Note that the RHS\ of (%
\ref{MultExp}) involves products over blocks of the partition $\pi $, in
which there is concentration over the diagonals associated with the blocks.
Thus, in view of (\ref{Triv1}) and (\ref{expD2}), one has that%
\begin{equation}
\left\langle \text{St}_{\pi }^{\varphi ,\left[ n\right] }\right\rangle
\left( A_{1}\times \cdot \cdot \cdot \times A_{n}\right) =\prod_{b\in \pi
}\left\langle \Delta _{\left\vert b\right\vert }^{\varphi }\right\rangle
\left( \cap _{j\in b}A_{j}\right) \text{,}  \label{Multexp1}
\end{equation}%
that is, we can express the LHS of (\ref{MultExp}) as a product of measures
involving sets in $\mathcal{Z}_{\nu }$. Observe that one can rewrite
relation (\ref{MultExp}) in the following (compact) way:%
\begin{equation}
\left\langle \text{St}_{\pi }^{\varphi ,\left[ n\right] }\right\rangle
=\bigotimes\limits_{b\in \pi }\left\langle \text{St}_{\hat{1}}^{\varphi ,%
\left[ \left\vert b\right\vert \right] }\right\rangle .  \label{multexp2}
\end{equation}

\bigskip

\textbf{Examples. }(i) When $\varphi $ is Gaussian with non-atomic control
measure $\nu $, relation (\ref{MultExp}) implies that $\left\langle \text{St}%
_{\pi }^{\varphi ,\left[ n\right] }\right\rangle $ is $0$ if $\pi $ contains
at least one block $b$ such that $\left\vert b\right\vert \neq 2$. If, on
the other hand, every block of $\pi $ contains exactly two elements, we
deduce from (\ref{GaussDia}) and (\ref{m-DG2}) that
\begin{equation}
\left\langle \text{St}_{\pi }^{\varphi ,\left[ n\right] }\right\rangle
\left( A_{1}\times \cdot \cdot \cdot \times A_{n}\right) =\prod_{b=\left\{
i,j\right\} \in \pi }\nu \left( A_{i}\cap A_{j}\right) ,
\label{gaussCombExp}
\end{equation}%
which is not atomic.

(ii) If $\varphi $ is a compensated Poisson measure with non-atomic control
measure $\nu $, then $\left\langle \text{St}_{\pi }^{\varphi ,\left[ n\right]
}\right\rangle $ is $0$ whenever $\pi $ contains at least one block $b$ such
that $\left\vert b\right\vert =1$ (indeed, that block would have measure 0,
since $\varphi $ is centered). If, on the other hand, every block of $\pi $
has more than two elements, then, by Corollary \ref{C : meanN}%
\begin{equation}
\left\langle \text{St}_{\pi }^{\varphi ,\left[ n\right] }\right\rangle
\left( A_{1}\times \cdot \cdot \cdot \times A_{n}\right)
=\prod_{k=2}^{n}\prod_{b=\left\{ j_{1},...,j_{k}\right\} \in \pi }\nu \left(
A_{j_{1}}\cap \cdot \cdot \cdot \cap A_{j_{k}}\right) ,  \label{PoissCombExp}
\end{equation}%
which is non-atomic. See \cite{RoWa} for (quite pathological) examples of
non-multiplicative measures.

\bigskip

\noindent\textbf{Notation. }In what follows, the notation%
\begin{equation}
\int_{Z^{n}}f\left( z_{1},...,z_{n}\right) \bigotimes\limits_{b\in \pi
}\left\langle \text{St}_{\hat{1}}^{\varphi ,\left[ \left\vert b\right\vert %
\right] }\right\rangle \left( dz_{1},...,dz_{n}\right) \triangleq
\bigotimes\limits_{b\in \pi }\left\langle \text{St}_{\hat{1}}^{\varphi ,%
\left[ \left\vert b\right\vert \right] }\right\rangle \left( f\right)
\label{ii}
\end{equation}%
will be used for every function $f\in \mathcal{E}\left( \nu ^{n}\right) .$
\footnote{%
The integral $\int_{Z^{n}}f\{d\bigotimes\limits_{b\in \pi }\left\langle
\text{St}_{\hat{1}}^{\varphi ,\left[ \left\vert b\right\vert \right]
}\right\rangle \}$ in (\ref{ii}) is well defined, since the set function $%
\bigotimes\limits_{b\in \pi }\left\langle \text{St}_{\hat{1}}^{\varphi ,%
\left[ \left\vert b\right\vert \right] }\right\rangle \left( \cdot \right) $
is a $\sigma $-additive signed measure (thanks to (\ref{sigmaEng})) on the
algebra generated by the products of the type $A_{1}\times \cdot \cdot \cdot
\times A_{n}$, where each $A_j$ is in $\mathcal{Z}_\nu$.  }

\bigskip

The next result gives a new universal combinatorial formula for the computation
of the cumulants and the moments associated with the multiple Wiener-It\^{o}
integrals with respect to a completely random \textsl{multiplicative}
measure.

\bigskip

\begin{theorem}[Diagram formulae]
\label{T : Diagrams}Let $\varphi $ be a good completely random measure, with
non-atomic control measure $\nu $, and suppose that $\varphi $ is also
multiplicative in the sense of Definition \ref{D : Multipl}. For every $%
n_{1},...,n_{k}\geq 1$, we write $n=n_{1}+\cdot \cdot \cdot +n_{k}$, and we
denote by $\pi ^{\ast }$ the partition of $\left[ n\right] =\left\{
1,...,n\right\} $ given by (\ref{pistar}). Then, for every collection of
kernels $f_{1},...,f_{k}$ such that $f_{j}\in \mathcal{E}_{s,0}\left( \nu
^{n_{j}}\right) $, one has that%
\begin{equation}
\mathbb{E}\left[ I_{n_{1}}^{\varphi }\left( f_{1}\right) \cdot \cdot \cdot
I_{n_{k}}^{\varphi }\left( f_{k}\right) \right] =\sum_{\left\{ \sigma \in
\mathcal{P}\left( \left[ n\right] \right) :\sigma \wedge \pi ^{\ast }=\hat{0}%
\right\} }\bigotimes\limits_{b\in \sigma }\left\langle {\rm St}\,\text{{}}%
_{\hat{1}}^{\varphi ,\left[ \left\vert b\right\vert \right] }\right\rangle
\left( f_{1}\otimes _{0}f_{2}\otimes _{0}\cdot \cdot \cdot \otimes
_{0}f_{k}\right) ,  \label{momMWI}
\end{equation}%
and%
\begin{equation}
\mathbb{\chi }\left( I_{n_{1}}^{\varphi }\left( f_{1}\right) ,\cdot \cdot
\cdot ,I_{n_{k}}^{\varphi }\left( f_{k}\right) \right) =\sum_{\substack{ %
\sigma \wedge \pi ^{\ast }=\hat{0}  \\ \sigma \vee \pi ^{\ast }=\hat{1}}}%
\bigotimes\limits_{b\in \sigma }\left\langle {\rm St}\,\text{{}}_{\hat{1}%
}^{\varphi ,\left[ \left\vert b\right\vert \right] }\right\rangle \left(
f_{1}\otimes _{0}f_{2}\otimes _{0}\cdot \cdot \cdot \otimes _{0}f_{k}\right)
,  \label{cumMWI}
\end{equation}
\end{theorem}

\begin{proof}
Formula (\ref{momMWI}) is a consequence of Theorem \ref{T : ProdRW} and (\ref%
{MultExp}). In order to prove (\ref{cumMWI}), we shall first show that the
following two equalities hold
\begin{eqnarray}
&&\mathbb{\chi }\left( I_{n_{1}}^{\varphi }\left( f_{1}\right) ,\cdot \cdot
\cdot ,I_{n_{k}}^{\varphi }\left( f_{k}\right) \right)  \notag \\
&=&\sum_{\pi ^{\ast }\leq \rho =\left( r_{1},...,r_{l}\right) \in \mathcal{P}%
\left( \left[ n\right] \right) }\mu \left( \rho ,\hat{1}\right)
\prod_{j=1}^{l}\mathbb{E}\left[ \prod_{a:\left\{ n_{1}+\cdot \cdot \cdot
+n_{a-1}+1,...,n_{1}+\cdot \cdot \cdot +n_{a}\right\} \subseteq
r_{j}}I_{n_{a}}^{\varphi }\left( f_{a}\right) \right]  \label{qo} \\
&=&\sum_{\pi ^{\ast }\leq \rho =\left( r_{1},...,r_{l}\right) \in \mathcal{P}%
\left( \left[ n\right] \right) }\mu \left( \rho ,\hat{1}\right) \sum
_{\substack{ \gamma \leq \rho  \\ \gamma \wedge \pi ^{\ast }=\hat{0}}}%
\bigotimes\limits_{b\in \gamma }\left\langle \text{St}_{\hat{1}}^{\varphi ,%
\left[ \left\vert b\right\vert \right] }\right\rangle \left( f_{1}\otimes
_{0}\cdot \cdot \cdot \otimes _{0}f_{k}\right) ,  \label{qoo}
\end{eqnarray}%
where $n_{1}+n_{0}=0$ by convention.

The proof of (\ref{qo}) uses arguments analogous to those in the proof of
Malyshev's formula (\ref{Maly}). Indeed, one can use relation (\ref{LS2}) to
deduce that%
\begin{equation}
\mathbb{\chi }\left( I_{n_{1}}^{\varphi }\left( f_{1}\right) ,\cdot \cdot
\cdot ,I_{n_{k}}^{\varphi }\left( f_{k}\right) \right) =\sum_{\sigma
=\left\{ x_{1},...,x_{l}\right\} \in \mathcal{P}\left( \left[ k\right]
\right) }\left( -1\right) ^{l-1}\left( l-1\right) !\prod_{j=1}^{l}\mathbb{E}%
\left[ \prod_{a\in x_{j}}I_{n_{a}}^{\varphi }\left( f_{a}\right) \right] .
\label{qqaqq}
\end{equation}%
Now observe that there exists a \textsl{bijection}%
\begin{equation*}
\mathcal{P}\left( \left[ k%
\right] \right) \rightarrow \left[ \pi ^{\ast },\hat{1}\right]:\sigma \mapsto \rho ^{\left( \sigma \right) } \text{,}
\end{equation*}%
between $\mathcal{P}\left( \left[ k\right] \right) $ and the segment $\left[
\pi ^{\ast },\hat{1}\right] $, which is defined as the set of those $\rho
\in \mathcal{P}\left( \left[ n\right] \right) $ such that $\pi ^{\ast }\leq
\rho $, where $\pi ^{\ast }$ is given by (\ref{pistar}). Such a bijection is
realized as follows: for every $\sigma =\left\{ x_{1},...,x_{l}\right\} \in
\mathcal{P}\left( \left[ k\right] \right) $, define $\rho ^{\left( \sigma
\right) }\in \left[ \pi ^{\ast },\hat{1}\right] $ by merging two blocks
\begin{equation*}
\left\{ n_{1}+\cdot \cdot \cdot +n_{a-1}+1,...,n_{1}+\cdot \cdot \cdot
+n_{a}\right\} \text{ \ and \ }\left\{ n_{1}+\cdot \cdot \cdot
+n_{b-1}+1,...,n_{1}+\cdot \cdot \cdot +n_{b}\right\}
\end{equation*}%
of $\pi ^{\ast }$ ($1\leq a\neq b\leq k$) if and only if $a\sim _{\sigma
}b $. Note that this construction implies that $\left\vert \sigma
\right\vert =\left\vert \rho ^{\left( \sigma \right) }\right\vert =l$, so
that (\ref{mp1}) yields
\begin{equation}
\left( -1\right) ^{l-1}\left( l-1\right) !=\mu \left( \sigma ,\hat{1}\right)
=\mu \left( \rho ^{\left( \sigma \right) },\hat{1}\right)  \label{momo}
\end{equation}%
(observe that the two M\"{o}bius functions appearing in (\ref{momo}) refer
to two different lattices of partitions). Now use the notation $\rho
^{\left( \sigma \right) }=\left\{ r_{1}^{\left( \sigma \right)
},...,r_{l}^{\left( \sigma \right) }\right\} $ to indicate the blocks of $%
\rho ^{\left( \sigma \right) }$: since, by construction,%
\begin{equation}
\prod_{j=1}^{l}\mathbb{E}\left[ \prod_{a\in x_{j}}I_{n_{a}}^{\varphi }\left(
f_{a}\right) \right] =\prod_{j=1}^{l}\mathbb{E}\left[ \prod_{a:\left\{
n_{1}+\cdot \cdot \cdot +n_{a-1}+1,...,n_{1}+\cdot \cdot \cdot
+n_{a}\right\} \subseteq r_{j}^{\left( \sigma \right) }}I_{n_{a}}^{\varphi
}\left( f_{a}\right) \right] ,  \label{qaq}
\end{equation}%
we immediately obtain (\ref{qo}) by plugging (\ref{momo}) and (\ref{qaq})
into (\ref{qqaqq}).

To prove (\ref{qoo}), fix $\rho =\left\{ r_{1},...,r_{l}\right\} $ such that
$\pi ^{\ast }\leq \rho .$ For $j=1,...,l$, we write $\pi ^{\ast }\left(
j\right) $ to indicate the partition of the block $r_{j}$ whose blocks are
the sets $\{n_{1}+\cdot \cdot \cdot +n_{a-1}+1,...,n_{1}+\cdot \cdot \cdot
+n_{a}\}$ such that
\begin{equation}
\left\{ n_{1}+\cdot \cdot \cdot +n_{a-1}+1,...,n_{1}+\cdot \cdot \cdot
+n_{a}\right\} \subseteq r_{j}\text{.}  \label{fgf}
\end{equation}%
According to (\ref{momMWI}),
\begin{equation*}
\mathbb{E}\left[ \prod_{a:\left\{ n_{1}+\cdot \cdot \cdot
+n_{a-1}+1,...,n_{1}+\cdot \cdot \cdot +n_{a}\right\} \subseteq
r_{j}}I_{n_{a}}^{\varphi }\left( f_{a}\right) \right] =\sum_{\left\{ \sigma
\in \mathcal{P}\left( r_{j}\right) :\sigma \wedge \pi ^{\ast }\left(
j\right) =\hat{0}\right\} }\bigotimes\limits_{b\in \sigma }\left\langle
{\rm St}\,\text{{}}_{\hat{1}}^{\varphi ,\left[ \left\vert b\right\vert %
\right] }\right\rangle \left( \left\{ \otimes _{r_{j},0}f\right\} \right)
\text{,}
\end{equation*}%
where the function $\left\{ \otimes _{r_{j},0}f\right\} $ is obtained by
juxtaposing the $\left\vert \pi ^{\ast }\left( j\right) \right\vert $
functions $f_{a}$ such that the index $a$ verifies (\ref{fgf}). Now observe
that $\gamma \in \mathcal{P}\left( \left[ n\right] \right) $ satisfies
\begin{equation*}
\gamma \leq \rho \text{ \ \ and \ \ }\gamma \wedge \pi ^{\ast }=\hat{0},
\end{equation*}%
if and only if $\gamma $ admits a (unique) representation as a union of
the type%
\begin{equation*}
\gamma =\bigcup\limits_{j=1}^{l}\sigma \left( j\right) \text{,}
\end{equation*}%
where each $\sigma \left( j\right) $ is an element of $\mathcal{P}\left(
r_{j}\right) $ such that $\sigma \left( j\right) \wedge \pi ^{\ast }\left(
j\right) =\hat{0}.$ This yields%
\begin{eqnarray*}
&&\prod_{j=1}^{l}\sum_{\left\{ \sigma \in \mathcal{P}\left( r_{j}\right)
:\sigma \wedge \pi ^{\ast }\left( j\right) =\hat{0}\right\}
}\bigotimes\limits_{b\in \sigma }\left\langle {\rm St}\,\text{{}}_{\hat{1}%
}^{\varphi ,\left[ \left\vert b\right\vert \right] }\right\rangle \left(
\left\{ \otimes _{r_{j},0}f\right\} \right) \\
&=&\sum_{\substack{ \gamma \leq \rho  \\ \gamma \wedge \pi ^{\ast }=\hat{0}}}%
\bigotimes\limits_{b\in \gamma }\left\langle \text{St}_{\hat{1}}^{\varphi ,%
\left[ \left\vert b\right\vert \right] }\right\rangle \left( f_{1}\otimes
_{0}\cdot \cdot \cdot \otimes _{0}f_{k}\right) \text{.}
\end{eqnarray*}%
This relation, together with (\ref{momo}) and (\ref{qaq}), shows that (\ref%
{qo}) implies (\ref{qoo}). To conclude the proof, just observe that, by
inverting the order of summation in (\ref{qoo}), one obtains that%
\begin{eqnarray*}
\mathbb{\chi }\left( I_{n_{1}}^{\varphi }\left( f_{1}\right) ,\cdot \cdot
\cdot ,I_{n_{k}}^{\varphi }\left( f_{k}\right) \right) &=&\sum_{\gamma
\wedge \pi ^{\ast }=\hat{0}}\bigotimes\limits_{b\in \gamma }\left\langle
\text{St}_{\hat{1}}^{\varphi ,\left[ \left\vert b\right\vert \right]
}\right\rangle \left( f_{1}\otimes _{0}\cdot \cdot \cdot \otimes
_{0}f_{k}\right) \sum_{\pi ^{\ast }\vee \gamma \leq \rho \leq \hat{1}}\mu
\left( \rho ,\hat{1}\right) \\
&=&\sum_{\substack{ \gamma \wedge \pi ^{\ast }=\hat{0}  \\ \pi ^{\ast }\vee
\gamma =\hat{1}}}\bigotimes\limits_{b\in \gamma }\left\langle \text{St}_{%
\hat{1}}^{\varphi ,\left[ \left\vert b\right\vert \right] }\right\rangle
\left( f_{1}\otimes _{0}\cdot \cdot \cdot \otimes _{0}f_{k}\right) ,
\end{eqnarray*}%
where the last equality is a consequence of the relation%
\begin{equation*}
\sum_{\pi ^{\ast }\vee \gamma \leq \rho \leq \hat{1}}\mu \left( \rho ,\hat{1}%
\right) =\left\{
\begin{array}{ll}
1 & \text{if }\pi ^{\ast }\vee \gamma =\hat{1} \\
0 & \text{otherwise,}%
\end{array}%
\right.
\end{equation*}%
which is in turn a special case of (\ref{maliMob}).
\end{proof}

\bigskip

\textbf{Remark. }Observe that the only difference between the moment formula
(\ref{momMWI}) and the cumulant formula (\ref{cumMWI}) is that in the first
case the sum is over all $\sigma $ such that $\sigma \wedge \pi ^{\ast }=%
\hat{0}=\left\{ \left\{ 1\right\} ,...,\left\{ n\right\} \right\} $, and
that in the second case $\sigma $ must satisfy in addition that $\sigma \vee
\pi ^{\ast }=\hat{1}=\left\{ \left[ n\right] \right\} $. Moreover, the
relations (\ref{momMWI}) and (\ref{cumMWI}) can be restated in terms of
diagrams by rewriting the sums as
\begin{equation*}
\sum_{\left\{ \sigma \in \mathcal{P}\left( \left[ n\right] \right) :\sigma
\wedge \pi ^{\ast }=\hat{0}\right\} }=\sum_{\sigma \in \mathcal{P}\left( %
\left[ n\right] \right) :\Gamma \left( \pi ^{\ast },\sigma \right) \text{ is
non-flat}}\text{ \ \ ; \ \ }\sum_{\substack{ \sigma \wedge \pi ^{\ast }=\hat{%
0}  \\ \sigma \vee \pi ^{\ast }=\hat{1}}}=\sum_{\substack{ \sigma \in
\mathcal{P}\left( \left[ n\right] \right) :\Gamma \left( \pi ^{\ast },\sigma
\right) \text{ is non-flat}  \\ \text{and connected}}}\text{,}
\end{equation*}%
where $\Gamma \left( \pi ^{\ast },\sigma \right) $ is the diagram of $\left(
\pi ^{\ast },\sigma \right) $, as defined in Section \ref{SS : diagrams}.

\subsection{The Gaussian case \label{SS : GaussCase}}

We shall now provide a version of Theorem \ref{T : Diagrams} in the case
where $\varphi $ is, respectively, Gaussian and Poisson. For convenience,
using the same notation as that in Theorem \ref{T : Diagrams}, let
\begin{eqnarray}
\mathcal{M}\left( \left[ n\right] ,\pi ^{\ast }\right) &\triangleq &\left\{
\sigma \in \mathcal{P}\left( \left[ n\right] \right) :\sigma \vee \pi ^{\ast
}=\hat{1}\text{ and }\sigma \wedge \pi ^{\ast }=\hat{0}\right\}  \label{form}
\\
\mathcal{M}^{0}\left( \left[ n\right] ,\pi ^{\ast }\right) &\triangleq
&\left\{ \sigma \in \mathcal{P}\left( \left[ n\right] \right) :\sigma \wedge
\pi ^{\ast }=\hat{0}\right\}  \label{form2}
\end{eqnarray}%
and%
\begin{eqnarray}
\mathcal{M}_{2}\left( \left[ n\right] ,\pi ^{\ast }\right) &\triangleq
&\left\{ \sigma \in \mathcal{M}\left( \left[ n\right] ,\pi ^{\ast }\right)
:\left\vert b\right\vert =2\text{, \ }\forall b\in \sigma \right\}
\label{M2} \\
\mathcal{M}_{2}^{0}\left( \left[ n\right] ,\pi ^{\ast }\right) &\triangleq
&\left\{ \sigma \in \mathcal{M}^{0}\left( \left[ n\right] ,\pi ^{\ast
}\right) :\left\vert b\right\vert =2\text{, \ }\forall b\in \sigma \right\}
\label{M2b} \\
\mathcal{M}_{\geq 2}\left( \left[ n\right] ,\pi ^{\ast }\right) &\triangleq
&\left\{ \sigma \in \mathcal{M}\left( \left[ n\right] ,\pi ^{\ast }\right)
:\left\vert b\right\vert \geq 2\text{, \ }\forall b\in \sigma \right\}
\label{Mplus2} \\
\mathcal{M}_{\geq 2}^{0}\left( \left[ n\right] ,\pi ^{\ast }\right)
&\triangleq &\left\{ \sigma \in \mathcal{M}^{0}\left( \left[ n\right] ,\pi
^{\ast }\right) :\left\vert b\right\vert \geq 2\text{, \ }\forall b\in
\sigma \right\}  \label{Mplus2b}
\end{eqnarray}%
where the partition $\pi ^{\ast }\in \mathcal{P}\left( \left[ n\right]
\right) $ is defined in (\ref{pistar}). The sets $\mathcal{M}_{2}\left( %
\left[ n\right] ,\pi ^{\ast }\right) $ and $\mathcal{M}_{2}^{0}\left( \left[
n\right] ,\pi ^{\ast }\right) $ appear in the case where $\varphi $ is
Gaussian. Note that, by using the formalism of diagrams $\Gamma $ and
multigraphs $\hat{\Gamma}$ introduced in Section \ref{S : DG}, one has that%
\begin{eqnarray}
\mathcal{M}\left( \left[ n\right] ,\pi ^{\ast }\right) &=&\left\{ \sigma \in
\mathcal{P}\left( \left[ n\right] \right) :\Gamma \left( \pi ^{\ast },\sigma
\right) \text{ is non-flat and connected}\right\}  \label{Lux1} \\
\mathcal{M}^{0}\left( \left[ n\right] ,\pi ^{\ast }\right) &=&\left\{ \sigma
\in \mathcal{P}\left( \left[ n\right] \right) :\Gamma \left( \pi ^{\ast
},\sigma \right) \text{ is connected}\right\}  \label{Lux2} \\
\mathcal{M}_{2}\left( \left[ n\right] ,\pi ^{\ast }\right) &=&\left\{ \sigma
\in \mathcal{P}\left( \left[ n\right] \right) :\Gamma \left( \pi ^{\ast
},\sigma \right) \text{ is Gaussian, non-flat and connected}\right\}
\label{Lux3} \\
&=&\left\{ \sigma \in \mathcal{P}\left( \left[ n\right] \right) :\hat{\Gamma}%
\left( \pi ^{\ast },\sigma \right) \text{ has no loops and is connected}%
\right\}  \notag \\
\mathcal{M}_{2}^{0}\left( \left[ n\right] ,\pi ^{\ast }\right) &=&\left\{
\sigma \in \mathcal{P}\left( \left[ n\right] \right) :\Gamma \left( \pi
^{\ast },\sigma \right) \text{ is Gaussian and non-flat}\right\}
\label{Lux5} \\
&=&\left\{ \sigma \in \mathcal{P}\left( \left[ n\right] \right) :\hat{\Gamma}%
\left( \pi ^{\ast },\sigma \right) \text{ has no loops}\right\} .  \notag
\end{eqnarray}%
Clearly, $\mathcal{M}_{2}\left( \left[ n\right] ,\pi ^{\ast }\right) \subset
\mathcal{M}_{2}^{0}\left( \left[ n\right] ,\pi ^{\ast }\right) $, $\mathcal{M%
}_{2}\left( \left[ n\right] ,\pi ^{\ast }\right) \subset \mathcal{M}_{\geq
2}\left( \left[ n\right] ,\pi ^{\ast }\right) $ and
\begin{equation*}
\mathcal{M}_{2}^{0}\left( \left[ n\right] ,\pi ^{\ast }\right) \subset
\mathcal{M}_{\geq 2}^{0}\left( \left[ n\right] ,\pi ^{\ast }\right) .
\end{equation*}%
The sets $\mathcal{M}_{\geq 2}\left( \left[ n\right] ,\pi ^{\ast }\right) $
and $\mathcal{M}_{\geq 2}^{0}\left( \left[ n\right] ,\pi ^{\ast }\right) $
appear when $\varphi $ is a compensated Poisson measure, namely $\varphi =%
\hat{N}$.

\bigskip

\begin{corollary}[Gaussian measures]
\label{C : DiaGauss}Suppose $\varphi =G$ is a centered Gaussian measure with
non-atomic control measure $\nu $, fix integers $n_{1},...,n_{k}\geq 1$ and
let $n=n_{1}+\cdot \cdot \cdot +n_{k}$. Write $\pi ^{\ast }$ for the
partition of $\left[ n\right] $ appearing in (\ref{pistar}). Then, for any
vector of functions $(f_{1},...,f_{k})$ such that $f_{j}\in L_{s}^{2}\left(
\nu ^{n_{j}}\right) $, $j=1,...,k$, the following relations hold:

\begin{enumerate}
\item If $\mathcal{M}_{2}\left( \left[ n\right] ,\pi ^{\ast }\right)
=\varnothing $ (in particular, if $n$ is odd), then $\mathbb{\chi }\left(
I_{n_{1}}^{G}\left( f_{1}\right) ,\cdot \cdot \cdot ,I_{n_{k}}^{G}\left(
f_{k}\right) \right) =0$;

\item If $\mathcal{M}_{2}\left( \left[ n\right] ,\pi ^{\ast }\right) \neq
\varnothing $, then
\begin{equation}
\mathbb{\chi }\left( I_{n_{1}}^{G}\left( f_{1}\right) ,\cdot \cdot \cdot
,I_{n_{k}}^{G}\left( f_{k}\right) \right) =\sum_{\sigma \in \mathcal{M}%
_{2}\left( \left[ n\right] ,\pi ^{\ast }\right) }\int_{Z^{n/2}}f_{\sigma
,k}d\nu ^{n/2},  \label{GaussDiagrammi}
\end{equation}%
where, for every $\sigma \in \mathcal{M}_{2}\left( \left[ n\right] ,\pi
^{\ast }\right) $, the function $f_{\sigma ,k}$, of $n/2$ variables, is
obtained by identifying the variables $x_{i}$ and $x_{j}$ in the argument of
$f_{1}\otimes _{0}\cdot \cdot \cdot \otimes _{0}f_{n_{k}}$ (as given in (\ref%
{0multContr})) if and only if $i\sim _{\sigma }j$;

\item If $\mathcal{M}_{2}^{0}\left( \left[ n\right] ,\pi ^{\ast }\right)
=\varnothing $ (in particular, if $n$ is odd), then $\mathbb{E}\left(
I_{n_{1}}^{G}\left( f_{1}\right) \cdot \cdot \cdot I_{n_{k}}^{G}\left(
f_{k}\right) \right) =0$;

\item If $\mathcal{M}_{2}^{0}\left( \left[ n\right] ,\pi ^{\ast }\right)
\neq \varnothing $,
\begin{equation}
\mathbb{E}\left( I_{n_{1}}^{G}\left( f_{1}\right) \cdot \cdot \cdot
I_{n_{k}}^{G}\left( f_{k}\right) \right) =\sum_{\sigma \in \mathcal{M}%
_{2}^{0}\left( \left[ n\right] ,\pi ^{\ast }\right) }\int_{Z^{n/2}}f_{\sigma
,k}d\nu ^{n/2}  \label{GaussDiagr}
\end{equation}
\end{enumerate}
\end{corollary}

\begin{proof}
First observe that, since $\varphi =G$ is Gaussian, then $\left\langle
{\rm St}\,\text{{}}_{\hat{1}}^{G,\left[ \left\vert b\right\vert \right]
}\right\rangle \equiv 0$ whenever $\left\vert b\right\vert \neq 2$. Assume
for the moment that $f_{j}\in \mathcal{E}_{s,0}\left( \nu ^{n_{j}}\right) $,
$j=1,...,k$. In this case, we can apply formula (\ref{cumMWI}) and obtain
that
\begin{eqnarray*}
&&\mathbb{\chi }\left( I_{n_{1}}^{G}\left( f_{1}\right) ,\cdot \cdot \cdot
,I_{n_{k}}^{G}\left( f_{k}\right) \right) \\
&=&\sum_{\substack{ \{\sigma :\sigma \wedge \pi ^{\ast }=\hat{0}\text{ ; }
\\ \sigma \vee \pi ^{\ast }=\hat{1}\text{ ; }\left\vert b\right\vert =2\text{%
\ }\forall b\in \sigma \}}}\bigotimes\limits_{b\in \sigma }\left\langle
{\rm St}\,\text{{}}_{\hat{1}}^{G,\left[ \left\vert b\right\vert \right]
}\right\rangle \left( f_{1}\otimes _{0}f_{2}\otimes _{0}\cdot \cdot \cdot
\otimes _{0}f_{k}\right) \\
&=&\sum_{\sigma \in \mathcal{M}_{2}\left( \left[ n\right] ,\pi ^{\ast
}\right) }\bigotimes\limits_{b\in \sigma }\left\langle {\rm St}\,\text{{}}%
_{\hat{1}}^{G,\left[ \left\vert b\right\vert \right] }\right\rangle \left(
f_{1}\otimes _{0}f_{2}\otimes _{0}\cdot \cdot \cdot \otimes _{0}f_{k}\right)
\text{,}
\end{eqnarray*}%
where we have used (\ref{M2}). The last relation trivially implies Point 1
in the statement. Moreover, since, for every $B,C\in \mathcal{Z}_{\nu }$,%
\begin{equation}
\left\langle {\rm St}\,\text{{}}_{\hat{1}}^{G,\left[ 2\right]
}\right\rangle \left( B\times C\right) =\left\langle {\rm St}\,\text{{}}_{%
\hat{1}}^{G,\left[ 2\right] }\right\rangle \left( \left( B\cap C\right)
\times \left( B\cap C\right) \right) =\left\langle \Delta
_{2}^{G}\right\rangle \left( B\cap C\right) \text{,}  \label{xx}
\end{equation}%
one deduces immediately that the support of the deterministic measure $%
\bigotimes\limits_{b\in \sigma }\left\langle {\rm St}\,\text{{}}_{\hat{1}%
}^{G,\left[ \left\vert b\right\vert \right] }\right\rangle $ is contained in
the set
\begin{equation*}
Z_{\geq \sigma }^{n}=\left\{ \left( z_{1},...,z_{n}\right) :z_{i}=z_{j}\text{
for every }i,j\text{ such that }i\sim _{\sigma }j\right\} .
\end{equation*}%
Since, by (\ref{GaussDia}) and (\ref{xx}),
\begin{equation}
\left\langle {\rm St}\,\text{{}}_{\hat{1}}^{G,\left[ \left\vert
b\right\vert \right] }\right\rangle \left( B\times C\right) =\nu \left(
B\cap C\right) \text{,}  \label{wq}
\end{equation}%
for every $B,C\in \mathcal{Z}_{\nu }$, we infer that%
\begin{equation}
\bigotimes\limits_{b\in \sigma }\left\langle {\rm St}\,\text{{}}_{\hat{1}%
}^{G,\left[ \left\vert b\right\vert \right] }\right\rangle \left(
f_{1}\otimes _{0}f_{2}\otimes _{0}\cdot \cdot \cdot \otimes _{0}f_{k}\right)
=\bigotimes\limits_{b\in \sigma }\left\langle {\rm St}\,\text{{}}_{\hat{1}%
}^{G,\left[ \left\vert b\right\vert \right] }\right\rangle \left( f_{\sigma
}\right) =\int_{Z^{n/2}}f_{\sigma ,k}d\nu ^{n/2}\text{.}  \label{fgg}
\end{equation}%
where the function $f_{\sigma ,k}$ is defined in the statement. To obtain
the last equality in (\ref{fgg}), one should start with functions $f_{j}$ of
type $f_{j}\left( z_{1},...,z_{n_{j}}\right) =\mathbf{1}_{C_{1}^{\left(
j\right) }\times \cdot \cdot \cdot \times C_{n_{j}}^{\left( j\right)
}}\left( z_{1},...,z_{n_{j}}\right) $, where the $C_{\ell }^{\left( j\right)
}\in \mathcal{Z}_{\nu }$ are disjoint, and then apply formula (\ref%
{gaussCombExp}), so that the extension to general functions $f_{j}\in
\mathcal{E}_{s,0}\left( \nu ^{n_{j}}\right) $ is obtained by the
multilinearity of the application%
\begin{equation*}
\left( f_{1},...,f_{k}\right) \mapsto \int_{Z^{n/2}}f_{\sigma ,k}d\nu ^{n/2}.
\end{equation*}%
To obtain (\ref{GaussDiagrammi}) for general functions $f_{1},...,f_{k}$
such that $f_{j}\in L_{s}^{2}\left( \nu ^{n_{j}}\right) $, start by
observing that $\mathcal{E}_{s,0}\left( \nu ^{n_{j}}\right) $ is dense in $%
L_{s}^{2}\left( \nu ^{n_{j}}\right) $, and then use the fact that, if a
sequence $f_{1}^{\left( r\right) },...,f_{k}^{\left( r\right) }$, $r\geq 1$,
is such that $f_{j}^{\left( r\right) }\in \mathcal{E}_{s,0}\left( \nu
^{n_{j}}\right) $ and $f_{j}^{\left( r\right) }\rightarrow f_{j}$ in $%
L_{s}^{2}\left( \nu ^{n_{j}}\right) $ ($j=1,...,k$), then
\begin{equation*}
\mathbb{\chi }\left( I_{n_{1}}^{G}\left( f_{1}^{\left( r\right) }\right)
,\cdot \cdot \cdot ,I_{n_{k}}^{G}\left( f_{k}^{\left( r\right) }\right)
\right) \rightarrow \mathbb{\chi }\left( I_{n_{1}}^{G}\left( f_{1}\right)
,\cdot \cdot \cdot ,I_{n_{k}}^{G}\left( f_{k}\right) \right) \text{,}
\end{equation*}%
by (\ref{GaussChaosCOntr}), and moreover%
\begin{equation*}
\int_{Z^{n/2}}f_{\sigma ,k}^{\left( r\right) }d\nu ^{n/2}\rightarrow
\int_{Z^{n/2}}f_{\sigma ,k}d\nu ^{n/2}\text{,}
\end{equation*}%
where $f_{\sigma ,k}^{\left( r\right) }$ is constructed from $f_{1}^{\left(
r\right) },...,f_{k}^{\left( r\right) }$, as specified in the statement (a
similar argument was needed in the proof of Proposition \ref{P :
2bleprodGauss}). Points 3 and 4 in the statement are obtained analogously,
by using the relations%
\begin{eqnarray*}
&&\mathbb{E}\left( I_{n_{1}}^{G}\left( f_{1}\right) ,\cdot \cdot \cdot
,I_{n_{k}}^{G}\left( f_{k}\right) \right) \\
&=&\sum_{\substack{ \{\sigma :\sigma \wedge \pi ^{\ast }=\hat{0}\text{ ; }
\\ \text{ }\left\vert b\right\vert =2\text{\ }\forall b\in \sigma \}}}%
\bigotimes\limits_{b\in \sigma }\left\langle {\rm St}\,\text{{}}_{\hat{1}%
}^{G,\left[ \left\vert b\right\vert \right] }\right\rangle \left(
f_{1}\otimes _{0}f_{2}\otimes _{0}\cdot \cdot \cdot \otimes _{0}f_{k}\right)
\\
&=&\sum_{\sigma \in \mathcal{M}_{2}^{0}\left( \left[ n\right] ,\pi ^{\ast
}\right) }\bigotimes\limits_{b\in \sigma }\left\langle {\rm St}\,\text{{}}%
_{\hat{1}}^{G,\left[ \left\vert b\right\vert \right] }\right\rangle \left(
f_{1}\otimes _{0}f_{2}\otimes _{0}\cdot \cdot \cdot \otimes _{0}f_{k}\right)
\text{,}
\end{eqnarray*}%
and then by applying the same line of reasoning as above.
\end{proof}

\bigskip

\textbf{Examples}. (i) We want to use Corollary \ref{C : DiaGauss} to
compute the cumulant of the two integrals%
\begin{eqnarray*}
I_{n_{1}}^{G}\left( f_{1}\right) &=&\int_{Z_{\hat{0}}^{n_{1}}}f_{1}\left(
z_{1},...,z_{n_{1}}\right) G\left( dz_{1}\right) \cdot \cdot \cdot G\left(
dz_{n_{1}}\right) \\
I_{n_{2}}^{G}\left( f_{2}\right) &=&\int_{Z_{\hat{0}}^{n_{2}}}f_{2}\left(
z_{1},...,z_{n_{2}}\right) G\left( dz_{1}\right) \cdot \cdot \cdot G\left(
dz_{n_{2}}\right) \text{,}
\end{eqnarray*}%
that is, the quantity
\begin{equation*}
\chi \left( I_{n_{1}}^{G}\left( f_{1}\right) ,I_{n_{2}}^{G}\left(
f_{2}\right) \right) =\mathbb{E}\left( I_{n_{1}}^{G}\left( f_{1}\right)
I_{n_{2}}^{G}\left( f_{2}\right) \right) \text{.}
\end{equation*}%
Here, $\pi ^{\ast }\in \mathcal{P}\left( \left[ n_{1}+n_{2}\right] \right) $
is given by $\pi ^{\ast }=\left\{ \left\{ 1,...,n_{1}\right\} ,\left\{
n_{1}+1,...,n_{1}+n_{2}\right\} \right\} $. It is easily seen that $\mathcal{%
M}_{2}\left( \left[ n_{1}+n_{2}\right] ,\pi ^{\ast }\right) \neq \varnothing
$ if and only if $n_{1}=n_{2}$. Indeed, each partition $\mathcal{M}%
_{2}\left( \left[ n_{1}+n_{2}\right] ,\pi ^{\ast }\right) $ is of the form
\begin{equation}
\sigma =\left\{ \left\{ i_{1},i_{2}\right\} :i_{1}\in \left\{
1,...,n_{1}\right\} ,i_{2}\in \left\{ n_{1}+1,...,n_{1}+n_{2}\right\}
\right\}  \label{sigma}
\end{equation}%
(this is the case because $\sigma $ must have blocks of size $\left\vert
b\right\vert =2$ only, and no blocks can be constructed using only the
indices $\left\{ 1,...,n_{1}\right\} $ or $\left\{
n_{1}+1,...,n_{1}+n_{2}\right\} $, since the corresponding diagram must be
non-flat). In the case where $n_{1}=n_{2}$, there are exactly $n_{1}!$
partitions as in (\ref{sigma}), since to each element in $\left\{
1,...,n_{1}\right\} $ one attaches one element of $\left\{
n_{1}+1,...,n_{1}+n_{2}\right\} $. Moreover, for any such $\sigma $ one has
that%
\begin{equation}
\int_{Z^{n/2}}f_{\sigma ,2}d\nu ^{n/2}=\int_{Z^{n_{1}}}f_{1}f_{2}d\nu
^{n_{1}}\text{,}  \label{from}
\end{equation}%
where $n=n_{1}+n_{2}$ and we have used the symmetry of $f_{1}$ and $f_{2}$
to obtain that
\begin{equation*}
f_{\sigma ,2}\left( z_{1},...,z_{\frac{n}{2}}\right) =f_{\sigma ,2}\left(
z_{1},...,z_{n_{1}}\right) =f_{1}\left( z_{1},...,z_{n_{1}}\right)
f_{2}\left( z_{1},...,z_{n_{1}}\right) \text{.}
\end{equation*}%
From (\ref{GaussDiagr}) and (\ref{from}), we deduce that%
\begin{equation*}
\mathbb{E}\left( I_{n_{1}}^{G}\left( f_{1}\right) I_{n_{2}}^{G}\left(
f_{2}\right) \right) =\mathbf{1}_{n_{1}=n_{2}}\times
n_{1}!\int_{Z^{n_{1}}}f_{1}f_{2}d\nu ^{n_{1}}\text{,}
\end{equation*}%
as expected (see (\ref{goo})). Note also that, since every diagram
associated with $\pi ^{\ast }$ has two rows, one also has
\begin{equation*}
\mathcal{M}_{2}\left( \left[ n_{1}+n_{2}\right] ,\pi ^{\ast }\right) =%
\mathcal{M}_{2}^{0}\left( \left[ n_{1}+n_{2}\right] ,\pi ^{\ast }\right) ,
\end{equation*}%
that is, every non-flat diagram is also connected, thus yielding (thanks to (%
\ref{GaussDiagrammi}) and (\ref{GaussDiagr}))%
\begin{equation*}
\chi \left( I_{n_{1}}^{G}\left( f_{1}\right) ,I_{n_{2}}^{G}\left(
f_{2}\right) \right) =\mathbb{E}\left( I_{n_{1}}^{G}\left( f_{1}\right)
I_{n_{2}}^{G}\left( f_{2}\right) \right) .
\end{equation*}

(ii) We fix an integer $k\geq 3$ and set $n_{1}=...=n_{k}=1$, that is, we
focus on functions $f_{j}$, $j=1,...,k$, of one variable, so that the
integral $I_{1}^{G}\left( f_{j}\right) $ is Gaussian for every $j$, and we
consider $\chi \left( I_{1}^{G}\left( f_{1}\right) ,...,I_{1}^{G}\left(
f_{k}\right) \right) $ and $\mathbb{E}\left[ I_{1}^{G}\left( f_{1}\right)
,...,I_{1}^{G}\left( f_{k}\right) \right] $. In this case, $n_{1}+\cdot
\cdot \cdot +n_{k}=k$, and $\pi ^{\ast }=\left\{ \left\{ 1\right\}
,...,\left\{ k\right\} \right\} =\hat{0}$. For instance, for $k=6$, $\pi
^{\ast }$ is represented in Fig. 19.

\begin{figure}[htbp]
\begin{center}
\psset{unit=0.7cm}
\begin{pspicture}(0,-2.0)(3.0,2.0)
\psframe[linewidth=0.02,dimen=outer](3.0,2.0)(0.0,-2.0)
\psdots[dotsize=0.15](1.48,1.5)
\psdots[dotsize=0.15](1.48,0.9)
\psdots[dotsize=0.15](1.48,0.32)
\psdots[dotsize=0.15](1.48,-0.28)
\psdots[dotsize=0.15](1.48,-0.9)
\psdots[dotsize=0.15](1.48,-1.5)
\end{pspicture}
\caption{\sl A representation of the partition $\hat{0}$}
\end{center}
\end{figure}
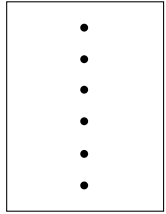

\newpage

\noindent In that case $\mathcal{M}%
_{2}\left( \left[ k\right] ,\pi ^{\ast }\right) =\varnothing $, because all
diagrams will be disconnected. One of such diagrams is represented in Fig. 20.

\begin{figure}[htbp]
\begin{center}
\psset{unit=0.7cm}
\begin{pspicture}(0,-2.0)(3.0,2.0)
\psframe[linewidth=0.02,dimen=outer](3.0,2.0)(0.0,-2.0)
\psdots[dotsize=0.15](1.48,1.5)
\psdots[dotsize=0.15](1.48,0.9)
\psdots[dotsize=0.15](1.48,0.32)
\psdots[dotsize=0.15](1.48,-0.28)
\psdots[dotsize=0.15](1.48,-0.9)
\psdots[dotsize=0.15](1.48,-1.5)
\rput{81.52395}(0.22706485,-2.6834104){\psarc[linewidth=0.02](1.67,-1.21){0.39}{50.194427}{146.54254}}
\rput{-90.0}(0.64,2.46){\psarc[linewidth=0.02](1.55,0.91){0.59}{0.0}{180.0}}
\rput{-90.0}(1.1,1.72){\psarc[linewidth=0.02](1.41,0.31){0.59}{180.0}{0.0}}
\end{pspicture}
\caption{\sl A disconnected Gaussian diagram}
\end{center}
\end{figure}
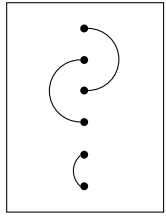
\noindent (\textbf{Exercise}: give an
algebraic proof of the fact that $\mathcal{M}_{2}\left( \left[ k\right] ,\pi
^{\ast }\right) =\varnothing )$. It follows from Point 1 of Corollary \ref{C : DiaGauss} that%
\begin{equation*}
\chi \left( I_{1}^{G}\left( f_{1}\right) ,...,I_{1}^{G}\left( f_{k}\right)
\right) =0
\end{equation*}%
(this is consistent with the properties of cumulants of Gaussian vectors
noted in Section \ref{S : cum}). Now focus on $\mathcal{M}_{2}^{0}\left( %
\left[ k\right] ,\pi ^{\ast }\right) $. If $k$ is odd the class $\mathcal{M}%
_{2}^{0}\left( \left[ k\right] ,\pi ^{\ast }\right) $ is empty, and, for $k$
even, $\mathcal{M}_{2}^{0}\left( \left[ k\right] ,\pi ^{\ast }\right) $
coincides with the collection of all partitions
\begin{equation}
\sigma =\left\{ \left\{ i_{1},j_{1}\right\} ,...,\left\{ i_{\frac{k}{2}},j_{%
\frac{k}{2}}\right\} \right\} \in \mathcal{P}\left( \left[ k\right] \right)
\label{match}
\end{equation}%
whose blocks have size two (that is, $\mathcal{M}_{2}^{0}\left( \left[ k%
\right] ,\pi ^{\ast }\right) $ is the class of all perfect matchings of the
first $k$ integers). For $\sigma $ as in (\ref{match}), we have%
\begin{equation*}
f_{\sigma ,k}\left( z_{1},...,z_{\frac{k}{2}}\right) =\prod_{\substack{ %
\left\{ i_{l},j_{l}\right\} \in \sigma  \\ l=1,...,k/2}}f_{i_{l}}\left(
z_{l}\right) f_{j_{l}}\left( z_{l}\right) .
\end{equation*}%
Points 3 and 4 of Corollary \ref{C : DiaGauss} yield therefore%
\begin{eqnarray*}
&&\mathbb{E}\left( I_{1}^{G}\left( f_{1}\right) \cdot \cdot \cdot
I_{1}^{G}\left( f_{k}\right) \right) \\
&=&\left\{
\begin{array}{ll}
\sum_{\sigma =\left\{ \left\{ i_{1},j_{1}\right\} ,...,\left\{
i_{k/2},j_{k/2}\right\} \right\} \in \mathcal{P}\left( \left[ k\right]
\right) }\int_{Z}f_{i_{1}}f_{j_{1}}d\nu \cdot \cdot \cdot
\int_{Z}f_{i_{k/2}}f_{j_{k/2}}d\nu , & k\text{ even} \\
0, & k\text{ odd.}%
\end{array}%
\right. ,
\end{eqnarray*}%
which is just a special case of (\ref{Feynman}), since $\mathbb{E}\left(
I_{1}^{G}\left( f_{i}\right) I_{1}^{G}\left( f_{j}\right) \right) $ $=$ $%
\int_{Z}f_{i_{1}}f_{j_{1}}d\nu $. For instance, if $k=4$, one has that
\begin{eqnarray*}
\mathbb{E}\left( I_{1}^{G}\left( f_{1}\right) \cdot \cdot \cdot
I_{1}^{G}\left( f_{4}\right) \right) &=&\int_{Z}f_{1}f_{2}d\nu \times
\int_{Z}f_{3}f_{4}d\nu +\int_{Z}f_{1}f_{3}d\nu \times \int_{Z}f_{2}f_{4}d\nu
\\
&&\int_{Z}f_{1}f_{4}d\nu \times \int_{Z}f_{2}f_{3}d\nu .
\end{eqnarray*}

(iii) Consider the case $k=3$, $n_{1}=2$, $n_{2}=n_{3}=1$. Here, $%
n=n_{1}+n_{2}+n_{3}=4$, and $\pi ^{\ast }=\left\{ \left\{ 1,2\right\}
,\left\{ 3\right\} ,\left\{ 4\right\} \right\} $. The partition $\pi ^{\ast
} $ is represented in Fig. 21.

\newpage

\begin{figure}[htbp]
\begin{center}
\psset{unit=0.7cm}
\begin{pspicture}(0,-1.5)(4.0,1.5)
\psframe[linewidth=0.02,dimen=outer](4.0,1.5)(0.0,-1.5)
\psdots[dotsize=0.15](0.6,0.88)
\psdots[dotsize=0.15](3.38,0.88)
\psdots[dotsize=0.15](0.6,-0.92)
\psdots[dotsize=0.15](0.6,-0.02)
\end{pspicture}
\caption{\sl A three-block partition}
\end{center}
\end{figure}
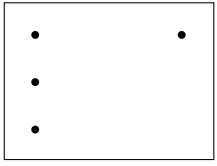

\noindent The class $\mathcal{M}%
_{2}\left( \left[ 4\right] ,\pi ^{\ast }\right) $ contains only two
elements, namely%
\begin{equation*}
\sigma _{1}=\left\{ \left\{ 1,3\right\} ,\left\{ 2,4\right\} \right\} \text{
\ and \ }\sigma _{2}=\left\{ \left\{ 1,4\right\} ,\left\{ 2,3\right\}
\right\} \text{,}
\end{equation*}%
whose diagrams are given in Fig. 22.

\begin{figure}[htbp]
\begin{center}
\psset{unit=0.7cm}
\begin{pspicture}(0,-1.51)(7.98,1.49)
\psframe[linewidth=0.02,dimen=outer](4.0,1.49)(0.0,-1.51)
\psdots[dotsize=0.15](0.6,0.87)
\psdots[dotsize=0.15](3.38,0.87)
\psdots[dotsize=0.15](0.6,-0.93)
\psdots[dotsize=0.15](0.6,-0.03)
\psframe[linewidth=0.02,dimen=outer](7.98,1.49)(3.98,-1.51)
\psdots[dotsize=0.15](4.5,0.87)
\psdots[dotsize=0.15](7.36,0.87)
\psdots[dotsize=0.15](4.5,-0.93)
\psdots[dotsize=0.15](4.6,-0.05)
\psline[linewidth=0.02cm](0.6,0.89)(0.6,-0.01)
\psline[linewidth=0.02cm](0.6,-0.93)(3.36,0.87)
\psline[linewidth=0.02cm](4.6,-0.03)(7.34,0.87)
\rput{-270.0}(5.99,-5.99){\psarc[linewidth=0.02](5.99,0.0){1.71}{59.470295}{120.96375}}
\end{pspicture}
\caption{\sl The two elements of $\mathcal{M}_2([4],\pi^*)$}
\end{center}
\end{figure}
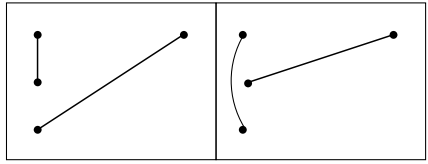
\noindent Since the rows of these diagrams
cannot be divided into two subsets (see Section \ref{SS : diagrams}), they
are connected, and one has $\mathcal{M}_{2}\left( \left[ 4\right] ,\pi
^{\ast }\right) =\mathcal{M}_{2}^{0}\left( \left[ 4\right] ,\pi ^{\ast
}\right) $, that is, cumulants equal moments by (\ref{GaussDiagrammi}) and (%
\ref{GaussDiagr}). Moreover,
\begin{eqnarray*}
f_{\sigma _{1},3}\left( z_{1},z_{2}\right) &=&f_{1}\left( z_{1},z_{2}\right)
f_{2}\left( z_{1}\right) f_{3}\left( z_{2}\right) \\
f_{\sigma _{2},3}\left( z_{1},z_{2}\right) &=&f_{1}\left( z_{1},z_{2}\right)
f_{2}\left( z_{2}\right) f_{3}\left( z_{1}\right) .
\end{eqnarray*}%
It follows that%
\begin{eqnarray*}
\chi \left( I_{2}^{G}\left( f_{1}\right) ,I_{1}^{G}\left( f_{2}\right)
,I_{1}^{G}\left( f_{3}\right) \right) &=&\mathbb{E}\left( I_{2}^{G}\left(
f_{1}\right) I_{1}^{G}\left( f_{2}\right) I_{1}^{G}\left( f_{3}\right)
\right) \\
&=&\int_{Z^{2}}\{f_{\sigma _{1},3}\left( z_{1},z_{2}\right) +f_{\sigma
_{2},3}\left( z_{1},z_{2}\right) \}\nu ^{2}\left( dz_{1},dz_{2}\right) \\
&=&2\int_{Z^{2}}f_{1}\left( z_{1},z_{2}\right) f_{2}\left( z_{1}\right)
f_{3}\left( z_{2}\right) \nu ^{2}\left( dz_{1},dz_{2}\right) \text{,}
\end{eqnarray*}%
where in the last equality we have used the symmetry of $f_{1}$.

(iv) We want to use Point 1 and 2 of Corollary \ref{C : DiaGauss} to compute
the $k$th cumulant%
\begin{equation*}
\chi _{k}\left( I_{2}^{G}\left( f\right) \right) =\chi \underset{k\text{
times.}}{(\underbrace{I_{2}^{G}\left( f\right) ,...,I_{2}^{G}\left( f\right)
})}\text{,}
\end{equation*}%
for every $k\geq 3$. This can be done by specializing formula (\ref%
{GaussDiagrammi}) to the case: $k\geq 3$ and $n_{1}=n_{2}=...=n_{k}=2$.
Here, $n=2k$ and $\pi ^{\ast }=\left\{ \left\{ 1,2\right\} ,\left\{
3,4\right\} ,...,\left\{ 2k-1,2k\right\} \right\} $; for instance, for $k=4$
the partition $\pi ^{\ast }$ can be represented as in Fig. 23.

\newpage

\begin{figure}[htbp]
\begin{center}
\psset{unit=0.7cm}
\begin{pspicture}(0,-1.99)(3.02,1.99)
\psframe[linewidth=0.02,dimen=outer](3.02,1.99)(0.0,-1.99)
\psdots[dotsize=0.15](0.44,1.57)
\psdots[dotsize=0.15](0.44,0.47)
\psdots[dotsize=0.15](0.44,-0.55)
\psdots[dotsize=0.15](0.44,-1.61)
\psdots[dotsize=0.15](2.6,1.57)
\psdots[dotsize=0.15](2.6,0.47)
\psdots[dotsize=0.15](2.6,-0.55)
\psdots[dotsize=0.15](2.6,-1.61)
\psline[linewidth=0.02cm](0.4,1.55)(0.38,1.55)
\end{pspicture}
\caption{\sl A four-block partition}
\end{center}
\end{figure}
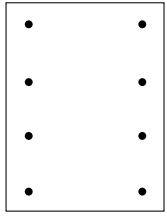 One element of the set $\mathcal{%
M}_{2}\left( \left[ 2k\right] ,\pi ^{\ast }\right) $ is given by the
partition%
\begin{equation*}
\sigma ^{\ast }=\left\{ \left\{ 1,2k\right\} ,\left\{ 2,3\right\} ,\left\{
4,5\right\} ,...,\left\{ 2k-2,2k\right\} \right\} \in \mathcal{P}\left( %
\left[ 2k\right] \right) \text{,}
\end{equation*}%

whose diagram (for $k=4$) appears in Fig. 24.

\begin{figure}[h]
\begin{center}
\psset{unit=0.7cm}
\begin{pspicture}(0,-2.22)(3.02,2.21)
\psframe[linewidth=0.02,dimen=outer](3.02,2.21)(0.0,-1.77)
\psdots[dotsize=0.15](0.44,1.79)
\psdots[dotsize=0.15](0.44,0.69)
\psdots[dotsize=0.15](0.44,-0.33)
\psdots[dotsize=0.15](0.44,-1.39)
\psdots[dotsize=0.15](2.6,1.79)
\psdots[dotsize=0.15](2.6,0.69)
\psdots[dotsize=0.15](2.6,-0.33)
\psdots[dotsize=0.15](2.6,-1.39)
\psline[linewidth=0.02cm](0.44,1.79)(0.38,1.77)
\psline[linewidth=0.02cm](0.44,0.69)(2.54,1.77)
\psline[linewidth=0.02cm](0.44,-0.33)(2.6,0.71)
\psline[linewidth=0.02cm](0.44,-1.39)(2.58,-0.33)
\psbezier[linewidth=0.02](0.38,1.81)(0.38,1.01)(2.58,-2.21)(2.58,-1.41)
\end{pspicture}
\caption{\sl A circular diagram with four rows}
\end{center}
\end{figure}
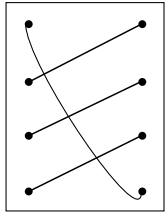

\noindent Note that such a diagram is circular, and that the corresponding multigraph looks like the one in Fig.
16. Therefore,
\begin{equation}
f_{\sigma ^{\ast },k}\left( z_{1},...,z_{k}\right) =f\left(
z_{1},z_{2}\right) f\left( z_{2},z_{3}\right) \cdot \cdot \cdot f\left(
z_{k-1},z_{k}\right) f\left( z_{k},z_{1}\right) .  \label{chose}
\end{equation}%
It is not difficult to see that $\mathcal{M}_{2}\left( \left[ 2k\right] ,\pi
^{\ast }\right) $ contains exactly $2^{k-1}\left( k-1\right) !$ elements and
that the diagram $\Gamma \left( \pi ^{\ast },\sigma \right) $ associated to
each $\sigma \in \mathcal{M}_{2}\left( \left[ 2k\right] ,\pi ^{\ast }\right)
$ is equivalent (up to a permutation, or equivalently to a renumbering, of the rows) to a circular diagram (see
Section \ref{SS : diagrams}). It follows that, for every $\sigma \in
\mathcal{M}_{2}\left( \left[ 2k\right] ,\pi ^{\ast }\right) $, one has that%
\begin{equation*}
f_{\sigma ,k}\left( z_{1},...,z_{k}\right) =f_{\sigma ^{\ast },k}\left(
z_{1},...,z_{k}\right) \text{,}
\end{equation*}%
where $f_{\sigma ^{\ast },2k}$ is given in (\ref{chose}). This yields the
classic formula (see e.g. \cite{FoxTaq})
\begin{eqnarray}
&&\chi _{k}\left( I_{2}^{G}\left( f\right) \right)  \label{qsq} \\
&=&2^{k-1}\left( k-1\right) !\int_{Z^{k}}f\left( z_{1},z_{2}\right) f\left(
z_{2},z_{3}\right) \cdot \cdot \cdot f\left( z_{k-1},z_{k}\right) f\left(
z_{k},z_{1}\right) \nu \left( dz_{1}\right) \cdot \cdot \cdot \nu \left(
dz_{k}\right) \text{.}  \notag
\end{eqnarray}%
For a non-combinatorial proof of (\ref{qsq}) see e.g. \cite[Section 2.2]%
{NouPec2007}.

\subsection{The Poisson case\label{SS : PossCase}}

The following result provides diagram formulae for Wiener-It\^{o} integrals
with respect to compensated Poisson measures. It is stated only for
elementary functions, so as not to have to deal with convergence issues. The
proof is similar to the one of Corollary \ref{C : DiaGauss}, and it is only
sketched. We let $\mathcal{M}_{\geq 2}\left( \left[ n\right] ,\pi ^{\ast
}\right) $ and $\mathcal{M}_{\geq 2}^{0}\left( \left[ n\right] ,\pi ^{\ast
}\right) $ be defined as in Section \ref{SS : GaussCase}.

\begin{corollary}[Poisson measures]
\label{C : DiaPoiss}Suppose $\varphi =\hat{N}$ is a centered Poisson measure
with non-atomic control measure $\nu $, fix integers $n_{1},...,n_{k}\geq 1$
and let $n=n_{1}+\cdot \cdot \cdot +n_{k}$. Write $\pi ^{\ast }$ for the
partition of $\left[ n\right] $ appearing in (\ref{pistar}). Then, for any
vector of functions $(f_{1},...,f_{k})$ such that $f_{j}\in \mathcal{E}%
_{s,0}\left( \nu ^{n_{j}}\right) $, $j=1,...,k$, the following relations
hold:

\begin{enumerate}
\item If $\mathcal{M}_{\geq 2}\left( \left[ n\right] ,\pi ^{\ast }\right)
=\varnothing $, then $\mathbb{\chi }\left( I_{n_{1}}^{\hat{N}}\left(
f_{1}\right) ,\cdot \cdot \cdot ,I_{n_{k}}^{\hat{N}}\left( f_{k}\right)
\right) =0$;

\item If $\mathcal{M}_{\geq 2}\left( \left[ n\right] ,\pi ^{\ast }\right)
\neq \varnothing $, then
\begin{equation}
\mathbb{\chi }\left( I_{n_{1}}^{\hat{N}}\left( f_{1}\right) ,\cdot \cdot
\cdot ,I_{n_{k}}^{\hat{N}}\left( f_{k}\right) \right) =\sum_{\sigma \in
\mathcal{M}_{\geq 2}\left( \left[ n\right] ,\pi ^{\ast }\right)
}\int_{Z^{\left\vert \sigma \right\vert }}f_{\sigma ,k}d\nu ^{\left\vert
\sigma \right\vert },  \label{poissDiagrammi}
\end{equation}%
where, for every $\sigma \in \mathcal{M}_{\geq 2}\left( \left[ n\right] ,\pi
^{\ast }\right) $, the function $f_{\sigma ,k}$, in $\left\vert \sigma
\right\vert $ variables, is obtained by identifying the variables $x_{i}$
and $x_{j}$ in the argument of $f_{1}\otimes _{0}\cdot \cdot \cdot \otimes
_{0}f_{n_{k}}$ (as defined in (\ref{0multContr})) if and only if $i\sim
_{\sigma }j$;

\item If $\mathcal{M}_{\geq 2}^{0}\left( \left[ n\right] ,\pi ^{\ast
}\right) =\varnothing $, then $\mathbb{E}\left( I_{n_{1}}^{\hat{N}}\left(
f_{1}\right) \cdot \cdot \cdot I_{n_{k}}^{\hat{N}}\left( f_{k}\right)
\right) =0$;

\item If $\mathcal{M}_{\geq 2}^{0}\left( \left[ n\right] ,\pi ^{\ast
}\right) \neq \varnothing $,
\begin{equation}
\mathbb{E}\left( I_{n_{1}}^{\hat{N}}\left( f_{1}\right) ,\cdot \cdot \cdot
,I_{n_{k}}^{\hat{N}}\left( f_{k}\right) \right) =\sum_{\sigma \in \mathcal{M}%
_{\geq 2}^{0}\left( \left[ n\right] ,\pi ^{\ast }\right)
}\int_{Z^{\left\vert \sigma \right\vert }}f_{\sigma ,k}d\nu ^{\left\vert
\sigma \right\vert }.  \label{poiss mom}
\end{equation}
\end{enumerate}
\end{corollary}

\noindent \begin{proof}[Sketch of the Proof]
The proof follows closely that of Corollary \ref{C : DiaGauss}. The only
difference is in evaluating (\ref{Multexp1}). Instead of having (\ref%
{gaussCombExp}) which requires considering $\mathcal{M}_{2}$ and $\mathcal{M}%
_{2}^{0}$, one has (\ref{PoissCombExp}), which implies that one must use $%
\mathcal{M}_{\geq 2}$ and $\mathcal{M}_{\geq 2}^{0}$.
\end{proof}

\bigskip

\textbf{Remark. }Corollaries \ref{C : DiaGauss} and \ref{C : DiaPoiss} are
quite similar. In the Poisson case, however, $f_{\sigma ,k}$ depends on $%
\left\vert \sigma \right\vert $ variables, whereas in the Gaussian case it
depends on $n/2$ variables.

\bigskip

\textbf{Examples. }All kernels appearing in the following examples are
symmetric, elementary and vanishing on diagonals (this ensures that multiple
integrals have moments of all orders).

(i) We apply Corollary \ref{C : DiaPoiss} in order to compute the cumulant
\begin{equation*}
\chi \left( I_{n_{1}}^{\hat{N}}\left( f_{1}\right) ,I_{n_{2}}^{\hat{N}%
}\left( f_{2}\right) \right) =\mathbb{E}\left( I_{n_{1}}^{\hat{N}}\left(
f_{1}\right) I_{n_{2}}^{\hat{N}}\left( f_{2}\right) \right) \text{,}
\end{equation*}%
where $n_{1},n_{2}\geq 1$ are arbitrary. In this case, $\pi ^{\ast }\in
\mathcal{P}\left( \left[ n_{1}+n_{2}\right] \right) $ is given by
\begin{equation*}
\pi ^{\ast }=\left\{ \left\{ 1,...,n_{1}\right\} ,\left\{
n_{1}+1,...,n_{1}+n_{2}\right\} \right\} .
\end{equation*}%
Moreover,
\begin{eqnarray*}
\mathcal{M}_{2}^{0}\left( \left[ n_{1}+n_{2}\right] ,\pi ^{\ast }\right) &=&%
\mathcal{M}_{2}\left( \left[ n_{1}+n_{2}\right] ,\pi ^{\ast }\right) \\
&=&\mathcal{M}_{\geq 2}\left( \left[ n_{1}+n_{2}\right] ,\pi ^{\ast }\right)
=\mathcal{M}_{\geq 2}^{0}\left( \left[ n_{1}+n_{2}\right] ,\pi ^{\ast
}\right)
\end{eqnarray*}%
(indeed, since any diagram of $\pi ^{\ast }$ is composed of two rows, every
non-flat diagram must be necessarily connected and Gaussian). This gives, in
particular, $\mathcal{M}_{\geq 2}\left( \left[ n_{1}+n_{2}\right] ,\pi
^{\ast }\right) \neq \varnothing $ if and only if $n_{1}=n_{2}$. The
computations performed in the Gaussian case thus apply and therefore yield%
\begin{equation*}
\chi \left( I_{n_{1}}^{\hat{N}}\left( f_{1}\right) ,I_{n_{2}}^{\hat{N}%
}\left( f_{2}\right) \right) =\mathbb{E}\left( I_{n_{1}}^{\hat{N}}\left(
f_{1}\right) I_{n_{2}}^{\hat{N}}\left( f_{2}\right) \right) =\mathbf{1}%
_{n_{1}=n_{2}}\times n_{1}!\int_{Z^{n_{1}}}f_{1}f_{2}d\nu ^{n_{1}}\text{,}
\end{equation*}%
which is once again consistent with (\ref{goo}).

(ii) Consider the case $k=3$, $n_{1}=n_{2}=2$, $n_{3}=1$. Here, $%
n=n_{1}+n_{2}+n_{3}=5$, and $\pi ^{\ast }=\left\{ \left\{ 1,2\right\}
,\left\{ 3,4\right\} ,\left\{ 5\right\} \right\} $. The partition $\pi
^{\ast }$ can be represented as in Fig. 25.

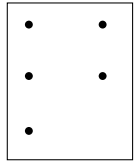
\begin{figure}[htbp]
\begin{center}
\psset{unit=0.7cm}
\begin{pspicture}(0,-1.5)(2.4,1.5)
\psframe[linewidth=0.02,dimen=outer](2.4,1.5)(0.0,-1.5)
\psdots[dotsize=0.15](0.42,1.08)
\psdots[dotsize=0.15](0.42,0.1)
\psdots[dotsize=0.15](0.42,-0.94)
\psdots[dotsize=0.15](1.82,1.08)
\psdots[dotsize=0.15](1.82,0.1)
\end{pspicture}
\caption{\sl A three-row partition}
\end{center}
\end{figure}

\noindent The class $\mathcal{M}_{\geq
2}^{0}\left( \left[ 5\right] ,\pi ^{\ast }\right) $, of $\sigma $'s such
that $\sigma \wedge \pi ^{\ast }=\hat{0}$, contains four elements, that is,%
\begin{eqnarray*}
\sigma _{1} &=&\left\{ \left\{ 1,3,5\right\} ,\left\{ 2,4\right\} \right\}
\text{, \ }\sigma _{2}=\left\{ \left\{ 1,4\right\} ,\left\{ 2,3,5\right\}
\right\} \\
\sigma _{3} &=&\left\{ \left\{ 1,3\right\} ,\left\{ 2,4,5\right\} \right\}
\text{ \ and \ }\sigma _{4}=\left\{ \left\{ 1,4,5\right\} ,\left\{
2,3\right\} \right\} \text{,}
\end{eqnarray*}%
whose diagrams are given in Fig. 26.

\begin{figure}[htbp]
\begin{center}
\psset{unit=0.7cm}
\begin{pspicture}(0,-1.7027255)(9.54,1.7127255)
\psframe[linewidth=0.02,dimen=outer](2.4,1.2972745)(0.0,-1.7027255)
\psdots[dotsize=0.15](0.42,0.8772745)
\psdots[dotsize=0.15](0.42,-0.102725506)
\psdots[dotsize=0.15](0.42,-1.1427255)
\psdots[dotsize=0.15](1.82,0.8772745)
\psdots[dotsize=0.15](1.82,-0.102725506)
\psframe[linewidth=0.02,dimen=outer](4.78,1.2972745)(2.38,-1.7027255)
\psdots[dotsize=0.15](2.8,0.8772745)
\psdots[dotsize=0.15](2.8,-0.102725506)
\psdots[dotsize=0.15](2.8,-1.1427255)
\psdots[dotsize=0.15](4.2,0.8772745)
\psdots[dotsize=0.15](4.2,-0.102725506)
\psframe[linewidth=0.02,dimen=outer](7.16,1.2972745)(4.76,-1.7027255)
\psdots[dotsize=0.15](5.18,0.8772745)
\psdots[dotsize=0.15](5.18,-0.102725506)
\psdots[dotsize=0.15](5.18,-1.1427255)
\psdots[dotsize=0.15](6.58,0.8772745)
\psdots[dotsize=0.15](6.58,-0.102725506)
\psframe[linewidth=0.02,dimen=outer](9.54,1.2972745)(7.14,-1.7027255)
\psdots[dotsize=0.15](7.56,0.8772745)
\psdots[dotsize=0.15](7.56,-0.102725506)
\psdots[dotsize=0.15](7.56,-1.1427255)
\psdots[dotsize=0.15](8.96,0.8772745)
\psdots[dotsize=0.15](8.96,-0.102725506)
\psellipse[linewidth=0.02,dimen=outer](0.42,-0.1727255)(0.24,1.29)
\psellipse[linewidth=0.02,dimen=outer](1.82,0.3572745)(0.26,0.78)
\psellipse[linewidth=0.02,dimen=outer](5.18,0.37727448)(0.26,0.78)
\psbezier[linewidth=0.02](2.5184243,0.8894871)(2.6222162,0.09624863)(4.6653676,-1.0281765)(4.561576,-0.23493807)(4.4577837,0.5583004)(2.4146323,1.6827255)(2.5184243,0.8894871)
\psline[linewidth=0.02](2.48,-0.042725503)(2.48,-1.5027255)(3.12,-1.5027255)(3.12,-0.3627255)(4.6,0.7972745)(4.24,1.2172745)(2.48,-0.042725503)(2.48,-0.042725503)
\psline[linewidth=0.02](6.94,-0.2627255)(6.94,1.1972744)(6.3,1.1972744)(6.3,0.057274498)(4.82,-1.1027255)(5.18,-1.5227255)(6.94,-0.2627255)(6.94,-0.2627255)
\psline[linewidth=0.02](7.32,1.0772744)(7.32,0.5772745)(8.56,-0.2027255)(7.34,-0.8627255)(7.34,-1.5027255)(9.42,-0.2027255)(7.32,1.2172745)(7.32,1.0372745)(7.32,1.0372745)
\psbezier[linewidth=0.02](9.321576,0.90948707)(9.217784,0.11624862)(7.174632,-1.0081766)(7.2784243,-0.21493807)(7.3822165,0.57830036)(9.425367,1.7027255)(9.321576,0.90948707)
\end{pspicture}
\caption{\sl The four elements of $\mathcal{M}_{\ge 2}([5],\pi^*)$}
\end{center}
\end{figure}

\noindent Since all these diagrams are
connected, the class $\mathcal{M}_{\geq 2}\left( \left[ 5\right] ,\pi ^{\ast
}\right) $ coincides with $\mathcal{M}_{\geq 2}^{0}\left( \left[ 5\right]
,\pi ^{\ast }\right) $. Note also that, since the above diagrams have an odd
number of vertices, $\mathcal{M}_{\geq 2}\left( \left[ 5\right] ,\pi ^{\ast
}\right) $ does not contain partitions $\sigma $ whose diagram is Gaussian.
Thus,
\begin{eqnarray*}
f_{\sigma _{1},3}\left( z_{1},z_{2}\right) &=&f_{1}\left( z_{1},z_{2}\right)
f_{2}\left( z_{1},z_{2}\right) f_{3}\left( z_{1}\right) \\
f_{\sigma _{2},3}\left( z_{1},z_{2}\right) &=&f_{1}\left( z_{1},z_{2}\right)
f_{2}\left( z_{2},z_{1}\right) f_{3}\left( z_{2}\right) \\
f_{\sigma _{3},3}\left( z_{1},z_{2}\right) &=&f_{1}\left( z_{1},z_{2}\right)
f_{2}\left( z_{1},z_{2}\right) f_{3}\left( z_{2}\right) \\
f_{\sigma _{4},3}\left( z_{1},z_{2}\right) &=&f_{1}\left( z_{1},z_{2}\right)
f_{2}\left( z_{2},z_{1}\right) f_{3}\left( z_{1}\right) .
\end{eqnarray*}%
For instance, $f_{\sigma _{1},3}\left( z_{1},z_{2}\right) $ has been
obtained by identifying the variables of $f_{1}\left( x_{1},x_{2}\right)
f_{2}\left( x_{3},x_{4}\right) f_{3}\left( x_{5}\right) $ as $%
x_{1}=x_{3}=x_{5}=z_{1}$ and $x_{2}=x_{4}=z_{2}$. By exploiting the symmetry
of $f_{1}$ and $f_{2}$, one deduces that the four quantities
\begin{equation*}
\int_{Z^{2}}f_{\sigma _{i},3}\left( z_{1},z_{2}\right) \nu ^{2}\left(
dz_{1},dz_{2}\right) \text{, \ \ }i=1,...,4\text{,}
\end{equation*}%
are equal. It follows from (\ref{poissDiagrammi}) and (\ref{poiss mom}) that
\begin{eqnarray*}
\chi \left( I_{2}^{\hat{N}}\left( f_{1}\right) ,I_{2}^{\hat{N}}\left(
f_{2}\right) ,I_{1}^{\hat{N}}\left( f_{3}\right) \right) &=&\mathbb{E}\left(
I_{2}^{\hat{N}}\left( f_{1}\right) I_{2}^{\hat{N}}\left( f_{2}\right) I_{1}^{%
\hat{N}}\left( f_{3}\right) \right) \\
&=&4\int_{Z^{2}}\left\{ f_{1}\left( z_{1},z_{2}\right) f_{2}\left(
z_{1},z_{2}\right) f_{3}\left( z_{1}\right) \right\} \nu ^{2}\left(
dz_{1},dz_{2}\right) .
\end{eqnarray*}

(iii) Consider the case $k=4$ and $n_{i}=1$, $i=1,...,4.$ Here, $\pi ^{\ast
}=\hat{0}=\left\{ \left\{ 1\right\} ,\left\{ 2\right\} ,\left\{ 3\right\}
,\left\{ 4\right\} \right\} $, and consequently $\pi ^{\ast }$ can be
represented as a single column of four vertices. The class $\mathcal{M}%
_{\geq 2}\left( \left[ 4\right] ,\pi ^{\ast }\right) $ contains only the
maximal partition $\hat{1}=\left\{ \left\{ 1,2,3,4\right\} \right\} $,
whereas $\mathcal{M}_{\geq 2}^{0}\left( \left[ 5\right] ,\pi ^{\ast }\right)
$ contains $\hat{1}$ and the three elements
\begin{eqnarray*}
\sigma _{1} &=&\left\{ \left\{ 1,2\right\} ,\left\{ 3,4\right\} \right\}
\text{, }\sigma _{2}=\left\{ \left\{ 1,3\right\} ,\left\{ 2,4\right\}
\right\} \text{, \ and} \\
\sigma _{3} &=&\left\{ \left\{ 1,4\right\} ,\left\{ 2,3\right\} \right\} .
\end{eqnarray*}%
The diagrams associated with the class $\mathcal{M}_{\geq 2}^{0}\left( \left[
5\right] ,\pi ^{\ast }\right) =\left\{ \hat{1},\sigma _{1},\sigma
_{2},\sigma _{3}\right\} $ are represented in Fig. 27.

\begin{figure}[htbp]
\begin{center}
\psset{unit=0.7cm}
\begin{pspicture}(0,-1.5)(9.54,1.5)
\psframe[linewidth=0.02,dimen=outer](2.4,1.5)(0.0,-1.5)
\psframe[linewidth=0.02,dimen=outer](4.78,1.5)(2.38,-1.5)
\psframe[linewidth=0.02,dimen=outer](7.16,1.5)(4.76,-1.5)
\psframe[linewidth=0.02,dimen=outer](9.54,1.5)(7.14,-1.5)
\psdots[dotsize=0.15](1.2,0.98)
\psdots[dotsize=0.15](1.2,0.38)
\psdots[dotsize=0.15](1.2,-0.22)
\psdots[dotsize=0.15](1.2,-0.82)
\psdots[dotsize=0.15](3.6,0.98)
\psdots[dotsize=0.15](3.6,0.38)
\psdots[dotsize=0.15](3.6,-0.22)
\psdots[dotsize=0.15](3.6,-0.82)
\psdots[dotsize=0.15](6.0,0.98)
\psdots[dotsize=0.15](6.0,0.38)
\psdots[dotsize=0.15](6.0,-0.22)
\psdots[dotsize=0.15](6.0,-0.82)
\psdots[dotsize=0.15](8.4,0.98)
\psdots[dotsize=0.15](8.4,0.38)
\psdots[dotsize=0.15](8.4,-0.22)
\psdots[dotsize=0.15](8.4,-0.82)
\psellipse[linewidth=0.02,dimen=outer](1.2,0.09)(0.24,1.13)
\psline[linewidth=0.02cm](3.6,1.02)(3.6,0.4)
\psline[linewidth=0.02cm](3.6,-0.16)(3.6,-0.82)
\rput{-270.0}(6.92,-6.18){\psarc[linewidth=0.02](6.55,0.37){0.85}{45.0}{133.21008}}
\rput{-270.0}(5.22,-5.68){\psarc[linewidth=0.02](5.45,-0.23){0.85}{226.78992}{315.0}}
\psline[linewidth=0.02cm](8.4,0.38)(8.4,-0.18)
\rput{-270.0}(9.14,-8.98){\psarc[linewidth=0.02](9.06,0.08){1.14}{37.568592}{140.71059}}
\end{pspicture}
\caption{\sl The elements of $\mathcal{M}^0_{\ge 2}([5],\pi^*)$}
\end{center}
\end{figure}
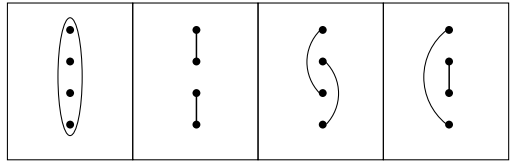

\noindent Now take $f_{i}=f$ for $i=1,...,4$%
, where $f$ is an elementary kernel. One has that
\begin{eqnarray*}
f_{\hat{1},4}\left( z\right) &=&f\left( z\right) ^{4} \\
f_{\sigma _{j},4}\left( z_{1},z_{2}\right) &=&f\left( z_{1}\right)
^{2}f\left( z_{2}\right) ^{2}\text{, \ \ }j=1,2,3\text{.}
\end{eqnarray*}%
It follows from (\ref{poissDiagrammi}) and (\ref{poiss mom}) that
\begin{eqnarray*}
\chi \left( I_{1}^{\hat{N}}\left( f\right) ,I_{1}^{\hat{N}}\left( f\right)
,I_{1}^{\hat{N}}\left( f\right) ,I_{1}^{\hat{N}}\left( f\right) \right)
&=&\chi _{4}\left( I_{1}^{\hat{N}}\left( f\right) \right) =\int_{Z}f\left(
z\right) ^{4}\nu \left( dz\right) \\
\mathbb{E}\left( I_{1}^{\hat{N}}\left( f\right) ^{4}\right)
&=&\int_{Z}f\left( z\right) ^{4}\nu \left( dz\right) \\
&&+3\int_{Z^{2}}f\left( z_{1}\right) ^{2}f\left( z_{2}\right) ^{2}\nu
^{2}\left( dz_{1},dz_{2}\right) \text{.}
\end{eqnarray*}

(iv) Let $Y$ be a centered random variable with finite moments of all
orders, and suppose that $Y$ is infinitely divisible and such that%
\begin{equation}
\mathbb{E}\left[ \exp \left( i\lambda Y\right) \right] =\exp \left[ \int_{%
\mathbb{R}}\left( e^{i\lambda u}-1-i\lambda u\right) \rho \left( du\right) %
\right] ,  \label{LKuu}
\end{equation}%
where the measure $\rho $ is such that $\rho \left( \left\{ 0\right\}
\right) =0$ and $\int_{\mathbb{R}}|u|^{k}\rho \left( du\right) <\infty $ for
every $k\geq 1$. Then, combining (\ref{LKuu}) and (\ref{cumDef}), one
deduces that
\begin{equation}
\chi _{k}\left( Y\right) =\int_{\mathbb{R}}u^{k}\rho \left( du\right) \text{%
, \ \ }k\geq 2,  \label{IDcim}
\end{equation}%
(note that $\chi _{1}\left( Y\right) =\mathbb{E}\left( Y\right) =0$). We
shall prove that (\ref{IDcim}) is consistent with (\ref{poissDiagrammi}).
Indeed, according to the discussion contained in Section \ref{SS : infDIV},
one has that
\begin{equation*}
Y\overset{law}{=}\int_{\mathbb{R}}\int_{0}^{1}u\hat{N}\left( du,dx\right)
=I_{1}^{\hat{N}}\left( f\right) \text{,}
\end{equation*}%
where $f\left( u,x\right) =u\mathbf{1}_{\left[ 0,1\right] }\left( x\right) $%
, and $\hat{N}$ is a centered Poisson measure on $\left[ 0,1\right] \times
\mathbb{R}$, with control $\rho \left( du\right) dx$. It follows that
\begin{equation}
\chi _{k}\left( Y\right) =\chi _{k}(I_{1}^{\hat{N}}\left( f\right) )=\chi (%
\underset{k\text{ times}}{\underbrace{I_{1}^{\hat{N}}\left( f\right)
,...,I_{1}^{\hat{N}}\left( f\right) }})\text{.}  \label{cumcum}
\end{equation}%
The RHS of (\ref{cumcum}) can be computed by means of Corollary \ref{C :
DiaPoiss} in the special case where $n_{j}=1$ ($\forall j=1,...,k$), $%
n=\Sigma _{j}n_{j}=k$, and $\pi ^{\ast }=\hat{0}=\left\{ \left\{ 1\right\}
,...,\left\{ k\right\} \right\} $. One has clearly that $\hat{1}$ is the
only partition such that the diagram $\Gamma \left( \hat{0},\hat{1}\right) $
is connected, so that
\begin{equation*}
\mathcal{M}_{\geq 2}\left( \left[ k\right] ,\pi ^{\ast }\right) =\left\{
\hat{1}\right\} =\left\{ \left\{ 1,...,k\right\} \right\} .
\end{equation*}%
Since%
\begin{equation*}
f_{\hat{1},k}\left( u,x\right) =u^{k}\mathbf{1}_{\left[ 0,1\right] }\left(
x\right) ,
\end{equation*}%
we can now use (\ref{poissDiagrammi}) to deduce that%
\begin{equation*}
\chi _{k}(I_{1}^{\hat{N}}\left( f\right) )=\int_{\mathbb{R}}\int_{0}^{1}f_{%
\hat{1},k}\left( u,x\right) \rho \left( du\right) dx=\int_{0}^{1}dx\int_{%
\mathbb{R}}u^{k}\rho \left( du\right) =\int_{\mathbb{R}}u^{k}\rho \left(
du\right) \text{.}
\end{equation*}

(v) As an explicit example of (\ref{IDcim}), consider the case where $Y$ is
a centered Gamma random variable with shape parameter $a>0$ and unitary
scale parameter, that is,%
\begin{equation*}
\mathbb{E}\left[ \exp \left( i\lambda Y\right) \right] =\frac{e^{-i\lambda a}%
}{\left( 1-i\lambda \right) ^{a}}=\exp \left[ a\int_{0}^{\infty }\left(
e^{i\lambda u}-1-i\lambda u\right) e^{-u}\frac{du}{u}\right] \text{.}
\end{equation*}%
Thus, $\rho(du) = a{\bf 1}_{\{u>0\}} u^{-1} e^{-u} du$. It follows that, $\chi _{1}\left( Y\right) =\mathbb{E}\left( Y\right) =0$
and, for $k\geq 2$,
\begin{equation*}
\chi _{k}\left( Y\right) =a\int_{\mathbb{R}}u^{k}e^{-u}\frac{du}{u}=a\Gamma
\left( k\right) =a\left( k-1\right) !\text{.}
\end{equation*}

\section{From Gaussian measures to isonormal Gaussian processes\label{S :
IsonormalGP}}

\setcounter{equation}{0}

For the sake of completeness, in this section we show how to generalize part
of the previous results to the case of an \textsl{isonormal Gaussian process}%
. These objects have been introduced by R.M.\ Dudley in \cite{Dudley}, and
are a natural generalization of the Gaussian measures introduced in Section %
\ref{SS : CRM}. In particular, the concept of isonormal Gaussian process can
be very useful in the study of fractional fields. See e.g. Pipiras and Taqqu
\cite{PipTaqPtrf, PipTaq, PipTaqSurv}, or the second edition of Nualart's book
\cite{Nualart}. For a general approach to Gaussian analysis by means of
Hilbert space techniques, and for further details on the subjects discussed
in this section, the reader is referred to Janson \cite{Janson}.

\subsection{General definitions and examples}

Let $\mathfrak{H}$ be a real separable Hilbert space with inner product $%
\left( \cdot ,\cdot \right) _{\mathfrak{H}}$. In what follows, we will
denote by%
\begin{equation*}
X=X\left( \mathfrak{H}\right) =\left\{ X\left( h\right) :h\in \mathfrak{H}%
\right\}
\end{equation*}%
an \textsl{isonormal Gaussian process }over $\mathfrak{H}$. This means that $%
X$ is a centered real-valued Gaussian family, indexed by the elements of $%
\mathfrak{H}$ and such that%
\begin{equation}
\mathbb{E}\left[ X\left( h\right) X\left( h^{\prime }\right) \right] =\left(
h,h^{\prime }\right) _{\mathfrak{H}}\text{, \ \ }\forall h,h^{\prime }\in
\mathfrak{H}.  \label{Inner}
\end{equation}%
In other words, relation (\ref{Inner}) means that $X$ is a centered Gaussian
Hilbert space (with respect to the inner product canonically induced by the
covariance) isomorphic to $\mathfrak{H}$.

\bigskip

\textbf{Example} (\textit{Euclidean spaces}). Fix an integer $d\geq 1$, set $%
\mathfrak{H}=\mathbb{R}^{d}$ and let $\left( e_{1},...,e_{d}\right) $ be an
orthonormal basis of $\mathbb{R}^{d}$ (with respect to the usual Euclidean
inner product). Let $\left( Z_{1},...,Z_{d}\right) $ be a Gaussian vector
whose components are i.i.d. $N\left( 0,1\right) $. For every $%
h=\sum_{j=1}^{d}c_{j}e_{j}$ (where the $c_{j}$ are real and uniquely
defined), set $X\left( h\right) =\sum_{j=1}^{d}c_{j}Z_{j}$ and define $%
X=\left\{ X\left( h\right) :h\in \mathbb{R}^{d}\right\} $. Then, $X$ is an
isonormal Gaussian process over $\mathbb{R}^{d}$.

\bigskip

\textbf{Example} (\textit{Gaussian measures}). Let $\left( Z,\mathcal{Z},\nu
\right) $ be a measure space, where $\nu $ is positive, $\sigma $-finite and
non atomic. Consider a completely random Gaussian measure $G=\left\{ G\left(
A\right) :A\in \mathcal{Z}_{\nu }\right\} $ (as defined in Section \ref{SS :
CRM}), where the class $\mathcal{Z}_{\nu }$ is given by (\ref{ZetaNu}). Set $%
\mathfrak{H}=L^{2}\left( Z,\mathcal{Z},\nu \right) $ (thus, for every $h,h'\in\mathfrak{H}$, $(h,h')_\mathfrak{H} = \int_Z h(z)h'(z)\nu(dz)$) and, for every $h\in
\mathfrak{H}$, define $X\left( h\right) =I_{1}^{G}\left( h\right) $ to be
the Wiener-It\^{o} integral of $h$ with respect to $G$, as defined in (\ref%
{WiIto1}). Recall that $X\left( h\right) $ is a centered Gaussian random
variable with variance given by $\left\Vert h\right\Vert _{\mathfrak{H}}^{2}$%
. Then, relation (\ref{ISO1}) implies that the collection $X=\left\{ X\left(
h\right) :h\in L^{2}\left( Z,\mathcal{Z},\nu \right) \right\} $ is an
isonormal Gaussian process over $L^{2}\left( Z,\mathcal{Z},\nu \right) $.

\bigskip

\textbf{Example} (\textit{Isonormal spaces built from covariances}). Let $%
Y=\left\{ Y_{t}:t\geq 0\right\} $ be a real-valued centered Gaussian\
process indexed by the positive axis, and set $R\left( s,t\right) =\mathbb{E}%
\left[ Y_{s}Y_{t}\right] $ to be the covariance function of $Y$. Then, one
can embed $Y$ into some isonormal Gaussian process as follows: (i) define $%
\mathcal{E}$\ as the collection of all finite linear combinations of
indicator functions of the type $\mathbf{1}_{\left[ 0,t\right] }$, $t\geq 0$%
; (ii) define $\mathfrak{H}=\mathfrak{H}_{R}$ to be the Hilbert space given
by the closure of $\mathcal{E}$ with respect to the inner product%
\begin{equation*}
\left( f,h\right) _{R}:=\sum_{i,j}a_{i}c_{j}R\left( s_{i},t_{j}\right) \text{%
,}
\end{equation*}%
where $f=\sum_{i}a_{i}\mathbf{1}_{\left[ 0,s_{i}\right] }$ and $%
h=\sum_{j}c_{j}\mathbf{1}_{\left[ 0,t_{j}\right] }\ $are two generic
elements of $\mathcal{E}$; (iii) for $h=\sum_{j}c_{j}\mathbf{1}_{\left[
0,t_{j}\right] }\in \mathcal{E}$, set $X\left( h\right)
=\sum_{j}c_{j}Y_{t_{j}}$; (iv) for $h\in \mathfrak{H}_{R}$, set $X\left(
h\right) $\ to be the $L^{2}\left( \mathbb{P}\right) $ limit\ of any
sequence of the type $X\left( h_{n}\right) $, where $\left\{ h_{n}\right\}
\subset \mathcal{E}$ converges to $h$ in $\mathfrak{H}_{R}$. Note that such
a sequence $\left\{ h_{n}\right\} $ necessarily exists and may not be unique
(however, the definition of $X\left( h\right) $ does not depend on the
choice of the sequence $\left\{ h_{n}\right\} $). Then, by construction, the
Gaussian space $\left\{ X\left( h\right) :h\in \mathfrak{H}\right\} $ is an
isonormal Gaussian process over $\mathfrak{H}_{R}$. See Janson \cite[Ch. 1]%
{Janson} or Nualart \cite{Nualart} for more details on this construction.

\bigskip

\textbf{Example} (\textit{Even functions and symmetric measures}).\textbf{\ }%
Other classic examples of isonormal Gaussian processes (see e.g., \cite%
{ChaSlud, GiSu, Major, Sur}) are given by objects of
the type $X_{\beta }=\left\{ X_{\beta }\left( \psi \right) :\psi \in
\mathfrak{H}_{\sl{E},\beta }\right\} $, where $\beta $ is a real
non-atomic symmetric measure on $\left( -\pi ,\pi \right] $ (that is, $\beta
\left( dx\right) =\beta \left( -dx\right) $), and
\begin{equation}
\mathfrak{H}_{\sl{E},\beta }=L_{\sl{E}}^{2}\left( \left( -\pi ,\pi %
\right] ,d\beta \right)  \label{evenN1}
\end{equation}%
stands for the collection of \textsl{real} linear combinations of
complex-valued \textsl{even} functions that are square-integrable with
respect to $\beta $ (recall that a function $\psi $ is even if $\overline{%
\psi \left( x\right) }=\psi \left( -x\right) $). The class $\mathfrak{H}_{%
\sl{E},\beta }$ is indeed a real Hilbert space, endowed with the inner
product
\begin{equation}
\left( \psi _{1},\psi _{2}\right) _{\beta }=\int_{-\pi }^{\pi }\psi
_{1}\left( x\right) \psi _{2}\left( -x\right) \beta \left( dx\right) \in
\mathbb{R}.  \label{evenN2}
\end{equation}
This type of construction is used in the spectral theory of time series.

\subsection{Hermite polynomials and Wiener chaos}

We shall now show how to extend the notion of \textsl{Wiener chaos} (as
defined in Section \ref{SS : CH R}) to the case of an isonormal Gaussian
process. The reader is referred to \cite[Ch. 1]{Nualart} for a complete
discussion of this subject. We need some further (standard) definitions.

\bigskip

\begin{definition}
\label{D : HP}The sequence of \textbf{Hermite polynomials }$\left\{
H_{q}:q\geq 0\right\} $ on $\mathbb{R}$, is defined via the following
relations: $H_{0}\equiv 1$ and, for $q\geq 1$,
\begin{equation*}
H_{q}\left( x\right) =\left( -1\right) ^{q}e^{\frac{x^{2}}{2}}\frac{d^{q}}{%
dx^{q}}e^{-\frac{x^{2}}{2}}\text{, \ \ }x\in \mathbb{R}\text{.}
\end{equation*}%
For instance, $H_{1}\left( x\right) =1$, $H_{2}\left( x\right) =x^{2}-1$ and
$H_{3}\left( x\right) =x^{3}-3x$.
\end{definition}

\bigskip

Recall that the sequence $\{\left( q!\right) ^{-1/2}H_{q}:q\geq 0\}$ is an
orthonormal basis of $L^{2}(\mathbb{R},\left( 2\pi \right) ^{-1/2}$ $%
e^{-x^{2}/2}dx).$

\bigskip

\begin{definition}
\label{D : multiindeces} From now on, the symbol $\mathcal{A}_{\infty }$
will denote the class of those sequences $\alpha =\left\{ \alpha _{i}:i\geq
1\right\} $ such that: (i) each $\alpha _{i}$ is a nonnegative integer, (ii)
$\alpha _{i}$ is different from zero only for a finite number of indices $i$%
. A sequence of this type is called a \textbf{multiindex.}\ For $\alpha \in
\mathcal{A}_{\infty }$, we use the notation $|\alpha |$ $=\sum_{i}\alpha
_{i} $. For $q\geq 1$, we also write%
\begin{equation*}
\mathcal{A}_{\infty ,q}=\left\{ \alpha \in \mathcal{A}_{\infty }:|\alpha |%
\text{ }=q\right\} .
\end{equation*}
\end{definition}

\bigskip

\textbf{Remark on notation. }Fix $q\geq 2$. Given a real separable Hilbert
space $\mathfrak{H}$, we denote by $\mathfrak{H}^{\otimes q}$ and $\mathfrak{%
H}^{\odot q}$, respectively, the $q$th \textsl{tensor power }of $\mathfrak{H}
$ and the $q$th \textsl{symmetric tensor power }of $\mathfrak{H}$ (see e.g.
\cite{Janson}). We conventionally set $\mathfrak{H}^{\otimes 1}=\mathfrak{H}%
^{\odot 1}=\mathfrak{H}$.

\bigskip

We recall four classic facts concerning tensors powers of Hilbert spaces
(see e.g. \cite{Janson}).
\begin{itemize}
\item[\textbf{(I)}] The spaces\textbf{\ }$\mathfrak{H}%
^{\otimes q}$ and $\mathfrak{H}^{\odot q}$ are real separable Hilbert
spaces, such that $\mathfrak{H}^{\odot q}\subset \mathfrak{H}^{\otimes q}$.

\item[(\textbf{II})] Let $\left\{ e_{j}:j\geq 1\right\} $ be an orthonormal basis of
$\mathfrak{H}$; then, an orthonormal basis of $\mathfrak{H}^{\otimes q}$ is
given by the collection of all tensors of the type%
\begin{equation*}
e_{j_{1}}\otimes \cdot \cdot \cdot \otimes e_{j_{q}}\text{, \ \ }%
j_{1},...,j_{d}\geq 1.
\end{equation*}%

\item[(\textbf{III})] Let $\left\{ e_{j}:j\geq 1\right\} $ be an orthonormal basis
of $\mathfrak{H}$ and endow $\mathfrak{H}^{\odot q}$ with the inner product $%
(\cdot ,\cdot )_{\mathfrak{H}^{\otimes q}}$; then, an orthogonal (and, in
general, \textsl{not} orthonormal) basis of $\mathfrak{H}^{\odot q}$ is
given by all elements of the type%
\begin{equation}
\mathbf{e}\left( j_{1},...,j_{q}\right) =\mathbf{sym}\left\{
e_{j_{1}}\otimes \cdot \cdot \cdot \otimes e_{j_{q}}\right\} \text{, \ \ }%
1\leq j_{1}\leq ...\leq j_{q}<\infty ,  \label{symbasis}
\end{equation}%
where $\mathbf{sym}\left\{ \cdot \right\} $ stands for a canonical
symmetrization. \textbf{Exercise}: find an orthonormal basis of $\mathfrak{H}%
^{\odot q}$.

\item[(\textbf{IV})] If $\mathfrak{H}=L^{2}\left( Z,\mathcal{Z},\nu
\right) $, where $\nu $ is $\sigma $-finite and non-atomic, then $\mathfrak{H%
}^{\otimes q}\ $can be identified with $L^{2}\left( Z^{q},\mathcal{Z}%
^{q},\nu ^{q}\right) $ and $\mathfrak{H}^{\odot q}$ can be identified with $%
L_{s}^{2}\left( Z^{q},\mathcal{Z}^{q},\nu ^{q}\right) $, where $%
L_{s}^{2}\left( Z^{q},\mathcal{Z}^{q},\nu ^{q}\right) $ is the subspace of $%
L^{2}\left( Z^{q},\mathcal{Z}^{q},\nu ^{q}\right) $ composed of symmetric
functions.
\end{itemize}
\bigskip

Now observe that, once an orthonormal basis of $\mathfrak{H}$ is fixed and
due to the symmetrization, each element $\mathbf{e}\left( j_{1},...,j_{q}\right) $
in (\ref{symbasis}) can be completely described in terms of a unique
multiindex $\alpha \in \mathcal{A}_{\infty ,q}$, as follows: (i) set $\alpha
_{i}=0$ if $i\neq j_{r}$ for every $r=1,...,q$, (ii) set $\alpha _{j}=k$ for
every $j\in \left\{ j_{1},...,j_{q}\right\} $ such that $j$ is repeated
exactly $k$ times in the vector $\left( j_{1},...,j_{q}\right) $ ($k\geq 1$).

\bigskip

\textbf{Examples. }(i) The multiindex $\left( 1,0,0,....\right) $ is
associated with the element of $\mathfrak{H}$ given by $e_{1}$.

(ii)\textbf{\ }Consider the element $\mathbf{e}\left( 1,7,7\right) $. In $(1,7,7)$ the number 1
is not repeated and 7 is repeated twice, hence $\mathbf{e}\left( 1,7,7\right) $ is associated
with the multiindex $\alpha \in \mathcal{A}_{\infty ,3}$ such that $\alpha
_{1}=1$, $\alpha _{7}=2$ and $\alpha _{j}=0$ for every $j\neq 1,7$, that is, $\alpha = (1,0,0,0,0,0,2,0,0,...)$.

(iii) The multindex $\alpha =\left( 1,2,2,0,5,0,0,0,...\right) $ is
associated with the element of $\mathfrak{H}^{\odot 10}$ given by $\mathbf{e}%
\left( 1,2,2,3,3,5,5,5,5,5\right) $.

\bigskip

In what follows, given $\alpha \in \mathcal{A}_{\infty ,q}$ ($q\geq 1$), we
shall write $\mathbf{e}\left( \alpha \right) $ in order to indicate the
element of $\mathfrak{H}^{\odot q}$ uniquely associated with $\alpha $.

\bigskip

\begin{definition}
For every $h\in \mathfrak{H}$, we set $I_{1}^{X}\left( h\right) =X\left(
h\right) $. Now fix an orthonormal basis $\left\{ e_{j}:j\geq 1\right\} $ of
$\mathfrak{H}$: for every $q\geq 2$ and every $h\in \mathfrak{H}^{\odot q}$
such that
\begin{equation*}
h=\sum_{\alpha \in \mathcal{A}_{\infty ,q}}c_{\alpha }e\left( \alpha \right)
\text{ }
\end{equation*}%
(with convergence in $\mathfrak{H}^{\odot q}$, endowed with the inner
product $\left( \cdot ,\cdot \right) _{\mathfrak{H}^{\otimes q}}$), we set
\begin{equation}
I_{q}^{X}\left( h\right) =\sum_{\alpha \in \mathcal{A}_{\infty ,q}}c_{\alpha
}\prod_{j}H_{\alpha _{j}}\left( X\left( e_{j}\right) \right) \text{,}
\label{DefGMWII}
\end{equation}%
where the products only involve the non-zero terms of each multiindex $%
\alpha $, and $H_{m}$ indicates the $m$th Hermite polynomial . For $q\geq 1$%
, the collection of all random variables of the type $I_{q}^{X}\left(
h\right) $, $h\in \mathfrak{H}^{\odot q}$, is called the $q$\textbf{th
Wiener chaos associated with }$X$ and is denoted by $C_{q}\left( X\right) $.
One sets conventionally $C_{0}\left( X\right) =\mathbb{R}$.
\end{definition}

\bigskip

\textbf{Examples.} (i) If $h=e\left( \alpha \right) $, where $\alpha =\left(
1,1,0,0,0,...\right) \in \mathcal{A}_{\infty ,2}$, then
\begin{equation*}
I_{2}^{X}\left( h\right) =H_{1}\left( X\left( e_{1}\right) \right)
H_{1}\left( X\left( e_{2}\right) \right) =X\left( e_{1}\right) X\left(
e_{2}\right) \text{.}
\end{equation*}

(ii) If $\alpha =\left( 1,0,1,2,0,...\right) \in \mathcal{A}_{\infty ,4}$,
then
\begin{eqnarray*}
I_{4}^{X}\left( h\right) &=&H_{1}\left( X\left( e_{1}\right) \right)
H_{1}\left( X\left( e_{3}\right) \right) H_{2}\left( X\left( e_{4}\right)
\right) \\
&=&X\left( e_{1}\right) X\left( e_{3}\right) \left( X\left( e_{4}\right)
^{2}-1\right) \\
&=&X\left( e_{1}\right) X\left( e_{3}\right) X\left( e_{4}\right)
^{2}-X\left( e_{1}\right) X\left( e_{3}\right) \text{.}
\end{eqnarray*}

(iii) If $\alpha =\left( 3,1,1,0,0,...\right) \in \mathcal{A}_{\infty ,5}$,
then
\begin{eqnarray*}
I_{5}^{X}\left( h\right) &=&H_{3}\left( X\left( e_{1}\right) \right)
H_{1}\left( X\left( e_{2}\right) \right) H_{1}\left( X\left( e_{3}\right)
\right) \\
&=&\left( X\left( e_{1}\right) ^{3}-3X\left( e_{1}\right) \right) X\left(
e_{2}\right) X\left( e_{3}\right) \\
&=&X\left( e_{1}\right) ^{3}X\left( e_{2}\right) X\left( e_{3}\right)
-3X\left( e_{1}\right) X\left( e_{2}\right) X\left( e_{3}\right) \text{.}
\end{eqnarray*}

\bigskip

The following result collects some well-known facts concerning Wiener chaos
and isonormal Gaussian processes. In particular: the first point
characterizes the operators $I_{q}^{X}$ as isomorphisms; the second point is
an equivalent of the chaotic representation property for Gaussian measures,
as stated in formula (\ref{CHAOS!}); the third point establishes a formal
relation between random variables of the type $I_{q}^{X}\left( h\right) $
and the multiple Wiener-It\^{o} integrals introduced in Section \ref{SS :
WISI} (see \cite[Ch. 1]{Nualart} for proofs and further discussions of all
these facts).

\begin{proposition}
\label{P : Iso---GM}
\begin{enumerate}
\item For every $q\geq 1$, the $q$th Wiener chaos $C_{q}\left( X\right) $ is
a Hilbert subspace of $L^{2}\left( \mathbb{P}\right) $, and the application%
\begin{equation*}
h\mapsto I_{q}^{X}\left( h\right) \text{, \ \ }h\in \mathfrak{H}^{\odot
q}\text{,}
\end{equation*}%
defines a Hilbert space isomorphism between $\mathfrak{H}^{\odot q}$,
endowed with the inner product $q!(\cdot ,\cdot )_{\mathfrak{H}%
^{\otimes q}}$, and $C_{q}\left( X\right) $.

\item For every $q,q^{\prime }\geq 0$ such that $q\neq q^{\prime }$, the
spaces $C_{q}\left( X\right) $ and $C_{q^{\prime }}\left( X\right) $ are
orthogonal in $L^{2}\left( \mathbb{P}\right) .$

\item Let $F$ be a functional of the isonormal Gaussian
process $X$ satisfying $\mathbb{E}[F(X)^2]<\infty$: then, there exists a unique sequence $\left\{ f_{q}:q\geq
1\right\} $ such that $f_{q}\in \mathfrak{H}^{\odot q}$, and
\begin{equation*}
F=\mathbb{E}\left( F\right) +\sum_{q=1}^{\infty }I_{q}^{X}\left(
f_{q}\right) \text{,}
\end{equation*}%
where the series converges in $L^{2}\left( \mathbb{P}\right) $.

\item Suppose that $\mathfrak{H}=L^{2}\left( Z,\mathcal{Z},\nu \right) $,
where $\nu $ is $\sigma $-finite and non-atomic. Then, for $q\geq 2$, the
symmetric power $\mathfrak{H}^{\odot q}$ can be identified with $%
L_{s}^{2}\left( Z^{q},\mathcal{Z}^{q},\nu ^{q}\right) $ and, for every $f\in
\mathfrak{H}^{\odot q}$, the random variable $I_{q}^{X}\left( f\right) $
coincides with the Wiener-It\^{o} integral (see Definition \ref{D :
MWIIdefinition}) of $f$ with respect to the Gaussian measure given by $%
A\rightarrow X\left( \mathbf{1}_{A}\right) $, $A\in \mathcal{Z}_{\nu }$.
\end{enumerate}
\end{proposition}

\bigskip

\textbf{Remark. }The combination of Point 1. anf Point 2. in the statement
of Proposition \ref{P : Iso---GM} implies that, for every $q,q^{\prime }\geq
1$,
\begin{equation*}
\mathbb{E}\left[ I_{q}^{X}\left( f\right) I_{q^{\prime }}^{X}\left(
f^{\prime }\right) \right] =\mathbf{1}_{q=q^{\prime }}q!\left( f,f^{\prime
}\right) _{\mathfrak{H}^{\otimes q}}
\end{equation*}%
(compare with (\ref{goo})).

\subsection{Contractions, products and some explicit formulae}

We start by introducing the notion of \textsl{contraction }in the context of
powers of Hilbert spaces.

\begin{definition}
Consider a real separable Hilbert space $\mathfrak{H}$, and let $%
\{e_{i}:i\geq 1\}$ be an orthonormal basis of $\mathfrak{H}$. For every $%
n,m\geq 1$, every $r=0,...,n\wedge m$ and every $f\in \mathfrak{H}^{\odot n}$
and $g\in \mathfrak{H}^{\odot m}$, we define the contraction of order $r$,
of $f$ and $g$, as the element of $\mathfrak{H}^{\otimes n+m-2r}$ given by%
\begin{equation}
f\otimes _{r}g=\sum_{i_{1},\ldots ,i_{r}=1}^{\infty }\ \left(
f,e_{i_{1}}\otimes \cdots \otimes e_{i_{r}}\right) _{\mathfrak{H}^{\otimes
r}}\otimes \left( g,e_{i_{1}}\otimes \cdots \otimes e_{i_{r}}\right) _{%
\mathfrak{H}^{\otimes r}},  \label{contrHilbert}
\end{equation}%
and we denote by $\widetilde{f\otimes _{r}g}$ its symmetrization.
\end{definition}

\bigskip

\textbf{Remark. }One can prove (\textbf{Exercise!}) the following result: if
$\mathfrak{H}=L^{2}\left( Z,\mathcal{Z},\nu \right) $, $f\in \mathfrak{H}%
^{\odot n}=L_{s}^{2}\left( Z^{n},\mathcal{Z}^{n},\nu ^{n}\right) $ and $g\in
\mathfrak{H}^{\odot m}=L_{s}^{2}\left( Z^{m},\mathcal{Z}^{m},\nu ^{m}\right)
$, then the definition of the contraction $f\otimes _{r}g$ given in (\ref%
{contrHilbert}) and the one given in (\ref{croceContr}) coincide.

\bigskip

The following result extends the product formula (\ref{ProdGauss}) to the
case of isonormal Gaussian processes. The proof (which is left to the
reader) can be obtained from Proposition \ref{P : 2bleprodGauss}, by using
the fact that every real separable Hilbert space is isomorphic to a space of
the type $L^{2}\left( Z,\mathcal{Z},\nu \right) $, where $\nu $ is $\sigma $%
-finite and non-atomic.\ An alternative proof (by induction) can be found in
\cite[Ch. 1]{Nualart}.

\begin{proposition}
Let $X$ be an isonormal Gaussian process over some real separable Hilbert
space $\mathfrak{H}$. Then, for every $n,m\geq 1$, $f\in \mathfrak{H}^{\odot
n}$ and $g\in \mathfrak{H}^{\odot m}$,
\begin{equation}
I_{n}^{X}\left( f\right) I_{m}^{X}\left( g\right) =\sum_{r=0}^{m\wedge n}r!%
\dbinom{m}{r}\dbinom{n}{r}I_{n+m-2r}^{X}\left( \widetilde{f\otimes _{r}g}%
\right) \text{,}  \label{ProdISO}
\end{equation}%
where the symbol ( $\widetilde{}$ ) indicates a symmetrization, the
contraction $f\otimes _{r}g$ is defined in (\ref{contrHilbert}), and for $%
m=n=r$, we write
\begin{equation*}
I_{0}^{X}\left( \widetilde{f\otimes _{n}g}\right) =\left( f,g\right) _{%
\mathfrak{H}^{\otimes n}}.
\end{equation*}
\end{proposition}

\bigskip

We stress that one can obtain a generalization of the cumulant formulae (\ref%
{GaussDiagrammi}) in the framework of isonormal Gaussian processes. To do
this, one should represent each integral of the type $\int_{Z^{n/2}}f_{%
\sigma }d\nu ^{n/2}$, appearing in (\ref{GaussDiagrammi}), as the inner
product between two iterated contractions of the kernels $\left\{
f_{n_{j}}\right\} $, and then use the canonical isomorphism between $%
\mathfrak{H}$ and a space of the form $L^{2}\left( Z,\mathcal{Z},\nu \right)
$. However, the formalism associated with this extension is rather heavy
(and not really useful for the discussion to follow), and is left to the
reader. Here, we will only state the following formula (proved in \cite%
{PNu05}) giving an explicit expression for the fourth cumulant of a random
variable of the type $I_{d}^{X}\left( f\right) $, $f\in \mathfrak{H}^{\odot
d}$, $d\geq 2$:%
\begin{eqnarray}
\chi _{4}\left( I_{d}^{X}\left( f\right) \right) &=&\mathbb{E}\left[
I_{d}^{X}\left( f\right) ^{4}\right] -3\left( d!\right) ^{2}\left\Vert
f\right\Vert _{\mathfrak{H}^{\otimes d}}^{4}  \label{4mecum} \\
&=&\sum_{p=1}^{d-1}\frac{\left( d!\right) ^{4}}{\left( p!\left( d-p\right)
!\right) ^{2}}\left\{ \left\Vert f\otimes _{p}f\right\Vert _{\mathfrak{H}%
^{\otimes 2\left( d-p\right) }}^{2}+\binom{2d-2p}{d-p}\left\Vert \widetilde{%
f\otimes _{p}f}\right\Vert _{\mathfrak{H}^{\otimes 2\left( d-p\right)
}}^{2}\right\} .  \label{4emecum2}
\end{eqnarray}

As pointed out in \cite[Corollary 2]{PNu05}, formula (\ref{4emecum2}) can be
used in order to prove that, for every isonormal Gaussian process $X$, every
$d\geq 2$ and every $f\in \mathfrak{H}^{\odot d}$, the random variable $%
I_{d}^{X}\left( f\right) $ \textit{cannot be Gaussian} (see also \cite[Ch. 6]%
{Janson}).

\bigskip

\textbf{Example. }We focus once again on isonormal Gaussian processes of the
type $X_{\beta }=\{X_{\beta }\left( \psi \right) :\psi \in \mathfrak{H}_{%
\sl{E},\beta }\}$, where the Hilbert space $\mathfrak{H}_{\sl{E}%
,\beta }$ is given in (\ref{evenN1}). In this case, for $d\geq 2$, the
symmetric power $\mathfrak{H}_{\sl{E},\beta }^{\odot d}$ can be
identified with the real Hilbert space of those functions $\psi _{d}$ that
are symmetric on $\left( -\pi ,\pi \right] ^{d}$, square integrable with
respect to $\beta ^{d}$, and such that $\overline{\psi _{d}\left(
x_{1},...,x_{d}\right) }=\psi _{d}\left( -x_{1},...,-x_{d}\right) .$ For
every $n_{1},...,n_{k}\geq 1$, one can write explicitly a diagram formula as
follows:%
\begin{equation*}
\mathbb{\chi }\left( I_{n_{1}}^{X_{\beta }}\left( \psi _{1}\right) ,\cdot
\cdot \cdot ,I_{n_{k}}^{X_{\beta }}\left( \psi _{k}\right) \right)
=\sum_{\sigma \in \mathcal{M}_{2}\left( \left[ n\right] ,\pi ^{\ast }\right)
}\int_{Z^{n/2}}\psi _{\sigma }d\beta ^{n/2}\text{,}
\end{equation*}%
where $\mathcal{M}_{2}\left( \left[ n\right] ,\pi ^{\ast }\right) $ is
defined in (\ref{M2}) and $\psi _{\sigma }$ is the function in $\left(
n_{1}+\cdot \cdot \cdot +n_{k}\right) /2$ variables obtained by setting $%
x_{i}=-x_{j}$ in $\psi _{1}\otimes _{0}\cdot \cdot \cdot \otimes _{0}\psi
_{d}$ if and only if $i\sim _{\sigma }j$. The field $X_{\beta }$ is often
defined in terms of a complex Gaussian measure (see \cite{ChaSlud, GiSu, Major}).

\section{Simplified CLTs, contractions and circular diagrams\label{S :
SImpliCLT}}

\setcounter{equation}{0}

In a recent series of papers
(see \cite{MaPeAb, NouNu, NouPeccPTRF, NouPecexact, NouPecReveillac, NuOrtiz, PNu05, Pecp, PeTaqMwi, PTu04} for the Gaussian
case, and \cite{PecSoleTaqUtz, PeTaqPoc, PeTaqMultAOP, PeTaq2bleP} for the Poisson case) a set of new results has been
established, allowing to obtain neat Central Limit Theorems (CLTs) for
sequences of random variables belonging to a fixed Wiener chaos of some
Gaussian or Poisson field. The techniques adopted in the above references
are quite varied, as they involve stochastic calculus (\cite{PNu05, Pecp, PTu04}),
Malliavin calculus (\cite{NouNu, NuOrtiz, PeTaqMwi}), Stein's method (\cite{NouPeccPTRF, NouPecexact, NouPecReveillac, PecSoleTaqUtz}) and decoupling (\cite%
{PeTaqMultAOP, PeTaqPoc, PeTaq2bleP}). However, all these
contributions may be described as \textquotedblleft drastic
simplifications\textquotedblright\ of the method of \textsl{moments and
cumulants} (see e.g. \cite{ChaSlud, Major}, as well as the discussion
below) which is a common tool for proving weak convergence results for non
linear functionals of random fields.

The aim of this section is to draw the connection between the above quoted
CLTs and the method of moments and cumulants into further light, by providing
a detailed discussion of the combinatorial implications of the former. This
discussion will involve the algebraic formalism introduced in Section \ref{S
: Lattice}--\ref{S : DG}, as well as the diagram formulae proved in Section %
\ref{S : DF}.

\subsection{A general problem}

In what follows, we will be interested in several variations of the
following problem.

\bigskip

\textbf{Problem A\ }\label{Probl : I}. \emph{Let }$\varphi $ \emph{be a
completely random Gaussian or Poisson measure over some space }$\left( Z,%
\mathcal{Z},\nu \right) $\emph{, where }$\nu $\emph{\ is }$\sigma $\emph{%
-finite and non-atomic. For }$m\geq 1$\emph{\ and }$d_{1},...,d_{m}\geq 1$%
\emph{, let }$\{f_{j}^{\left( k\right) }:j=1,...,m$, \ $k\geq 1\}$\emph{\ be
a collection of kernels such that }$f_{j}^{\left( k\right) }\in
L_{s}^{2}\left( Z^{d_{j}},\mathcal{Z}^{d_{j}},\nu ^{d_{j}}\right) $\emph{\
(the vector }$\left( d_{1},...,d_{m}\right) $\emph{\ does not depend on }$k$%
\emph{), and }%
\begin{equation}
\lim_{k\rightarrow \infty }\mathbb{E}\left[ I_{d_{i}}^{\varphi }\left(
f_{i}^{\left( k\right) }\right) I_{d_{j}}^{\varphi }\left( f_{j}^{\left(
k\right) }\right) \right] =C\left( i,j\right) \text{, \ \ }1\leq i,j\leq m%
\text{,}  \label{NormMWI}
\end{equation}%
\emph{where the integrals }$I_{d_{i}}^{\varphi }\left( f_{i}^{\left(
k\right) }\right) $\emph{\ are defined via (\ref{MWIdef}) and }$\mathbf{C}%
=\left\{ C\left( i,j\right) \right\} $\emph{\ is a }$m\times m$\emph{\
positive definite matrix. We denote by }$\mathbf{N}_{m}\left( 0\mathbf{,C}%
\right) $\emph{\ a }$m$\emph{-dimensional centered Gaussian vector with
covariance matrix }$\mathbf{C}$\emph{.\ Find conditions on the sequence }$%
\left( f_{1}^{\left( k\right) },...,f_{m}^{\left( k\right) }\right) $\emph{,
}$k\geq 1$\emph{, in order to have the CLT}%
\begin{equation}
\mathbf{F}_{k}\triangleq \left( I_{d_{1}}^{\varphi }\left( f_{1}^{\left(
k\right) }\right) ,...,I_{d_{m}}^{\varphi }\left( f_{m}^{\left( k\right)
}\right) \right) \overset{law}{\longrightarrow }\mathbf{N}_{m}\left( 0,%
\mathbf{C}\right) \text{, \ \ }k\rightarrow \infty \text{.}  \label{TLCmwi}
\end{equation}

\bigskip

We observe that, if $d_{i}\neq d_{j}$ in (\ref{NormMWI}), then necessarily $%
C\left( i,j\right) =0$ by Point 2 in Proposition \ref{P : Iso---GM}. The relevance of Problem A comes from the chaotic
representation (\ref{CHAOS!}), implying that a result such as (\ref{TLCmwi})
may be a key tool in order to establish CLTs for more general functionals of
the random measure $\varphi $.\ We recall (see e.g. \cite[Ch.\ 6]{Janson})
that, if $\varphi $ is Gaussian and $d\geq 2$, then a random variable of the
type $I_{d}^{\varphi }\left( f\right) $ cannot be Gaussian.

\bigskip

Plainly, when $\varphi $ is Gaussian, a solution of Problem A can be
immediately deduced from the results discussed in Section \ref{SS :
GaussCase}. Indeed, if the normalization condition (\ref{NormMWI}) is
satisfied, then the moments of the sequence $\left\{ \mathbf{F}_{k}\right\} $
are uniformly bounded (to see this, one can use (\ref{GaussChaosCOntr})),
and the CLT (\ref{TLCmwi}) takes place if and only if every cumulant of
order $\geq 3$ associated with $\mathbf{F}_{k}$ converges to zero when $%
k\rightarrow \infty $. Moreover, an explicit expression for the cumulants
can be deduced from (\ref{GaussDiagrammi}). This method of proving the CLT (%
\ref{TLCmwi}) (which is known as the \textsl{method of cumulants}) has been
used e.g.\ in the references \cite{BrMa, ChaSlud, GiSu, Marinucci, Maruyama82, Maruy}, where the authors proved
CLTs for non-linear functionals of Gaussian fields with a non trivial
covariance structure (for instance, sequences with long memory or isotropic
spherical fields). However, such an approach (e.g. in the study of
fractional Gaussian processes) may be technically quite demanding, since it
involves an infinity of asymptotic relations (one for every cumulant of
order $\geq 3$). If one uses the diagram formulae (\ref{GaussDiagrammi}),
the method of cumulants requires that one explicitly computes and controls
an infinity of expressions of the type $\int_{Z^{n/2}}f_{\sigma ,k}d\nu
^{n/2}$, where the partition $\sigma $ is associated with a non-flat,
Gaussian and connected diagram (see Section \ref{SS : diagrams}).

\bigskip

\textbf{Remarks. }(i)\textbf{\ }We recall that (except for trivial cases),
when $\varphi $ is Gaussian the explicit expression of the characteristic function
of a random variable of the type $I_{d}^{\varphi }\left( f\right) $, $d\geq
3 $, is unknown (for $d=2$ see e.g. \cite[p. 185]{PNu05}).

(ii) Thanks to the results discussed in Section \ref{SS : PossCase} (in
particular, formula (\ref{poissDiagrammi})), the method of cumulants and
diagrams can be also used when $\varphi $ is a completely random Poisson
measure. Clearly, since (\ref{poissDiagrammi}) also involves non-Gaussian
diagrams, the use of this approach in the Poisson case is even more
technically demanding.

\bigskip

In the forthcoming sections we will show how one can successfully bypass the
method of moments and cumulants when dealing with CLTs on a fixed Wiener
chaos.

\subsection{\label{S : TCLmwi}One-dimensional CLTs in the Gaussian case}

We now consider an isonormal Gaussian process $X=\left\{ X\left( h\right)
:h\in \mathfrak{H}\right\} $ over some real separable Hilbert space $%
\mathfrak{H}$. Recall (see Section \ref{S : IsonormalGP}) that the notion of
isonormal Gaussian process is more general than the one of Gaussian measure.
The following result involves one-dimensional sequences of multiple
stochastic integrals, and collects the main findings of \cite{PNu05} and
\cite{NouPeccPTRF}. We recall that the \textsl{total variation distance},
between the law of two general real-valued random variables $Y$ and $Z$, is
given by%
\begin{equation*}
d_{TV}\left( Y,Z\right) =\sup \left\vert \mathbb{P}\left( Y\in B\right) -%
\mathbb{P}\left( Z\in B\right) \right\vert \text{,}
\end{equation*}%
where the supremum is taken over all Borel sets $B\in \mathcal{B}\left(
\mathbb{R}\right) $. Observe that the topology induced by $d_{TV}$, on the
class of probability measures on $\mathbb{R}$, is strictly stronger than the
topology of weak convergence and thus $\lim_{n\rightarrow\infty} d_{TV} (B_n,B)= 0$
is a stronger result than $B_n \stackrel{law}{\rightarrow}B$ (see e.g. Dudley \cite[Ch. 11]{Dudley} for a
discussion of other relevant properties of $d_{TV}$).

\begin{theorem}[See \protect\cite{PNu05} and \protect\cite{NouPeccPTRF}]
\label{T : Nu&Pe05}Fix an integer $d\geq 2$, and define the operator $%
I_{d}^{X}$ according to (\ref{DefGMWII}). Then,for every sequence $\left\{
f^{\left( k\right) }:k\geq 1\right\} $ such that $f^{\left( k\right) }\in
\mathfrak{H}^{\odot d}$ for every $k$, and
\begin{equation}
\lim_{k\rightarrow \infty }d!\left\Vert f^{\left( k\right) }\right\Vert _{%
\mathfrak{H}^{\otimes d}}^{2}=\lim_{k\rightarrow \infty }\mathbb{E}\left[
I_{d}^{X}\left( f^{\left( k\right) }\right) ^{2}\right] =1\text{,}
\label{convtoone1}
\end{equation}%
the following three conditions are equivalent

\begin{enumerate}
\item $\lim_{k\rightarrow \infty }\chi _{4}\left( I_{d}^{X}\left( f^{\left(
k\right) }\right) \right) =0$ $;$

\item for every $r=1,...,d-1$,
\begin{equation}
\lim_{k\rightarrow \infty }\left\Vert f^{\left( k\right) }\otimes
_{r}f^{\left( k\right) }\right\Vert _{\mathfrak{H}^{\otimes 2(d-r)}}^{2}=0%
\text{ },  \label{conCOntrTnp05}
\end{equation}%
where the contraction $f^{\left( k\right) }\otimes _{r}f^{\left( k\right) }$
is defined according to (\ref{contrHilbert});

\item as $k\rightarrow \infty $, the sequence $\left\{ I_{d}^{X}\left(
f^{\left( k\right) }\right) :k\geq 1\right\} $ converges towards a centered
standard Gaussian random variable $Z\sim N\left( 0,1\right) $.
\end{enumerate}

Moreover, the following bound holds for every fixed $k$:%
\begin{eqnarray}
&&d_{TV}\left( I_{d}^{X}\left( f^{\left( k\right) }\right) ,Z\right)
^{2}\leq \left( 1-d!\left\Vert f^{\left( k\right) }\right\Vert _{\mathfrak{H}%
^{\otimes d}}^{2}\right) ^{2}  \label{BTV} \\
&&+d^{2}\sum_{r=1}^{d-1}\left( 2q-2r\right) !\left( r-1\right) !^{2}\binom{%
q-1}{r-1}^{2}\left\Vert f^{\left( k\right) }\otimes _{r}f^{\left( k\right)
}\right\Vert _{\mathfrak{H}^{\otimes 2(d-r)}}^{2}  \notag
\end{eqnarray}
\end{theorem}

\bigskip

Observe that condition (1.) in the previous statement holds if and only if%
\begin{equation*}
\lim_{k\rightarrow \infty }\mathbb{E}\left[ I_{d}^{X}\left( f^{\left(
k\right) }\right) ^{4}\right] =3.
\end{equation*}%
The equivalence of (1.), (2.) and (3.) has been first proved in \cite{PNu05}
by means of stochastic calculus techniques. The paper \cite{NuOrtiz}
contains an alternate proof with additional necessary and sufficient
conditions, as well as several crucial connections with Malliavin calculus
operators (see e.g. \cite{Nualart}). The upper bound (\ref{BTV}) is proved
in \cite{NouPeccPTRF}, by means of Malliavin calculus and the so-called
\textsl{Stein's method }for normal approximation (see e.g. \cite{STeinbook}).

\bigskip

\textbf{Remark. }Theorem \ref{T : Nu&Pe05}, as well as its multidimensional
generalizations (see Section \ref{S : TLCcnj} below), has been applied to a
variery of frameworks, such as: quadratic functionals of bivariate Gaussian
processes (see \cite{DPY}), quadratic functionals of fractional processes
(see \cite{PNu05}), high-frequency limit theorems on homogeneous spaces (see
\cite{MaPeAb, MaPeSphere}), self-intersection local times of
fractional Brownian motion (see \cite{HuNu, NuOrtiz}), Berry-Ess\'{e}%
en bounds in CLTs for Gaussian subordinated sequences (see \cite{NouPeccPTRF, NouPecexact, NouPecReveillac}),
needleets analysis on the sphere (see \cite{BKMP}), power variations of iterated processes (see \cite%
{NouPecIBM}), weighted variations of fractional processes (see \cite{NNTud, NouREv}) and of related random functions (see \cite{Barndorff..., CorNuWoe}).

\subsection{Combinatorial implications of Theorem \protect\ref{T : Nu&Pe05}}

The implication (1.) $\Longrightarrow $ (3.) in Theorem \ref{T : Nu&Pe05}
provides the announced \textquotedblleft drastic
simplification\textquotedblright\ of the methods of moments and cumulants.
However, as demonstrated by the applications of Theorem \ref{T : Nu&Pe05}
listed above, condition (\ref{conCOntrTnp05}) is often much easier to
verify. Indeed, it turns out that the implication (2.) $\Longrightarrow $
(3.) has an interesting combinatorial interpretation.

\bigskip

To see this, we shall fix $d\geq 2$ and suppose that $\mathfrak{H}%
=L^{2}\left( Z,\mathcal{Z},\nu \right) $, with $\nu $ a $\sigma $-finite and
non-atomic measure. According to Proposition \ref{P : Iso---GM}, in this
case the random variable $I_{d}^{X}\left( f\right) $, where $f\in
L_{s}^{2}\left( Z^{d},\mathcal{Z}^{d},\nu ^{d}\right) =L_{s}^{2}\left( \nu
^{d}\right) $, is the multiple Wiener-It\^{o} integral of $f$ with respect
to the Gaussian measure $A\rightarrow X\left( \mathbf{1}_{A}\right) $, as
defined in Definition \ref{D : MWIIdefinition}. We shall also use some
notation from Sections \ref{S : Lattice}--\ref{S : DF}, in particular:

\begin{itemize}
\item For every $n\geq 2$, the symbol $\pi ^{\ast }\left( \left[ nd\right]
\right) \in \mathcal{P}\left( \left[ nd\right] \right) $ stands for the
partition of $\left[ nd\right] =\left\{ 1,2,....,nd\right\} $ obtained by
taking $n$ consecutive blocks of size $d$, that is:%
\begin{equation*}
\pi ^{\ast }\left( \left[ nd\right] \right) =\left\{ \left\{ 1,...,d\right\}
,\left\{ d+1,....,2d\right\} ,...,\left\{ \left( n-1\right)
d+1,...,nd\right\} \right\} .
\end{equation*}

\item The class of partitions $\mathcal{M}_{2}\left( \left[ nd\right] ,\pi
^{\ast }\left( \left[ nd\right] \right) \right) $ is defined according to
formula (\ref{M2}). Recall that, according to (\ref{Lux3}), a partition $%
\sigma \in \mathcal{P}\left( \left[ nd\right] \right) $ is an element of $%
\mathcal{M}_{2}\left( \left[ nd\right] ,\pi ^{\ast }\left( \left[ nd\right]
\right) \right) $ if and only if the diagram $\Gamma \left( \pi ^{\ast
}\left( \left[ nd\right] \right) ,\sigma \right) $ (see Section \ref{SS :
diagrams}) is Gaussian, non-flat and connected, which is equivalent to
saying that the graph $\hat{\Gamma}\left( \pi ^{\ast }\left( \left[ nd\right]
\right) ,\sigma \right) $ (see Section \ref{ss : mgr}) is connected and has
no loops.

\item As in formula (\ref{GaussDiagrammi}), for every $f\in L_{s}^{2}\left(
\nu ^{d}\right) $, every $n$ such that $nd$ is even, and every $\sigma \in
\mathcal{M}_{2}(\left[ nd\right] ,$ $\pi ^{\ast }\left( \left[ nd\right]
\right) )$, we denote by $f_{\sigma ,n}$ the function in $dn/2$ variables,
obtained by identifying two variables $x_{i}$ and $x_{j}$ in the argument of
\begin{equation}
f\underset{n\text{ times}}{\underbrace{\otimes _{0}\cdot \cdot \cdot \otimes
_{0}}}f  \label{aqqa}
\end{equation}%
if and only if $i\sim _{\sigma }j$.
\end{itemize}

\bigskip

We will also denote by
\begin{equation*}
\mathcal{M}_{2}^{c}\left( \left[ nd\right] ,\pi ^{\ast }\left( \left[ nd%
\right] \right) \right)
\end{equation*}%
the subset of $\mathcal{M}_{2}\left( \left[ nd\right] ,\pi ^{\ast }\left( %
\left[ nd\right] \right) \right) $ composed of those partitions $\sigma $
such that the diagram%
\begin{equation*}
\Gamma (\pi ^{\ast }\left( \left[ nd\right] \right) ,\sigma )
\end{equation*}%
is \textsl{circular }(see Section \ref{SS : diagrams}). We also say that a
partition $\sigma \in \mathcal{M}_{2}^{c}\left( \left[ nd\right] ,\pi ^{\ast
}\left( \left[ nd\right] \right) \right) $ has \textsl{rank} $r$ ($%
r=1,...,d-1$) if the diagram $\Gamma \left( \pi ^{\ast }\left( \left[ nd%
\right] \right) ,\sigma \right) $ has exactly $r$ edges linking the first
and the second row.

\bigskip

\textbf{Examples. }(i) The partition whose diagram is given in Fig. 24 (Section \ref{SS : GaussCase}) is an
element of%
\begin{equation*}
\mathcal{M}_{2}^{c}\left( \left[ 8\right] ,\pi ^{\ast }\left( \left[ 9\right]
\right) \right)
\end{equation*}
and has rank $r=1$.

(ii) Consider the case $d=3$ and $n=4$, as well as the partition $\sigma \in
\mathcal{M}_{2}\left( \left[ 12\right] ,\pi ^{\ast }\left( \left[ 12\right]
\right) \right) $ given by
\begin{equation*}
\sigma =\left\{ \left\{ 1,4\right\} ,\left\{ 2,5\right\} ,\left\{
3,12\right\} ,\left\{ 6,9\right\} ,\left\{ 7,10\right\} ,\left\{
8,11\right\} \right\} \text{.}
\end{equation*}%
Then, the diagram $\Gamma \left( \pi ^{\ast }\left( \left[ 12\right] \right)
,\sigma \right) $ is the one in Fig. 28, and therefore $\sigma \in \mathcal{M}_{2}^{c}\left( \left[ 12\right] ,\pi
^{\ast }\left( \left[ 12\right] \right) \right) $ and $\sigma $ has rank $%
r=2 $.

\begin{figure}[htbp]
\begin{center}
\psset{unit=0.7cm}
\begin{pspicture}(0,-1.51)(4.02,1.51)
\psframe[linewidth=0.02,dimen=outer](4.02,1.51)(0.0,-1.51)
\psdots[dotsize=0.15](0.46,0.91)
\psdots[dotsize=0.15](3.66,0.91)
\psdots[dotsize=0.15](2.0,0.91)
\psdots[dotsize=0.15](0.46,0.31)
\psdots[dotsize=0.15](3.66,0.31)
\psdots[dotsize=0.15](2.0,0.31)
\psdots[dotsize=0.15](0.46,-0.27)
\psdots[dotsize=0.15](3.66,-0.27)
\psdots[dotsize=0.15](2.0,-0.27)
\psdots[dotsize=0.15](0.46,-0.87)
\psdots[dotsize=0.15](3.66,-0.87)
\psdots[dotsize=0.15](2.0,-0.87)
\psline[linewidth=0.02cm](0.46,0.93)(0.46,0.33)
\psline[linewidth=0.02cm](2.0,0.93)(2.0,0.35)
\psline[linewidth=0.02cm](0.46,-0.25)(0.46,-0.81)
\psline[linewidth=0.02cm](2.0,-0.25)(2.0,-0.87)
\psline[linewidth=0.02cm](3.64,0.33)(3.62,0.31)
\psline[linewidth=0.02cm](3.64,0.33)(3.64,0.31)
\psline[linewidth=0.02cm](3.66,0.33)(3.66,-0.23)
\rput{-270.0}(4.31,-4.27){\psarc[linewidth=0.02](4.29,0.02){1.11}{36.869896}{142.43141}}
\end{pspicture}
\caption{\sl A circular diagram}
\end{center}
\end{figure}
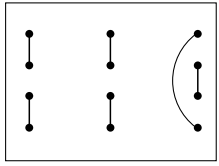

\bigskip

The following technical result links the notions of circular diagram, rank
and contraction. For $d\ge 2$ and $\sigma \in \mathcal{M}_{2}^{c}\left( \left[ 4d\right]
,\pi ^{\ast }\left( \left[ 4d\right] \right) \right) $, let $f_{\sigma ,4}$
be the function in $2d$ variables obtained by identifying $x_{i}$ and $x_{j}$
in the argument of the tensor product (\ref{aqqa}) (with $n=4$) if and only if $i\sim _{\sigma }j$. For instance, if $d=3$ and $\sigma \in \mathcal{M}_{2}^{c}\left( \left[ 12\right] ,\pi^{\ast }\left( \left[ 12\right] \right) \right) $ is associated with the diagram in Fig. 28, then
\begin{equation*}
f_{\sigma,4} (x_1 ,x_2,x_3,x_4,x_5,x_6) = f(x_1,x_2,x_3)f(x_1,x_2,x_4)f(x_5,x_6,x_4) f(x_5,x_6,x_3).
\end{equation*}

\begin{lemma}
\label{L : CP}Fix $f\in L_{s}^{2}\left( \nu ^{d}\right) $, $d\geq 2$, and,
for $r=1,...,d-1$, define the contraction $f\otimes _{r}f$ according to (\ref%
{croceContr}). Then, for every $\sigma \in \mathcal{M}%
_{2}^{c}\left( \left[ 4d\right] ,\pi ^{\ast }\left( \left[ 4d\right] \right)
\right) $ with rank $r\in \left\{ 1,...,d-1\right\} $,
\begin{equation}
\int_{Z^{2d}}f_{\sigma ,4}d\nu ^{2d}=\left\Vert f\otimes _{r}f\right\Vert
_{L^{2}\left( \nu ^{2\left( d-r\right) }\right) }^{2}=\left\Vert f\otimes
_{d-r}f\right\Vert _{L^{2}\left( \nu ^{2r}\right) }^{2}  \label{uff}
\end{equation}
\end{lemma}

\begin{proof}
It is sufficient to observe that $f$ is symmetric by definition, and then to
use the relation%
\begin{eqnarray*}
&&\left\Vert f\otimes _{r}f\right\Vert _{L^{2}\left( \nu ^{2\left(
d-r\right) }\right) }^{2} \\
&=&\int_{Z^{d-r}}\int_{Z^{d-r}}\int_{Z^{r}}\int_{Z^{r}}f\left( \mathbf{a}%
_{d-r},\mathbf{b}_{r}\right) f\left( \mathbf{b}_{r},\mathbf{a}_{d-r}^{\prime
}\right) \times \\
&&\text{ \ \ \ \ \ \ \ \ \ \ \ \ \ }\times f\left( \mathbf{a}_{d-r}^{\prime
},\mathbf{b}_{r}^{\prime }\right) f\left( \mathbf{b}_{r}^{\prime },\mathbf{a}%
_{d-r}\right) \nu ^{d-r}\left( d\mathbf{a}_{d-r}\right) \nu ^{d-r}\left( d%
\mathbf{a}_{d-r}^{\prime }\right) \nu ^{r}\left( d\mathbf{b}_{r}\right) \nu
^{r}\left( d\mathbf{b}_{r}^{\prime }\right) .
\end{eqnarray*}
\end{proof}

\bigskip

\textbf{Remark. }Formula (\ref{uff}) implies that, for a fixed $f$ and for
every $\sigma \in \mathcal{M}_{2}^{c}\left( \left[ 4d\right] ,\pi ^{\ast
}\left( \left[ 4d\right] \right) \right) $, the value of the integral $%
\int_{Z^{2d}}f_{\sigma }d\nu ^{2d}$ \textsl{depends on $\sigma$ uniquely through} $r$ (or $d-r$),
where $r$ is the rank of $\sigma $.

\bigskip

By using Lemma \ref{L : CP}, one obtains immediately the following result,
which provides a combinatorial description of the implication (2.) $%
\Longrightarrow $ (3.) in Theorem \ref{T : Nu&Pe05}.

\begin{proposition}
\label{P : CombCons}For every $d\geq 2$ and every sequence $\left\{
f^{\left( k\right) }:k\geq 1\right\} \subset L_{s}^{2}\left( \nu ^{d}\right)
$ such that $d!\left\Vert f^{\left( k\right) }\right\Vert _{L^{2}\left( \nu
^{d}\right) }^{2}\rightarrow 1$ ($k\rightarrow \infty $), the following
relations are equivalent:

\begin{enumerate}
\item as $k\rightarrow \infty $
\begin{equation}
\sum_{\sigma \in \mathcal{M}_{2}\left( \left[ nd\right] ,\pi ^{\ast }\left( %
\left[ nd\right] \right) \right) }\int_{Z^{nd/2}}f_{\sigma }^{\left(
k\right) }d\nu ^{nd/2}\rightarrow 0\text{, \ \ }\forall n\geq 3;
\label{comb1}
\end{equation}

\item for every partition $\sigma \in \mathcal{M}_{2}^{c}\left( \left[ 4d%
\right] ,\pi ^{\ast }\left( \left[ 4d\right] \right) \right) $, as $%
k\rightarrow \infty $,%
\begin{equation}
\int_{Z^{2d}}f_{\sigma }^{\left( k\right) }d\nu ^{2d}\rightarrow 0\text{.}
\label{Comb2}
\end{equation}
\end{enumerate}
\end{proposition}

\begin{proof}
Thanks to formula \ref{GaussDiagrammi}, one deduces that
\begin{equation*}
\sum_{\sigma \in \mathcal{M}_{2}\left( \left[ nd\right] ,\pi ^{\ast }\left( %
\left[ nd\right] \right) \right) }\int_{Z^{nd/2}}f_{\sigma }^{\left(
k\right) }d\nu ^{nd/2}=\chi _{n}\left( I_{d}^{X}\left( f^{\left( k\right)
}\right) \right) \text{,}
\end{equation*}%
where $I_{d}^{X}\left( f^{\left( k\right) }\right) $ is the multiple
Wiener-It\^{o} integral of $f^{\left( k\right) }$ with respect to the
Gaussian measure induced by $X$, and $\chi _{n}$ indicates the $n$th
cumulant. It follows that, since $d!\left\Vert f^{\left( k\right)
}\right\Vert _{L^{2}\left( \nu ^{d}\right) }^{2}=\mathbb{E}\left[
I_{d}^{X}\left( f^{\left( k\right) }\right) \right] \rightarrow 1$, relation
(\ref{comb1}) is equivalent to $I_{d}^{X}\left( f^{\left( k\right) }\right)
\overset{law}{\rightarrow }Z\sim N\left( 0,1\right) $. On the other hand,
one deduces from Lemma \ref{L : CP} that (\ref{Comb2}) takes place if and only if (\ref{conCOntrTnp05}) holds. Since, according to Theorem \ref{T :
Nu&Pe05}, condition (\ref{conCOntrTnp05}) is necessary and sufficient in
order to have $I_{d}^{X}\left( f^{\left( k\right) }\right) \overset{law}{%
\rightarrow }Z$, we immediately obtain the desired conclusion.
\end{proof}

\medskip

\begin{corollary}\label{C : OHIO}
Fix $d\ge 2$ and suppose that the sequence $\left\{
f^{\left( k\right) }:k\geq 1\right\} \subset L_{s}^{2}\left( \nu ^{d}\right)
$ is such that $d!\left\Vert f^{\left( k\right) }\right\Vert _{L^{2}\left( \nu
^{d}\right) }^{2}\rightarrow 1$ ($k\rightarrow \infty $). Then, (\ref{Comb2}) takes places if and only if
$I^X_d(f^{(k)})\stackrel{law}{\rightarrow}Z\sim N(0,1)$.
\end{corollary}
\begin{proof}
As pointed out in the proof of Proposition \ref{P :
CombCons}, since the normalization condition $d!\left\Vert f^{\left(
k\right) }\right\Vert _{L^{2}\left( \nu ^{d}\right) }^{2}\rightarrow 1$ is
in order, relation (\ref{comb1}) is equivalent to the fact that the sequence
$I_{d}^{X}\left( f^{\left( k\right) }\right) $, $k\geq 1$, converges in law
to a standard Gaussian random variables. The implication (2.)
$\Longrightarrow $ (1.) in the statement of Proposition \ref{P : CombCons}
yields the desired result.
\end{proof}

\bigskip

\textbf{Remarks. }(1) Corollary \ref{C : OHIO} implies that,
in order to prove a CLT on a fixed Wiener chaos, \textsl{it is
sufficient to compute and control a finite number of expressions of the type
}$\int_{Z^{2d}}f_{\sigma }^{\left( k\right) }d\nu ^{2d}$\textsl{, where }$%
\sigma $\textsl{\ is associated with a connected Gaussian circular diagram
with four rows. }Moreover, these expressions determine the speed of
convergence in total variation, via the upper bound given in (\ref{BTV}).

(2) Relation (\ref{uff}) also implies that: (i) for $d$ even, (\ref{Comb2})
takes place for every
\begin{equation*}
\sigma \in \mathcal{M}_{2}^{c}\left( \left[ 4d\right] ,\pi _{d}^{\ast
}\left( \left[ 4d\right] \right) \right)
\end{equation*}
if and only if for every $r=1,...,d/2$, there exists a partition $\sigma
\in \mathcal{M}_{2}^{c}\left( \left[ 4d\right] ,\pi _{d}^{\ast }\left( \left[
4d\right] \right) \right) $ with rank $r$ and such that (\ref{Comb2}) holds;
(ii) for $d$ odd, (\ref{Comb2}) takes place for every $\sigma \in \mathcal{M}%
_{2}^{c}\left( \left[ 4d\right] ,\pi _{d}^{\ast }\left( \left[ 4d\right]
\right) \right) $ if and only if for every $r=1,...,(d+1)/2$, there exists
a partition $\sigma \in \mathcal{M}_{2}^{c}\left( \left[ 4d\right] ,\pi
_{d}^{\ast }\left( \left[ 4d\right] \right) \right) $ with rank $r$ and such
that (\ref{Comb2}) holds.

(3) When $d=2$, the implication (\ref{Comb2}) $\Rightarrow $ (\ref{comb1})
is a consequence of the fact that, for every $n\geq 3$ and up to a
permutation of the rows, the diagram associated with any element of $%
\mathcal{M}_{2}\left( \left[ 2n\right] ,\pi _{2}^{\ast }\left( \left[ 2n%
\right] \right) \right) \ $is equivalent to a circular diagram (this fact
has been already pointed out at the end of Section \ref{SS : GaussCase}).
For instance, it is always possible to permute the blocks of $\pi _{2}^{\ast
}\left( \left[ 10\right] \right) $ in such a way that the diagram $\Gamma
(\pi _{2}^{\ast }\left( \left[ 10\right] \right) ,\sigma )$, associated with
some $\sigma \in \mathcal{M}_{2}\left( \left[ 10\right] ,\pi _{2}^{\ast
}\left( \left[ 10\right] \right) \right) $, has the form of the diagram in Fig. 29.
By using this fact, one can prove that (\ref{Comb2}) $\Rightarrow $ (\ref%
{comb1}) by means of the Cauchy-Schwarz inequality and of a recurrence
argument (for another proof of Theorem \ref{T : Nu&Pe05} in the case $d=2$,
by means of an explicit expression of the Fourier transform of $%
I_{2}^{X}\left( f^{\left( k\right) }\right) $, see \cite[p. 185]{PNu05}).

\begin{figure}[htbp]
\begin{center}
\psset{unit=0.7cm}
\begin{pspicture}(0,-2.19)(3.02,2.18)
\psframe[linewidth=0.02,dimen=outer](3.02,2.18)(0.0,-2.0)
\psdots[dotsize=0.15](0.6,1.78)
\psdots[dotsize=0.15](2.4,1.78)
\psdots[dotsize=0.15](0.6,1.0)
\psdots[dotsize=0.15](2.4,1.0)
\psdots[dotsize=0.15](0.6,0.16)
\psdots[dotsize=0.15](2.4,0.16)
\psdots[dotsize=0.15](0.6,-0.6)
\psdots[dotsize=0.15](2.4,-0.6)
\psdots[dotsize=0.15](0.6,-1.38)
\psdots[dotsize=0.15](2.4,-1.38)
\psline[linewidth=0.02cm](0.6,1.78)(0.6,1.02)
\psline[linewidth=0.02cm](2.38,0.98)(2.4,1.0)
\psline[linewidth=0.02cm](2.4,1.0)(2.4,0.22)
\psline[linewidth=0.02cm](0.6,0.14)(0.6,-0.56)
\psline[linewidth=0.02cm](2.4,-0.6)(2.4,-0.62)
\psline[linewidth=0.02cm](2.4,-0.6)(2.4,-1.34)
\psbezier[linewidth=0.02](2.4,1.78)(2.4,0.98)(0.6,-2.18)(0.6,-1.38)
\end{pspicture}
\caption{\sl A circular diagram with five rows}
\end{center}
\end{figure}
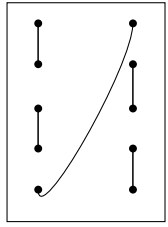

\subsection{\label{S : TLCcnj}A multidimensional CLT}

The paper \cite{PTu04} (but see also \cite{NouPecReveillac, NuOrtizSPA, Pecp}) contains a complete solution of Problem A in
the Gaussian case, for every $m\geq 2$. For such an index $m$, we denote by $%
V_{m}$ the set of all vectors $\left( i_{1},i_{2},i_{3},i_{4}\right) \in
\left( 1,...,m\right) ^{4}$ such that at least one of the following three
properties is verified: (a) $i_{1}\neq i_{2}=i_{3}=i_{4}$, (b) $i_{1}\neq
i_{2}=i_{3}\neq i_{4}$ and $i_{4}\neq i_{1}$, (c) the elements of $\left(
i_{1},...,i_{4}\right) $ are all different. In what follows, $X=\left\{
X\left( h\right) :h\in \mathfrak{H}\right\} $ indicates an isonormal
Gaussian process over a real separable Hilbert space $\mathfrak{H}$.

\begin{theorem}
\label{T : IoCipCVcong}Let $m\geq 2$ and $d_{1},...,d_{m}\geq 1$ be fixed
and let
\begin{equation*}
\left\{ f_{j}^{\left( k\right) }:j=1,...,m\text{, \ }k\geq 1\right\}
\end{equation*}%
be a collection of kernels such that $f_{j}^{\left( k\right) }\in \mathfrak{H%
}^{\odot d_{j}}$ and the normalization condition (\ref{NormMWI}) is
verified. Then, the following conditions are equivalent:

\begin{enumerate}
\item as $k\rightarrow \infty $, the vector $\mathbf{F}_{k}=\left(
I_{d_{1}}^{X}\left( f_{1}^{\left( k\right) }\right) ,...,I_{d_{m}}^{X}\left(
f_{m}^{\left( k\right) }\right) \right) $ converges in law towards a $m$%
-dimensionnel Gaussian vector $\mathbf{N}_{m}\left( 0,\mathbf{C}\right)
=\left( N_{1},...,N_{m}\right) $ with covariance matrix $\mathbf{C}=\left\{
C\left( i,j\right) \right\} $;

\item
\begin{eqnarray*}
&&\lim_{k\rightarrow \infty }\mathbb{E}\left[ \left(
\sum_{i=1,...,m}I_{d_{i}}^{X}\left( f_{i}^{\left( k\right) }\right) \right)
^{4}\right] \\
&=&3\left( \sum_{i=1}^{m}C\left( i,i\right) +2\sum_{1\leq i<j\leq m}C\left(
i,j\right) \right) ^{2}=\mathbb{E}\left[ \left( \sum_{i=1}^{m}N_{i}\right)
^{4}\right] ,
\end{eqnarray*}%
and
\begin{equation*}
\lim_{k\rightarrow \infty }\mathbb{E}\left[
\prod_{l=1}^{4}I_{d_{i_{l}}}^{X}\left( f_{i_{l}}^{\left( k\right) }\right) %
\right] =\mathbb{E}\left[ \prod_{l=1}^{4}N_{i_{l}}\right]
\end{equation*}%
$\forall \left( i_{1},i_{2},i_{3},i_{4}\right) \in V_{m}$;

\item for every $j=1,...,m$, the sequence $I_{d_{j}}^{X}\left( f_{j}^{\left(
k\right) }\right) $ $k\geq 1$, converges in law towards $N_{j}$, that is,
towards a centered Gaussian variable with variance $C\left( j,j\right) ;$

\item $\forall j=1,...,m$, $\lim_{k\rightarrow \infty }\chi _{4}\left(
I_{d_{j}}^{X}\left( f_{j}^{\left( k\right) }\right) \right) =0$;

\item $\forall j=1,...,m$%
\begin{equation}
\lim_{k\rightarrow \infty }\left\Vert f_{j}^{\left( k\right) }\otimes
_{r}f_{j}^{\left( k\right) }\right\Vert _{\mathfrak{H}^{\otimes 2\left(
d_{j}-r\right) }}=0,  \label{contrConj}
\end{equation}%
$\forall r=1,...,d_{j}-1.$
\end{enumerate}
\end{theorem}

\bigskip

The original proof of Theorem \ref{T : IoCipCVcong} uses arguments from
stochastic calculus. See \cite{NouPecReveillac} and \cite{NuOrtizSPA},
respectively, for alternate proofs based on Malliavin calculus and Stein's
method.\ In particular, in \cite{NouPecReveillac} one can find bounds
analogous to (\ref{BTV}), concerning the multidimensional Gaussian
approximation of $\mathbf{F}_{k}$ in the Wasserstein distance. The crucial
element in the statement of Theorem \ref{T : IoCipCVcong} is the implication
(3.) $\Rightarrow $ (1.), which yields the following result.

\begin{corollary}
Let the vectors $\mathbf{F}_{k}$, $k\ge 1$, be as in the statement of Theorem \ref{T : IoCipCVcong},
and suppose that (\ref{NormMWI}) is satisfied. Then, the convergence in law of each component of the vectors $\mathbf{F%
}_{k}$, towards a Gaussian random variable, always implies the joint
convergence of $\mathbf{F}_{k}$ towards a Gaussian vector with covariance $%
\mathbf{C}$.
\end{corollary}

Thanks to Theorem \ref{T : Nu&Pe05}, it follows that a CLT such
as (\ref{TLCmwi}) can be uniquely deduced from (\ref{NormMWI}) and from the
relations (\ref{contrConj}), involving the contractions of the kernels $%
f_{j}^{\left( k\right) }$.

\bigskip

When $\mathfrak{H}=L^{2}\left( Z,\mathcal{Z},\nu \right) $ (with $\nu $ non
atomic), the combinatorial implications of Theorem \ref{T : IoCipCVcong} are
similar to those of Theorem \ref{T : Nu&Pe05}. Indeed, thanks to the
implication (5.) $\Rightarrow $ (1.), one deduces that, for a sequence $%
\left( f_{1}^{\left( k\right) },...,f_{m}^{\left( k\right) }\right) $,%
\textit{\ }$k\geq 1$, as in (\ref{NormMWI}), if
\begin{equation*}
\int_{Z^{2d}}\left( f_{j}^{\left( k\right) }\right) _{\sigma }d\nu
^{2d}\rightarrow 0\text{, \ \ }\forall \sigma \in \mathcal{M}_{2}^{c}\left( %
\left[ 4d\right] ,\pi ^{\ast }\left( \left[ 4d\right] \right) \right) \text{,%
}
\end{equation*}%
then
\begin{equation*}
\sum_{\sigma \in \mathcal{M}_{2}\left( \left[ n\right] ,\pi ^{\ast }\right)
}\int_{Z^{n/2}}f_{\sigma ,\ell }^{\left( k\right) }d\nu ^{n/2}\rightarrow 0%
\text{,}
\end{equation*}%
for every integer $n$ which is the sum of $\ell \geq 3$ components $\left(
d_{i_{1}},d_{i_{2}},...,d_{i_{\ell }}\right) $ of the vector $\left(
d_{1},...,d_{m}\right) $ (with possible repetitions of the indices $%
i_{1},...,i_{\ell }$), with%
\begin{equation*}
\pi ^{\ast }=\left\{ \left\{ 1,...,d_{i_{1}}\right\} ,...,\left\{
d_{1}+...+d_{i_{\ell -1}}+1,...,n\right\} \right\} \in \mathcal{P}\left( %
\left[ n\right] \right) ,
\end{equation*}%
and every function $f_{\sigma ,\ell }^{\left( k\right) }$, in $n/2$
variables, is obtained by identifying two variables $x_{k}$ and $x_{j}$ in
the argument of $f_{i_{1}}\otimes _{0}\cdot \cdot \cdot \otimes
_{0}f_{i_{\ell }}$ if and only if $k\sim _{\sigma }j$.

\bigskip

As already pointed out, the chaotic representation property (\ref{CHAOS!})
allows to use Theorem \ref{T : IoCipCVcong} in order to obtain CLTs for
general functionals of an isonormal Gaussian process $X$. We now present a
result in this direction, obtained in \cite{HuNu}, whose proof can be
deduced from Theorem \ref{T : IoCipCVcong}.

\begin{theorem}[See \cite{HuNu}]
\label{T : H&Nclt}We consider a sequence $\left\{ F_{k}:k\geq 1\right\} $ of
centered and square-integrable functionals of an isonormal Gaussian process $%
X $, admitting the chaotic decomposition%
\begin{equation*}
F_{k}=\sum_{d=1}^{\infty }I_{d}^{X}\left( f_{d}^{\left( k\right) }\right)
\text{, \ \ }k\geq 1\text{.}
\end{equation*}%
Assume that

\begin{itemize}
\item $\lim_{N\rightarrow \infty }\lim \sup_{k\rightarrow \infty
}\sum_{d\geq N+1}d!\left\Vert f_{d}^{\left( k\right) }\right\Vert _{%
\mathfrak{H}^{\otimes d}}^{2}\rightarrow 0,$

\item for every $d\geq 1$, $\lim_{k\rightarrow \infty }d!\left\Vert
f_{d}^{\left( k\right) }\right\Vert _{\mathfrak{H}^{\otimes d}}^{2}=\sigma
_{d}^{2},$

\item $\sum_{d=1}^{\infty }\sigma _{d}^{2}\triangleq \sigma ^{2}<\infty ,$

\item for every $d\geq 1$, $\lim_{k\rightarrow \infty }\left\Vert
f_{d}^{\left( k\right) }\otimes _{r}f_{d}^{\left( k\right) }\right\Vert _{%
\mathfrak{H}^{\otimes 2\left( d-r\right) }}=0$, $\forall r=1,...,d-1.$
\end{itemize}

Then, as $k\rightarrow \infty $, $F_{k}\overset{law}{\rightarrow }N\left(
0,\sigma ^{2}\right) $, where $N\left( 0,\sigma ^{2}\right) $ is a centered
Gaussian random variable with variance $\sigma ^{2}.$
\end{theorem}

\subsection{Simplified CLTs in the Poisson case: the case of double integrals%
}

We conclude this survey by discussing a simplified CLT for sequences of
double integrals with respect to a Poisson random measure. Note that this
result (originally obtained in \cite{PeTaq2bleP}) has been
generalized in \cite{PecSoleTaqUtz}, where one can find CLTs for sequences
of multiple integrals of arbitrary orders -- with explicit Berry-Ess\'{e}en
bounds in the Wasserstein distance obtained once again via Stein's method.

\bigskip

In this section, $\left( Z,\mathcal{Z},\nu \right) $ is a measure space,
with $\nu $ $\sigma $-finite and non-atomic. Also, $\hat{N}=\{\hat{N}\left(
B\right) :B\in \mathcal{Z}_{\nu }\}$ is a compensated Poisson measure with
control measure given by $\nu $. In \cite{PeTaq2bleP}, we have used some
decoupling techniques developed in \cite{PeTaqPoc} in order to prove CLTs
for sequences of random variables of the type:%
\begin{equation}
F_{k}=I_{2}^{\hat{N}}\left( f^{\left( k\right) }\right) \text{, \ \ }k\geq 1%
\text{,}  \label{seq}
\end{equation}%
where $f^{\left( k\right) }\in L_{s}^{2}\left( \nu ^{2}\right) $. In
particular, we focus on sequences $\left\{ F_{k}\right\} $ satisfying the
following assumption

\bigskip

\textbf{Assumption N}\textsc{.}\textbf{\ }The sequence $f^{\left( k\right) }$%
, $k\geq 1$, in (\ref{seq}) verifies :

\begin{description}
\item[\textbf{N.i }(\textit{integrability}) ] $\forall k\geq 1$,
\begin{equation}
\int_{Z}f^{\left( k\right) }\left( z,\cdot \right) ^{2}\nu \left( dz\right)
\in L^{2}\left( \nu \right) \text{ \ \ and \ \ }\left\{ \int_{Z}f^{\left(
k\right) }\left( z,\cdot \right) ^{4}\nu \left( dz\right) \right\} ^{\frac{1%
}{2}}\in L^{1}\left( \nu \right) ;  \label{N-1}
\end{equation}

\item[\textbf{N.ii }(\textit{normalization}) ] As $k\rightarrow \infty $,
\begin{equation}
2\int_{Z}\int_{Z}f^{\left( k\right) }\left( z,z^{\prime }\right) ^{2}\nu
\left( dz\right) \nu \left( dz^{\prime }\right) \rightarrow 1;  \label{N0}
\end{equation}

\item[\textbf{N.iii }(\textit{fourth power}) ] As $k\rightarrow \infty $,
\begin{equation}
\int_{Z}\int_{Z}f^{\left( k\right) }\left( z,z^{\prime }\right) ^{4}\nu
\left( dz\right) \nu \left( dz^{\prime }\right) \rightarrow 0  \label{N1}
\end{equation}%
(in particular, this implies that $f^{\left( k\right) }\in L^{4}\left( \nu
^{2}\right) $).
\end{description}

\bigskip

\textbf{Remarks.\ }(1) The conditions in (\ref{N-1}) are technical : the
first ensures the existence of the stochastic integral of $\int_{Z}f^{\left(
k\right) }\left( z,\cdot \right) ^{2}\nu \left( dz\right) $ with respect to $%
\hat{N}$; the second allows to use some Fubini arguments in the proof of the
results to follow.

(2) Suppose that there exists a set $B$, independent of $n$, such that $\nu
\left( B\right) <\infty $ and $f^{\left( k\right) }=f^{\left( k\right) }%
\mathbf{1}_{B}$, a.e.--$d\nu ^{2}$, $\forall k\geq 1$ (this holds, in
particular, when $\nu $ is finite). Then, by the Cauchy-Schwarz inequality,
if (\ref{N1}) is true, then $\left( f^{\left( k\right) }\right) $ converges
necessarily to zero$.$ Therefore, in order to study more general sequences $\left(
f^{\left( k\right) }\right) $, we must assume that $\nu \left( Z\right)
=+\infty $.

\bigskip

The next theorem is the main result of \cite{PeTaq2bleP}.

\begin{theorem}
\label{T : PecTaq2blePclt}Let $F_{k}=I_{2}^{\hat{N}}(f^{\left( k\right) })$
with $f^{\left( k\right) }\in L_{s}^{2}(\nu ^{2})$, $k\geq 1$, and suppose
that Assumption N is verified. Then, $f^{\left( k\right) }\star
_{1}^{0}f^{\left( k\right) }\in L^{2}(\nu ^{3})$ and $f^{\left( k\right)
}\star _{1}^{1}f^{\left( k\right) }\in L_{s}^{2}(\nu ^{2})$, $\forall k\geq
1 $, and also :

\begin{enumerate}
\item if
\begin{equation}
\left\Vert \text{\ }f^{\left( k\right) }\star _{2}^{1}f^{\left( k\right)
}\right\Vert _{L^{2}(\nu )}\rightarrow 0\text{\ and }\left\Vert f^{\left(
k\right) }\star _{1}^{1}f^{\left( k\right) }\right\Vert _{L^{2}(\nu
^{2})}\rightarrow 0\text{ },  \label{G*}
\end{equation}%
then
\begin{equation}
F_{k}\overset{\text{law}}{\rightarrow }N\left( 0,1\right) \text{,}
\label{GG}
\end{equation}%
where $N\left( 0,1\right) $ is a centered Gaussian random variable with
unitary variance.

\item if $F_{k}\in L^{4}\left( \mathbb{P}\right) $, $\forall k$, a
sufficient condition in order to have (\ref{G*}) is%
\begin{equation}
\chi _{4}\left( F_{n}\right) \rightarrow 0;  \label{GGG}
\end{equation}

\item if the sequence $\left\{ F_{k}^{4}:k\geq 1\right\} $ is uniformly
integrable, then the three conditions (\ref{G*}), (\ref{GG}) and (\ref{GGG})
are equivalent.
\end{enumerate}
\end{theorem}

\bigskip

\textbf{Remark. }See \cite{DBPP} and \cite{PecPru} for several applications
of Theorem \ref{T : PecTaq2blePclt} to Bayesian non-parametric survival
analysis.

\bigskip

We now give a combinatorial interpretation (in terms of diagrams) of the
three asymptotic conditions appearing in formulae (\ref{N1}) and (\ref{G*}).
To do this, consider the set $\left[ 8\right] =\left\{ 1,...,8\right\} $, as
well as the partition $\pi ^{\ast }=\left\{ \left\{ 1,2\right\} ,\left\{
3,4\right\} ,\left\{ 5,6\right\} ,\left\{ 7,8\right\} \right\} \in \mathcal{P%
}\left( \left[ 8\right] \right) $. We define the set of partitions $\mathcal{%
M}_{\geq 2}\left( \left[ 8\right] ,\pi ^{\ast }\right) \subset \mathcal{P}%
\left( \left[ 8\right] \right) $ according to (\ref{Mplus2}). Given an
element $\sigma \in \mathcal{M}_{\geq 2}\left( \left[ 8\right] ,\pi ^{\ast
}\right) $ and given $f\in L_{s}^{2}\left( \nu ^{2}\right) $, the function $%
f_{\sigma ,4}$, in $\left\vert \sigma \right\vert $ variables, is obtained
by identifying the variables $x_{i}$ and $x_{j}$ in the argument of $%
f\otimes _{0}f\otimes _{0}f\otimes _{0}f$ (as defined in (\ref{0multContr}))
if and only if $i\sim _{\sigma }j.$ We define three partitions $\sigma
_{1},\sigma _{2},\sigma _{3}\in \mathcal{M}_{\geq 2}\left( \left[ 8\right]
,\pi ^{\ast }\right) $ as follows:%
\begin{eqnarray*}
\sigma _{1} &=&\left\{ \left\{ 1,3,5,7\right\} ,\left\{ 2,4,6,8\right\}
\right\} \\
\sigma _{2} &=&\left\{ \left\{ 1,3,5,7\right\} ,\left\{ 2,4\right\} ,\left\{
6,8\right\} \right\} \\
\sigma _{3} &=&\left\{ \left\{ 1,3\right\} ,\left\{ 4,6\right\} ,\left\{
5,7\right\} ,\left\{ 2,8\right\} \right\} .
\end{eqnarray*}%
The diagrams $\Gamma \left( \pi ^{\ast },\sigma _{1}\right) $, $\Gamma
\left( \pi ^{\ast },\sigma _{2}\right) $ and $\Gamma \left( \pi ^{\ast
},\sigma _{3}\right) $ are represented (in order) in Fig. 30.

\begin{figure}[htbp]
\begin{center}
\psset{unit=0.7cm}
\begin{pspicture}(0,-1.41)(8.6,1.42)
\psframe[linewidth=0.02,dimen=outer](8.6,1.41)(0.0,-1.41)
\psline[linewidth=0.02cm](2.72,1.41)(2.72,-1.37)
\psline[linewidth=0.02cm](5.6,1.41)(5.6,-1.37)
\psdots[dotsize=0.15](0.62,0.99)
\psdots[dotsize=0.15](0.62,0.39)
\psdots[dotsize=0.15](0.62,-0.23)
\psdots[dotsize=0.15](0.62,-0.81)
\psdots[dotsize=0.15](2.0,0.99)
\psdots[dotsize=0.15](2.0,0.39)
\psdots[dotsize=0.15](2.0,-0.23)
\psdots[dotsize=0.15](2.0,-0.81)
\psdots[dotsize=0.15](3.44,0.99)
\psdots[dotsize=0.15](3.44,0.39)
\psdots[dotsize=0.15](3.44,-0.23)
\psdots[dotsize=0.15](3.44,-0.81)
\psdots[dotsize=0.15](4.82,0.99)
\psdots[dotsize=0.15](4.82,0.39)
\psdots[dotsize=0.15](4.82,-0.23)
\psdots[dotsize=0.15](4.82,-0.81)
\psdots[dotsize=0.15](6.36,0.97)
\psdots[dotsize=0.15](6.36,0.37)
\psdots[dotsize=0.15](6.36,-0.25)
\psdots[dotsize=0.15](6.36,-0.83)
\psdots[dotsize=0.15](7.74,0.97)
\psdots[dotsize=0.15](7.74,0.37)
\psdots[dotsize=0.15](7.74,-0.25)
\psdots[dotsize=0.15](7.74,-0.83)
\psellipse[linewidth=0.02,dimen=outer](0.63,0.06)(0.25,1.19)
\psellipse[linewidth=0.02,dimen=outer](1.99,0.06)(0.25,1.19)
\psellipse[linewidth=0.02,dimen=outer](3.43,0.06)(0.25,1.19)
\psellipse[linewidth=0.02,dimen=outer](4.82,0.68)(0.32,0.55)
\psellipse[linewidth=0.02,dimen=outer](4.82,-0.5)(0.32,0.55)
\psline[linewidth=0.02cm](6.36,0.99)(6.36,0.37)
\psline[linewidth=0.02cm](6.36,-0.23)(6.36,-0.81)
\psline[linewidth=0.02cm](7.74,0.41)(7.74,-0.27)
\rput{-270.0}(8.45,-8.29){\psarc[linewidth=0.02](8.37,0.08){1.13}{37.234833}{142.76517}}
\end{pspicture}
\caption{\sl Three diagrams associated with contractions}
\end{center}
\end{figure}
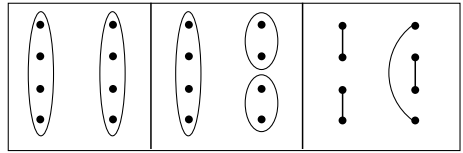

One has therefore the following combinatorial representation of the three
norms appearing in formulae (\ref{N1}) and (\ref{G*}) (the proof is
elementary, and left to the reader).

\begin{proposition}
\label{P : Last}For every $f\in L_{s}^{2}\left( \nu ^{2}\right) $, one has
that
\begin{eqnarray*}
\int_{Z}\int_{Z}f\left( z,z^{\prime }\right) ^{4}\nu \left( dz\right) \nu
\left( dz^{\prime }\right) &=&\left\Vert f\right\Vert _{L^{4}\left( \nu
^{2}\right) }^{4}=\int_{Z^{2}}f_{\sigma _{1},4}\left( z,z^{\prime }\right)
\nu \left( dz\right) \nu \left( dz^{\prime }\right) \\
\int_{Z}\left[ \int_{Z}f\left( z,z^{\prime }\right) ^{2}\nu \left( dz\right) %
\right] ^{2}\nu \left( dz^{\prime }\right) &=&\left\Vert \text{\ }f\star
_{2}^{1}f\right\Vert _{L^{2}(\nu )}^{2}=\int_{Z}f_{\sigma _{2},4}\left(
z\right) \nu \left( dz\right) \\
\int_{Z^{2}}\left[ \int_{Z}f\left( z,z^{\prime }\right) f\left( z,z^{\prime
\prime }\right) \nu \left( dz\right) \right] ^{2}\nu \left( dz^{\prime
}\right) \nu \left( dz^{\prime \prime }\right) &=&\left\Vert \text{\ }f\star
_{1}^{1}f\right\Vert _{L^{2}(\nu )}^{2}=\int_{Z}f_{\sigma _{3},4}\left(
z\right) \nu \left( dz\right) \text{.}
\end{eqnarray*}
\end{proposition}

\bigskip

In particular, Proposition \ref{P : Last} implies that, on the second Wiener
chaos of a Poisson measure, one can establish CLTs by focusing uniquely on
expressions related to three connected diagrams with four rows. Similar
characterizations for sequences belonging to chaoses of higher orders can be
deduced from the main findings of \cite{PecSoleTaqUtz}.

\bigskip

\textbf{Acknowledgements. }Part of this survey has been written while the
authors were visiting the Departement of Mathematics and Applied Statistics
of Turin University, in june 2007. The authors heartily thank Massimo\
Marinacci and Igor\ Pr\"{u}nster for their kind hospitality and support. Giovanni
Peccati acknowledges support from ISI Foundation--Lagrange Project. Murad S.
Taqqu acknowledges support by NSF under grant DMS-0706786 at Boston University. We thank
Florent Benaych-Georges, Domenico Marinucci and Marc Yor for a careful
reading of an earlier draft of this manuscript, as well as for valuable
suggestions.


\bigskip

\begin{thebibliography}{999}
\bibitem{Adler} R.J. Adler (1990). \textit{An Introduction to Continuity,
Extrema, and Related Topics for General Gaussian Processes}. Lecture
Notes-Monograph Series \textbf{12}, Institut of Mathematical Statistics,
Hayward, California.

\bibitem{Aigner} M. Aigner (1979). \textit{Combinatorial theory. }%
Springer-Verlag, Berlin Heidelberg New York.

\bibitem{Ans1} M. Anshelevich (2001). Partition-dependent stochastic
measures and $q$-deformed cumulants. \textit{Documenta Mathematica} \textbf{6%
}, 343-384.

\bibitem{Ans2} M. Anshelevich (2005). Linearization coefficients for
orthogonal polynomials using stochastic processes. \textit{The Annals of
Probabability} \textbf{33}(1), 114-136.

\bibitem{BKMP} P. Baldi, G. Kerkyacharian, D. Marinucci and D. Picard
(2006). High-frequency asymptotics for wavelet-based tests for Gaussianity
and isotropy on the torus. Preprint.

\bibitem{Barndorff...} O. Barndorff-Nielsen, J.\ Corcuera, M.\ Podolskij and
J.\ Woerner (2008). Bipower variations for Gaussian processes with
stationary increments. Preprint.

\bibitem{BillBook} P. Billingsley (1995). {\it Probability and Measure}, $3^{\rm rd}$ Edition. Wiley, New York.

\bibitem{Bitcheler} K.\ Bitcheler (1980). Stochastic integration and $L^{p}$
theory of semimartingales. \textit{The Annals of Probability }\textbf{9}(1), 49-89.

\bibitem{BrMa} P.\ Breuer et P. Major (1983). Central limit theorems for
non-linear functionals of Gaussian fields. \textit{Journal of Multivariate
Analysis} \textbf{13}, 425-441.

\bibitem{Brod} M.S. Brodskii (1971), \textit{Triangular and Jordan
Representations of Linear Operators}. Transl. Math. Monographs \textbf{32},
AMS, Providence.

\bibitem{ChaSlud} D. Chambers et E. Slud (1989). Central limit theorems for
nonlinear functionals of stationary Gaussian processes. \textit{Probability
Theory and Related Fields} \textbf{80}, 323-349.

\bibitem{CohTaq} S. Cohen et M.S. Taqqu (2004). Small and large scale
behavior of the Poissonized Telecom process. \textit{Methodology and
Computing in Applied Probability }\textbf{6}, 363-379.

\bibitem{CorNuWoe} J.M. Corcuera, D. Nualart et J.H.C. Woerner (2006). Power
variation of some integral long memory process. \textit{Bernoulli} \textbf{12%
}(4), 713-735.

\bibitem{DBPP} P.\ de Blasi, G.\ Peccati and I. Pr\"{u}nster (2008).
Asymptotics for posterior hazards. To appear in: \textit{The Annals of
Statistics.}

\bibitem{DPY} P. Deheuvels, G. Peccati et M. Yor (2006) On quadratic
functionals of the Brownian sheet and related processes. \textit{Stochastic
Processes and their Applications} \textbf{116}, 493-538.

\bibitem{DMM5} C. Dellacherie, B. Maisonneuve et P.-A. Meyer (1992). \textit{%
Probabilit\'{e}s et Potentiel (Chapitres XVII \`{a} XXIV)}. Hermann, Paris.

\bibitem{Dudley} R.M. Dudley (1967). The sizes of compact subsets of Hilbert
space and continuity of Gaussian processes. \textit{Journal of Functional
Analysis} \textbf{1}, 290-330.

\bibitem{Engel} D.D. Engel (1982). The multiple stochastic integral. \textit{%
Memoirs of the AMS }\textbf{38}, 1-82.

\bibitem{FarreJoUtz} M. Farr\'{e}, M.\ Jolis and F.\ Utzet (2008). Multiple
Stratonovich integral and Hu-Meyer formula for L\'{e}vy processes. Preprint.

\bibitem{Feigin} P.D. Feigin (1985). Stable convergence of semimartingales.
\textit{Stochastic Processes and their Applications }\textbf{19},
125-134.T.S.

\bibitem{Ferg1973} Ferguson (1973). A Bayesian analysis of some
non-parametric problems.\ \textit{The Annals of Statistics }\textbf{1} (2),
209-230.

\bibitem{FoxTaq} R. Fox et M.S. Taqqu (1987). Multiple stochastic integrals
with dependent integrators. \textit{Journal of Multivariate Analysis}
\textbf{21}(1), 105-127.

\bibitem{GineDelaPena} E. Gin\'{e} and V.H. de la Pe\~{n}a (1999). \textit{%
Decoupling}. Springer-Verlag. Berlin Heidelberg New York.

\bibitem{GiSu} L. Giraitis and D. Surgailis (1985). CLT and other limit
theorems for functionals of Gaussian processes. \textit{Zeitschrift f\"{u}r
Wahrsch. verw. Gebiete} \textbf{70}, 191-212.

\bibitem{Go Taqqu} J. Goldberg and M.S. Taqqu (1982). Regular multigraphs and their applications to the Monte Carlo evaluation of moments of non-linear functions of Gaussian processes. \textit{Stochastic Processes and their Applications} \textbf{13}, 121-138.

\bibitem{GripNorros} G. Gripenberg and I. Norros (1996). On the prediction
of fractional Brownian motion. \textit{Journal of Applied Probability}
\textbf{33}, 400-410.

\bibitem{Handa1} K.\ Handa (2005). Sampling formulae for symmetric
selection. \textit{Electronic Communications in Probabability} \textbf{10}, 223-234 (Electronic).

\bibitem{Handa2} K.\ Handa (2007). The two-parameter Poisson-Dirichlet point
process. Preprint.

\bibitem{HuNu} Y. Hu and D. Nualart (2005). Renormalized self-intersection
local time for fractional Brownian motion. \textit{The Annals of
Probabability} \textbf{33}(3), 948-983.

\bibitem{Ito} K. It\^{o} (1951). Multiple Wiener integral. \textit{J. Math.
Soc. Japan} \textbf{3}, 157--169

\bibitem{JKM} J. Jacod, A. Klopotowski et J. M\'{e}min (1982). Th\'{e}or\`{e}%
me de la limite centrale et convergence fonctionnelle vers un processus \`{a}
accroissements ind\'{e}pendants : la m\'{e}thode des martingales. \textit{%
Annales de l'Institut H. Poincar\'{e} }(\textit{PR}) \textbf{1}, 1-45.

\bibitem{JacSh} J. Jacod et A.N. Shiryaev (1987). \textit{Limit Theorems for
Stochastic Processes. }Springer, Berlin Heidelberg New York.

\bibitem{Jakubowki} A. Jakubowski (1986). Principle of conditioning in limit
theorems for sums of random variables. \textit{The Annals of Probability }%
\textbf{11}(3), 902-915.

\bibitem{JamesLijPru} L.F.\ James, A. Lijoi and I. Pr\"{u}nster (2005).
Conjugacy as a distinctive feature of the Dirichlet process. \textit{%
Scandinavian Journal of Statistics} \textbf{33}, 105-120.

\bibitem{JamesRYsurvey} L.\ James, B.\ Roynette and M.\ Yor (2008).
Generalized Gamma Convolutions, Dirichlet means, Thorin measures, with
explicit examples. \textit{Probability Surveys} \textbf{5}, 346-415.

\bibitem{Janson} S. Janson (1997). \textit{Gaussian Hilbert Spaces. }%
Cambridge University Press, Cambridge.

\bibitem{JuliaNualart} O. Juli\`{a} et D. Nualart (1988). The distribution
of a double stochastic integral with respect to two independent Brownian
Sheets. \textit{Stochastics\ }\textbf{25}, 171-182.

\bibitem{Kab} {Y. Kabanov (1975). }On extended stochastic integrals.\emph{\ }%
\textit{Theory of Probability and its applications}{\ \textbf{20}},{\
710-722. }

\bibitem{KaSz} O. Kallenberg et J. Szulga (1991). Multiple integration with
respect to Poisson and L\'{e}vy processes. \textit{Probability Theory and
Related Fields }\textbf{83}, 101-134.

\bibitem{kingman67} J.F.C. Kingman (1967). Completely random measures.
\textit{Pacific Journal of Mathematics} \textbf{21}, 59-78.

\bibitem{Kuo} H.-H. Kuo (1975). \textit{Gaussian measures in Banach spaces. }%
LNM \textbf{463}.\textit{\ }Springer-Verlag, Berlin Heidelberg New-York.

\bibitem{Kussmaul} A.U.\ Kussmaul (1977). \textit{Stochastic integration and
generalized martingales. }Pitman research notes in mathematic, \textbf{11}.
London.

\bibitem{KW91} S. Kwapie\'{n} and W.A. Woyczy\'{n}ski (1991). Semimartingale
integrals via decoupling inequalities and tangent processes. \textit{%
Probability and Mathematical Statisitics} \textbf{12}(2), 165-200.

\bibitem{KW} S. Kwapie\'{n} and W.A. Woyczy\'{n}ski (1992). \textit{Random
Series and Stochastic Integrals: Single and Multiple}. Birkh\"{a}user, Basel.

\bibitem{LeoShy} V.P. Leonov and A.N. Shiryaev (1959). On a method of
calculations of semi-invariants. \textit{Theory of Probability and its
Applications} \textbf{4}, 319-329.

\bibitem{M.A. Lifshits} M.A. Lifshits (1995). \textit{Gaussian Random
Functions. }Kluwer, Dordrecht.

\bibitem{WLindeBook} W.\ Linde (1986). \textit{Probability in Banach spaces:
stable and infinitely divisible ditributions. }Wiley, New York.

\bibitem{MPS} {J. Ma, Ph. Protter and J. San Martin (1998) }Anticipating
integrals for a class of martingales.{\ \emph{Bernoulli}, \textbf{4}},{\
81-114. }

\bibitem{Major} P. Major (1981). \textit{Multiple Wiener-It\^{o} integrals}.
LNM \textbf{849}. Springer-Verlag, Berlin Heidelberg New York.

\bibitem{Maly} V.A. Malyshev (1980). Cluster expansion in lattice models of
statistical physics and quantum fields theory. \textit{Uspehi Mat. Nauk}
\textbf{35}, 3-53.

\bibitem{Marinucci} D. Marinucci (2006). High resolution asymptotics for the
angular bispectrum of spherical random fields. \textit{The Annals of
Statistics}\emph{\ }\textbf{34}, 1-41.

\bibitem{MaPeAb} D. Marinucci and G. Peccati (2007). High-frequency
asymptotics for subordinated stationary fields on an Abelian compact group.
\textit{Stochastic Processes and their Applications }\textbf{118}(4),
585-613.

\bibitem{MaPeSphere} D. Marinucci and G. Peccati (2007). Group
representations and high-frequency central limit theorems for subordinated
random fields on a sphere. Preprint.

\bibitem{Maruyama82} G. Maruyama (1982). Applications of the multiplication
of the It\^{o}-Wiener expansions to limit theorems. \textit{Proc. Japan
Acad.\ }\textbf{58}, 388-390.

\bibitem{Maruy} G. Maruyama (1985). Wiener functionals and probability limit
theorems, I: the central limit theorem. \textit{Osaka Journal of Mathematics}
\textbf{22}, 697-732.

\bibitem{Masani} P.R. Masani (1995) The homogeneous chaos from the
standpoint of vector measures, \textit{Phil. Trans. R. Soc. Lond, A} \textbf{%
355}, 1099-1258

\bibitem{Mcc} P. McCullagh (1987). \textit{Tensor Methods in Statistics.}
Chapman and Hall. London.

\bibitem{MSW} R.D. Mauldin, W.D. Sudderth and S.C. Williams (1992). P\'{o}%
lya trees and random distributions. \textit{The Annals of Statistics},%
\textit{\ }\textbf{20} (3), 1203-1221

\bibitem{Mey78} P.-A. Meyer (1976). Un cours sur les int\'{e}grales
stochastiques. \textit{S\'{e}minaire de Probabilit\'{e}s X}, LNM \textbf{511}%
. Springer-Verlag, Berlin Heidelberg New York, pp. 245-400.

\bibitem{Mey92} P.-A. Meyer (1992). \textit{Quantum probability for
probabilists. }LNM \textbf{1538}.\textit{\ }Springer-Verlag, Berlin
Heidelberg New York.

\bibitem{NeuNourdin} A. Neuenkirch and I. Nourdin (2006). Exact rate of
convergence of some approximation schemes associated to SDEs driven by a
fractional Brownian motion. Pr\'{e}publication.

\bibitem{Neveu-1968} J. Neveu (1968). \textit{Processus Al\'{e}atoires
Gaussiens. }Presses de l'Universit\'{e} de Montr\'{e}al.

\bibitem{Nica Speicher} A.\ Nica and R.\ Speciher (2006). \textit{Lectures
on the combinatorics of free probability. }London Mathematical Society
Lecture Notes Series \textbf{335. }Cambridge University Press, Cambridge.

\bibitem{Nourdin05} I. Nourdin (2005). Sch\'{e}mas d'approximation associ%
\'{e}s \`{a} une \'{e}quation diff\'{e}rentielle dirig\'{e}e par une
fonction h\"{o}ld\'{e}rienne; cas du mouvement Brownien fractionnaire.
\textit{C.R.A.S.} \textbf{340}(8), 611-614.

\bibitem{NouNu} I.\ Nourdin and D.\ Nualart (2007). Central limit theorems
for multiple stable integrals. Preprint.

\bibitem{NNTud} I. Nourdin, D. Nualart and C.A. Tudor (2007). Central and
non-central limit theorems for weighted power variations of fractional
Brownian motion. Preprint.

\bibitem{NouPec2007} I. Nourdin and G. Peccati (2007). Non-central
convergence of multiple integrals. To appear in: \textit{The Annals of
Probability.}

\bibitem{NouPecIBM} I.\ Nourdin and G.\ Peccati (2008). Weighted power
variations of iterated Brownian motion. \textit{The Electronic Journal of
Probability.} {\bf 13}, n. 43, 1229-1256 (Electronic).

\bibitem{NouPeccPTRF} I. Nourdin and G. Peccati (2008). Stein's method on
Wiener chaos. To appear in: \textit{Probability Theory and Related Fields.}

\bibitem{NouPecexact} I. Nourdin and G. Peccati (2008). Stein's method and
exact Berry-Ess\'{e}en asymptotics for functionals of Gaussian fields.
Preprint.

\bibitem{NouPecReveillac} I. Nourdin, G. Peccati and A.\ R\'{e}veillac
(2008). Multivariate normal approximation using Stein's method and Malliavin
calculus. To appear in: \textit{Annales de l'Institut H.\ Poincar\'{e}.}

\bibitem{NouREv} I. Nourdin and A. R\'{e}veillac (2008). Asymptotic behavior
of weighted quadratic variations of fractional Brownian motion: the critical
case $H=1/4$. Preprint.

\bibitem{Nualart83} D. Nualart (1983), On the distribution of a double
stochastic integral. \textit{Z. Wahr\-schein\-lichkeit verw. Gebiete }%
\textbf{65}, 49-60

\bibitem{Nualart2} D. Nualart\ (1998). Analysis on Wiener space and
anticipating stochastic calculus. \textit{Lectures on Probability Theory and
Statistics. \'{E}cole de probabilit\'{e}s de St. Flour XXV (1995)},\textit{\
}LNM\textit{\textbf{\ }}\textbf{1690}. Springer-Verlag, Berlin Heidelberg
New York, pp. 123-227.

\bibitem{Nualart} D. Nualart (2006). \textit{The Malliavin Calculus and
related topics }(2$^{\text{\`{e}me}}$ \'{e}dition).\textit{\ }%
Springer-Verlag, Berlin Heidelberg New York.

\bibitem{NuOrtiz} D. Nualart and S. Ortiz-Latorre\ (2007). Intersection
local times for two independent fractional Brownian motions. \textit{Journal of
Theoretical Probabability }\textbf{20}(4), 759-767.

\bibitem{NuOrtizSPA} D. Nualart and S. Ortiz-Latorre\ (2008). Central limit
theorems for multiple stochastic integrals and Malliavin calculus. \textit{%
Stochastic Processes and their Applications} \textbf{118}(4), 614-628

\bibitem{NuaParx} D. Nualart and E. Pardoux (1988). Stochastic calculus with
anticipating integrands. \textit{Probabability Theory Related Fields},
\textbf{78}, 535-581.

\bibitem{PNu05} D. Nualart and G. Peccati (2005). Central limit theorems for
sequences of multiple stochastic integrals. \textit{The Annals of Probability%
}, \textbf{33}(1), 177-193.

\bibitem{NuSch} D. Nualart and W. Schoutens (2000). Chaotic and predictable
representation for L\'{e}vy processes. \textit{Stochastic Processes and
their Applications }\textbf{90}, 109-122.

\bibitem{NV} D. {Nualart and J. Vives (1990). }Anticipative calculus for the
Poisson space based on the Fock space.{\ \textit{S\'{e}minaire de Probabilit%
\'{e}s XXIV}, LNM\ \textbf{1426}. Springer-Verlag, Berlin Heidelberg New
York, pp. 154-165.}

\bibitem{gphUN} G. Peccati (2001). On the convergence of multiple random
integrals. \textit{Studia Sc. Mat. Hungarica}, \textbf{37}, 429-470.

\bibitem{Pecp} G. Peccati (2007). Gaussian approximations of multiple
integrals. \textit{Electronic Communications in Probability} \textbf{12},
350-364 (electronic).

\bibitem{PecBer2008dir} G. Peccati (2008). Multiple integral representation
for functionals of Dirichlet processes. \textit{Bernoulli }\textbf{14}(1),
91-124

\bibitem{PecPru} G.\ Peccati and I. Pr\"{u}nster (2008). Linear and
quadratic functionals of random hazard rates: an asymptotic analysis.
\textit{The Annals of Applied Probability }\textbf{18}(5), 1910-1943

\bibitem{PecSoleTaqUtz} G.\ Peccati, J.-L. Sol\'{e}, F.\ Utzet and M.S.
Taqqu (2008). Stein's method and Gaussian approximation of Poisson
functionals. Preprint.

\bibitem{PeTaqPoc} G. Peccati and M.S.\ Taqqu (2007). Stable convergence of
generalized $L^{2}$ stochastic integrals and the principle of conditioning.
\textit{The Electronic Journal of Probability}, \textbf{12}, 447-480, n. 15
(electronic).

\bibitem{PeTaqMultAOP} G. Peccati and M.S.\ Taqqu (2008). Limit theorems for
multiple stochastic integrals. To appear in: \textit{ALEA.}

\bibitem{PeTaq2bleP} G. Peccati and M.S.\ Taqqu (2008). Central limit
theorems for double Poisson integrals. \textit{Bernoulli }\textbf{14}(3),
791-821.

\bibitem{PeTaqMwi} G. Peccati and M.S. Taqqu (2007). Stable convergence of
multiple Wiener-It\^{o} integrals. To appear in: \textit{The Journal of
Theoretical Probability.}

\bibitem{PTu04} G.\ Peccati and C.A. Tudor (2005). Gaussian limits for
vector-valued multiple stochastic integrals. \textit{S\'{e}minaire de
Probabilit\'{e}s XXXVIII}, LNM \textbf{1857}. Springer-Verlag, Berlin
Heidelberg New York, pp. 247-262.

\bibitem{GPMY04a} G. Peccati and M. Yor (2004). Hardy's inequality in $%
L^{2}\left( \left[ 0,1\right] \right) $ and principal values of Brownian
local times. \textit{Asymptotic Methods in Stochastics}, AMS, Fields
Institute Communications Series, 49-74.

\bibitem{GPMY04b} G. Peccati and M. Yor (2004). Four limit theorems for
quadratic functionals of Brownian motion and Brownian bridge. \textit{%
Asymptotic Methods in Stochastics}, AMS, Fields Institute Communication
Series, 75-87.

\bibitem{PipTaqPtrf} V. Pipiras and M.S. Taqqu (2000). Integration questions
related to fractional Brownian motion. \textit{Probability Theory and
Related Fields }\textbf{118}(2), 251-291.

\bibitem{PipTaq} V. Pipiras and M.S. Taqqu (2001). Are classes of
deterministic integrands for fractional Brownian motion complete? \textit{%
Bernoulli }\textbf{7}(6), 873-897

\bibitem{PipTaqSurv} V. Pipiras and M.S. Taqqu (2003). Fractional calculus
and its connection to fractional Brownian motion. In: \textit{Long Range
Dependence}, 166-201, Birkh\"{a}user, Basel.

\bibitem{Pit} J. Pitman (2006). \textit{Combinatorial Stochastic Processes. }%
LNM\textit{\ }\textbf{1875}.\textit{\ }Springer-Verlag, Berlin Heidelberg
New York.

\bibitem{P1} N. Privault (1994). Chaotic and variational calculus in
discrete and continuous time for the Poisson process. \textit{Stochastics
and Stochastics Reports} \textbf{51}, 83-109.

\bibitem{P2} N. Privault (1994). In\'{e}galit\'{e}s de Meyer sur l'espace de
Poisson. \textit{C.R.A.S.} \textbf{318}, 559-562.

\bibitem{PrivaultSoleVives} N. Privault, J.L. Sol\'{e} and J.\ Vives (2001).
Chaotic Kabanov formula for the Az\'{e}ma martingales. \textit{Bernoulli}
\textbf{6}(4), 633-651.

\bibitem{PrivaultWu} N. Privault et J.-L. Wu (1998). Poisson stochastic
integration in Hilbert spaces. \textit{Ann. Math. Blaise Pascal} \textbf{6}%
(2), 41-61.

\bibitem{Protter} P. Protter (2005). \textit{Stochastic Integration and
Differential Equations} (2$^{\text{\`{e}me}}$ \'{e}dition). Springer-Verlag,
Berlin Heidelberg New York.

\bibitem{Raj Ros} B.S. Rajput and J. Rosinski (1989). Spectral
representation of infinitely divisible processes. \textit{Probability Theory
and Related Fields} \textbf{82}, 451-487.

\bibitem{RY} D. Revuz and M. Yor (1999). \textit{Continuous martingales and
Brownian motion. }Springer-Verlag, Berlin Heidelberg New York.

\bibitem{Ro SHen} G.-C. Rota and J. Shen (2000). On the combinatorics of
cumulants. \textit{Journal of Combinatorial Theory Series A} \textbf{91},
283-304.

\bibitem{RoWa} G.-C. Rota and C. Wallstrom (1997). Stochastic integrals: a
combinatorial approach. \textit{The Annals of Probability} \textbf{25}(3),
1257-1283.

\bibitem{RosWoy} J. Rosi\'{n}sky and W.A. Woyczy\'{n}ski (1984). Products of
random measures, multilinear random forms and multiple stochastic integrals.
\textit{Proc. Conference of Measure Theory}, Oberwolfach 1983, LNM 1089.
Springer-Verlag, Berlin Heidelberg New York, pp. 294-315.

\bibitem{RussoVal1997} F. Russo and P. Vallois (1998). Product of two
multiple stochastic integrals with respect to a normal martingale. \textit{%
Stochastic Processes and their Applications} \textbf{73}(1), 47-68.

\bibitem{SamoTaqqu} G. Samorodnitsky and M.S. Taqqu (1994). \textit{Stable
Non-Gaussian Random Processes. }Chapman and Hall, New York.

\bibitem{Sato} K.-I. Sato (1999). \textit{L\'{e}vy Processes and Infinitely
Divisible Distributions. }Cambridge Studies in Advanced Mathematics \textbf{%
68}. Cambridge University Press.

\bibitem{SchoutBook} W. Schoutens (2000). Stochastic processes and orthogonal polynomials. Lecture Notes in Staistics 146. Springer-Verlag, Berlin Heidelberg New York.

\bibitem{Sch} M. Schreiber (1969). Fermeture en probabilit\'{e} de certains
sous-espaces d'un espace $L^{2}$. \textit{Zeitschrift Warsch. verw.\ Gebiete
}\textbf{14}, 36-48.

\bibitem{Shir} A.N. Shyryayev (1984). \textit{Probability}. Springer-Verlag.
Berlin Heidelberg New York.

\bibitem{Slud} E.V. Slud (1993). The moment problem for polynomial forms in
normal random variables. \textit{The Annals of Probabability} \textbf{21}%
(4), 2200-2214.

\bibitem{SolUtzAOP} J.-L. Sol\'{e} and F. Utzet (2008). Time-space harmonic polynomials associated with a L\'{e}vy process. \textit{Bernoulli} {\bf 14}(1), 1-13.

\bibitem{SolUtzBer} J.-L. Sol\'{e} and F. Utzet (2008). On the orthogonal polynomials associated to a L\'{e}vy process. To appear in: \textit{The Annals of Probability}.

\bibitem{Speed} T. Speed (1983). Cumulants and partitions lattices. \textit{%
Australian Journal of Statistics\ }\textbf{25}(2), 378-388.

\bibitem{Stanley} R. Stanley (1997). \textit{Enumerative combinatorics, Vol.
1. }Cambridge University Press.

\bibitem{STeinbook} Ch. Stein (1986). \textit{Approximate computation of
expectations.} Institute of Mathematical Statistics Lecture Notes -
Monograph Series, \textbf{7}. Institute of Mathematical Statistics, Hayward,
CA.

\bibitem{Stroock} D.W. Stroock (1987). Homogeneous chaos revisited. \textit{S%
\'{e}minaire de Probabilit\'{e}s XXI}, LNM \textbf{1247}, Springer-Verlag,
Berlin Heidelberg New York, pp. 1-8.

\bibitem{Surg1984} D. Surgailis (1984). On multiple Poisson stochastic
integrals and associated Markov semigroups. \textit{Probab. Math. Statist.}
\textbf{3}(2), 217-239.

\bibitem{Sur} D. Surgailis (2000). CLTs for Polynomials of Linear Sequences:
Diagram Formulae with Applications. Dans : \textit{Long Range Dependence}.
Birkh\"{a}user, Basel, pp. 111-128.

\bibitem{SUr2} D. Surgailis (2000). Non-CLT's: U-Statistics, Multinomial
Formula and Approximations of Multiple Wiener-It\^{o} integrals. Dans :
\textit{Long Range Dependence}. Birkh\"{a}user, Basel, pp. 129-142.

\bibitem{Versh} N.V. Tsilevich and A.M. Vershik (2003). Fock factorizations
and decompositions of the $L^{2}$ spaces over general L\'{e}vy processes.
\textit{Russian Math. Surveys} \textbf{58}(3), 427-472.

\bibitem{TsVeYorJFA} N.\ Tsilevich, A.M.\ Vershik and M.\ Yor (2001). An
infinite-dimensional analogue of the Lebesgue measure and distinguished
properties of the gamma process. \textit{J. Funct. Anal.} \textbf{185}(1),
274-296.

\bibitem{CT} C. Tudor (1997). Product formula for multiple Poisson-It\^{o}
integrals. \textit{Revue Roumaine de Math. Pures et Appliqu\'{e}es}\emph{\ }%
\textbf{42}(3-4), 339-345.

\bibitem{TV} {C.A. Tudor and J. Vives (2002) }The indefinite Skorohod
integral as integrator on the Poisson space.\emph{\ }\textit{Random
Operators and Stochastic Equations}{\ }\textbf{10}, 29-46.

\bibitem{UZ} A.S. \"{U}st\"{u}nel and M. Zakai (1997). The construction of
filtrations on Abstract Wiener Space. \textit{Journal of Functional
Analysis\ }\textbf{143}, 10-32.

\bibitem{Vitale} R.A. Vitale (1990). Covariances of symmetric statistics.
\textit{Journal of Multivariate Analysis} \textbf{41},\textbf{\ }14-26.

\bibitem{Wiener1938} N. Wiener (1938). The homogeneous chaos. \textit{Amer.
J. Math.} \textbf{60}, 879-936.

\bibitem{WolTaq} R.L. Wolpert and M.S. Taqqu (2005). Fractional
Ornstein-Uhlenbeck L\'{e}vy Processes and the Telecom Process: Upstairs and
Downstairs. \textit{Signal Processing }\textbf{85}(8), 1523-1545.

\bibitem{Wu} L.M. Wu (1990). Un traitement unifi\'{e} de la repr\'{e}%
sentation des fonctionnelles de Wiener. \textit{S\'{e}minaire de Probabilit%
\'{e}s XXIV}, LNM \textbf{1426}, Springer-Verlag, Berlin Heidelberg New
York, pp.\ 166-187.

\bibitem{Xue} X.-H. Xue (1991). On the principle of conditioning and
convergence to mixtures of distributions for sums of dependent random
variables. \textit{Stochastic Processes and their Applications }\textbf{37}%
(2), 175-186.

\bibitem{Yoshida} K. Yosida (1980). \textit{Functional analysis}.
Springer-Verlag, Berlin Heidelberg New York.
\end{thebibliography}
\end{document}